\documentclass[10pt]{article}
\usepackage{amssymb}
\usepackage[mathscr]{eucal}
\newcommand{\bsf}[1]{{\textbf{\textsf{#1}}}}
\newcommand{\Gup}[2]{{\displaystyle \bigcup_{#1}^{#2}}}
\DeclareMathAlphabet{\mathpzc}{OT1}{pzc}{m}{it}
\usepackage{bbm}
\usepackage{dsfont}
\usepackage{color}
\usepackage{xcolor}
\usepackage{mathrsfs}
\usepackage{amssymb}
\usepackage{calligra}
\usepackage{fontenc}
\usepackage{amsbsy}
\usepackage{amsmath}

\newcommand{\bfpzc}[1]{{\pmb{\mathpzc{#1}}}}
\usepackage{graphicx}
\graphicspath{%
    {converted_graphics/}
    {/}
}
\usepackage{longtable}
\begin{document}

\setlength\topmargin{-0.4in}
\newcommand{\simge}{\ba{cc}\vspace*{-2.4mm}>\\ \sim\ea }
\newcommand{\simle}{\ba{cc}\vspace*{-2.4mm}<\\ \sim\ea }
\newcommand{\Cdot}{\!\cdot\!}
\newcommand{\sq}{{$\sqcap\!\!\!\!\sqcup$}}
\newcommand{\Eu}{{\rm I\,\!\! E}}
\newcommand{\Io}{\Int{\Omega}{}}
\newcommand{\Id}{\Int{\cald}{}}
\newcommand{\Div}{\mbox{\rm div}\,}
\newcommand{\tr}{\mbox{\rm tr}\,}
\newcommand{\grad}{\mbox{\rm grad}\,}
\newcommand{\supp}{\mbox{\rm supp}\,}
\newcommand{\curl}{\mbox{\rm curl}\,}
\newcommand{\Ido}{\Int{\partial\Omega}{}}
\newcommand{\IdS}{\Int{\Sigma}{}}
\newcommand{\Oint}[2]{{\displaystyle \oint_{#1}^{#2}}}
\newcommand{\Int}[2]{{\displaystyle \int_{ #1}^{ #2}}}
\newcommand{\Lim}[1]{{\displaystyle \lim_{ #1}}}
\newcommand{\Limsup}[1]{{\displaystyle \limsup_{\footnotesize #1}}}
\newcommand{\Liminf}[1]{{\displaystyle \liminf_{\footnotesize #1}}}
\newcommand{\Sup}[1]{{\displaystyle \sup_{#1}}}
\newcommand{\Inf}[1]{{\displaystyle \inf_{#1}}}
\newcommand{\Max}[1]{{\displaystyle \max_{#1}}}
\newcommand{\Min}[1]{{\displaystyle \min_{#1}}}
\newcommand{\Sum}[2]{{\displaystyle \sum_{#1}^{#2}}}
\newcommand{\Prod}[2]{{\displaystyle \prod_{#1}^{#2}}}
\newcommand{\BCup}[2]{{\displaystyle \bigcup_{#1}^{#2}}}
\newcommand{\BCap}[2]{{\displaystyle \bigcap_{#1}^{#2}}}
\newcommand{\Ccup}[2]{{\displaystyle \cup_{#1}^{#2}}}
\newcommand{\Ccap}[2]{{\displaystyle \cap_{#1}^{#2}}}
\newcommand{\Frac}[2]{\displaystyle{\frac{\displaystyle{#1}}{\displaystyle{#2}}}}
\newcommand{\norm}[1]{\left\|{#1}\right\|}
\newcommand{\Norm}[1]{\langle\langle{#1}\rangle\rangle_q}
\newcommand{\No}[1]{\langle\!\langle{#1}\rangle\!\rangle}
\newcommand{\NO}[1]{{\langle{#1}\rangle}_{\lambda,q}}
\newcommand{\beea}{\begin{eqnarray}}
\newcommand{\eeea}{\end{eqnarray}}
\newcommand{\ms}{\medskip\smallskip}
\newcommand{\bs}{\bigskip}
\newcommand{\ps}{\par\smallskip}
\newcommand{\bfe}{{\mbox{\boldmath $e$}} }
\newcommand{\pni}{{\par\noindent}}
\newcommand{\bfq}{{\mbox{\boldmath $q$}} }
\newcommand{\bfz}{{\mbox{\boldmath $z$}} }
\newcommand{\0}{{\mbox{\boldmath $0$}} }
\newcommand{\LE}{\!\!\!&\le&\!\!\!}
\newcommand{\BL}[1]{{\par\smallskip{\bf Lemma #1.}}}
\newcommand{\BT}[1]{{\par\smallskip{\bf Theorem #1.}}}
\newcommand{\Ln}{[\!|}
\newcommand{\Rn}{|\!]}
\newcommand{\n}[1]{{\Ln{#1}\Rn}} 
\newcommand{\nq}[1]{{\Ln{#1}\Rn}_{q}} 
\newcommand{\nqr}[1]{{\Ln{#1}\Rn}_{q,r}} 
\newcommand{\Nq}[1]{{\langle{#1}\rangle}_{q}} 
\newcommand{\Nql}[1]{{\langle{#1}\rangle}_{\lambda,q}} 
\newcommand{\Nqr}[1]{{\langle{#1}\rangle}_{q,r}}
\newcommand{\N}[1]{{|\!\!|\!\!|\,{#1}\,|\!|\!\!|_2}}
\newcommand{\EA}[2]{$$#1$$%
\vspace{-6.mm}
\begin{equation}
\end{equation}
\vspace{-6.mm}
$$
#2
\setlength{\belowdisplayskip}{3mm}
\setlength{\belowdisplayshortskip}{3mm}
$$
}
\newcommand{\A}[2]{$$#1$$%
\vspace{-4.mm}
$$
#2
\setlength{\belowdisplayskip}{3mm}
\setlength{\belowdisplayshortskip}{3mm}
$$
}
\newcommand{\BF}{\begin{footnotesize}}
\newcommand{\EF}{\end{footnotesize}}
\setlength{\jot}{.15in}
\newcommand{\pde}[2]{{\displaystyle \frac{\mbox{$\partial #1$}}{\mbox{$\partial #2$}}}}
\newcommand{\ode}[2]{{\displaystyle \frac{\mbox{$d #1$}}{\mbox{$d #2$}}}}
\newcommand{\f}[2]{\frac{\mbox{$#1$}}{\mbox{$ #2$}}}
\newcommand{\bi}{\begin{itemize}}
\newcommand{\ei}{\end{itemize}}
\newcommand{\ed}{\end{document}}
\newcommand{\be}{\begin{equation}}
\newcommand{\ba}{\begin{array}}
\newcommand{\ea}{\end{array}}
\newcommand{\ee}{\end{equation}}
\newcommand{\eeq}[1]{\label{eq:#1}\end{equation}}
\newcommand{\real}{{\mathbb R}}
\newcommand{\compl}{{\mathbb C}}
\def\Id{\mbox{\boldmath $1$}}
\def\zero{\mbox{\boldmath $0$}}
\newcommand{\PP}{{\rm I\!\!\,P}}
\newcommand{\nat}{{\mathbb N}}
\newcommand{\bfpsi}{\mbox{\boldmath $\psi$}}
\newcommand{\bfchi}{\mbox{\boldmath $\chi$}}
\newcommand{\bfomega}{\mbox{\boldmath $\omega$}}
\newcommand{\bfvaromega}{\mbox{\boldmath $\varpi$}}
\newcommand{\bfOmega}{\mbox{\boldmath $\Omega$}}
\newcommand{\bfTheta}{\mbox{\boldmath $\Theta$}}
\newcommand{\bfxi}{\mbox{\boldmath $\xi$}}
\newcommand{\bfmu}{\mbox{\boldmath $\mu$}}
\newcommand{\bfx}{\mbox{\boldmath $x$}}
\newcommand{\bfy}{\mbox{\boldmath $y$}}
\newcommand{\bfPsi}{\mbox{\boldmath $\Psi$}}
\newcommand{\bfphi}{\mbox{\boldmath $\varphi$}}
\newcommand{\bfhi}{\mbox{\boldmath $\phi$}}
\newcommand{\bfPhi}{\mbox{\boldmath $\Phi$}}
\newcommand{\bfv}{{\mbox{\boldmath $v$}} }
\newcommand{\bfu}{{\mbox{\boldmath $u$}} }
\newcommand{\bfsf}{{\mbox{\footnotesize\boldmath $s$}} }
\newcommand{\bfuf}{{\mbox{\footnotesize\boldmath $u$}} }
\newcommand{\bfw}{{\mbox{\boldmath $w$}} }
\newcommand{\bff}{{\mbox{\boldmath $f$}} }
\newcommand{\bfa}{{\mbox{\boldmath $a$}} }
\newcommand{\bfi}{{\mbox{\boldmath $i$}} }
\newcommand{\bfj}{{\mbox{\boldmath $j$}} }
\newcommand{\bfc}{{\mbox{\boldmath $c$}} }
\newcommand{\bfo}{{\mbox{\boldmath $o$}} }
\newcommand{\bfp}{{\mbox{\boldmath $p$}} }
\newcommand{\bft}{{\mbox{\boldmath $t$}} }
\newcommand{\bfd}{{\mbox{\boldmath $d$}} }
\newcommand{\bfl}{{\mbox{\boldmath $l$}} }
\newcommand{\bfr}{{\mbox{\boldmath $r$}} }
\newcommand{\bfk}{{\mbox{\boldmath $k$}} }
\newcommand{\bfA}{{\mbox{\boldmath $A$}} }
\newcommand{\bfS}{{\mbox{\boldmath $S$}} }
\newcommand{\bfO}{{\mbox{\boldmath $O$}} }
\newcommand{\bfM}{{\mbox{\boldmath $M$}} }
\newcommand{\bfP}{{\mbox{\boldmath $P$}} }
\newcommand{\bfB}{{\mbox{\boldmath $B$}} }
\newcommand{\bfR}{{\mbox{\boldmath $R$}} }
\newcommand{\bfC}{{\mbox{\boldmath $C$}} }
\newcommand{\bfD}{{\mbox{\boldmath $D$}} }
\newcommand{\bfQ}{{\mbox{\boldmath $Q$}} }
\newcommand{\bfZ}{{\mbox{\boldmath $Z$}} }
\newcommand{\bfG}{{\mbox{\boldmath $G$}} }
\newcommand{\bfE}{{\mbox{\boldmath $E$}} }
\newcommand{\bfX}{{\mbox{\boldmath $X$}} }
\newcommand{\bfY}{{\mbox{\boldmath $Y$}} }
\newcommand{\bfH}{{\mbox{\boldmath $H$}} }
\newcommand{\bfI}{{\mbox{\boldmath $I$}} }
\newcommand{\bfJ}{{\mbox{\boldmath $J$}} }
\newcommand{\bfN}{{\mbox{\boldmath $N$}} }
\newcommand{\bfh}{{\mbox{\boldmath $h$}} }
\newcommand{\bfm}{{\mbox{\boldmath $m$}} }
\newcommand{\bfone}{{\mbox{\boldmath $1$}} }
\newcommand{\hs}{{\rm I}\!\!\,{\rm R}^3_+}
\newcommand{\cala}{{\cal A}}
\newcommand{\calb}{{\cal B}}
\newcommand{\calc}{{\cal C}}
\newcommand{\cald}{{\cal D}}
\newcommand{\cale}{{\cal E}}
\newcommand{\calf}{{\cal F}}
\newcommand{\calg}{{\cal G}}
\newcommand{\calh}{{\cal H}}
\newcommand{\cali}{{\cal I}}
\newcommand{\calj}{{\cal J}}
\newcommand{\calk}{{\cal K}}
\newcommand{\call}{{\cal L}}
\newcommand{\calm}{{\cal M}}
\newcommand{\caln}{{\cal N}}
\newcommand{\calo}{{\cal O}}
\newcommand{\calp}{{\cal P}}
\newcommand{\calq}{{\cal Q}}
\newcommand{\calr}{{\cal R}}
\newcommand{\cals}{{\cal S}}
\newcommand{\calt}{{\cal T}}
\newcommand{\calu}{{\cal U}}
\newcommand{\calv}{{\cal V}}
\newcommand{\calx}{{\cal X}}
\newcommand{\caly}{{\cal Y}}
\newcommand{\calw}{{\cal W}}
\newcommand{\calz}{{\cal Z}}
\newcommand{\bfsigma}{\mbox{\boldmath $\sigma$}}
\newcommand{\bfSigma}{\mbox{\boldmath $\Sigma$}}
\newcommand{\bftau}{\mbox{\boldmath $\tau$}}
\newcommand{\bfeta}{\mbox{\boldmath $\eta$}}
\newcommand{\bfT}{{\mbox{\boldmath $T$}} }
\newcommand{\bfV}{{\mbox{\boldmath $V$}} }
\newcommand{\bfU}{{\mbox{\boldmath $U$}} }
\newcommand{\bfW}{{\mbox{\boldmath $W$}} }
\newcommand{\bfF}{{\mbox{\boldmath $F$}} }
\newcommand{\bfK}{{\mbox{\boldmath $K$}} }
\newcommand{\bfL}{{\mbox{\boldmath $L$}} }
\newcommand{\bfb}{{\mbox{\boldmath $b$}} }
\newcommand{\bfg}{{\mbox{\boldmath $g$}} }
\newcommand{\bfn}{{\mbox{\boldmath $n$}} }
\newcommand{\bfs}{{\mbox{\boldmath $s$}} }
\newcommand{\cf}{{\it cf.} }
\newcommand{\io}{\int_\Omega}
\newcommand{\1}{\item[({\it i})]}
\newcommand{\2}{\item[({\it ii})]}
\newcommand{\3}{\item[({\it iii})]}
\newcommand{\4}{\item[({\it iv})]}
\newcommand{\5}{\item[({\it v})]}
\newcommand{\6}{\item[({\it vi})]}
\newcommand{\7}{\item[({\it vii})]}
\newcommand{\8}{\item[({\it viii})]}
\newcommand{\9}{\item[({\it xi})]}
\newcommand{\ido}{\int_{\partial\Omega}}
\newcommand{\half}{\mbox{$\frac{1}{2}$}}
\def\parallel{\|}
\def\mid{|}
\def\Bbb R{\real}
\def\hat{\widehat}
\def\tilde{\widetilde}
\def\bar{\overline}
\newcommand{\threehalves}{3\over 2}
\newcommand{\bfPi}{\mbox{\boldmath $\Pi$}}
\newcommand{\bfXi}{\mbox{\boldmath $\Xi$}}
\newcommand{\bfalpha}{\mbox{\boldmath $\alpha$}}
\newcommand{\bfbeta}{\mbox{\boldmath $\beta$}}
\newcommand{\bfgamma}{\mbox{\boldmath $\gamma$}}
\newcommand{\bfdelta}{\mbox{\boldmath $\delta$}}
\newcommand{\bfzeta}{\mbox{\boldmath $\zeta$}}
\newcommand{\bfUpsilon}{\mbox{\boldmath $\Upsilon$}}
\newcommand{\bfGamma}{\mbox{\boldmath $\Gamma$}}
\newcommand{\bfvGamma}{\mbox{\boldmath $\varGamma$}}
\newcommand{\bfcala}{\mbox{\boldmath ${\cal A}$}}
\newcommand{\bfcalm}{\mbox{\boldmath ${\cal M}$}}
\newcommand{\bfcaln}{\mbox{\boldmath ${\cal N}$}}
\newcommand{\bfcalq}{\mbox{\boldmath ${\cal Q}$}}
\newcommand{\bfcalb}{\mbox{\boldmath ${\cal B}$}}
\newcommand{\bfcalc}{\mbox{\boldmath ${\cal C}$}}
\newcommand{\bfcali}{\mbox{\boldmath ${\cal I}$}}
\newcommand{\bfcalg}{\mbox{\boldmath ${\cal G}$}}
\newcommand{\bfcalh}{\mbox{\boldmath ${\cal H}$}}
\newcommand{\bfcalk}{\mbox{\boldmath ${\cal K}$}}
\newcommand{\bfcalt}{\mbox{\boldmath ${\cal T}$}}
\newcommand{\bfcalx}{\mbox{\boldmath ${\cal X}$}}
\newcommand{\bfcall}{\mbox{\boldmath ${\cal L}$}}
\newcommand{\bfcalf}{\mbox{\boldmath ${\cal F}$}}
\newcommand{\bfcalr}{\mbox{\boldmath ${\cal R}$}}
\newcommand{\bfcals}{\mbox{\boldmath ${\cal S}$}}
\newcommand{\bfcalw}{\mbox{\boldmath ${\cal W}$}}
\newcommand{\bfcalu}{\mbox{\boldmath ${\cal U}$}}
\newcommand{\bfcalv}{\mbox{\boldmath ${\cal V}$}}
\newcommand{\bfcalz}{\mbox{\boldmath ${\cal Z}$}}
\pagenumbering{roman}
\newcommand{\art}[6]{{\I[{\sc #1,}] {#2}, {\it #3}, {\bf #4}, {#5} {[#6]}}}
\newcommand{\ED}{\end{description}}
\newcommand{\I}{\item }
\newcommand{\ra}{\rm a}
\newcommand{\rb}{\rm b}
\newcommand{\rc}{\rm c}
\newcommand{\Hsp}{{\rm I}\!\!\,{\rm R}^n_+}
\newcommand{\Hsn}{{\rm I}\!\!\,{\rm R}^n_-}
\newcommand{\po}[1]{\mbox{$\displaystyle \frac{\mbox{$\partial #1$}}
{\mbox{$\partial x_{1}$}}$}}
\newcommand{\PO}[1]{\mbox{$\displaystyle \frac{\mbox{$\partial #1$}}
{\mbox{$\partial y_{1}$}}$}}
\newcommand{\OP}{\left(\Delta+2\lambda\PO{}\right)}
\newcommand{\op}{\left(\Delta+2\lambda\po{}\right)}
\newcommand{\ft}[1]{
\Frac{1}{(2\pi)^{n/2}}\Int{{\Bbb R}^{n}}{}e^{i{\bf x}\cdot \bfxi}
#1(\xi)d\xi}
\newcommand{\Ft}[1]{
\Frac{1}{2\pi}\Int{{\Bbb R}^{2}}{}e^{i{x}\cdot \xi}
#1(\xi)d\xi}
\newcommand{\Z}{\item[({\it a})]}
\newcommand{\B}{\item[({\it b})]}
\newcommand{\C}{\item[({\it c})]}
\newcommand{\D}{\item[({\it d})]}
\newcommand{\E}{\item[({\it e})]}
\newcommand{\G}{\item[({\it g})]}
\def\tag{\renewcommand{\theequation}}
\newcommand{\Footnote}{~\footnote}
\newcommand{\ie}{{\it i.e.}}
\newcommand{\dist}{\mbox{\rm dist\,}}
\newcommand{\const}{\mbox{\rm const}}
\newcommand{\trace}{\mbox{\rm trace}}
\newcommand{\Bo}{\par\hfill{$\Box$}\par\noindent}
\newcommand{\Nor}[1]{\langle{#1}\rangle_q}
\newcommand{\vs}{\vspace*{.5cm}\par\noindent}
\newcommand{\Vs}{\vspace*{.6cm}\par\noindent}
\newcommand{\Vvs}{\vspace*{.7cm}\par\noindent}
\newcommand{\VVs}{\vspace*{.8cm}\par\noindent}
\newtheorem{definition}{Definition}[section]
\newcommand{\Bd}{\begin{definition}\begin{it}}
\newcommand{\Ed}{\end{it}\end{definition}}
\newtheorem{remark}{Remark}[section]
\newcommand{\Br}{\begin{remark}\begin{rm}}
\newcommand{\Er}{\end{rm}\end{remark}}
\newcommand{\EED}[1]{\end{it}\label{definition:#1}\end{definition}}
\newcommand{\defref}[1]{{\rm Definition \ref{definition:#1}}}
\newtheorem{proposition}{Proposition}[section]
\newcommand{\Bp}{\begin{proposition}\begin{sl}}
\newcommand{\EP}[1]{\end{sl}\label{proposition:#1}\end{proposition}}
\newcommand{\propref}[1]{{\rm Proposition \ref{proposition:#1}}}
\newcommand{\Bt}{\begin{theorem}\begin{sl}}
\newcommand{\Et}{\end{sl}\end{theorem}}
\newcommand{\Bl}{\begin{lemma}\begin{sl}}
\newcommand{\El}{\end{sl}\end{lemma}}
\newtheorem{theorem}{Theorem}[section]
\newtheorem{lemma}{Lemma}[section]
\newtheorem{corollary}{Corollary}[section]
\newcommand{\Eqref}[1]{{\rm (\ref{eq:#1})}}
\newcommand{\Bc}{\begin{corollary}\begin{sl}}
\newcommand{\Ec}{\end{sl}\end{corollary}}
\newcommand{\ET}[1]{\end{sl}\label{theorem:#1}\end{theorem}}
\newcommand{\EL}[1]{\end{sl}\label{lemma:#1}\end{lemma}}
\newcommand{\exempref}[1]{{\rm Example \ref{example:#1}}}
\newcommand{\theoref}[1]{{\rm Theorem \ref{theorem:#1}}}
\newcommand{\ER}[1]{\end{rm}\label{remark:#1}\end{remark}}
\newcommand{\EC}[1]{\end{sl}\label{corollary:#1}\end{corollary}}
\newcommand{\remref}[1]{{\rm Remark \ref{remark:#1}}}
\newcommand{\cororef}[1]{{\rm Corollary \ref{corollary:#1}}}
\newcommand{\lemmref}[1]{{\rm Lemma \ref{lemma:#1}}}
\newcommand{\essup}[1]{{\rm ess}\,{{\displaystyle \sup_{\hspace*{-5mm}{#1}}}}}

\pagenumbering{arabic}
\newcommand{\DIV}{{\rm Div}\,}
\newcommand{\Rey}{{\rm Re}\,}
\newcommand{\Grad}{{\rm Grad}\,}
\renewcommand{\thefootnote}{(\arabic{footnote})}
\title{On Self-Propulsion by Oscillations in a Viscous Liquid } 
\medskip\bigskip
\author{Giovanni P. Galdi~\thanks{Department of Mechanical Engineering and Materials Science, University of  Pittsburgh, U.S.A.}\quad       Boris Muha~\thanks{Department of Mathematics, Faculty of Science, University of Zagreb, Croatia}\quad  \& Ana  Rado\v{s}evi\'c~\thanks{Department of Mathematics, Faculty of Economics and Business, University of Zagreb, Croatia}
}
\date{}
\maketitle
\begin{abstract}
Suppose that a body $\mathscr B$ can move by translatory motion with velocity $\bfgamma$ in an otherwise quiescent Navier-Stokes liquid, $\mathscr L$, filling the entire
 space outside $\mathscr B$. 
 Denote by $\Omega = \Omega(t)$, $t\in\real$, the one-parameter family of bounded, sufficiently smooth
 domains of $\real^3$, each one representing the configuration of $\mathscr B$ at time $t$ with respect to a
 frame with the origin at the center of mass $G$ and axes parallel to those of an inertial frame. We assume
 that there are no external forces acting on the coupled system $\mathscr S := \mathscr B +\mathscr L$ and that the
 only driving mechanism is a prescribed change in shape of $\Omega$ with time.
 The self-propulsion problem that we would like to address can be thus qualitatively
 formulated as follows. Suppose that $\mathscr B$ changes its shape in a given time-periodic fashion,
 namely, $\Omega(t+T) = \Omega(t)$, for some $T > 0$ and all $t \in \real$. Then, find necessary and sufficient
 conditions on the map $t\mapsto \Omega(t)$ securing that $\mathscr B$ self-propels, that is, 
 $G$ covers any given finite distance in a finite time. We show that this problem is solvable, in a suitable function class, provided the amplitude of the oscillations is below a given constant. Moreover, we provide examples where the propelling velocity of $\mathscr B$ is explicitly evaluated in terms of the physical parameters  and the frequency of oscillations.
\ms\par\noindent
\smallskip\par\noindent
{\bf Primary AMS code.} 76D05, 35A01, 35B10, 35Q70, 35Q74.

\end{abstract}
\renewcommand{\theequation}{\arabic{section}.\arabic{equation}}
\section{Introduction}
Self-propulsion of bodies in water or air has always been an intriguing topic of research.
The fundamental problem that one wants to investigate can be roughly formulated as follows\vspace*{-2mm}
\begin{quote}
{\em How can living creatures or mechanical devices move in a fluid by changing the shape of their bodies?}\end{quote}\vspace*{-2mm}
\par\noindent
This is one of the two questions (the other being the nature of turbulence) that tormented Leonardo da Vinci throughout his life, and to which he also dedicated a short essay (Codex ``On the Flight of Birds") now preserved in the Royal Library of Turin \cite{LdV}.
\par
The first systematic study of the locomotion of aquatic and aerial animals dates back to 1681, with the famous treatise of the Neapolitan physiologist and physicist Giovanni Alfonso Borelli \cite{Bore}. Borelli, in fact, is also credited with designing the first submarine \cite{DHR}. In modern times, the topic has been further studied by James Gray, in particular with his work on the motion of eels \cite{Gray2} and the paradoxical conclusion regarding the swimming efficiency of dolphins \cite{Gray1}; see also \cite{Gray}.
\par
The use of a mathematical approach to the study of self-propulsion is, however, relatively recent and began with the seminal work of G.I. Taylor \cite{Ta}. Taylor's analysis concerns the intriguing problem of the motion of microorganisms in a liquid at zero Reynolds number, that is, in the absence of inertia of the liquid. The remarkable feature of such motion is that it cannot be generated by a reversible periodic oscillation of the body, since whatever the creature would accomplish with one flap of a part of his body will be immediately lost with the next reverted flap. This argument was later made precise in the well-known ``scallop theorem" of E.M. Purcell \cite{EP}.
Since Taylor's pioneering paper, several notable contributions have been made to the mathematical modeling and corresponding quantitative analysis of  self-propulsion of both rigid and deformable bodies in a fluid, both at small and large scales. The list of such contributions is too long to be mentioned in full here.\Footnote{We are interested in self-propulsion in a viscous liquid. The same problem in a {\em inviscid} liquid has also received a large amount of significant contributions. For this, we refer, e.g., to the recent work \cite{NaU} and the literature cited there.} In addition to the monographs \cite{Chil,Ligh,WBB} and the references therein, it is worth mentioning, in particular, the studies performed in \cite{Ke,Pu,Se} in the case of rigid bodies,   and, for deformable bodies, those in \cite{SW1,SW2} on the relation between shape change and propulsion velocity and \cite{FJ1,FJ2,FJ3} on the swimming of spheres in a Navier-Stokes liquid.  
It must be emphasized that, probably due to the inherent difficulty of the problem, all  works mentioned above are based on different degrees of approximations, including linearization of the relevant equations and formal amplitude expansions combined with quasi-steady assumptions. 
\par
A {\em rigorous} and {\em quantitative} mathematical analysis of the self-propulsion of a (finite) body, $\mathscr B$, in a Navier-Stokes liquid began with the work \cite{Gasp}; see also \cite{Gah}. There, $\mathscr B$ is assumed {\em rigid} and moving with time-independent motion. The propulsion is generated by a non-zero momentum flux across its boundary, $\varSigma$, or by portions of it in tangential motion (or by a combination of both mechanisms). Thus, $\varSigma$ acts as the engine of $\mathscr B$ and the velocity distribution $\bfv_*$ on $\varSigma$ as its thrust. In addition to the well-posedness of the problem 
--~and more importantly~-- in \cite{Gasp}  necessary and sufficient conditions were provided on $\bfv_*$ to generate thrust. Further, it was shown a {\em quantitative} relationship between thrust and propulsion velocity in the range of medium, small and also zero Reynolds number, along with several applications.  The analysis developed in \cite{Gasp} has been subsequently completed and extended  by several authors in different directions, including attainability from rest \cite{ALS,AS2,HiMa}, motion of several bodies \cite{Staro}, controllability \cite{Hishida,Hishida1} and unsteady self-propulsion \cite{GaspU}, always in the case of rigid bodies.  
\par
At this point, the next, natural question to ask is: what can a {\em rigorous} mathematical analysis predict in the more difficult and intriguing case where propulsion is generated just by a {\em shape-changing} mechanism. Investigation of this problem  has received the attention of several mathematicians; see e.g. \cite{Court,Court1,Court2, Kapa,MS,Macha,Necasova,Nau,Raymond,Tuc}.  The model used for this study is rather general and a classic in fluid-structure interaction theory. Specifically, the body $\mathscr B$ is completely immersed in a Navier-Stokes liquid, $\mathscr L$, and its configuration (shape) is a prescribed  and sufficiently smooth  function of time. The only force acting on $\mathscr B$ is that exerted by the liquid on the surface of the body $\varSigma$. Therefore, the propulsive thrust comes only from the interaction of $\varSigma$ with $\mathscr L$. The unknowns are velocity and pressure fields of $\mathscr L$, and translational and angular velocity of $\mathscr B$; see Section \ref{Section: ProblemFormulation} for the precise formulation. It should be emphasized, though, that the main goal of the works cited above is limited to prove the well-posedness of the corresponding initial-boundary value problem (not an easy task!) in classes of functions of different regularity, 2D or 3D situations, covering also the case when the fluid is compressible. However, the important aspect completely omitted from all of them is to be able to ascertain whether the solutions to the self-propulsion problem given therein actually allow for a non-zero net motion of the body, let alone the relation between the change in shape of the body and its (net) translational velocity. In other words, through pure mathematical analysis, we cannot (yet)  guarantee that a ``fish'' actually moves! Some might argue that a rigorous analysis  on the full model, i.e.  without the approximations previously mentioned, could lead to overly general results that, ultimately, will furnish little or no quantitative information. However, this is not the case. Indeed, one of the objectives of this paper is to show not only that a rigorous analysis can produce accurate and remarkable quantitative results, but also that such results exhibit interesting features that are not  captured in the linearized or approximate approaches referenced above; see Section \ref{section:Application}.
\par
Precisely, suppose that the body $\mathscr B$ (a compact subset of $\real^3$) moves in a Navier-Stokes liquid that fills all the space outside $\mathscr B$,\footnote{This assumption is made in order to avoid  ``wall effects" that could blur the true cause of propulsion.}  periodically changing its shape over time. Our final goal is to establish conditions on the periodic deformation that ensure a net motion of $\mathscr B$, that is, that its center of mass, $G$, can travel any given distance in a finite time. Moreover, we furnish a {\em quantitative}  relation between such deformations and the net velocity of $G$. In order to accomplish these objectives, we assume that the motion of $\mathscr B$ is translatory. Indeed, allowing $\mathscr B$ to rotate as well will introduce a number of technical difficulties that could obscure the clarity of our results. Therefore, we prefer to defer the investigation of this more general case to future work. 
Before presenting our approach and the corresponding challenges, however, we would like to emphasize two important aspects of the problem. The first is, as expected, that not every periodic deformation of $\mathscr B$ can produce self-propulsion; see \remref{1.2}. The second is that the problem of self-propulsion is genuinely {\em nonlinear}; see \remref{3.1}.  
\par
The general strategy  that we will employ to reach our goals  develops according to the following steps. 
\begin{itemize}
\item[1.] We reformulate the original problem in a  reference configuration of $\mathscr B$, say $\Omega_0$, by means of a suitable diffeomorphism, $\bfchi=\bfchi(x,t)$, that is time-periodic of period $T$ ($T${\em-periodic}). 
In this Reformulated Problem (RP), the  domain occupied by the liquid becomes time independent,  given by $\varOmega:=\real^3\backslash \bar{\Omega_0}$. Likewise, the ``leading" datum becomes the transformed boundary velocity, $\partial_t\bfchi=:\bfu_*$; see Section \ref{Section:ProblemInReferenceConfiguration}. 
\item[2.] Since $\bfu_*$ is {$T$-periodic} (and, analogously, all coefficients in the RP equations are), it is natural to look for  $T$-periodic solutions. Existence of such solutions is accomplished provided the amplitude, $\delta$, of the oscillations of $\mathscr B$ is suitably restricted; see Section \ref{Section:NonlinearProblem}. Of course, the solutions depend on the given deformation, $\nabla\bfchi$, and $\bfu_*$.  
\item[3.] The final step is to find conditions on $\nabla\bfchi$ and $\bfu_*$ ensuring that $\mathscr B$ indeed self-propels, that is, $G$ performs a non-zero net motion.  This crucial step is {\em equivalent}  to show that the average over a period $[0,T]$ of the velocity of $G$ is nonzero, that is,  denoted by $\bfgamma$ such a velocity,  
$$
\bar{\bfgamma}:=\frac1T\int_0^{T}{\bfgamma}(s)\textrm{d}s\neq\0\,.
$$
\end{itemize} 
While completing the first step is fairly routine (see Section \ref{Section:ProblemInReferenceConfiguration}), completing the other two is anything but trivial, as we are going to explain. 
\par
For Step 2, discussed in Sections \ref{Section:LinearProblem} through \ref{nonl}, the leading idea is to construct a solution around the ($T$-periodic) one, $(\bfV,{\sf p},\bfzeta)$, to the linear problem obtained by neglecting all the nonlinear terms involving the velocity field; see \Eqref{3.1_0}. The role of such a solution, whose existence and uniqueness is proved in \propref{O_k},  is to ``lift" the boundary data $\bfu_*$. However, it also depends on the deformation $\nabla\bfchi$ and, in fact, it contributes to propulsion at the order of $\delta$, provided $\nabla\bfchi\not\equiv\0$, i.e., $\mathscr B$ is not rigid; see \propref{5.1} and \remref{5.1}. Precisely, we show
\be
\bar{\bfzeta}=\delta\,\mathbb F\cdot \bfpzc G_1+O(\delta^2)
\eeq{0._0}
where $\mathbb F$ is a symmetric positive-definite matrix depending only on $\Omega_0$, and $\bfpzc G_1$ is a vector that depends on $\Omega_0$, $\bfu_*$ and $\nabla\bfchi$; see \Eqref{5.5}--\Eqref{zizi}.   
 With this result in hand, we then look for a $T$-periodic solution $(\bfu,p,\bfgamma)$ to RP in the form 
\be
\bfu=\bsf u+\delta\,\bfV\,,\ \ p={\sf q}+\delta\,{\sf p}\,,\ \ \bfgamma=\bfxi+\delta\,\bfzeta\,, 
\eeq{0._01}
with $(\bsf u,{\sf q},\bfxi)$ satisfying a ``perturbed" problem, RPP, where the boundary data $\bfu_*$  are now replaced by ``driving forces" acting on both $\mathscr B$ and $\mathscr L$; see \Eqref{6.1}--\Eqref{6.3}. Our next task is to prove the existence  of a $T$-periodic solution $(\bsf u,{\sf q},\bfxi)$ to RPP.
Though the method employed  is simple in its formulation, it is  quite challenging in its implementation. 
The approach we use is the classical ``invading domain technique"  introduced in \cite{GaSi1} and successfully tested in several analogous circumstances \cite{GaKyH,6,7,GaSi}. In our case, it consists in redefining RPP in a sequence of bounded domains $\{\varOmega_{k}\}$ such that
$$
\bar{\varOmega_{k}}\subset \varOmega_{{k+1}}\,,\ \ \mbox{for all $k\in\nat$}\,;\ \ 
\cup_{k\in\nat}\varOmega_{k}=\varOmega\,;
$$
see \Eqref{6.10_lin}. 
In each of these domains one establishes the existence of a $T$-periodic solution $\bfs_k$, say, with corresponding estimates in term of the data, so that one can eventually pass to the limit $k\to\infty$ and prove  convergence of $\{\bfs_k\}$ to a solution of the original RPP in $\varOmega$.
For the success of this procedure it is essential that the constants entering the estimates are independent of $k$. The Galerkin method is appropriate for this purpose. However, adapting the method to the present situation is not an easy task, as RPP presents three challenging features. 
\par
First, the field, $\bsf u$ is not solenoidal. Rather, we have
\be
\Div \bsf u= -\mathbb C^\top:\nabla\bsf u
\eeq{0.1}
where $\mathbb C$ is a tensor field depending on $\nabla\bfchi$; see \Eqref{2.12}.  
However, the method requires to reformulate RPP in terms of a {\em solenoidal} vector field $\bsf v$ (say), i.e., with $\Div \bsf v=0$. By \Eqref{0.1}, this means to write $\bsf u=\bsf v+{\sf B}(\bsf u)$ where ${\sf B}$ is the Bogovskii operator such that 
$$
\Div {\sf B}(\bsf u)=\mathbb C^\top:\nabla\bsf u\,. 
$$ 
However, according to the classical Bogovskii theory \cite[Section III.3]{Gab}, this entails 
$$
\partial_t{\sf B}(\bsf u)\simeq \partial_t\nabla\bsf u\,,
$$
which would break the fundamental ``parabolic" character of the problem. To this end, in Section \ref{Section:DivOperator} we introduce and study the properties of a ``generalized" Bogovskii operator, $\bfB$, that allows us to overcome the above difficulty; see \propref{Boo}.  Roughly speaking, $\bfB$ satisfies $\partial_t\bfB(\bsf u)\simeq \partial_t\bsf u+\bfsigma(\partial_t\bsf u)$, where $\bfsigma(\partial_t\bsf u)$ involves only negative Sobolev trace norms at $\partial\Omega_0$; see \Eqref{Trunz}. 
\par  
Second, the linear momentum equation in RPP has an {\em extra} perturbative  term containing first- and second-order derivatives of $\bsf u$ that prevent us from obtaining a uniform bound in time of the kinetic energy norm of the solution, which, in the ``classical'' approach, is crucial to establish the existence of a fixed point for the Poincar\'e map ${\sf M}$; see \cite{GaKyH}. One way to overcome this problem would be to employ ``higher-order'' energy estimates. However, as is well known also for the ``simpler'' Navier-Stokes case, this necessarily entails a restriction on the size of the initial data, and this would be incompatible with finding a fixed point for ${\sf M}$; see \cite{GaKyH}. To resolve this deadlock, we resort to an appropriate linearization of the equations combined with careful use of the Schauder fixed point theorem; see Section \ref{linear}. \par Third, the linear momentum equations for both the liquid and the body in RPP include terms dependent on the pressure field ${\sf q}$, which may not coincide with the pressure field ${\sf q}_*$ (say) recovered {\em a posteriori} by the Galerkin method. This requires an appropriate perturbation argument to show ${\sf q}\equiv {\sf q}_*$.
\par
Once all of the above is accomplished, we are then able to prove the existence of a $T$-periodic  solution $(\bsf u,{\sf q},\bfxi)$ to RPP, at least for "small" $\delta$ (see \theoref{7.1}), which then leads to the completion of step 2.
\par 
The completion of the last Step 3 is obtained in Section \ref{asy}.  We proceed as follows. In view of \Eqref{0.1}, in order to provide conditions for self-propulsion, we  evaluate the average of the solutions found in \theoref{7.1} over the interval $[0,T]$, namely, $(\bar{\bsf u},\bar{\sf q},\bar{\bfxi})$. As expected, we show that, at order $\delta$, it is $\bar{\bfxi}=\0$, confirming that the phenomenon is nonlinear. We thus scale $(\bar{\bsf u},\bar{\sf q},\bar{\bfxi})$ by $\delta^2$ and then pass to the limit $\delta\to 0$ in the equations for the {\em scaled} fields; see \lemmref{11.3}. In such a way, we prove  that the limiting fields $(\bsf v_0,{\sf r}_0,\bfsigma_0)$ obey a suitable time-independent, non-homogeneous Stokes problem corresponding to the body $\mathscr B$ moving with  constant velocity $\bfsigma_0$, subject to a force $\bar{\bsf G}_0$,  while a body force $\bar{\bsf g}_0$ acts on the liquid; see \Eqref{11.2}. It is important to emphasize that both $\bar{\bsf g}_0$ and $\bar{\bsf G}_0$  depend only on $\bfu_*$, $\nabla\bfchi$, the physical parameters and the reference configuration; see \Eqref{11.1}. The other crucial point is that, at order of $\delta^2$, we show that $\bar{\bfxi}\neq\0$ {\em if and only if} $\bfsigma_0\neq\0$, precisely, $\bar{\bfxi}=\delta^2\bfsigma_0+o(\delta^2)$. By using a classical procedure (an adaptation of Lorentz reciprocity theorem) we are able to express $\bfsigma_0$ in terms of the known quantities $\bar{\bsf g}_0$ and $\bar{\bsf G}_0$ and show that
\be
\bar{\bfxi}=\delta^2\mathbb F\cdot \bfpzc G_0+o(\delta^2)\,,
\eeq{00.2}
with $\bfpzc G_0$ a vector whose components are functionals of $\bar{\bsf g}_0$ and $\bar{\bsf G}_0$ and, as such, dependent only on $\bfu_*$, $\nabla\bfchi$, the physical parameters and the reference configuration; see \Eqref{G0}. 
Thus, recalling \Eqref{0._0}, \Eqref{0._01} and setting $\bfpzc G:=\bfpzc G_0+\bfpzc G_1$, we conclude
\be
\bar{\bfgamma}=\delta^2\mathbb F\cdot \bfpzc G+o(\delta^2)\,,
\eeq{00.3}
from which we infer that,  at order $\delta^2$, self-propulsion can occur {\em if and only if}  the ``thrust" $\bfpzc G\neq\0$. It must be emphasized that, in principle, once  $\bfu_*$, $\nabla\bfchi$, $\Omega_0$ and the physical constants  are prescribed,
both $\mathbb F$ and $\bfpzc G$ can be explicitly evaluated. This is precisely done in Section \ref{section:Application}, where we apply our results to the classical benchmark problem where $\Omega_0$ is a ball of radius $a$; see  \cite{FJ1,FJ2,GiW,SW1,SW2}. We also choose as  deformation field the one generated by a $T$-periodic  dipole flow pattern combined with a rigid oscillation around the $\bfe_3$ axis; see \Eqref{00}. After a suitable non-dimensionalization, we then show that $\mathbb F$ is a multiple of the identity while $\bfpzc G$ is a function only of the Stokes number $h$; see \Eqref{STOK}, \Eqref{th} and \Eqref{th1}. This vector function can be easily computed with MATLAB. It turns out that both components $\mathpzc G_{1}
$ and $\mathpzc G_{3}$ are zero, so that $\bfpzc G=\mathpzc G(h)\bfe_2$. 
The graph of $\mathpzc G=\mathpzc G(h)$
 is reported in Figure \ref{graph_h_G}. Among other interesting features, it  shows the existence of an ``optimal" frequency that maximizes the speed of the body. It should be pointed out that most of these features are not present in any of analogous researches performed on similar but  {\em linearized} models; see, e.g., \cite{FJ1,FJ2}. This is because the functions $\bar{\bsf g_0}$ and $\bar{\bsf G_0}$ characterizing the thrust, involve a certain number of terms that may be missing in a linearization procedure imposed {\em directly} on the starting equations.  
\par
The plan of the paper is as follows. In Section \ref{Section: ProblemFormulation} we  formulate the self-propulsion  problem, along with some relevant remarks, while in Section \ref{Section:Notation_and_FunctionSpaces} we recall and introduce the main function spaces and their related  properties. In Section \ref{Section:ProblemInReferenceConfiguration}
we present a suitable ``deformation" map that allows us to reformulate the original problem in a reference configuration. Section \ref{Section:DivOperator} is entirely dedicated to the construction of a ``generalized" Bogovskii operator, whose properties are collected in \propref{Boo}. In the following Section \ref{Section:LinearProblem} we give the proof of existence of the $T$-periodic solution $(\bfV,{\sf p},\bfzeta)$ to the linear problem mentioned earlier on; see \propref{O_k}. The crucial feature of this solution is that it  decays sufficiently fast at large spatial distances (see \Eqref{AMPP_0}), a property that we show with the help of the method given in \cite{GaZ}. We also show that, in this approximation, $\mathscr B$ can propel, and find the corresponding velocity at the order of $\delta$ in \propref{5.1}.  In Section \ref{Section:NonlinearProblem}, we state our main result concerning the existence of  $T$-periodic solution to the nonlinear problem (see \theoref{7.1}) and outline the strategy for its proof. As mentioned above, the proof of the theorem is quite complex and  is split in several parts, developed in Sections \ref{linear}--\ref{nonl}. In order not to obscure the main ideas, we have preferred to postpone the demonstration of some technical results to  Appendices A--C. Employing the results established in \theoref{7.1}, in Section \ref{asy} we are able to give necessary and sufficient conditions for self-propulsion at the order of $\delta^2$, as well as identify the thrust, and provide in \theoref{11.1} a precise {\em quantitative}  relation between the thrust and $\bar{\bfgamma}$. In the Final Section \ref{section:Application} we give an application of \theoref{11.1} in the way described earlier on. Some related technical details are furnished in Appendix D.
\par 
To facilitate reading the paper, at the end we have inserted a table containing all the frequently used symbols, their description and the indication of the page on which they are defined.

\section{Formulation of the Problem}
\label{Section: ProblemFormulation}
Let $\mathscr B$ be a body moving in an otherwise quiescent Navier-Stokes liquid, $\mathscr L$, that fills the entire space outside $\mathscr B$. We will consider the case where $\mathscr B$ is  prevented from performing rigid rotations around its center of mass $G$, a condition that can be realized by applying a suitable torque on $\mathscr B$.   
\par
Denote by $\Omega=\Omega(t)$, $t\in\real$, a one-parameter family of bounded, sufficiently smooth domains of $\real^3$, each one representing the configuration of $\mathscr B$ at time $t$ with respect to a frame, $\calf$, with the origin at $G$ and axes parallel to those of an inertial frame. We assume that there are no external forces acting on the coupled system $\mathscr S:=\mathscr B\cup\mathscr L$ and that the only driving mechanism is  a {\em prescribed} change in shape of $\Omega$ with time, in a way that will be made precise later on.
\par
The self-propulsion problem that we would like to address can be thus qualitatively formulated as follows. Suppose that $\mathscr B$ changes its shape in a given time-periodic fashion, so that, for some $T>0$ and all $t\in \real$, $\Omega(t+T)=\Omega(t)$. Then, the goal is to find sufficient conditions on the map $t\mapsto \Omega(t)$ securing that $\mathscr B$ self-propels, namely, the center of mass $G$ covers any given finite distance in a finite time. 
\par
We begin to observe that the position of $G$ may vary with time for two reasons. The first, due to the prescribed deformation of the body,   and the second because of the interaction of $\mathscr B$ with $\mathscr L$. Of course,  only the latter  is responsible for propulsion, whereas the former is irrelevant, since it occurs even in absence of the liquid.  This situation can be accounted for by the following  simple  rescaling. 
Indicating by $\bfgamma=\bfgamma(t)$  the velocity of $G$ in $\mathcal F$, and by $M$ the total mass of $\mathscr B$, we have
\be
\bfgamma(t)=\frac1M\int_{\Omega(t)}\varrho(y,t)\bfw(y,t)\,{\rm d}y\,,
\eeq{lM}
where $\varrho=\varrho(y,t)$ is the material density of $\mathscr B$ and $\bfw(y,t)$ is the velocity at the point $y\in \Omega(t)$ at time $t$. Let $\bfv_*:=\bfv_*(y,t)$ be the velocity at $y\in\Omega(t)$ just due to the change in shape of $\Omega(t)$. Setting $\hat{\bfw}=\bfw-\bfv_*$, \Eqref{lM} entails   
$$
\bfgamma(t)=\frac1M\int_{\Omega(t)}\varrho(y,t)\hat{\bfw}(y,t)\,{\rm d}y+\bfgamma_{\rm d}\,,
$$
where
$$
\bfgamma_{\rm d}(t):=\frac1M\int_{\Omega(t)}\varrho(y,t)\bfv_*(y,t)\,{\rm d}y\,,
$$
stands for the component of the velocity of $G$ merely due to the rate of deformation of $\mathscr B$.
Thus, in what follows we shall tacitly understand that the velocity of $G$ and  that of points of $\mathscr B$ are rescaled by $\bfgamma_{\rm d}$. For simplicity, we continue to denote by  $\bfgamma$ and $\bfv_*$  the quantity   $\bfgamma-\bfgamma_{\rm d}$ and $\bfv_*+\bfgamma_{\rm d}$, respectively.

\begin{figure}[tbp] 
  \centering
  \includegraphics[bb=3 0 123 77,width=1.62in,height=1.04in,keepaspectratio]{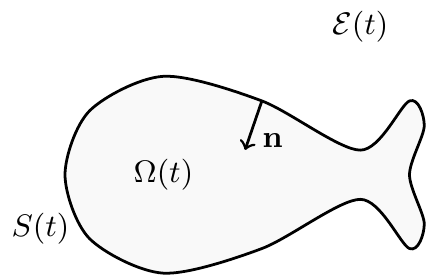}
  \caption{Domain}
  \label{fig:domain}
\end{figure}

\par
 In this wise, setting $\mathscr E(t):=\real^3\backslash \Omega(t)$ and $S(t):=\partial\Omega(t)\, (\equiv\partial\mathcal E(t))$ {(see Figure \ref{fig:domain})}, the equations governing the motion of the coupled system $\mathscr S$ in $\calf$   are given by
\be\ba{cc}\medskip\left.\ba{ll}\medskip
\partial_t\bfv(y,t)+\bfv(y,t)\cdot\nabla\bfv(y,t)=\Div\mathbb T(\bfv(y,t),{\sf p}(y,t))\\
\Div\bfv(y,t)=0\ea\right\}\ \ (y,t)\in \Gup{t\in\real}{}\,\mathcal E(t)\times\{t\}\,,\\ \medskip
\bfv(y,t)=\bfv_*(y,t)+\bfgamma\,, \ \ (y,t)\in \Gup {t\in\real}{}\, S(t)\times \{t\}\,;\ \ \Lim{|y|\to\infty}\bfv(y,t)=\0\,,\ t\in\real\,;\\ 
{\sf m}\dot{\bfgamma}(t)=-\Int{S(t)}{}\mathbb T(\bfv,{\sf p})\cdot\bsf{n}{\rm d}S\,,\ \ t\in\real\,.
\ea
\eeq{Eu}
In these equations, $\bfv$ and $\rho\,{\sf p}$ are velocity and pressure field  of $\mathscr L$,  $\rho$, $\nu$  its density and coefficient of kinematic viscosity, and ${\sf m}:=M/\rho$. (We assume no-slip boundary conditions.) Furthermore, the tensor field 
\be
\mathbb T(\bfz,\phi)=2\nu\mathbb D(\bfz)-\phi\,\mathbb I\,,\ \ \mathbb D(\bfz):=\half[\nabla\bfz+(\nabla\bfz)^\top]\,,
\eeq{Cauchy}
with $\mathbb I$ identity matrix, is the Cauchy stress tensor, and  ${\sf n}={\sf n}(y,t)$ is the unit outer (to $\cale$) normal at the point $y\in S(t)$ at time $t$.
\par
We shall now outline the strategy  we shall employ to give an answer to the self-propulsion problem. The first step is to reformulate  \Eqref{Eu} in a fixed reference configuration, say $\Omega(0)=:\Omega_0$, by means of a suitable time dependent diffeomorphism. In this Reformulated Problem (RP), the  domain occupied by the liquid becomes time independent, and is given by $\real^3\backslash \bar{\Omega_0}$. Likewise, the ``leading" datum becomes the transformed boundary velocity, say $\bfu_*$. Since $\bfu_*$ is time-periodic of period $T$ (and, likewise, all coefficients in the equations in RP), we look for the existence of time-periodic solutions to RP. We accomplish the latter,  provided the magnitude of $\bfu_*$ is appropriately restricted. The final, and more challenging, step is to find conditions ensuring that $\mathscr B$ indeed self-propels, that is, $G$ performs a non-zero net motion. Taking into account that the position vector of $G$, $\bfeta=\bfeta(t)$ counted from time $t=0$ (say), is given by 
\be
\bfeta(t)=\int_0^t\bfgamma(s){\rm d}s+\bfeta(0)\,,
\eeq{eta}
and that $\bfgamma$ is $T$-periodic, we deduce that the distance $d_T$ covered by $G$ in any interval of length $T$ is given by
$$
d_T=\left|\int_0^T\bfgamma(s){\rm d}s\right|
\,.
$$
As a result, self-propulsion is  equivalent to show that the time-average of $\gamma$ over a period is non-zero:
\be
\bar{\bfgamma}:=\frac1T\int_0^T\bfgamma(t){\rm d}t\neq \0\,.     
\eeq{gamma}
\par
The following remarks provide as many reasons as to why  the resolution of the problem \Eqref{Eu}, \Eqref{gamma} is far from being obvious.
\Br 
The main difficulty in investigating the validity of \Eqref{gamma} lies in the fact that the propulsive mechanism, constituted by the oscillation of the body, produces data with zero temporal average. In other words, the given boundary velocity distribution, $\bfv_*=\bfv_*(y,t)$, is such that
$$
\bar{\bfv_*}(x):=\frac1T\int_0^T\bfv_*(x,t){\rm d}t=\0\,;
$$
see \remref{av}. 
As actually demonstrated by {\sc G.I. Taylor} with his famous mechanical fish experiment \cite{Tay}, this circumstance implies that propulsion can only occur if the inertia of the liquid is taken into account. That is, in mathematical terms, problem \Eqref{Eu}, \Eqref{gamma} can only be solved if nonlinear effects are factored in. This circumstance will be made precise in \remref{3.1}.   
\ER{1.1}
\Br
As expected, not every $T$-periodic deformation can produce self-propulsion, as shown by the following example. Take $\Omega_0=B_1$, with $B_1$ unit ball centered at the origin, and let
$$
\Omega(t):=\{\bfy\in\real^3:\, \bfy=R(t)\bfx\,,\ \bfx\in B_1\}\,,
$$
where $R=R(t)$ is a smooth, positive and bounded function of time only such that
$$
R(t)=R(t+T)\,, \  R(0)=1\,,\ \mbox{for all}\ t\in\real\,. 
$$
The map $t\mapsto \Omega(t)$ represents then a $T$-periodic expansion and contraction of the unit ball. 
Correspondingly, we have
$$
S(t)=\{y\in\real^3:\ |y|=R(t)\}\,;\ \ \ \bfv_*(y,t)=\frac{\dot{R}(t)}{R(t)}\,\bfy\,,\ y\in S(t)\,.
$$
Consider the pair of fields
$$
\bfv(y,t)=a(t)\nabla\psi(y)\,,\ \ p(y,t)=-\dot{a}(t)\psi(y) -\half a^2(t)(\nabla\psi(y))^2\,;\ \ \psi(y):=\frac{1}{|\bfy|}\,,\ \ a(t):=-R^2(t)\dot{R}(t)\,.
$$
It is at once checked that  $(\bfv,p)$ satisfies \Eqref{Eu}$_{1,2,4}$. Moreover,
$$
\int_{S(t)}\left(\pde{v_i}{y_j}+\pde{v_j}{y_i}\right){\sf n}_j{\rm d}S=-\frac{a(t)}{R(t)}\int_{S^2}\left(2\delta_{ij}x_j-6x_i\right){\rm d}\sigma=4\frac{a(t)}{R(t)}\int_{S^2}x_i\,{\rm d}\sigma=0\,,
$$ 
and, likewise,
$$
\int_{S(t)}p(y,t)\bsf n\,{\rm d}S=-\left(\dot{a}(t)R(t)+\frac{a^2(t)}{R^2(t)}\right)\int_{S^2}\bfx\,{\rm d}\sigma=\0\,.
$$
The last two displayed equations then show that  \Eqref{Eu}$_{5}$ is satisfied with $\bfgamma\equiv\0$.
Finally, observing that  
$$
\bfv(y,t)=\big(\dot{R}(t)/R\big) \bfy=\bfv_*(y,t)\,,\ \  y\in S(t)\,,
$$ 
we deduce that also \Eqref{Eu}$_{4}$ is satisfied  with $\bfgamma\equiv\0$, which completes the proof of our statement.  
\ER{1.2}
\section{Notation and Relevant Function Spaces}
\label{Section:Notation_and_FunctionSpaces}
We begin by recalling some basic notation. By $\{\bfe_1,\bfe_2,\bfe_3\}$ we denote the canonical base in $\real^3$. Let $\mathbb A=\mathbb A_{ik}\bfe_i\otimes\bfe_k$ and $\mathbb B=\mathbb B_{lj}\bfe_l\otimes\bfe_j$ be second-order tensors, where   $\otimes$ denotes dyadic product. We set
$$
\mathbb A\cdot\mathbb B:=\mathbb A_{il}\mathbb B_{lj}\bfe_i\otimes\bfe_j\,,\ \ \mathbb A:\mathbb B:={\rm trace}\,\big(\mathbb A\cdot\mathbb B^\top\big)=\mathbb A_{il}\mathbb B_{il}\,,
$$
with ${}^\top$ denoting transpose.
Moreover, if $\bfa$ is a vector with components $a_i$, $i=1,2,3$, we set $\mathbb A\cdot\bfa:=\mathbb A_{ik}a_k\bfe_i$ and $\bfa\cdot\mathbb A:=\mathbb A^\top\cdot \bfa$. We also define
$$
\Div\mathbb A=\pde{\mathbb A_{ij}}{x_i}\bfe_j\,, \ \ \nabla\bfa=\pde{a_{k}}{x_i}\bfe_i\otimes\bfe_k\,,\ \ \nabla \phi=\pde{\phi}{x_i}\bfe_i\,,
$$
with $\phi$ scalar field.
\par
If $\varOmega$ is the complement of the closure of the bounded domain $\Omega_0$, we set
$$
\varOmega_R:=\varOmega\cap B_R\,,\ \ \varOmega^R:=\real^3\backslash\bar{B_R}
{ \,,\ \ \varOmega_{R_1,R_2}:=B_{R_2}\backslash\bar{B_{R_1}},}
$$
where $B_R$ stands for the open ball  with the origin in $\Omega_0$, and radius $R>R_*:={\rm diam}\,(\Omega_0)$, and the bar denotes closure. 
\par
As customary, for $A$ a domain of $\real^3$,  $L^q=L^q(A)$, $W^{m,2}=W^{m,2}(A)$, $q\in[1,\infty]$, $m\in\nat$, are Lebesgue and  Sobolev spaces with norm  $\|\cdot\|_{q,A}$, and $\|\cdot\|_{m,2,A}$. By $(\,\ ,\,\ )_A$ we mean the $L^2(A)$-scalar product.  Furthermore, $D^{m,q}=D^{m,q}(A)$ is the homogeneous Sobolev space with semi-norm $\sum_{|l|=m}\|D^lu\|_{q,A}$. 
In all the above notation we shall typically omit the subscript ``$A$", unless confusion arises. If $X$ is Banach space, we may  occasionally indicate its norm by $\|\cdot\|_X$.
\par 
Let ${\sf A}\subseteq \real^3$ be a domain with ${\sf A}\supset \bar{\Omega_0}$. We set
\be\ba{ll}\medskip
\calk=\mathcal K({\sf A}):=\big\{\bfphi\in C_0^\infty({\sf A}) : 
\exists\,\hat{\bfphi}\in\real^3 \mbox{ s.t. }\bfphi(x)=\hat{\bfphi} \mbox{ in a neighborhood of }\Omega_0\big\}\,,
\\ \medskip
\calc=\mathcal C({\sf A}):=\{\bfphi\in\calk({\sf A}):\ \Div\bfphi=0\ \mbox{in ${\sf A}$}\}\,,\\
\calc_0=\mathcal C_0({\sf A}):=\{\bfphi\in\calc({\sf A}): \hat{\bfphi}=\0\}\,.
\ea
\eeq{kcc}
In $\calk({\sf A})$ we introduce the scalar product
\be
\langle \bfphi,\bfpsi\rangle_{\sf A}:=
{\sf m}\,\hat{\bfphi}\cdot\hat{\bfpsi}+(\bfphi,\bfpsi)_{{\sf A}\cap{\varOmega}}\,,\ \ \bfphi,\bfpsi\in\calk\,,
\eeq{0.0}
and define 
\be\ba{ll}\medskip
\call^2(\real^3):= \,\big\{\mbox{completion of $\calk(\real^3)$ in the norm induced by \Eqref{0.0}}\big\}\,,\\ \medskip
\calh(\real^3):=\,\big\{\mbox{completion of $\calc(\real^3)$ in the norm induced by \Eqref{0.0}}\big\}\,,\\ 
\calg(\real^3):=\big\{\bfh\in \call^2(\real^3): \, \exists\, p\in D^{1,2}({\varOmega}) \mbox{ s.t. } \bfh=\nabla p \mbox{ in } {\varOmega},
\mbox{and}\ \bfh=-\varpi\int_{\partial\Omega}p\,\bfn   \ \mbox{in}\ \Omega_0\,\big\}\,.

\ea
\eeq{spazi}
It is shown in \cite[Theorem 3.1 and Lemma 3.2]{ALS} that
$$\ba{ll}\medskip
\call^2(\real^3)=\{\bfu\in L^2(\real^3): \ \bfu=\hat{\bfu}\ \mbox{in}\ \Omega_0, \ \mbox{for some $\hat{\bfu}\in \real^3$}\}\\
\calh(\real^3)=\{\bfu\in \call^2(\real^3): \ \Div\bfu=0\,\}\,,
\ea
$$
along with the following orthogonal decomposition \cite[Theorem 3.2]{ALS}
\be
\call^2(\real^3)=\calh(\real^3)\oplus\calg(\real^3)\,.
\eeq{Helm}

We next define the space
$$
\cald^{1,2}=\cald^{1,2}(\real^3):=\,\big\{\mbox{completion of $\calc(\real^3)$ in the norm $\|\mathbb D(\cdot)\|_2$}\big\}
$$
whose basic properties are collected in the next lemma; see \cite[Lemmas 9--11]{Gah}.
\Bl  ${\cald^{1,2}}$ is a separable Hilbert space when equipped with the scalar product
$$
(\mathbb D(\bfu_1),\mathbb D(\bfu_2))\,,\ \ \bfu_i\in {\cald^{1,2}}\,, \, \ i=1,2\,.
$$
Moreover, we have the characterization:
\be
{\cald^{1,2}}=\big\{\bfu\in L^6(\real^3)\cap D^{1,2}(\real^3)\,;\ \Div\bfu=0\,;\,
\bfu=\hat{\bfu} \ \mbox{in $\Omega_0$}\,,\ \mbox{for some $\hat{\bfu}\in\real^3$} \big\}\,.
\eeq{1.7777}
Also, for each $\bfu\in{\cald^{1,2}}$, it holds
\be
\|\nabla\bfu\|_2=\sqrt{2}\|\mathbb D(\bfu)\|_2\,,
\eeq{1.8}
and
\be
\|\bfu\|_6\le \kappa_0\,\|\mathbb D(\bfu)\|_2\,,
\eeq{1.9}
for some numerical constant $\kappa_0>0$. Finally, 
there is another positive constant $\kappa_1$ such that
\be
|\hat{\bfu}|\le \kappa_1\,\|\mathbb D(\bfu)\|_2\,.
\eeq{1.10}
\EL{1.1_1}
\par
\Br A relevant consequence of the above lemma is that 
$$
\|\mathbb D(\bfu)\|_2+\|\bfu\|_6+|\hat{\bfu}|\ \ \bfu\in\cald^{1,2}
$$
is 
an equivalent norm in $\cald^{1,2}$.
\ER{norm}
Along with the spaces $\call^2,\calh$, and  $\cald^{1,2}$ defined above, we introduce suitable ``local" versions of these spaces. Precisely, we set
$$\ba{ll}\medskip
\call^2(B_R):= \{\bfphi\in L^2(B_R): \, \bfphi|_{\Omega_0}=\hat{\bfphi}\,\ \mbox{for some $\hat{\bfphi}\in\real^3$}\}\,,\\ \medskip
\calh(B_R):=\{\bfphi\in \call^2(B_R):\, \Div\bfphi=0\,,\ \bfphi\cdot\bfn|_{\partial B_R}=0\}\,,\\ \medskip
\cald^{1,2}(B_R):= \{\bfphi\in W^{1,2}(B_R): \, \Div\bfphi=0\,,\ \bfphi|_{\Omega_0}=\hat{\bfphi}\,\ \mbox{for some $\hat{\bfphi}\in\real^3$}\,,\ \bfphi|_{\partial B_R}=\0 \}\,,\\
\cald^{1,2}_0({B_R}):= \{\bfphi\in \cald^{1,2}(B_R): \, \hat{\bfphi}=\0 \}\,.
\ea
$$
Then $\calh(B_R)$ and $\cald^{1,2}(B_R)$ are Hilbert spaces with scalar products
$$
\langle\bfphi_1,\bfphi_2\rangle_{B_R}\,,\ \ \bfphi_i\in\calh(B_R)\,;\ \ \ (\mathbb D(\bfpsi_1),\mathbb D(\bfpsi_2))_{B_R}\,,\ \bfpsi_i\in \cald^{1,2}(B_R)\,,\ \ i=1,2.
$$
Moreover, the following decomposition holds, analogous to \Eqref{Helm}  \cite[Theorem 3.1 and Lemma 3.2]{ALS}
\be
\call^2(B_R)=\calh(B_R)\oplus\calg(B_R)\,,
\eeq{Helm1}
where $\calg(B_R)$ is defined as in \Eqref{spazi}$_3$, by replacing {$\varOmega$} with ${\varOmega}_R$.
\par
Finally, the dual spaces of $\cald^{1,2}(B_R)$ and $\cald^{1,2}_0({B_R})$ will be denoted by $\cald^{-1,2}(B_R)$ and $\cald^{-1,2}_0({B_R})$, respectively.
\Br The space $\cald^{1,2}(B_R)$ can be viewed as a subspace of $W^{1,2}(\Omega)\cap\cald^{1,2}(\real^3)$, by extending its generic element to 0 in $\real^3\backslash B_R$. Therefore, all the properties mentioned in \lemmref{1.1_1} and \remref{norm} continue to hold for $\cald^{1,2}(B_R)$. 
\ER{2.1}
\par
We conclude this section by introducing certain spaces of time-periodic functions. With $A\subseteq\real^d$, a function $u:A\times \real\mapsto \real^3$ is 
{\em $T$-periodic}, $T>0$, if $u(\cdot,t+T)=u(\cdot\,t)$, for a.a. $t\in \real$,
 and we denote by
$$
{\bar u}:=\frac{1}{T}\int_{0}^{T}u(t){\rm d}t\,,
$$
its {\em average}. 
Let $B$ be a function space over $A$, endowed with seminorm $\|\cdot\|_B$, $r=[1,\infty]$, and $I$ an interval in $\real$. Then, $L^r(I;B)$ is the class of functions
$u:I\rightarrow B$ such that 
$$
\|u\|_{L^r(B)}\equiv\left\{\ba{ll}\smallskip\big( \Int{I}{}\|u(t)\|_B^r \big)^{\frac 1r}<\infty, \ \ \mbox{if 
$r\in [1,\infty)\,;$}\\   
\essup{t\in {I}}\,\|u(t)\|_B <\infty, \ \ \mbox{if $r=\infty.$}
\ea\right.
$$
Likewise, we put ($m\ge 0$)
$$
W^{m,r}(I;B)=\Big\{u\in L^{r}(I;B): {\partial_t^ku\in L^{r}(I;B), \, k=1,\ldots,m}\Big\}\,.
$$
We shall simply write $L^r(B)$ for $L^r(I;B)$, etc. unless otherwise stated, and when $B\equiv\real^d$, $d\ge1$, we set $L^r(I;B)=L^r(I)$, etc. We also define
$$
W_T^{m,r}(\real;B):=\{u\in W^{m,r}(I,B),\ \mbox{for all $I\subset \real$}:\ u\, \,\mbox{is $T$-periodic}\}\,, 
$$
endowed with norm
$$
\|\cdot\|_{W_T^{m,r}(\real;B)}:=\|\cdot\|_{W^{m,r}(0,T;B)}\,,
$$
and often use the abbreviation $W_T^{m,r}(B):=W_T^{m,r}(\real;B)$. In the case $B=\real^3$, for simplicity, the norm $ \|\cdot\|_{W_T^{m,r}(\real;\real^3)}$ will be denoted by $\|\cdot\|_{W^{m,r}(0,T)}$. 
\par 
\setcounter{equation}{0}
\section{Formulation of the Problem in the Reference Configuration}
\label{Section:ProblemInReferenceConfiguration}
We begin to introduce a time-dependent map  between the reference configuration $\Omega_0$ ($:=\Omega(0)$) and the configuration $\Omega(t)$, characterizing the ``deformation" of the body $\mathscr B$. To this end, 
for some $K>0$ and fixed $q>3$ let $\hat{\bsf{s}}: (x,t)\in\Omega_0\times\real \mapsto \real^3$ be such that 
\be
\hat{\bsf{s}}\in W_{T}^{3,\infty}(\real;W^{4,q}(\Omega_0)),\,q>3;\ \ \essup{t\in\real}\|\partial_t^k\hat{\bsf{s}}(t)\|_{4,q}\le K,\ k=0,1,2,3\,.
\eeq{Esse}
In view of the continuous embedding $ W^{2,q}(\Omega_0)\subset C^{1}(\bar{\Omega_0})$,  we have
\be
\sup_{t\in\real}\left(\max_{x\in\Omega_0}|\nabla\hat{\bsf{s}}(x,t)|\right)\le C_0\,K\,,
\eeq{emb}
with $C_0=C_0(\Omega_0,q)$. We then assume that the ``deformation" of $\mathscr B$ is described by the map
\be
\hat{\bfchi}:(\bfx,t)\in \Omega_0\times\real\mapsto \bfx+\delta\,\hat{\bsf{s}}(x,t)\equiv\bfy\in \Omega(t)\,,
\eeq{map}
where $\Omega(t)$ is the range of $\hat{\bfchi}$.
The following extension lemma holds.
\Bl
There is $\delta_1=\delta_1(\Omega_0)>0$  such that if $\delta\in(0,\delta_1)$, there exists   $\bfchi: (x,t)\in\real^3\times\real\mapsto\real^3$, satisfying the listed properties:
\begin{itemize}
\item[{\rm (i)}]
 $\bfchi(\bfx,t)=\hat{\bfchi}(x,t)$ for all $(\bfx,t)\in{\Omega_0\times\real}$\,;
\item[{\rm (ii)}] $\bfchi(\bfx,t)=\bfx$, for all $(\bfx,t)\in \left(\real^3\backslash {B_R}\right)\times \real$, and some $R>R_*$\,;
\item[{\rm (iii)}] $\bfchi(\cdot,t)$ is a $C^3$-diffeomorphism of $\real^3$ onto itself, for all $t\in\real$.  
\end{itemize}
Thus, in particular, $\bfchi(\cdot,t)$ is a $C^3$-diffeomorphism from $\real^3\backslash\Omega_0$ onto $\bfchi(\real^3\backslash\Omega_0,t)$, for all $t\in\real$.
\EL{1}
{\em Proof.} By well-known result, we can extend the function $\hat{\bsf s}$ to a function $\check{\bsf s}$ defined in $\real^3\times\real$ and such that
$$
\check{\bsf{s}}\in W_T^{3,\infty}(\real;W^{4,q}(\real^3))\,,\ \ \check{\bsf s}\in S_{c_0K}\,,\ \ c_0>0\,.
$$
Let $R$ be  such that $\bar{\Omega_0}\subset  B_R$,   and
let $\beta=\beta(x)$ be a smooth cut-off function that is equal to 1 for $x\in B_R$ and is zero for $x\in \real^3\backslash {B_{2R}}$. This function can be chosen in such a way that 
\be
|\nabla\beta(x)|\le c\,R^{-1}\,\ \  |D^2\beta(x)|\le c\,R^{-2}\,,
\eeq{beta} 
where $c$ is a constant independent of  $R$. Set $\bsf{s}:=\beta\,\check{\bsf{s}}$. Clearly, $\bsf{s}$ meets the following conditions
\be\ba{ll}\medskip
{\rm(a)}\,\ 
\bsf{s}\in W_T^{3,\infty}(\real;W^{4,q}(\real^3))\,;\\ \medskip
{\rm(b)}\ 
\, \bsf{s}(\bfx,t)=\left\{\ba{ll}\medskip \hat{\bsf{s}}(x,t)\ &   \mbox{if $(\bfx,t)\in\Omega_0\times\real$}\\ \medskip
\0\ &   \mbox{if $(\bfx,t)\in(\real^3\backslash {B_{2R}})\times\real$}
\ea\right.
\,;\\
{\rm(c)}\,\
 \|\bsf{s}\|_{W^{3,\infty}(W^{4,q}(\real^3))}\le c\,\|\hat{\bsf{s}}\|_{W^{3,\infty}(W^{4,q}(\Omega_0))}\,,   \, c=c(\Omega_0)>0\,. 
\ea
\eeq{Esse1}
Define, now, the following map:
\be
\bfchi:(\bfx,t)\in\real^3\times\real\mapsto \bfx+\delta\,\bsf{s}(x,t)\equiv \bfz\in \real^3\,.
\eeq{7_}
In view of (b), we at once deduce the validity of the stated properties (i) and (ii). Moreover,
since $q>3$, by   \Eqref{emb} and (c) and the Sobolev embedding  
\be
W^{4,q}(\real^3)\subset C^3(\real^3)\,,
\eeq{embr}
we get, on one hand, that $\bfchi$ is of class $C^3$ for each $t\in\real$, and, on the other hand,\be
\sup_{t\in\real}\left(\max_{x\in\real^3}|\nabla\bsf{s}(x,t)|\right)\le c_1K
\eeq{8}
where $c_1=c_1(\Omega_0,q)$. 
From this and well known results, it follows that $\bfchi$ is injective, provided we choose 
\be
0<\delta<1/(c_1K)\,.
\eeq{delta}  
It is easy to see that, in fact,  under the assumption \Eqref{delta} $\bfchi$ is also surjective, which, since $\bfchi(\cdot,t)=\hat{\bfchi}(\cdot,t)$ on $\bar\Omega_0$, would in turn imply that $\bfchi(\cdot,t)$ is a diffeomorphism of class $C^3$ of $\real^3\backslash \Omega_0$ onto $\bfchi(\real^3\backslash \Omega_0,t)$ for all $t\in\real$. To show surjectivity, for any fixed $\bfz\in\real^3$ and $t\in\real$ consider the map 
$$
\bsf{P}:\bfx\in\real^3\mapsto \big(\bfz-\delta\,\bsf{s}(\bfx,t)\big)\in\real^3\,.
$$ 
By (a) and \Eqref{embr}, we deduce  $\bsf{P}\in C^1(\real^3)$ and moreover, by \Eqref{8} and  \Eqref{delta},  
$$\Sup{{\footnotesize \bfx}\in\real^3}|\nabla\bsf{P}(\bfx)|<1\,.$$
Therefore, by classical results (see, e.g., \cite[Chapter XVII, Theorem 1]{KaAk}) it follows that $\bsf{P}$ has a fixed point, namely,  $\bfchi(\cdot,t)$ is surjective for all $t\in\real$. The proof of the lemma is  completed.\par\hfill$\square$\par
Let
\be
\mathbb{A}:=\nabla\bfchi^{-1}\,,\ \mathbb B:=\mathbb A-\mathbb I\,,\ \ \mathbb C:=J\mathbb A-\mathbb I\,,\ \ J:={\rm det}\,\nabla\bfchi\,.
\eeq{2.9}
The following result holds.
\Bl Suppose \Eqref{Esse} holds. Then, there exists $\delta_2>0$ such that, for all $\delta\in (0,\delta_2)$,  $\mathbb{A}$, $\mathbb B$ and $\mathbb C$ are $T$-periodic functions satisfying
\be
\mathbb{A}\in W^{3,\infty}_T(W^{2,\infty}(\varOmega))\,;\ \mathbb B,\,\mathbb C\in W^{3,\infty}_{T}(W^{3,q}(\varOmega))\,.
\eeq{Si0}
Moreover,  $\mathbb B$ and $\mathbb C$ have bounded support, $B_0$,  and,  
for any $k=0,1,2,3$ and $|\alpha|=0,1,2$, satisfy the following conditions 
\be\ba{ll}\medskip
\partial_t^kD^\alpha \mathbb B(x,t)=-\delta\,\partial_t^kD^\alpha\nabla \bsf s(x,t)+\mathbb Z_1(\delta,\bsf s)\,,\\
\partial_t^kD^\alpha \mathbb C(x,t)=-\delta\,\partial_t^kD^\alpha\left(\nabla \bsf s(x,t)-\Div\bsf s\,\mathbb I\right)+\mathbb Z_2(\delta,\bsf s)\,,
\ea
\eeq{Si1}
where 
$$
|\mathbb Z_i(\delta,\bsf s)|\le c\,\delta^2\,,\ \ i=1,2.
$$
Thus, in particular,
\be
\essup{(x,t)\in\varOmega\times\real}\left(|\partial_t^kD^\alpha \mathbb B(x,t)|+|\partial_t^kD^\alpha \mathbb C(x,t)|\right)\le c\,\delta\,.
\eeq{Si2}
\EL{Def}
{\em Proof.} In view of the definition of $\bfchi$ given in \Eqref{Esse1}--{\Eqref{7_}}, $T$-periodicity follows at once.    
From {\Eqref{7_}}, \Eqref{8} and \Eqref{delta}, we deduce  that  $\mathbb A$ can be expressed as the Neumann series
\be
\mathbb A=\sum_{n=0}^\infty(-1)^n(\delta\,\nabla\bsf s)^n\,,
\eeq{Si3}
which, by \Eqref{Esse1} and the embedding $W^{3,q}(\varOmega)\subset W^{2,\infty}(\varOmega)$,  furnishes property \Eqref{Si0} for $\mathbb A$, provided $\delta>0$ is chosen below a certain constant, depending on $\varOmega$ and $K$. From \Eqref{Si3} we also have 
\be
\mathbb B=-\delta\,\nabla\bsf s+\sum_{n=2}^\infty(-1)^{n}(\delta\,\nabla\bsf s)^n\,.
\eeq{JuSa}
Thus, by \Eqref{Esse1}, $\supp(\mathbb B)\subset \supp(\nabla\bsf s)\subset B_R$. Moreover, properties \Eqref{Si0} and \Eqref{Si1} for $\mathbb B$  follow from \Eqref{JuSa} by taking $\delta$ sufficiently small and using again \Eqref{Esse1} along with the  embedding $W^{4,q}(\varOmega)\subset C^3(\bar{\varOmega})$, $q>3$.
Furthermore, by a straightforward calculation, we get
$$
J=1+\delta\,\Div\bsf s+\delta^2\,P(\nabla\bsf s)+\delta^3{\rm det}\,\nabla\bsf s:=1+\delta\,\Div\bsf s+f(\delta\nabla \bsf s)\,,
$$
where $P$ is a second order homogenous polynomial of components of $\nabla\bsf s$.
As a result,
$$
\mathbb C=\mathbb B+\delta\,\Div\bsf s\,\mathbb I + \delta\,\Div\bsf s\,\mathbb B+ f(\delta\nabla\bsf s)\mathbb A\,.
$$
The stated properties for $\mathbb C$ then follow from the latter in conjunction with \Eqref{Si1}$_1$ and \Eqref{Si3}.
\par\hfill$\square$\par
We shall now reformulate problem \Eqref{Eu} in the reference configuration $\Omega_0$. To this end, let
\be
\bfpsi(x,t):=\bfchi(x,t)+\bfeta(t)
\eeq{2.7}
and define the following fields
\be\left.\ba{ll}\medskip
\bfu(x,t):=\bfv(\bfpsi(x,t),t)\,,\ \  p(x,t)={\sf p}(\bfpsi(x,t),t)\,,\\ 
\bfu_*(x,t):=\delta^{-1}\bfv_*(\bfchi(x,t),t)=  \partial_t{\bsf{s}}(x,t)\,,
\ea \right\}\ \ (x,t)\in{\varOmega}\times\real\,.
\eeq{2.8}
We have
$$
\pde{{p}}{x_i}=\pde{\chi_l}{x_i}\pde{{\sf p}}{y_l}\,,
$$
which gives
$$
\mathbb A_{ji}\pde{{p}}{x_i}=\delta_{jl}\pde{\sf p}{y_l}=\pde{\sf p}{y_j}\,,
$$
that is
\be
\nabla_y{\sf p}=\mathbb A\cdot\nabla_x p\,,
\eeq{000}
or also
\be
\nabla_y{\sf p}=\nabla_x p+\mathbb B\cdot\nabla_x p\,,
\eeq{0}
Likewise, from 
$$
\pde{u_k}{x_i}=\pde{\chi_l}{x_i}\pde{v_k}{y_l}
$$
we get
\be
\mathbb A_{ji}\pde{u_k}{x_i}=\mathbb A_{ji}\pde{\chi_l}{x_i}\pde{v_k}{y_l}=\delta_{jl}\pde{v_k}{y_l}=\pde{v_k}{y_j}\,,
\eeq{1}
namely,
\be
\nabla_y\bfv=\mathbb{A}\cdot\nabla_x\bfu\,,
\eeq{2}
or also
\be
\nabla_y\bfv=\nabla_x\bfu+\mathbb{B}\cdot\nabla_x\bfu\,.
\eeq{2_1}
Taking $k=j$ in \Eqref{1}, and multiplying both side of the resulting equation by $J$, we get
\be
J\mathbb A_{ji}\pde{u_j}{x_i}=J\pde{v_j}{y_j}=0\,.
\eeq{3}
Thanks to  Piola identity
\be
\Div(J\mathbb A^\top)=\Div(\mathbb C^\top)=\0\,,
\eeq{PId}
it follows that 
\be
(J\mathbb A_{ji}-\delta_{ji})\pde{u_j}{x_i}=\pde{}{x_i}[(J\mathbb A_{ji}-\delta_{ji})u_j]-u_j\pde{}{x_i}(J\mathbb A_{ji})=\pde{}{x_i}[(J\mathbb A_{ji}-\delta_{ji})u_j]\,.
\eeq{4}
Thus, collecting \Eqref{3} and \Eqref{4} we get
\be
\Div\bfu=-\Div(\mathbb C^\top\cdot\bfu)\,.
\eeq{5}
Taking into account \Eqref{2} and the definition of $\mathbb B$, we show
$$
\mathbb T(\bfv,{\sf p})=\nu (\nabla_y\bfv+(\nabla_y\bfv)^\top)-{\sf p}\,\mathbb I=\nu \big(\nabla_x\bfu+(\nabla_x\bfu)^\top\big)-{p}\,\mathbb I+\nu\big(\mathbb B\cdot\nabla_x\bfu+(\mathbb B\cdot\nabla_x\bfu)^\top\big)\,.
$$
Hence, from the latter and the well-known identity 
$$
\bsf{n}\,{\rm d}S=J\,\mathbb A\cdot\bfn \,{\rm d}\varSigma\,,
$$
we deduce, in particular,
\be
\Int S{}\mathbb T(\bfv,{\sf p})\cdot\bsf n\,{\rm d}S=\Int \varSigma{}
\mathbb T(\bfu,{p})\cdot\bfn\,{\rm d}\varSigma +\Int \varSigma{}
\mathbb T(\bfu,{p})\cdot\mathbb C\cdot\bfn\,{\rm d}\varSigma+\nu\Int \varSigma{}
\big(\mathbb B\cdot\nabla_x\bfu+(\mathbb B\cdot\nabla_x\bfu)^\top\big)\,J\,\mathbb A\cdot\bfn\,{\rm d}\varSigma\,.
\eeq{6_1}
Moreover, again from \Eqref{000}, \Eqref{1} and the definition of $\mathbb B$, we have (with $\partial_s\equiv\partial/\partial x_s$, $s=1,2,3$)
$$\ba{rl}\medskip
[\Div_y\mathbb T(\bfv,{\sf p})]_k&\!\!\!\!=
\pde{}{y_i}\left[\nu\left(\pde{v_k}{y_i}+\pde {v_i}{y_k}\right)-{\sf p}\delta_{ik}\right]=\mathbb A_{il}\partial_l\left[\nu\left(\mathbb A_{im}\partial_m{u_k}+\mathbb A_{km}\partial_mu_i\right)-{p}\delta_{ik}\right]\\ 
&\!\!\!\!=[\Div_x\mathbb T(\bfu,{p})]_k+\mathbb B_{il}\partial_l[\nu(\partial_i u_k+\partial_k u_i)-{p}\delta_{ik}]+\nu\mathbb A_{il}\partial_l(\mathbb B_{im}\partial_m u_k+\mathbb B_{km}\partial_m u_i)
\ea
$$
namely, in intrinsic form, 
\be
\Div_y\mathbb T(\bfv,{\sf p})=\Div\mathbb T(\bfu,{p})+\mathbb B^\top:\nabla\,\mathbb T(\bfu,p)+\nu\,\mathbb A^\top:\nabla\big(\mathbb B\cdot\nabla\bfu+(\mathbb B\cdot\nabla\bfu)^\top\big)\,,
\eeq{8_0}
where all the differential operators on the right-hand side act on the $x$-variables. 
Finally, also using \Eqref{2}, we show 
\be
\partial_t{\bfu}=\partial_t{\bfpsi}\cdot\nabla_y\bfv+\partial_t{\bfv}=(\delta\bfu_*+\bfgamma)\cdot\mathbb{A}\cdot\nabla_x\bfu+\partial_t{\bfv}\,.
\eeq{2.11_1}
Therefore, 
we conclude that the original problem is reformulated,  in the reference configuration, 
 as follows:
\be\ba{cc}\medskip\left.\ba{ll}\medskip
\partial_t{\bfu}+(\bfu-\bfgamma)\cdot\nabla\bfu=\Div\mathbb T(\bfu,p)+\bff_{\ell}(\bsf{s},\bfu,p)+\bff_{n\ell}(\bsf{s},\bfu,\bfgamma)\\
\Div\bfu=-\Div(\mathbb C^\top\cdot\bfu)\ea\right\}\ \ \mbox{in $\varOmega\times\real$}\,,\\ \medskip
\bfu=\delta\bfu_*+\bfgamma\,, \ \ \mbox{at $\varSigma\times \real$}\,;\ \ \Lim{|x|\to\infty}\bfu(x,t)=\0\,,\ t\in\real\,;\\ 
{\sf m}\dot{\bfgamma}(t)=-\Int{\varSigma}{}\mathbb T(\bfu,{p})\cdot\bfn{\rm d}\varSigma+\bfF(\bsf{s},\bfu,p)\,,\ \ t\in\real\,,
\ea
\eeq{2.12}
where $\bfn$ is the unit outer (to $\varOmega$) normal at $\varSigma$, and
\be\ba{rl}\medskip
\bff_{\ell}({\bsf{s},\bfu,p})&\!\!\!\!:=\delta\bfu_*\!\cdot\!\mathbb{A}\cdot\nabla\bfu
+\mathbb B^\top:\nabla\,\mathbb T(\bfu,p)+\nu\,\mathbb A^\top:\nabla\big(\mathbb B\cdot\nabla\bfu+(\mathbb B\cdot\nabla\bfu)^\top\big)\\ \medskip
\bff_{n\ell}({\bsf{s},\bfu,\bfgamma})&\!\!\!\!:=-(\bfu-\bfgamma)\cdot\mathbb{B}\!\cdot\!\nabla\bfu\\
\bfF(\bsf{s},\bfu,p)&\!\!\!\!:=-\Int \varSigma{}
\mathbb T(\bfu,{p})\cdot\mathbb C\cdot\bfn\,{\rm d}\varSigma-\nu\Int \varSigma{}J\,
\big(\mathbb B\cdot\nabla\bfu+(\mathbb B\cdot\nabla\bfu)^\top\big)\cdot\mathbb A\cdot\bfn\,{\rm d}\varSigma\,.
\ea
\eeq{2.13}
\Br In the reformulated problem \Eqref{2.12}--\Eqref{2.13}, the (only) driving mechanism becomes the boundary velocity distribution $\delta\,\bfu_*$. It is worth emphasizing that, since
$$
\bfu_*(x,t)=\partial_t{\bsf s}(x,t)\,,
$$
and $\bsf s$ is  $T$-periodic, it follows that
$$
\bar{\bfu_*}=\0\,.
$$ 
As explained later on in \remref{3.1}, this circumstance implies that, in our framework, self-propulsion is a strictly {\em nonlinear}  problem.   
\ER{av}
\Br Since $\bfu$ is {\em not} solenoidal, 
we have
$$
\Div\mathbb D(\bfu)\neq \nu\Delta\bfu\,.
$$
Instead, 
from \Eqref{Cauchy} and  \Eqref{2.12}$_2$ it follows that
$$
\Div\mathbb T(\bfu,p)=\nu\left(\Delta\bfu+\nabla\Div\bfu\right)-\nabla p=\nu\left(\Delta\bfu-\nabla\Div(\mathbb C^\top\cdot\bfu)\right)-\nabla p\,.
$$
\ER{div}
\Br 
In the context of problem  \Eqref{2.12}--\Eqref{2.13}, the parameter $\delta$ $(\ge 0)$   serves as the magnitude of the driving mechanism. It can be viewed as the ratio of the largest displacement of $\mathscr B$ to its diameter $L$:
$$
\delta:=\Frac{\max_{(x,t)\in \Omega_0\times[0,T]}|\bsf s(x,t)|}L\,.   
$$
It is clear (and immediately checked) that for $\delta=0$ a corresponding solution is the identically vanishing one, which trivially implies no self-propulsion. However, self-propulsion does not occur also at $O(\delta)$. In fact, suppose $(\bfu,p,\bfgamma)$ is a $T$-periodic solution to  \Eqref{2.12}--\Eqref{2.13}. If we write $\bfu=\delta\,\bfu_0$, $p=\delta\,p_0$, $\bfgamma=\delta\,\bfgamma_0$, and disregard terms in $\delta^2$ and higher, also with the help of \remref{div} and \lemmref{fF} we then deduce that $(\bfu_0,p_0,\bfgamma_0)$ obeys the following problem
$$
\ba{cc}\medskip\left.\ba{ll}\medskip
\partial_t{\bfu}_0=\nu\Delta\bfu_0-\nabla{p}_0\\
\Div\bfu_0=0\ea\right\}\ \ \mbox{in $\varOmega\times\real$}\,,\\ \medskip
\bfu_0=\bfu_*+\bfgamma_0\,, \ \ \mbox{at $\varSigma\times \real$}\,;\ \ \Lim{|x|\to\infty}\bfu_0(x,t)=\0\,,\ t\in\real\,;\\ 
{\sf m}\dot{\bfgamma}_0(t)=-\Int{\varSigma}{}\mathbb T(\bfu_0,{p}_0)\cdot\bfn{\rm d}\varSigma\,,\ \ t\in\real\,,
\ea
$$   
Henceforth, the averaged fields $(\bar{\bfu_0},\bar{p_0},\bar{\bfgamma_0})$ and  $\bar{\bfu_*}$ obey the boundary-value problem:
$$
\ba{cc}\medskip\left.\ba{ll}\medskip
\nu\Delta\bar{\bfu_0}-\nabla\bar{{p}_0}=\0\\
\Div\bar{\bfu_0}=0\ea\right\}\ \ \mbox{in $\varOmega$}\,,\\ \medskip
\bar{\bfu_0}=\bar{\bfu_*}+\bar{\bfgamma_0}\,, \ \ \mbox{at $\varSigma$}\,;\ \ \Lim{|x|\to\infty}\bar{\bfu_0}(x)=\0\,;\\ 
\Int{\varSigma}{}\mathbb T(\bar{\bfu_0},\bar{{p}_0})\cdot\bfn{\rm d}\varSigma=\0\,.
\ea
$$
In view of classical results on the Stokes problem \cite[Section V.7]{Gab}, it follows that, in a very large class of solutions, we may have $(\bar{\bfu_0},\nabla\bar{p_0},\bar{\bfgamma_0})\neq(\0,\0,\0)$ {\em if and only if} $\bar{\bfu_*}\not\equiv\0$. However, as pointed out in \remref{av}, in our case, it is $\bar{\bfu_*}\equiv\0$, intimating that  self-propulsion can only occur at the order of $\delta^2$ or higher, that is, only when nonlinear effects are taken into account. The main objective of this paper is to provide a characterization for the  this circumstance to take place. 
\ER{3.1}
\medskip\par
We conclude this section with the following simple but useful result, which is a direct consequence of \lemmref{Def}, of the classical trace inequalities, and, possibly, of the addition to $p$ of a suitable function of time only.
\Bl Let $\delta_2, B_0$ be as in \lemmref{Def}. Then, for all $\delta\in (0,\delta_2)$ and $q\in [1,\infty]$ the following estimates hold (with $(\cdot)_t:=\partial_t(\cdot)$)
$$\ba{ll}\medskip
\|\bff_{\ell}\|_{q}\le c\,\delta\,\left(\|\nabla\bfu\|_{1,q,B_0}+\|\nabla p\|_{q,B_0}\right)\\ 
\|(\bff_{\ell})_t\|_{q}\le c\,\delta\,\left(\|\nabla\bfu\|_{1,q,B_0}+\|\nabla p\|_{q,B_0}+\|\nabla\bfu_t\|_{1,q,B_0}+\|\nabla p_t\|_{q,B_0}\right)\,,
\ea$$ 
and
$$
|\bfF|\le c\,\delta\,\left(\|\nabla\bfu\|_{1,2,B_0}+\|\nabla p\|_{2,B_0}\right)\,;\ \ |\bfF_t|\le c\,\delta\,\left(\|\nabla\bfu\|_{1,2,B_0}+\|\nabla\bfu_t\|_{1,2,B_0}+\|\nabla p_t\|_{2,B_0}\right)\,.
$$
\EL{fF}

\setcounter{equation}{0}
\section{Suitable Inversion of the ``{\bf div}" Operator}
\label{Section:DivOperator}
A notable feature of the reformulated problem \Eqref{2.12}--\Eqref{2.13} is that, unlike the original one, the velocity field ${\bfu}$ is no longer solenoidal. To address this point we are therefore naturally led to the study of some relevant properties of the ``$\Div$" operator, with particular regard to its inversion.
Seemingly, these properties are not directly amenable to the classical theory. 
\Bl Let $k=0,1,2,3$, and $m=1,2$. Suppose $\bfg$ satisfy the following conditions:
\begin{itemize}
\item[{\rm (a)}] $\bfg\in W^{k,2}_{{T}}(W^{m,2}(\varOmega))$\,; 
\item[{\rm (b)}] There is $\rho_0>R_*$  such that $\bfg(x,t)=\0$,  for all $(x,t)\in \varOmega^{\rho_0}\times \real$;
\item[{\rm (c)}]\, $\Int{\varSigma}{}\bfg(x,t)\cdot\bfn\,{\rm d}\varSigma=0$\,,\ for a.a. $t\in\real$\,.   
\end{itemize} 
Then, we can find a field $\bfw=\bfw(\bfg)$  for which the following properties hold, for any fixed $\rho>\rho_0$.
\begin{itemize}
\item[{\rm (i)}] $\Div\bfw=\Div\bfg $\, in $\varOmega_\rho$\,;
\item[{\rm (ii)}] $\bfw(x,t)=\0$, for all $(x,t)\in\varSigma\times \real$; 
\item[{\rm (iii)}] $\bfw(x,t)=\0$,  for all $(x,t)\in\varOmega^{\rho}\times \real$; 
\item[{\rm (iv)}] $\bfw\in W^{k,2}_{{T}}(W^{m,2}(\varOmega))$\,;
\item[{\rm (v)}] For a.a. $t\in\real$,  the following inequalities hold
$$\ba{c}\medskip\|\partial^k_t\bfw(t)\|_{m,2}\le c\,\|\partial^k_t\bfg(t)\|_{m,2}\,,\\ \|\partial^k_t\bfw(t)\|_{2}\le c\,\big(\|\partial^k_t\bfg(t)\|_{2}+\|\partial^k_t\bfg(t)\cdot\bfn\|_{-\frac12,2(\Sigma)}\big)\,,
\ea$$
where the positive constant $c$ depends only on $\rho$.
\end{itemize} 
\EL{g}
{\em Proof.} For each fixed $t\in\real$, consider the following Stokes problem
\be\ba{cc}\medskip\left.\ba{ll}\medskip
\Delta\bfv-\nabla p=0\\
\Div\bfv=\Div\bfg\ea\right\}\ \ \mbox{in ${\varOmega}_\rho$}\\
\bfv=\0\ \ \mbox{at $\partial{\varOmega}_\rho$}\,.
\ea
\eeq{3.1}
Since, by (b) and (c),
$$
\Int{{\varOmega}_\rho}{}\Div\bfv=\Int{{\varOmega}_\rho}{}\Div\bfg=\Int{\varSigma}{}\bfg\cdot\bfn\,{\rm d\varSigma}=0=\int_\varSigma\bfv\cdot\bfn\,{\rm d}\varSigma\,,
$$
from classical results we show the existence of a unique solution $\bfv(t)\in W^{m,2}({\varOmega}_\rho)$ such that
\be
\|\bfv(t)\|_{m,2}\le c\,\|\Div\bfg(t)\|_{m-1,2}\le c\,\|\bfg(t)\|_{m,2}\,.  
\eeq{sfal}
In view of the assumptions made on $\bfg$, we can also easily show that $(\partial_t^k\bfv,\partial_t^k p)$  solves \Eqref{3.1} with $\bfg$ replaced by $\partial^k_t\bfg$, with estimate analogous to \Eqref{sfal}. We thus conclude
\be 
\|\partial_t^k\bfv\|_{m,2}\le c\,\|\partial^k_t\bfg\|_{m,2}\,.
\eeq{3.2}
Next, for a given $\bff\in C_0^\infty({\varOmega}_\rho)$, denote by $(\bfphi,\tau)\in W^{2,2}({\varOmega}_\rho)\times W^{1,2}({\varOmega}_\rho)$ the solution to the problem
\be\ba{cc}\medskip\left.\ba{ll}\medskip
\Delta\bfphi-\nabla \tau=\bff\\
\Div\bfphi=0\ea\right\}\ \ \mbox{in ${\varOmega}_\rho$}\\
\bfphi=\0\ \ \mbox{at $\partial{\varOmega}_\rho$}\,, \ \ \Int{{\varOmega}_\rho}{}\tau=0\,.
\ea
\eeq{3.3}
Again by classical results, such a (unique) solution exists and satisfies
\be 
\|\bfphi\|_{2,2}+\|\tau\|_{1,2}\le c\,\|\bff\|_2\,.
\eeq{3.4}
By testing \Eqref{3.1}$_1$ with $\bfphi$, integrating by parts, and  using the boundary conditions along with the divergence equations,  we get
$$
0=(\Delta\bfv-\nabla p,\bfphi)=(\bfv,\bff)+(\nabla\tau,\bfv)=(\bfv,\bff)-(\tau,\Div\bfv)=(\bfv,\bff)+(\nabla\tau,\bfg)-\int_\varSigma\tau\,\bfg\cdot\bfn\,{\rm d}\varSigma\,.
$$
Thus, from Schwarz inequality, and   trace theorem we deduce
$$
|(\bfv,\bff)|\le c\big(\|\bfg\|_2+\|\bfg\cdot\bfn\|_{-\frac12,2 (\varSigma)}\big)\|\tau\|_{1,2}
$$
which, in turn, by  \Eqref{3.4}  and the arbitrariness of $\bff\in C_0^\infty({\varOmega}_\rho)$, allows us to conclude
\be
\|\bfv\|_2\le c\,\big(\|\bfg\|_2+\|\bfg\cdot\bfn\|_{-\frac12,2 (\varSigma)}\big)\,.
\eeq{3.5}
In a similar way, we prove
\be
\|\partial^k_t\bfv\|_2\le c\,\big(\|\partial^k_t\bfg\|_2+\|\partial^k_t\bfg\cdot\bfn\|_{-\frac12,2 (\varSigma)}\big)\,.
\eeq{3.6}
Pick $r\in(\rho_0,\rho)$ and let $\psi=\psi(|x|)$ be a smooth non-increasing ``cut-off" function such that $\psi(|x|)=1$, for $|x|\le r$, and $\psi(|x|)=0$, for $|x|\ge \rho$. Notice that, in view of assumption (b) on $\bfg$, we have
\be\ba{ll}\medskip
\psi\,\Div\bfg=\Div\bfg \,,\  \mbox{in ${\varOmega}_\rho$}\,,\\ 
\Div\bfv=\0\,,\  \mbox{in $\varOmega_{\rho_0,\rho}$}\,.
\ea
\eeq{3.7}
Set
\be
\bfw:=\psi\,\bfv-\bsf{v}\,\ \ \mbox{in ${\varOmega}_\rho$}
\eeq{3.8}
where $\bsf{v}$ is a solution to the problem
\be\ba{cc}\ms\left\{\ba{ll}\medskip
\Div\bsf{v}=\nabla\psi\cdot\bfv\,,\ \ \mbox{in $\varOmega_{\rho_0,\rho}$}\,,\\
\bsf{v}\in W^{k,2}(W_0^{m,2}(\varOmega_{\rho_0,\rho}))\,,\ea\right.
\\
\|\partial_t^k\bsf{v}\|_{m,2}\le c\, \|\partial_t^k\bfv\|_{m-1,2}\,.
\ea
\eeq{3.9}
Since  the properties of $\psi$ and  \Eqref{3.2}  entail $\nabla\psi\cdot\partial_t^k\bfv\in W_0^{1,2}({\varOmega}_\rho)$, the existence of such a $\bsf{v}$ follows from \cite[Exercise III.3.7]{Gab} provided
\be
\Phi:=\int_{\varOmega_{\rho_0,\rho}}\nabla\psi\cdot\bfv=0\,.
\eeq{3.10}
Now, with the help of \Eqref{3.7}$_2$, \Eqref{3.1}$_{1,2}$ and (b) we show
$$
\Phi=\int_{\varOmega_{\rho_0,\rho}}\Div(\psi\bfv)=\int_{\partial{\varOmega}_{\rho_0}}\bfv\cdot\bfn=\int_{{\varOmega}_{\rho_0}}\Div\bfv=\int_{{\varOmega}_{\rho_0}}\Div\bfg=\int_{\varSigma}\bfg\cdot\bfn\,{\rm d}\varSigma\,,
$$
so that \Eqref{3.10} follows from the latter and the assumption (c). Our next objective is to show estimates for the time derivatives of $\bsf{v}$ in the $L^2(\varOmega_{\rho_0,\rho})$-norm:
\be
\|\partial^k_t\bsf{v}\|_2\le c\,\|\partial^k_t\bfv\|_2\,.
\eeq{3.11}
These inequalities, in turn, follow from \cite[Theorem III.3.4 and Exercise III.3.7]{Gab}, if we prove that there exists a field $\bfz$ such that
$$
\Div\bfz=\nabla\psi\cdot\bfv\,,\ \ \bfz\cdot\bfn|_{\partial \varOmega_{\rho_0,\rho}}=0\,,\ \ \|\partial^k_t\bfz\|_2\le c\,\|\partial^k_t\bfv\|_2\,.
$$
To this end, for each $t\in\real$ consider the Neumann problem
\be
\Delta w(t)=\nabla\psi\cdot\bfv(t)\ \ \mbox{in $\varOmega_{\rho_0,\rho}$}\,,\ \ \nabla w(t)\cdot\bfn|_{\partial \varOmega_{\rho_0,\rho}}=0\,.
\eeq{3.12}
Because of \Eqref{3.10} this problem is uniquely solvable up to a constant in the space variables. Operating with $\partial_t^k$ in both equations in \Eqref{3.12} and testing the resulting equation in \Eqref{3.12}$_1$ with $\partial_t^kw$, after integration by parts  we get
\be
\|\nabla (\partial_t^kw)\|_2^2=(\nabla\psi\cdot\partial_t^k\bfv,\partial_t^kw)_{\varOmega_{\rho_0,\rho}}\,.
\eeq{3.13}
Again in view of \Eqref{3.10}, the right-hand side of \Eqref{3.13} remains unaltered if we replace $\partial_t^kw$ with $\partial_t^kw-\tilde{\partial_t^kw}$, where $\tilde{\textcolor{white}{abc}}$ denotes integral average over $\varOmega_{\rho_0,\rho}$. Therefore, after using Schwarz and Poincar\'e inequality, we deduce
$$
\|\nabla (\partial_t^k w)\|_2^2\le c\,\|\partial_t^k\bfv\|_2\,.
$$ 
Therefore, \Eqref{3.11} follows by choosing $\bfz=\nabla w$. It is now easy to see that the field $\bfw$ defined in \Eqref{3.8} satisfies all requirements  stated in the lemma. In fact, by extending $\bsf{v}$ to 0 outside $\varOmega_{\rho_0,\rho}$, from \Eqref{3.1}$_{2,3}$ and \Eqref{3.7} we infer properties (i)--(iii). Likewise, \Eqref{3.2}, \Eqref{3.6} and \Eqref{3.11} along with the assumptions on $\bfg$ secure the validity of (iv) and (v). The proof of the lemma is thus completed. 
\par\hfill$\square$\par
\Br We immediately see that the result just proved continues to hold, in fact, for any $k\in\nat$. However, its extension to $m>2$ would require a subsequent increase in regularity for $\Omega_0$.
\ER{reg}
\par
From \lemmref{g} we shall now derive a simple but crucial corollary. To this end, pick $\rho_0$ such that ${\varOmega}_{\rho_0}\supseteq B_0$, with $B_0$ as in \lemmref{Def},  define for $k=0,1,2$, $m=1,2$,
$$
\cals^{k,m}:=\left\{\bfg\in W^{k,2}_{T}(W^{m,2}(\varOmega)):\, \int_\varSigma \bfg(t)\cdot\bfn=0\,,\ \mbox{all $t\in\real$}; \ \bfg(x,t)=\0 \,\ \mbox{for all $|x|\ge \rho_0$}\ \mbox{and}\ t\in\real\right\}\,, 
$$ 
and also, for fixed $\rho>\rho_0$, set
$$
\cals_0^{k,m}:=\left\{\bfw\in W^{k,2}_{T}(W^{m,2}(\varOmega)):\, \bfw|_{\varSigma}=\0\,; \ \bfw(x,t)=\0 \,\ \mbox{for all $|x|\ge \rho$}\right\}\,. 
$$
We endow both spaces with the norm of $W^{k,2}(0,T;W^{m,2}(\varOmega))$.
The following result holds.
\Bp There exists $\delta_3>0$ such that for any fixed $\delta\in(0,\delta_3)$ there is a linear bounded operator
$$
\bfB:\bfg\in\cals^{k,m}\mapsto \bfB(\bfg)\in\cals_0^{k,m}
$$
such that
$$
\Div\bfB(g)+\Div(\mathbb C^\top\cdot\bfB(g))=\Div\bfg\ \ \mbox{in ${\varOmega}_\rho$}\,,
$$
and satisfying the following inequalities for $k=0,1,2$ and $m=1,2,$
\be\ba{c}\medskip\|\partial^k_t\bfB(\bfg)\|_{m,2}\le c \Sum{\ell=0}{k}\|\partial^\ell_t\bfg\|_{m,2}\,,\\ 
\|\partial^k_t\bfB(\bfg)\|_{2}\le c \Sum{\ell=0}{k}\left(\|\partial^\ell_t\bfg\|_{2}+\|\partial^\ell_t\bfg\cdot\bfn\|_{-\frac12,2(\Sigma)}\right)\,,
\ea\eeq{Trunz}
\EP{Boo}
{\em Proof.} Let $\bfg\in\cals^{k,m}$
and consider the map
$$
M:\hat{\bfw}\in \cals_0^{k,m}\mapsto\bfw
$$
where
\be
\Div\bfw=-\Div(\mathbb C^\top\cdot\hat{\bfw}-\bfg):=\Div\bfG\ \ \mbox{in ${\varOmega}_\rho$}\,.
\eeq{dw}
It is at once checked that $\bfG\in\cals^{k,m}$, which means that $\bfG$ satisfies all the requests made on the function $\bfg$ in \lemmref{g}. Therefore, by that lemma, there exists  $\bfw\in\cals_0^{k,m}$ satisfying \Eqref{dw}. Thus, also with the help of \lemmref{Def}, we show for $k=0,1,2,$ and $m=1,2$ 
\be\ba{c}\medskip\|\partial^k_t\bfw\|_{m,2}\le c\,\|\partial^k_t\bfG\|_{m,2}\le c\left[\delta\left(\Sum{\ell=0}{k}\|\partial_t^\ell\hat{\bfw}\|_{m,2}\right)+\|\partial^k_t\bfg\|_{m,2}\right]\,,\\ \|\partial^k_t{\bfw}\|_{2}\le c\big(\|\partial^k_t\bfG\|_{2}+\|\partial^k_t\bfG\cdot\bfn\|_{-\frac12,2(\Sigma)}\big)\le c\,\left[\delta\left(\Sum{\ell=0}{k}\|\partial_t^\ell\hat{\bfw}\|_2\right)+\|\partial^k_t\bfg\|_{2}+\|\partial^k_t\bfg\cdot\bfn\|_{-\frac12,2(\Sigma)}\right]\,.
\ea\eeq{dw1}
From \Eqref{dw1}$_1$ it readily follows that, for sufficiently small $\delta$, there exists $R>0$ such that $M$ maps the ball of $\cals_0^{k,m}$ of radius $R$ into itself and, in fact, $M$ is contracting. This immediately leads to the existence of the operator $\bfB$ with the stated properties, which concludes the proof of the lemma. 
\par\hfill$\square$\par
We conclude this section with the following result.
\Bl Let 
$$
\bfg\in W^{3,2}_T(L^2({\varOmega}))\cap W^{2,2}_T(W^{2,2}({\varOmega}))\cap W^{3,2}_T(W^{-\frac12,2}(\varSigma))\,. 
$$ 
Then, for each $t\in (0,T)$ the problem
\be\ba{cc}\medskip
\Div\bfv(t)=\Div\bfg(t)\,,\ \ \mbox{in $\varOmega$}\\
\bfv(t)=\0\,,\ \ \mbox{at $\varSigma$}\,,
\ea
\eeq{Bo1}
has at least one solution $\bfv \in W^{3,2}_T(L^2({\varOmega}))\cap W^{2,2}_T(W^{2,2}({\varOmega}))$ such that
\be\ba{c}\medskip\|\partial_t^\ell D^2\bfv(t)\|_{2}+\|\partial_t^\ell\nabla\bfv(t)\|_{2}\le c\,\|\partial_t\bfg(t)\|_{2,2}\,,\ \ \ell=0,1,2\,,\\ \|\partial^k_t\bfv(t)\|_{2}\le c\,\big(\|\partial^k_t\bfg(t)\|_{2}+\|\partial^k_t\bfg(t)\cdot\bfn\|_{-\frac12,2(\Sigma)}\big)\,,\ \ k=0,1,2,3\,.
\ea\eeq{Bo11}
Moreover, if $\bar{\bfg}=\0$ we can take $\bar{\bfv}=\0$ as well. 
\EL{Bog}
{\em Proof.} We begin to show the last statement. Suppose we find $\bfv$ satisfying \Eqref{Bo1} and \Eqref{Bo11}. If $\bar{\bfg}=\0$, then $\Div\bar{\bfv}=0$, which implies that  the function $\bsf v:=\bfv-\bar{\bfv}$ is also a solution to  \Eqref{Bo1} and \Eqref{Bo11} with $\bar{\bsf v}=\0$, which proves the claim. Consider now the Neumann problem
\be
\ba{cc}\medskip
\Delta\psi=\Div\bfg\,,\ \ \mbox{in $\varOmega$}\\
\bfn\cdot\nabla\psi=\0\,,\ \ \mbox{at $\varSigma$}\,.
\ea
\eeq{Bo2}
By standard methods  we show the existence of a unique solution satisfying
\be\ba{cc}\medskip
\|\partial_t^\ell\nabla\psi(t)\|_{2,2}\le c\,\|\partial_t^\ell\bfg(t)\|_{2,2}\,,\\ \medskip \|\partial^k_t\psi(t)\|_{6}+\|\partial^k_t\nabla\psi(t)\|_{2}\le c\,\big(\|\partial^k_t\bfg(t)\|_{2}+\|\partial^k_t\bfg(t)\cdot\bfn\|_{-\frac12,2(\Sigma)}\big)\,;
\ea\eeq{Bo3} 
see, e.g., \cite[Sections 7.29-7.30]{Medk}. Since $\bfg$ is $T$-periodic,   by the uniqueness property so is $\psi$.
Next, consider the problem
\be\ba{cc}\medskip
\Div\bfw=0\,,\ \ \mbox{in $\varOmega$}\\
\bfw=-\nabla\psi\,,\ \ \mbox{at $\varSigma$}\,.
\ea
\eeq{Bo4}
We look for a solution to \Eqref{Bo4} of the form
\be
\bfw:=-\nabla(\zeta\psi)+\bfu\,,
\eeq{Bo5}
where $\zeta$ is a ``cut-off" function that is 1 in a neighborhood of $\varSigma$ and  0 for $|x|\ge R/2$, $R>2R_*$, while $\bfu$ solves the problem 
\be\ba{cc}\medskip
\Div\bfu=\Div(\nabla(\zeta\psi))\,,\ \ \mbox{in $\varOmega_R$}\\
\bfu=0\,,\ \ \mbox{at $\partial\varOmega_R$}\,.
\ea
\eeq{Bo6}
Since
$$
\supp(\nabla(\zeta\psi))\subset {\varOmega}_{R/2}\,,\ \ \int_\varSigma\nabla(\zeta\psi)\cdot\bfn=\int_\varSigma\nabla\psi\cdot\bfn\,,
$$
by \lemmref{g}, \Eqref{Bo2}$_2$, and \Eqref{Bo3} it follows that \Eqref{Bo6} has a solution $\bfu\in W^{k,2}_{T}(W^{2,2}(\varOmega))$ with $\supp(\bfu)\subset \varOmega_R$ and such that
\be\ba{c}\medskip
\|\partial_t^\ell\bfu(t)\|_{2,2}\le c\,(\|\partial_t^\ell\nabla\psi(t)\|_{2,2}+\|\partial_t^\ell\psi\|_6)\,,\\ \|\partial^k_t\bfu(t)\|_{2}\le c\,(\|\partial^k_t\nabla\psi(t)\|_{2}+\|\partial_t^k\psi\|_6)\,.
\ea
\eeq{Bo7}
In view of \Eqref{Bo2}--\Eqref{Bo7}, it is then readily checked that the vector field $\bfv:=\nabla\psi+\bfw$ satisfies all the properties stated in the lemma.
\par\hfill$\square$\par

\setcounter{equation}{0}
\section{Existence of $T$-Periodic Solutions to a Linear Problem} \label{Section:LinearProblem}
Objective of this section is the study of some relevant properties of $T$-periodic solutions to the following linear problem 
\be\ba{cc}\medskip\left.\ba{ll}\medskip
\partial_t{\bfV}=\Div\mathbb T(\bfV,{\sf p})\\
\Div\bfV=-\Div(\mathbb C^\top\cdot\bfV)\ea\right\}\ \ \mbox{in $\varOmega\times\real$}\,,\\ \medskip
\bfV=\bfu_*+\bfzeta\,, \ \, \mbox{at $\varSigma\times \real$}\,;\\ 
{\sf m}\dot{\bfzeta}(t)=-\Int{\varSigma}{}\mathbb T(\bfV,{\sf p})\cdot\bfn{\rm d}\varSigma\,,\ \ t\in\real\,.
\ea
\eeq{3.1_0}
More specifically, we shall first show the existence and uniqueness of such solutions in a suitable function class, and, successively,  will provide an explicit relation between the time-averaged velocity, $\bar{\bfzeta}$, and the data, $\bfu_*$.  The current section will be devoted to the first problem, while the second one will be considered in the next section. Note that in these two sections, we do not explicitly state the conditions at infinity. However, such conditions are inherently embedded within the function spaces employed.
\par
We begin to split the solution into its average and oscillatory components:
\be
\bfV=\bar{\bfV}+\bsf V\,,\ \ {\sf p}=\bar{\sf p}+{\sf q}\,,\ \ \bfzeta=\bar{\bfzeta}+\bfchi\,.
\eeq{Split}
Also, for $s\in (1,\frac32)$, we define
\be\ba{ll}\medskip
\cald^s(\varOmega):=L^{\frac{3s}{3-2s}}(\varOmega)\cap D^{1,\frac{3s}{3-s}}(\varOmega)\cap D^{2,s}(\varOmega)\cap D^{2,2}(\varOmega)\,,\ \calp^s(\varOmega):= L^{\frac{3s}{3-s}}(\varOmega)\cap D^{1,s}(\varOmega)\cap D^{1,2}(\varOmega)\,,
\\ \medskip
\calw({\varOmega}):=W^{3,2}_T(0,T;L^2({\varOmega}))\cap W^{2,2}_T(0,T;W^{2,2}({\varOmega}))\,,\ \ \calq(\varOmega):=W^{2,2}_T(0,T;D^{1,2}({\varOmega})\cap L^6({\varOmega}))\,;\\  \calv:=W_T^{3,2}(W^{\frac32,2}(\varSigma))\,.
\ea
\eeq{calw}
Our goal is thus to prove the following result.
\Bp Suppose 
$$
\bfu_*\in \calv\,.
$$
Then, there exists $\delta^*>0$ such that, for all $\delta\in(0,\delta^*)$, problem \Eqref{3.1_0}
has one and only one $T$-periodic solution  \Eqref{Split} where 
\be
(\bar{\bfV},\bar{ \sf p},\bar{\bfzeta})\in \cald^{s}(\varOmega)\times \calp^{s}(\varOmega)\times\real^3\,,\ \ (\bsf{V},{\sf q},\bfzeta) \in\calw(\varOmega)\times \calq(\varOmega)\times W^{3,2}_T(\real^3)\,,
\eeq{AMPP}
for arbitrary  $s\in(1,\frac32)$.
This solution satisfies the estimate
\be
\|\bar{\bfV}\|_{\cald^s(\varOmega)}+\|\bar{\sf p}\|_{\calp^s(\varOmega)}+|\bar{\bfzeta}|+\|\bsf V\|_{\calw(\varOmega)}+\|{\sf q}\|_{\calq(\varOmega)}+\|\bfchi\|_{W^{3,2}(0,T)}\le c\,\|\bfu_*\|_\calv\,.
\eeq{AMPP_0}
Moreover, 
there exist  $\rho_*, C>0$  such that
\be
\Sup{(x,t)\in \varOmega^{\rho_*}\times [0,T]}\left(|x|^{m+2}|D^k{\bfV}(x)|\right)\le C\,,\ \ |k|=m\ge0\,.
\eeq{west1_1}\label{O_k}
\EP{O_k}

The proof of \propref{O_k} will be given in Subsection \ref{sub:proof}, after showing a  number of auxiliary results.
\subsection{Unique Solvability of Auxiliary Linear Problems}
This section is devoted to the well-posedness of steady-state and time-periodic linear problems, relevant to the proof of \propref{O_k}. 
\Bl Let $s\in (1,\frac32)$, and suppose
$$
\bfG\in L^s(\varOmega)\cap L^2(\varOmega)\,,\ \ g\in W^{1,s}(\varOmega)\cap W^{1,2}(\varOmega)\,. 
$$
Then, the problem
\be\ba{cc}\medskip\left.\ba{ll}\medskip
\nu\Delta\bfv-\nabla{p}=\bfG\\
\Div\bfv=g\ea\right\}\ \ \mbox{in $\varOmega$}\,,\\ \medskip
\bfv=\bfgamma\,, \ \, \mbox{at $\varSigma$}\,;\\ 
\Int{\varSigma}{}\mathbb T(\bfv,{ p})\cdot\bfn{\rm d}\varSigma=0\,,
\ea
\eeq{4.27}
has one and only one solution 
\be(\bfv,p,\bfgamma)\in \cald^s(\varOmega)\times\calp^s(\varOmega)\times\real^3\,.
\eeq{4.28} 
Moreover, this solution satisfies
\be
\|\bfv\|_{\cald^s(\varOmega)}+\|p\|_{\calp^s(\varOmega)}+|\bfgamma|\le C\,(\|\bfG\|_2+\|\bfG\|_s+\|g\|_{1,s}+\|g\|_{1,2})\,.
\eeq{4.29}
\EL{bar}
{\em Proof.} We write $\bfv=\bfv_1+\bfv_2$, $p=p_1+p_2$, where
\be\ba{cc}\medskip\left.\ba{ll}\medskip
\nu\Delta\bfv_1-\nabla{p}_1=\0\\
\Div\bfv_1=g\ea\right\}\ \ \mbox{in $\varOmega$}\,,\\ 
\bfv_1=0\,, \ \, \mbox{at $\varSigma$}\,;\ea
\eeq{4.30}
and
\be\ba{cc}\medskip\left.\ba{ll}\medskip
\nu\Delta\bfv_2-\nabla{p}_2=\bfG\\
\Div\bfv_2=0\ea\right\}\ \ \mbox{in $\varOmega$}\,,\\ \medskip
\bfv_2=\bfgamma\,, \ \, \mbox{at $\varSigma$}\,;\\ 
\Int{\varSigma}{}\mathbb T({\bfv_2,{ p_2}})\cdot\bfn{\rm d}\varSigma=-\Int{\varSigma}{}\mathbb T(\bfv_1,{ p_1})\cdot\bfn{\rm d}\varSigma:=\bfF\,.
\ea
\eeq{4.31}
Under the stated assumption on $g$, from \cite[Theorem V.4.8 and Exercise V.4.9]{Gab} it follows the existence of $(\bfv_1,p_1)$ in the class specified in \Eqref{4.28}, and satisfying the estimate \Eqref{4.29} (with $\bfgamma=\0$ {and \bfG=0}). Next, for $i=1,2,3$, consider the following Stokes problems
\be\ba{cc}\medskip\left.\ba{ll}\medskip
\Div\mathbb T(\bfh^{(i)},{p}^{(i)})=\0\\
\Div\bfh^{(i)}=0\ea\right\}\ \ \mbox{in $\varOmega$}\,,\\ 
\bfh^{(i)}=\bfe_i\,, \ \, \mbox{at $\varSigma$}\,.\ea
\eeq{4.32}
Again by  \cite[Theorem V.4.8]{Gab} we deduce that \Eqref{4.32} has a unique solution in the class \Eqref{4.28}, and 
\be
\|\bfh^{(i)}\|_{\cald^s(\varOmega)}+\|p^{(i)}\|_{\calp^s(\varOmega)}\le C\,.
\eeq{4.33}
Then, a solution to \Eqref{4.31} is given by
\be
\bfv_2=\sum_{i=1}^3\gamma_i\bfh^{(i)}\,,\ \ p_2=\sum_{i=1}^3\gamma_ip^{(i)}\,,
\eeq{4.34}
where the $\gamma_i$, $i=1,2,3$, solve the linear algebraic system
$$
\gamma_i\Int{\varSigma}{}\mathbb T(\bfh^{(i)},{ p}^{(i)})\cdot\bfn{\rm d}\varSigma=\bfF\,.
$$
Since the matrix
\be
\mathbb M_{ji}=\Int{\varSigma}{}\mathbb T_{jk}(\bfh^{(i)},{ p}^{(i)})\,n_k{\rm d}\varSigma\,,\ \ i,j=1,2,3,
\eeq{Matrix}
is invertible \cite[Lemma 4.1]{Gah}, the above system has one and only one solution $\bfgamma$ satisfying
$$
|\bfgamma|\le c\,|\bfF|\,.
$$
The latter, in combination with \Eqref{4.33} and \Eqref{4.34}, furnishes
\be
\|\bfv_2\|_{\cald^s(\varOmega)}+\|p_2\|_{\calp^s(\varOmega)}+|\bfgamma|\le C\,|\bfF|\,.
\eeq{4.35}
However, by classical trace theorems, 
$$
|\bfF|\le c\,(\|\bfv_1\|_{2,s,\varOmega_\rho}+\|p_1\|_{1,s,\varOmega_\rho})
$$
for any fixed $\rho>R_*$. Therefore, the stated existence along with the validity of \Eqref{4.29} follows from the latter, \Eqref{4.35} and the previously proved properties of $(\bfv_1,p_1)$. Consider, now, \Eqref{4.27} with $g\equiv0$. Then, testing \Eqref{4.27}$_1$ with $\bfv$, integrating by parts over $\varOmega$ and using \Eqref{4.27}$_4$ gives $\|\mathbb D(\bfv)\|_2=0$, which, since $\bfv$ is in the class \Eqref{4.28}, furnishes at once $\bfv\equiv\bfgamma\equiv\0$,  $p\equiv0$, which shows uniqueness.
\par\hfill$\square$\par
\Bl Suppose 
$$
\bsf f\in  W^{2,2}_T(L^2)\,, \   {\bsf g}\in \calw(\varOmega)\cap W^{3,2}_T(W^{-\frac12,2}(\varSigma))\,,\ \bsf F\in W^{2,2}_T(\real^3)\,,\ \ \bar{\bsf f}=\bar{\bsf g}=\bar{\bsf F}=\0\,.
$$
Then, the problem
\be\ba{cc}\medskip\left.\ba{ll}\medskip
\partial_t{\bsf v}=\nu\Delta\bsf v-\nabla{\sf r}+\bsf{f}\\
\Div\bsf v=\Div{\bsf g}\ea\right\}\ \ \mbox{in $\varOmega\times\real$}\,,\\ \medskip
\bsf v=\bfmu\,, \ \, \mbox{at $\varSigma\times \real$}\,;\\ 
{\sf m}\dot{\bfmu}(t)=-\Int{\varSigma}{}\mathbb T(\bsf v,{\sf r})\cdot\bfn{\rm d}\varSigma+\bsf F\,,\ \ t\in\real\,.
\ea
\eeq{AMP_0}
has one and only one $T$-periodic solution 
\be
(\bsf{v},{\sf r},\bfmu)\in \calw(\varOmega)\times \calq(\varOmega)\times W^{3,2}_T(\real^3)\,.
\eeq{AMP_1}
Moreover, 
\be
\|\bsf{v}\|_{\calw({\varOmega})}+\|{\sf r}\|_{\calq({\varOmega})}+\|\bfmu\|_{W^{3,2}(0,T)}\le c\,\left(\|\bsf{f}\|_{W^{2,2}(L^2)}+\|{\bsf g}\|_{\calw({\varOmega})}+\|\bsf g\|_{W^{3,2}(W^{-\frac12,2}(\varSigma))}+\|\bsf{F}\|_{W^{2,2}(0,T)}\right)\,.
\eeq{est_1}
\EL{Os}
{\em Proof.} We begin to operate a  lift of $\Div\bsf g$. 
From assumption and \lemmref{Bog}, there exists $\bfv\in \calw(\varOmega)$, $\bar{\bfv}=\0$,  satisfying \Eqref{Bo1} and \Eqref{Bo11} with $\bfg\equiv\bsf g$. 
Thus, setting
\be
\bff_1:=-\partial_t\bfv+\nu\Delta\bfv,\ \ \bfF_1:=-2\Int{\varSigma}{}\mathbb D(\bfv)\cdot\bfn\,{\rm d}\varSigma\,,
\eeq{AMP_3}
because of  \Eqref{Bo11} and classical trace theorems, we get
\be
\|\bff_1\|_{W^{2,2}(L^2(\varOmega))}+\|\bfF_1\|_{W^{2,2}(0,T)}\le c_1 \|\bfv\|_{\calw(\varOmega)}\le c_2\,\left(\|{\bsf g}\|_{\calw({\varOmega})}+\|\bsf g\|_{W^{3,2}(W^{-\frac12,2}(\varSigma))}\right)\,.
\eeq{AMP_4}
If we now write
\be
\bsf v=\bsf u+\bfv\,,
\eeq{SpLiT}
and define \be\bff:=\bff_1+\bsf f\,,\ \  \bfF:=\bfF_1+\bsf F\,,\eeq{AMP_8} 
from the properties of $\bfv$, we infer that  
\Eqref{AMP_0} entails
\be\ba{cc}\medskip\left.\ba{ll}\medskip
\partial_t{\bsf u}=\nu\Delta\bsf u-\nabla{\sf r }+\bff\\
\Div\bsf u=0\ea\right\}\ \ \mbox{in $\varOmega\times\real$}\,,\\ \medskip
\bsf u=\bfmu\,, \ \, \mbox{at $\varSigma\times \real$}\,;\\ 
{\sf m}\dot{\bfmu}(t)=-\Int{\varSigma}{}\mathbb T(\bsf u,{\sf r})\cdot\bfn{\rm d}\varSigma+\bfF\,,\ \ t\in\real\,.
\ea
\eeq{AMP_9}
Since ${\bff},{\bfF}$ are $T$-periodic with $\bar{\bff}=\bar{\bfF}=\0$, by \cite[Lemma 5.2]{GaspU}, there exists one and only one $T$-periodic solution to \Eqref{AMP_9} such that
\be
(\bsf u,{\sf r},\bfmu)\in [W^{1,2}(L^2(\varOmega))\cap L^2(W^{2,2}(\varOmega))]\times L^2(D^{1,2}(\varOmega))\times W^{1,2}(0,T).
\eeq{AMP_10}
This solution satisfies also the estimate
\be
\|\bsf u\|_{W^{1,2}(L^2)\cap L^2(W^{2,2})}+\|{\sf r}\|_{L^2(D^{1,2})}+\|\bfmu\|_{W^{1,2}(0,T)}\le c\,\left(\|\bff\|_{L^{2}(L^2)}+\|\bfF\|_{L^{2}(0,T)}\right)\,.
\eeq{AMP_11}
Let $h=h(x,t)$,  be a $T$-periodic (scalar or vector) function, and denote by $h_\eta$ its Friederichs mollifier in time, namely,
$$
h_\eta(x,t)=\int_\real j_\eta(s)h(x,t-s){\rm d}s\,,  
$$
where $j_\eta(\tau):=\eta^{-1}j(\tau/\eta)$, with $j\in C_0^\infty(-1,1)$. It then readily follows that the functions
$$
\bsf u^\prime:=\partial_t\bsf u_\eta\,,\ \ {\sf r}^\prime:=\partial_t{\sf r}_\eta\,,\ \ \bfmu^\prime:=\dot{\bfmu}_\eta
$$
belong to the class \Eqref{AMP_10}, and 
satisfy \Eqref{AMP_9}  with $\bff,\bfF$ replaced by $\partial_t\bff_\eta$ and $\dot{\bfF}_\eta$, respectively. Again by \cite[Lemma 5.2]{GaspU}, and since  $\bff\in W^{1,2}(L^2(\varOmega))$,  $\bfF\in W^{1,2}(0,T)$, from \Eqref{AMP_11} it also follows that  
\be
\|\bsf u^\prime\|_{W^{1,2}(L^2)\cap L^2(W^{2,2})}+\|{\sf r}^\prime\|_{L^2(D^{1,2})}+\|\bfmu^\prime\|_{W^{1,2}(0,T)}\le c\,\left(\|\bff\|_{W^{1,2}(L^2)}+\|\bfF\|_{W^{1,2}(0,T)}\right)\,.
\eeq{AMP_12}
Thus, letting $\eta\to0$ in \Eqref{AMP_12} we show that
$$
(\bsf u,{\sf r},\bfmu)\in [W^{2,2}(L^2(\varOmega))\cap W^{1,2}(W^{2,2}(\varOmega))]\times W^{1,2}(D^{1,2}(\varOmega))\times W^{2,2}(0,T).
$$
along with the estimate
$$
\|\bsf u\|_{W^{2,2}(L^2)\cap W^{1,2}(W^{2,2})}+\|{\sf r}\|_{W^{1,2}(D^{1,2})}+\|\bfmu\|_{W^{2,2}(0,T)}\le c\,\left(\|\bff\|_{W^{1,2}(L^2)}+\|\bfF\|_{W^{1,2}(0,T)}\right)\,.
$$
Employing the mollifying procedure one more time, we can thus prove
$$
(\bsf u,{\sf r},\bfmu)\in \calw(\varOmega)\times \calq(\varOmega)\times W^{3,2}(0,T)\,,
$$
as well as
$$
\|\bsf u\|_{\calw(\varOmega)}+\|{\sf r}\|_{\calq(\varOmega)}+\|\bfmu\|_{W^{3,2}(0,T)}\le c\,\left(\|\bff\|_{W^{2,2}(L^2)}+\|\bfF\|_{W^{2,2}(0,T)}\right)\,.
$$
Therefore, combining the latter with \Eqref{AMP_4}--\Eqref{AMP_8}, we conclude the proof of the lemma. 
\par\hfill$\square$\par
\subsection{Proof of Proposition \ref{O_k}}\label{sub:proof}
We commence by rewriting problem \Eqref{3.1_0} in a form that is amenable to a fixed-point argument. To this end,
we construct a suitable lift ofthe boundary data  $\bfu_*$ as follows. For each $t\in [0,T]$, consider the problem\Footnote{Throughout the proof, we set $\nu=1$  for simplicity, since its real value is irrelevant.} 
\be\ba{cc}\medskip\left.\ba{ll}\medskip
\bsf U(t)=\Delta\bsf U(t)-\nabla\tau(t)\\
\Div\bsf U(t)=0\ea\right\}\ \ \mbox{in $\varOmega$}\,,\\ 
\bsf U(t)=\bfu_*(t)\,, \ \, \mbox{at $\varSigma$}\,.
\ea
\eeq{AMP_5}
It is then known that, under the given assumptions on $\bfu_*$, problem \Eqref{AMP_5} has a unique $T$-periodic solution $(\bsf U,\tau)\in W^{2,2}(W^{2,2}(\varOmega))\times W^{2,2}(D^{1,2}(\varOmega))$ with $\bar{\bsf U}=\0$, such that
\be
\|\bsf U\|_{W^{3,2}(W^{2,2})}+\|\nabla\tau\|_{W^{3,2}(L^{2})}\le c\,\|\bfu_*\|_\calv\,.
\eeq{AMP_6}
Moreover, setting
$$
\bff:=-\partial_t\bsf U+\Delta\bsf U,\ \ \bfg:=\mathbb C^\top\cdot\bsf U\,,\ \ \bfF:=-2\Int{\varSigma}{}\mathbb D(\bsf U)\cdot\bfn\,{\rm d}\varSigma\,, 
$$
from \Eqref{AMP_6},  classical trace theorems, and \lemmref{Def}, we deduce 
\be\ba{cc}\medskip
\|\bff\|_{W^{2,2}(L^2(\varOmega))}+\|\bfg\|_{\calw(\varOmega)}+\|\bfg\|_{\calv}+\|\bfF\|_{W^{2,2}(0,T)}\le c\,\|\bfu_*\|_\calv\\
\|\Div\bar{\bfg}\|_{1,q}\le c\,\|\bfu_*\|_\calv\,,\ \ \mbox{all $q\in(1,2]$}\,.
\ea
\eeq{AMP_7}
We also recall that, see \remref{div},
\be
\Div\mathbb T(\bfV,{\sf p})=\Delta\bfV-\nabla\Div(\mathbb C^\top\cdot\bfV)-\nabla{\sf p}\,.
\eeq{divuno}
Thus, with
\be
\bfV:=\bfU+\bsf U\,,
\eeq{pupa}
problem \Eqref{3.1_0} can be equivalently rewritten as follows
\be\ba{cc}\medskip\left.\ba{ll}\medskip
\partial_t{\bfU}=\Delta\bfU-\nabla\Div(\mathbb C^\top\cdot\bfU+\bfg)-\nabla{\sf p}+\bff\\
\Div\bfU=-\Div(\mathbb C^\top\cdot\bfU+\bfg)\ea\right\}\ \ \mbox{in $\varOmega\times\real$}\,,\\ \medskip
\bfU=\bfzeta\,, \ \, \mbox{at $\varSigma\times \real$}\,;\\ 
{\sf m}\dot{\bfzeta}(t)=-\Int{\varSigma}{}\mathbb T(\bfU,{\sf p})\cdot\bfn{\rm d}\varSigma+\bfF\,,\ \ t\in\real\,.
\ea
\eeq{3.1_00}
Notice that, by the property of the extension $\bsf U$, we have
$$
\bar{\bff}=\bar{\bfF}=\0\,.
$$
Let 
$$\ba{ll}\medskip
\mathscr D^s:=\left\{(\bfv,p,\bfgamma)\in\cald^s(\varOmega)\times\calp^s(\varOmega)\times\real^3:\, \bfv=\bfgamma\ \mbox{at $\varSigma$}\right\}\,,\ \ s\in(1,\frac32)\,,\\ \medskip
\mathscr U:=\left\{(\bfu,{\sf q},\bfmu)\in \calw(\varOmega)\times \calq(\varOmega)\times W^{3,2}_T(\real^3)\, : \ \bfu=\bfmu \ \mbox{at $\varSigma$}\,;\ \bar{\bfu}=\bar{\sf q}=\bar{\bfmu}=0\right\}\\
\mathscr H:=\{(\bfW,q,\bfxi): \, (\bar{\bfW},\bar{q},\bar{\bfxi})\in \mathscr D^s\,,\, (\bfw:=\bfW-\bar{\bfW},{\sf q}:=q-\bar{q},\bfchi:=\bfxi-\bar{\bfxi})\in \mathscr U\}\,,
\ea
$$
and set
$$
\|(\bfW,q,\bfxi)\|_{\mathscr H}:=\|\bar{\bfW}\|_{\cald^s(\varOmega)}+\|\bar{q}\|_{\calp^s(\varOmega)}+|\bar{\bfxi}| +\|\bfw\|_{\calw(\varOmega)}+\|{\sf q}\|_{\calq(\varOmega)}+\|\bfchi\|_{W^{3,2}(0,T)}\,.
$$
Next, define the map
$$
M: (\bfW,q,\bfxi)\in \mathscr H\mapsto (\bfU=\bar{\bfU}+\bsf S\,, {\sf p}=\bar{\sf p}+{\sf q}\,, \bfzeta=\bar{\bfzeta}+\bfsigma)
$$
where
\be\ba{cc}\medskip\left.\ba{ll}\medskip
\Delta\bar{\bfU}=\nabla\bar{\sf p}+\nabla\Div(\bar{\mathbb C^\top\cdot\bfW+\bfg})\\
\Div\bar{\bfU}=-\Div(\bar{\mathbb C^\top\cdot\bfW+\bfg})\ea\right\}\ \ \mbox{in $\varOmega$}\,,\\ \medskip
\bar{\bfU}=\bar{\bfzeta}\,, \ \, \mbox{at $\varSigma$}\,;\\ 
\Int{\varSigma}{}\mathbb T(\bar{\bfU},\bar{\sf p})\cdot\bfn{\rm d}\varSigma=\0\,,
\ea
\eeq{AMP_14}
and
\be\ba{cc}\medskip\left.\ba{ll}\medskip
\partial_t{\bsf S}=\Delta\bsf S-\nabla{\sf q}-\nabla\Div(\mathbb C^\top\cdot\bfW+\bfg-\bar{\mathbb C^\top\cdot\bfW+\bfg})+\bff\\
\Div\bsf S=-\Div(\mathbb C^\top\cdot\bfW+\bfg-\bar{\mathbb C^\top\cdot\bfW+\bfg})\ea\right\}\ \ \mbox{in $\varOmega\times\real$}\,,\\ \medskip
\bsf S=\bfsigma\,, \ \, \mbox{at $\varSigma\times \real$}\,;\\ 
{\sf m}\dot{\bfsigma}(t)=-\Int{\varSigma}{}\mathbb T(\bsf S,{\sf q})\cdot\bfn{\rm d}\varSigma+\bfF\,,\ \ t\in\real\,.
\ea
\eeq{AMP_15}
Clearly, if $M$ has a fixed point then, in view of \Eqref{divuno}, the triple $(\bfU,{\sf p},\bfzeta)$ is a solution to \Eqref{3.1_00} in the class $\mathscr H$. Observing that, by \Eqref{PId},
$$
g:=\Div(\bar{\mathbb C^\top\cdot\bfW+\bfg})=\bar{\mathbb C^\top:\nabla\bfW}+\Div\bar{\bfg}\,,
$$
with the help of  \lemmref{Def} and  \Eqref{AMP_7}$_2$ we get (recall that $\mathbb C$ is of bounded support)
\be
\|g\|_{1,s}+\|g\|_{1,2}\le c\left[\delta\left(\|\bar{\bfW}\|_{\cald^s(\varOmega)}+\|\bfw\|_{\calw(\varOmega)}\right)+\|\bfu_*\|_\calv\right]\,.
\eeq{AMP_16}
Then, from \Eqref{AMP_16} and \lemmref{bar}, we conclude that \Eqref{AMP_14} has one and only one solution: 
\be(\bar{\bfU},\bar{\sf p},\bar{\bfzeta})\in \cald^s(\varOmega)\times\calp^s(\varOmega)\times\real^3\,,
\eeq{Ak}
such that
\be
\|\bar{\bfU}\|_{\cald^s(\varOmega)}+\|\bar{\sf p}\|_{\calp^s(\varOmega)}+|\bar{\bfzeta}|\le  c\,\left[\delta\,\left(\|\bar{\bfW}\|_{\cald^s(\varOmega)}+\|\bfw\|_{\calw(\varOmega)}\right)+\|\bfu_*\|_\calv\right]\,.
\eeq{AMP_17}
Next, setting
$$
\bsf G:=-\Div\bsf g\,,\ \ 
\bsf g:=\mathbb C^\top\cdot\bfW+\bfg-\bar{\mathbb C^\top\cdot\bfW+\bfg}\,,
$$
again by employing the properties of $\mathbb C$ proved in \lemmref{Def} and \Eqref{AMP_7}$_1$, we show $\bsf g\in \calw(\varOmega)\cap \calv$, and that
\be
\|\nabla \bsf G\|_{W^{2,2}(L^2)}+\|\bsf g\|_{\calw(\varOmega)}+\|\bsf g\|_{\calv}\le c\,\left[\delta\left(\|\bar{\bfW}\|_{\cald^s(\varOmega)}+\|\bfw\|_{\calw(\varOmega)}+|\bar\bfxi|+\|\bfchi\|_{W^{3,2}(0,T)}\right)+\|\bfu_*\|_\calv\right]\,.
\eeq{KTM}
Therefore, combining the property of $\bff$ given in \Eqref{AMP_7}$_1$ with \Eqref{KTM}, and \lemmref{Os}, we infer that \Eqref{AMP_15} has one and only one solution 
\be
(\bsf{S},{\sf q},\bfsigma)\in \calw(\varOmega)\times \calq(\varOmega)\times W^{3,2}(0,T)\,,
\eeq{AMP_18}
and, in addition,
\be
\|\bsf{S}\|_{\calw({\varOmega})}+\|{\sf q}\|_{\calq({\varOmega})}+\|\bfsigma\|_{W^{3,2}(0,T)}\le c\,\left[\delta\left(\|\bar{\bfW}\|_{\cald^s(\varOmega)}+\|\bfw\|_{\calw(\varOmega)}+|\bar{\bfxi}|+\|\bfxi\|_{W^{3,2}(0,T)}\right)+\|\bsf{u}_*\|_\calv\right]\,.
\eeq{AMP_19}
Thus, we conclude that, on one hand, by \Eqref{Ak} and \Eqref{AMP_18},  $M$ maps $\mathscr H$ into itself and, on the other hand, 
\be   
\|M(\bfW,q,\bfxi)\|_{\mathscr H}\le c\,(\delta\|(\bfW,q,\bfxi)\|_{\mathscr H}+\|\bfu_*\|_\calv)\,.
\eeq{AMP_20}
From the latter, it readily follows that if we choose $\delta$ less than a suitable constant, we can find  $\rho>0$ such that $M$ maps the ball of radius $\rho$ of $\mathscr H$ into itself. 
Moreover, in the same manner, we obtain
\begin{equation*}
\|M(\bfW_1-\bfW_2,q_1-q_2,\bfxi_1-\bfxi_2)\|_{\mathscr H}\le c\delta\|(\bfW_1-\bfW_2,q_1-q_2,\bfxi_1-\bfxi_2)\|_{\mathscr H}.
\end{equation*}
Consequently, $M$ has a fixed point $(\bfU,{\sf p},\bfzeta)$ there. Additionally, by \Eqref{AMP_20},
$$
\|(\bfU,{\sf p},\bfzeta)\|_{\mathscr H}\le c\,\|\bfu_*\|_\calv\,.
$$
The latter, along with \Eqref{pupa} and \Eqref{AMP_6}, concludes the existence and uniqueness proof. We shall next show the asymptotic properties stated in \Eqref{west1_1}. By averaging \Eqref{3.1_0}, we get
\be\ba{cc}\medskip\left.\ba{ll}\medskip
\Div\mathbb T(\bar{\bfV},\bar{\sf p})=\0\\
\Div\bar{\bfV}=-\Div(\bar{\mathbb C^\top\cdot\bfV})\ea\right\}\ \ \mbox{in $\varOmega$}\,,\\ \medskip
\bar{\bfV}=\bar{\bfzeta}\,, \ \, \mbox{at $\varSigma$}\,;\\ 
\Int{\varSigma}{}\mathbb T(\bar{\bfV},\bar{\sf p})\cdot\bfn{\rm d}\varSigma=\0\,.
\ea
\eeq{AMP_22}
Pick sufficiently large $R$ such that $B_R$ contains the support of $\mathbb C$ and integrate both sides of \Eqref{AMP_22}$_1$ in ${\varOmega}_R$. Taking into account \Eqref{AMP_22}$_4$, we deduce that $(\bar{\bfV},\bar{\sf p})$ obeys the following equations
\be\ba{cc}\medskip\left.\ba{ll}\medskip
\Div\mathbb T(\bar{\bfV},\bar{\sf p})=\0\\
\Div\bar{\bfV}=0\ea\right\}\ \ \mbox{in $\varOmega^R$}\,,\\  
\Int{\partial\varOmega^R}{}\mathbb T(\bar{\bfV},\bar{\sf p})\cdot\bfn{\rm d}\varSigma=\0\,,
\ea
\eeq{AMP_23}
so that the estimate \Eqref{west1_1} for $\bar{\bfV}$ follows from classical results \cite[Theorem V.3.2]{Gab}. It remains to  show the validity of \Eqref{west1_1} for $\bsf V$. Let $\psi=\psi_{\rho,\rho_1}(|x|)$ be a ``cut-off" function that is 0 in $\varOmega_{\rho}$ and 1 in $\varOmega^{\rho_1}$, for  fixed $\rho_1>\rho>R$, and set
\be
\hat{\bsf V}=\psi\,{\bsf V}\,,\ \ \hat{\sf q}=\psi\,{\sf q}. 
\eeq{hat} 
From \Eqref{3.1_0}$_{1,2}$ it then follows that $\hat{\bsf V},\hat{\sf q}$ is a $T$-periodic solution to
\be\left.\ba{ll}\medskip
\partial_t\hat{\bsf V}=\nu\Delta\hat{\bsf{V}}-\nabla \hat{\sf q}+\bsf f\\
\Div\hat{\bsf{V}}={\sf h}\ea\right\}\ \ \mbox{in $\real^3\times(0,T)$}\,,
\eeq{4.50}
where
\be
\bsf f:=-\Delta\psi\,\bsf V-2\nabla\psi\cdot\nabla\bsf V+{\sf q}\,\nabla\psi\,;\ \ \ {\sf h}:=\nabla\psi\cdot\bsf V\,.
\eeq{4.51}
From \Eqref{AMPP} and the embedding $W^{1,2}\subset L^\infty$, 
we easily prove that
\be
\|\nabla\partial_t{\bsf V}\|_{L^\infty(L^2)}+\|\bsf V\|_{L^\infty(W^{2,2})}+\|{\sf q}\|_{L^\infty(L^6)}+\|\nabla{\sf q}\|_{L^\infty(L^2)}\le C\,, 
\eeq{car}
which, by the properties of $\psi$, implies
$$
\|\nabla\partial_t\hat{\bsf V}\|_{L^\infty(L^2)}+\|\nabla\bff\|_{L^\infty(L^2)}+\|D^2{\sf h}\|_{L^\infty(L^2)}\le C\,.
$$
Thus, \cite[Theorem IV.2.1]{Gab} applied to \Eqref{4.50}--\Eqref{4.51}, along with \Eqref{car} ensures that
\be
\|{\bsf V}\|_{L^\infty(W^{3,2}(\varOmega^{\rho_1}))}+\|\nabla {\sf q}\|_{L^\infty(W^{1,2}(\varOmega^{\rho_1}))}\le C\,.
\eeq{car1}
We next consider problem \Eqref{4.50}--\Eqref{4.51} with $\psi=\psi_{\rho_1,\rho_2}(|x|)$, with $\rho_2>\rho_1$. Thus, from \Eqref{car1}, classical embedding theorems 
 and the $T$-periodicity of $\bsf V$, we have 
\be
\supp(\bsf f)\subset {B_{\rho_2}}\,,\ \ \int_{B_{\rho_2}}\sup_{t\in\real} |\bsf f(y,t)|\le C\,.
\eeq{4.52}
In order to lift the divergence in \Eqref{4.50}, we introduce the field
\be
{\bsf V}_1(x,t):=\nabla \int_{\real^3}\cale(x-y){\sf h}(y,t)\,{\rm d} y\,,
\eeq{Vi1} 
where $\cale=\cale(\xi)$, $\bfxi\in\real^3\backslash\{\0\}$, is the Laplace fundamental solution. Clearly, $\Div\bsf V_1={\sf h}$.  Furthermore, as is well known,
\be
|D^k(\nabla\cale(\xi))|\le c\, |\xi|^{-m-2}\,,\ \ |k|=m\ge 0\,.
\eeq{LapEs}
By \Eqref{car},  embedding theorems and $T$-periodicity, we obtain ${\sf h}\in L^\infty({\varOmega}\times \real)$, and since  $\supp ({\sf h})\subset B_{\rho_2}$. Thus, from \Eqref{Vi1} it follows
\be
\sup_{(x,t)\in \varOmega^r\times[0,T]}(|x|^{m+2}|D^k\bsf V_1(x,t)|)\le C\,,\ \ |k|=m\ge0\,,
\eeq{estv}
for $r>2\rho_2$ and some constant $C>0$. 
Setting
\be
\bsf v=\hat{\bsf V}-\bsf V_1\,,
\eeq{4.54}
collecting \Eqref{4.50}--\Eqref{4.51}, it then follows
\be\left.\ba{ll}\medskip
\partial_t{\bsf v}=\Delta {\bsf{v}}-\nabla {\sf r}+\bsf f\\
\Div{\bsf{v}}=0\ea\right\}\ \ \mbox{in $\real^3\times(0,T)$}\,,
\eeq{4.55}
with 
\be
{\sf r}(x,t):=\hat{\sf q}(x,t)+\int_{\real^3}\cale (x-y)[\partial_t{\sf h}(y,t)-\nu\Delta{\sf h}(y,t)]\,{\rm d}y.
\eeq{4.56}
By \Eqref{hat}, \Eqref{estv} and \Eqref{4.54}, it is clear that  to prove the claimed spatial asymptotic property of $\bsf V$, it suffices to show the same for $\bsf v$. 
To establish the latter, we follow the approach introduced in \cite[p. 1249-1258]{GaZ}. This approach ensures the validity of the claimed property, provided the same  is satisfied, uniformly in $t\in\real$, by the solution, $\bsf u=\bsf u(x,t)$, say, to the {\em Cauchy problem} for \Eqref{4.55} corresponding to the initial condition $\bsf u(0)=\0$. As is well known this solution is given by
\be
{\sf u}_i(x,t)=\int_0^t\int_{\real^3}\varGamma_{ij}(x-y,s){\sf f}_{j}(y,t-s){\rm d}y\,{\rm d}s\,,\ \ i=1,2,3\,,
\eeq{4.57} 
where $\bfvGamma=\bfvGamma(\xi,\tau)$ is the Oseen fundamental solution to the Stokes problem \cite[Theorem VIII.4.2]{Gab}, for which  the following estimates hold:
\be
\int_0^\infty|D^k\bfvGamma(\xi,\tau)|{\rm d}\tau\le C_k|\xi|^{-m-1}\,,\ \ |k|=m\ge 0\,,\ \ \bfxi\neq\0\,;
\eeq{4.58}
see \cite[Lemma VIII.3.3 and Exercise VIII.3.1]{Gab}.
From \Eqref{4.52}, \Eqref{4.57} and and \Eqref{4.58} we deduce, as before,
$$
|x|^{m+1}|D^k\bsf u(x,t)|\le C\,,\ \ |k|=m\ge 0\,,\ \ |x|>2\rho_2\,, 
$$
which, by what we said above, implies the analogous property for $\bsf V$, that is,
\be
|x|^{m+1}|D^k\bsf V(x,t)|\le C\,,\ \ |k|=m\ge 0\,,\ \ |x|>2\rho_2\,. 
\eeq{4.60}
Next, from \Eqref{4.55}--\Eqref{4.56}, we deduce
\be
\hat{\sf q}(x,t)=-\int_{\real^3}\cale (x-y){\sf H}(y,t)\,{\rm d}y +\int_{\real^3}\nabla_y\cale(x-y){\bsf f}(y,t)\,{\rm d}y\,,
\eeq{hq}
with
$$
{\sf H}(y,t):=\partial_t{\sf h}(y,t)-\nu\Delta{\sf h}(y,t)\,.
$$
From \Eqref{car} and the fact that both $\bsf f$ and {\sf h} have compact support in $B_{\rho_2}$, we infer
\be
\|{\sf H}\|_{L^\infty(L^1)}+\|{\bsf f}\|_{L^\infty(L^1)}\le C\,.
\eeq{carro}
Thus, using \Eqref{carro} in \Eqref{hq} along with the estimates \Eqref{LapEs}, and recalling that $\hat{\sf q}={\sf q}$ in ${\varOmega}^{\rho_2}$, we can show
\be
\sup_{(x,t)\in \varOmega^r\times[0,T]}(|x|^{m+1}|D^k{\sf q}(x,t)|)\le C\,,\ \ |k|=m\ge0\,,\ \ r>2\rho_2\,.
\eeq{estq}
Since  $\bar{D^k\bsf V}(x)\equiv\0$, from the Poincar\'e inequality we get,  for all $x\in\varOmega$:
$$
\int_0^T|D^k\bsf V(x,t)|{\rm d}t\le T\int_0^T|\partial_tD^k\bsf V(x,t)|{\rm d}t\,,\ \ |k|\ge0\,. 
$$
which once combined with an elementary embedding inequality and $T$-periodicity, implies
\be
\sup_{s\in \real}|D^k\bsf V(x,s)|\le c\int_0^T|\partial_tD^k\bsf V(x,t)|{\rm d}t\,.
\eeq{4.63}
Now, from \Eqref{3.1_0}$_1$, we deduce
$$
\partial_t{\bsf V}=\Delta\bsf V-\nabla{\sf q}\,,\ \ \mbox{in $\varOmega^{2\rho_2}\times\real$}\,,
$$
so that, with the help of \Eqref{4.60} and \Eqref{estq}, it follows that
$$
|x|^{m+2}|\partial_tD^k\bsf V(x,t)|\le C\,,\ \ \ |x|>2\rho_2\,,\ \ |k|=m\ge 0\,.
$$
The latter, in conjunction with  \Eqref{4.63},  
implies the stated property on the spatial decay of $\bsf V$.\par\hfill$\square$\par

\setcounter{equation}{0}
\section{Further Properties of Solutions to Problem \Eqref{3.1_0}}
The main goal of this section is to find an explicit relation between the average velocity, $\bar{\bfzeta}$, and the data, $\bfu_*$, for the solutions to \Eqref{3.1_0} given in Proposition \ref{O_k}.
\Bl Let $\bfu_*\in\calv$. Then, the problem
\be\ba{cc}\medskip\left.\ba{ll}\medskip
\partial_t{\bfV}_0=\nu\Delta\bfV_0-\nabla{\sf p}_0\\
\Div\bfV_0=0\ea\right\}\ \ \mbox{in $\varOmega\times(0,T)$}\,,\\ \medskip
{\bfV_{0}}=\bfu_*+\bfzeta_0\,, \ \, \mbox{at $\varSigma\times (0,T)$}\,;\\ 
{\sf m}\dot{\bfzeta}_0(t)=-\Int{\varSigma}{}\mathbb T(\bfV_0,{\sf p}_0)\cdot\bfn{\rm d}\varSigma\,,\ \ t\in(0,T)\,.
\ea
\eeq{5.1}
has one and only one $T$-periodic solution
$$
(\bfV_0,{\sf p}_0,\bfzeta_0) \in\calw(\varOmega)\times \calq(\varOmega)\times W^{2,2}(0,T)\,,\ \ \bar{\bfV}_0=\bar{\sf p}_0=\bar{\bfzeta_0}=0.
$$
\EL{5.1}
{\em Proof.} From \cite[Lemma 6.2]{GaspU}, we know, in particular, that there exists a unique $T$-periodic solution  to \Eqref{AMP_9} such that, for all $q\in(1,2]$,
\be
(\bfV_0,{\sf p}_0,\bfzeta_0)\in [W^{1,q}(L^q(\varOmega))\cap L^q(W^{2,q}(\varOmega))]\times L^q(D^{1,q}(\varOmega)\cap L^{\frac{3q}{3-q}}(\varOmega))\times W^{1,q}(0,T)\,;\ \  \bar{\bfV}_0=\bar{\sf p}_0=\bar{\bfzeta_0}=0.
\eeq{5.2}
This solution satisfies also the estimate
\be
\|\bfV_0\|_{W^{1,q}(L^q(\varOmega))\cap L^q(W^{2,q}(\varOmega))}+\|{\sf p}_0\|_{L^q(D^{1,q}(\varOmega)\cap L^{\frac{3q}{3-q}}(\varOmega))}+\|\bfzeta_0\|_{W^{1,q}(0,T)}\le c\,\|\bfu_*\|_\calv\,.
\eeq{5.3}
Combining \Eqref{5.2}--\Eqref{5.3} with the regularization argument employed in the proof of \lemmref{Os} and the assumption on $\bfu_*$ leads directly to the property stated in the lemma.\par\hfill$\square$\par
\Bl Let $(\bfV,{\sf p},\bfzeta)$ and $(\bfV_0,{\sf p}_0,\bfzeta_0)$ be the solutions to \Eqref{3.1_0} and \Eqref{5.1} given in {\rm Proposition \ref{O_k}} and \lemmref{5.1}, respectively, and set
$$
\bfU:=\bfV-\bfV_0\,,\ \ p={\sf p}-{\sf p}_0\,,\ \ \bfz=\bfzeta-\bfzeta_0\,.
$$
Then, the following estimate holds
$$
\|\bar{\bfV}\|_{\cald^s(\varOmega)}+\|\bar{\sf p}\|_{\calp^s(\varOmega)}+|\bar{\bfzeta}|+\|\bfU\|_{\calw(\varOmega)}+\|{p}\|_{\calq(\varOmega)}+\|\bfxi\|_{W^{2,2}(0,T)}\le c\,\delta\,\|\bfu_*\|_\calv\,.
$$
\EL{5.2}
{\em Proof.} Setting
$$
g:=-\Div(\mathbb C^\top\cdot\bfV)\,,
$$ 
from \Eqref{3.1_0} and \Eqref{5.1}
we get
\be\ba{cc}\medskip\left.\ba{ll}\medskip
\partial_t{\bfU}=\nu\Delta\bfU-\nabla{p}+\nu\nabla g\\
\Div\bfU=g\ea\right\}\ \ \mbox{in $\varOmega\times(0,T)$}\,,\\ \medskip
\bfU=\bfz\,, \ \, \mbox{at $\varSigma\times (0,T)$}\,;\\ 
{\sf m}\,\dot{\bfz}(t)=-\Int{\varSigma}{}\mathbb T(\bfU,{p})\cdot\bfn{\rm d}\varSigma\,,\ \ t\in(0,T)\,.
\ea
\eeq{5.4}
Using \lemmref{bar}, \lemmref{Os},  arguing as in \Eqref{KTM}, \Eqref{AMP_16} and recalling that $\bar{\bfV}_0=\bar{\sf p}_0=\bar{\bfzeta_0}=0$, we show
$$
\|\bar{\bfV}\|_{\cald^s(\varOmega)}+\|\bar{\sf p}\|_{\calp^s(\varOmega)}+|\bar{\bfzeta}|+\|\bfU\|_{\calw(\varOmega)}+\|{p}\|_{\calq(\varOmega)}+\|\bfz\|_{W^{2,2}(0,T)}\le c\,\delta\left (\|\bar{\bfV}\|_{\cald^s(\varOmega)}+\|\bsf V\|_{\calw(\varOmega)}\right)\,.
$$
Thus, the desired property follows from this inequality and \Eqref{AMPP_0}.\par\hfill$\square$\par

We are now in a position to  show the main result of this section.
%
\Bp Let $(\bfV,{\sf p},\bfzeta)$ be the solution to \Eqref{3.1_0} determined in {\rm Proposition \ref{O_k}}. 
Then, the following relation holds:
\be
\bar{\bfzeta}=\mathbb M^{-1}\cdot\sum_{i=1}^3\big[\int_\varOmega p^{(i)}\Div(\bar{\mathbb C^\top\cdot\bfV})\big]\bfe_i\,, 
\eeq{5.5}
where $\mathbb M$ is the matrix defined in \Eqref{Matrix}.
This relation implies, in particular,
\be
\bar{\bfzeta}=\delta\, \mathbb M^{-1}\cdot  \bfpzc G_1+\bsf z(\delta,\bsf s)\,,
\eeq{5.6}
where
\be
\bfpzc G_1:=\sum_{i=1}^3\big[\int_\varOmega p^{(i)}\bar{\mathbb H(\bsf s):\nabla\bfV_0}\big]\bfe_i
\eeq{5.6_1}
with
\be
\mathbb H(\bsf s):= \Div\bsf s\,\mathbb I-(\nabla\bsf s)^\top\,, 
\eeq{Hh}
and
\be
|\bsf z(\delta,\bsf s)|\le c\,\delta^2\,.
\eeq{zizi}
\EP{5.1}
{\em Proof.} We observe that, once \Eqref{5.5} is proved, \Eqref{5.6} is immediately established  by using in \Eqref{5.5} the relation \Eqref{PId} along with the results proved in \lemmref{Def} and \lemmref{5.2}. To show  \Eqref{5.5}, we average \Eqref{3.1_0} to get
\be\ba{cc}\medskip\left.\ba{ll}\medskip
\Div\mathbb T(\bar{\bfV},\bar{\sf p})=\0\\
\Div\bar{\bfV}=-\Div(\bar{\mathbb C^\top\cdot\bfV})\ea\right\}\ \ \mbox{in $\varOmega$}\,,\\ \medskip
\bar{\bfV}=\bar{\bfzeta}\,, \ \, \mbox{at $\varSigma$}\,;\\ 
\Int{\varSigma}{}\mathbb T(\bar{\bfV},\bar{\sf p})\cdot\bfn{\rm d}\varSigma=\0\,.
\ea
\eeq{5.7}
Testing \Eqref{5.7}$_1$ with $\bfh^{(i)}$,  integrating by parts and using \Eqref{5.7}$_4$, we deduce
\be
\Int{\varOmega}{}\mathbb D(\bar{\bfV}):\mathbb D(\bfh^{(i)})=\0\,.
\eeq{5.8}
Likewise, testing \Eqref{4.32}$_1$ with $\bar{\bfV}$ and using \Eqref{5.7}$_{2,3}$ and \Eqref{5.8}, entails
$$
\bar{\bfzeta}\cdot\int_\varSigma \mathbb T(\bfh^{(i)},p^{(i)})\cdot\bfn=\int_\varOmega p^{(i)}\Div(\bar{\mathbb C^\top\cdot\bfV})\,,
$$
which completes the proof.
\par\hfill$\square$\par
\Br It is worth emphasizing that, in view of \Eqref{5.5}, \Eqref{5.6}, for the linearized problem \Eqref{3.1_0} the propulsion velocity $\bar{\bfzeta}$ depends only on the choice of the reference configuration, the physical parameters and  the displacement function $\bsf s(x,t)$. Moreover, $\bar{\bfzeta}\neq\0$ only if $\nabla\bsf s\not\equiv\0$. In other words, within this approximation, the body can self-propel {\em only} by changing its shape.
\ER{5.1}
\setcounter{equation}{0}
\section{On the Resolution of the Nonlinear Problem}
\label{Section:NonlinearProblem}
The starting point is  to rewrite \Eqref{2.12}--\Eqref{2.13} in an equivalent form, obtained by lifting the boundary data $\delta\,\bfu_*$ appropriately. To this end, we shall use the solution $(\bfV,{\sf p},\bfzeta)$ constructed in Proposition \ref{O_k}. More precisely, setting
\be
\bfu=\bsf u+\delta\,\bfV\,,\,\ p= {\sf q}+\delta\,{\sf p}\,,\,\ \bfgamma=\bfxi+\delta\,\bfzeta\,. 
\eeq{6.L}
the resolution of \Eqref{2.12}--\Eqref{2.13} becomes formally equivalent to that of the following problem in the unknowns $(\bsf u, {\sf q},\bfxi)$
\be\ba{cc}\medskip\left.\ba{ll}\medskip
\partial_t{\bsf u}+(\bsf u-\bfxi)\cdot\nabla\bsf u=\Div\mathbb T(\bsf u,{\sf q})+\bsf f_{\ell}(\bsf{s},\bsf u,{\sf q},\bfxi)+\bff_{n\ell}(\bsf s,\bsf u,\bfxi)+\delta^2\bff_{\mbox{\tiny $\bfV$}}\\
\Div\bsf u=-\Div(\mathbb C^\top\cdot\bsf u)\ea\right\}\,  \mbox{in $\varOmega\times\real$}\,,\\ \medskip
\bsf u=\bfxi\,, \ \ \mbox{at $\varSigma\times (0,T)$}\,;\ \ \Lim{|x|\to\infty}\bsf u(x,t)=\0\,,\ t\in\real\,;\\ 
{\sf m}\dot{\bfxi}(t)=-\Int{\varSigma}{}\mathbb T(\bsf u,{\sf q})\cdot\bfn{\rm d}\varSigma+\bsf F(\bsf{s},\bsf u, {{\sf q}})+\delta^2\bfF_{\mbox{\tiny $\bfV$}}\,,\ \ t\in\real\,,
\ea
\eeq{6.1}
where
\be\ba{rl}\medskip
\bsf f_{\ell}(\bsf{s},\bsf u,{\sf q},\bfxi)&\!\!\!\!\!:=\bff_{\ell}(\bsf{s},\bsf u,0)-\mathbb B\cdot\nabla{\sf q}-\delta(\bfV\!-\!\bfzeta)\cdot\mathbb A\cdot\nabla\bsf u-\delta\,(\bsf u-\bfxi)\cdot\mathbb A\cdot\nabla\bfV\\ \medskip
\bff_{\mbox{\tiny $\bfV$}}&\!\!\!\!\!:=\delta^{-1}\,\bff_\ell({\bsf s},\bfV,{ \sf p})-
(\bfV\!-\!\bfzeta)\cdot\mathbb A\cdot\nabla\bfV\,,\\ \medskip 
\bsf F(\bsf{s},\bsf u,{\sf q})&\!\!\!\!\!:=\bfF(\bsf{s},\bsf u,0)-\Int\varSigma{}{\sf q}\,\mathbb  C\cdot\bfn\,{\rm d}\varSigma\\
\bfF_{\mbox{\tiny $\bfV$}}&\!\!\!\!\!:=\delta^{-1}\bfF(\bsf{s},\bfV,{\sf p})
\,.
\ea
\eeq{6.2}
In this reformulation, $\bff_{\mbox{\tiny$\bfV$}}$ and $\bfF_{\mbox{\tiny$\bfV$}}$ are {\em given}  ``driving forces" acting on the liquid and on the body, respectively. 
\par
We shall address the study of the nonlinear problem  \Eqref{6.1}--\Eqref{6.2}, with a 
dual objective. Firstly, we show the existence of $T$-periodic solutions in an appropriate function class, and, successively, provide sufficient conditions ensuring that $\bar{\bfxi}\neq\0$. These results, combined with those established for the linearized problem in \propref{5.1}, will eventually allow us to provide a complete characterization of the averaged velocity of the center of mass at the order of $\delta^2$.
The accomplishment of the first objective will be the purpose of this and the following three sections, whereas the realization of the second one is postponed until Section \ref{asy}.  
\par
Our first goal is thus to prove that, for a given $(\bfV,{\sf p},\bfzeta)$, namely, by \propref{O_k}, for a given $\bfu_*\in\calv$, there is a corresponding  
$T$-periodic solution to \Eqref{6.1}--\Eqref{6.2}, provided $\delta$ is taken sufficiently ``small." 
\par
Precisely, let ${\sf A}$  be a subdomain of ${\varOmega}$,  $T>0$, and define 
$$
\mathscr W_T({\sf A}):=\left\{ \bfu\in W^{1,2}_T(L^6({\sf A})\cap D^{1,2}({\sf A})\cap D^{2,2}({\sf A}));\,\ \partial_t\bfu\in W^{1,2}_T(L^2({\sf A}))\right\}\,,
$$
with corresponding norm
$$
\|\bfu\|_{\mathscr W_T({\sf A})}:=\|\bfu\|_{W^{1,2}(0,T;L^6\cap D^{1,2}\cap D^{2,2})}+\|\partial_t\bfu\|_{W^{1,2}(0,T:L^2)}\,,
$$
and recall that the space $\calq(\varOmega)$ is defined in \Eqref{calw}. 
Our first goal is to prove the following result.
\Bt Suppose 
$$
\bfu_*\in \calv\,.
$$
Then, there exists $\delta_0>0$ such that, for all $\delta\in(0,\delta_0)$, problem \Eqref{2.12}--\Eqref{2.13}
has at least one corresponding $T$-periodic solution $(\bfu,p,\bfgamma)\in \mathscr W_T(\varOmega)\times\calq(\varOmega)\times W_T^{2,2}(\real^3)$. Moreover, this solution satisfies the estimate 
\be
\|\bsf u\|_{\mathscr W_T(\varOmega)}+\|{\sf q}\|_{\calq(\varOmega)}+\|\bfxi\|_{W^{2,2}(0,T)}\le c\,\delta^2\,,
\eeq{goest}
with $(\bsf u, {\sf q},\bfxi)$ defined in \Eqref{6.L}.
\label{theo:6.1}
\ET{7.1}
\par
The proof of this theorem will be achieved through a number of  intermediate steps that we are going to describe next.
\subsection*{Step 1: Solenoidal Formulation} \label{sub:7.1}
We begin to put \Eqref{6.1} in an equivalent form, by a suitable lifting of the divergence equation in \Eqref{6.1}$_2$ that enables us to reformulate \Eqref{6.1}-\Eqref{6.2} in terms of a {\em solenoidal} velocity field. To this end, we formally write
\be
\bsf u=\bfv +\bfB(\bfv)\,,
\eeq{6.8}
where $\bfB=\bfB(\bfv)$ is (formally) the operator $\bfB$ introduced in \propref{Boo} and, in order to simplify notation we set $\bfB(\bfv)\equiv \bfB(-\mathbb C^\top\cdot\bfv).$

Replacing \Eqref{6.8} in \Eqref{6.1} and taking into account the above definition of $\bfB(\bfv)$ we then show that $\bfv$ must solve the following set of equations
\be\ba{cc}\medskip\left.\ba{ll}\medskip
\partial_t{\bfv}+\partial_t\bfB(\bfv)+(\bfv+\bfB(\bfv)-\bfxi)\cdot\nabla\bfv+(\bfv-\bfxi)\cdot\nabla\bfB(\bfv)\\ 
\medskip
\hspace*{3.3cm}
=\nu\Delta\bfv-\nabla{\sf q}+\bfpzc{f}_\ell(\bsf{s},\bfv,{\sf q},\bfxi)+
\bfpzc{f}_{n\ell}(\bsf{s},\bfv,\bfxi)+\delta^2\bff_{\mbox{\tiny $\bfV$}}
\\
\Div\bfv=0\ea\right\}\,  \mbox{in $\varOmega\times\real$}\,,\\ \medskip
\bfv=\bfxi\,, \ \ \mbox{at $\varSigma\times \real$}\,;\ \ \Lim{|x|\to\infty}\bfv(x,t)=\0\,,\ t\in\real\,;\\ 
{\sf m}\dot{\bfxi}(t)=-\Int{\varSigma}{}\mathbb T(\bfv,{\sf q})\cdot\bfn{\rm d}\varSigma+ \bfpzc{F}(\bsf{s},\bfv,{\sf q})+\delta^2\bfF_{\mbox{\tiny $\bfV$}}\,,\ \ t\in\real\,,
\ea
\eeq{6.10}
where
\be\ba{rl}\medskip
\bfpzc{f}_\ell(\bsf s,\bfv,{\sf q},\bfxi)&\!\!\!\!:=\bsf f_\ell(\bsf s,\bfv+\bfB(\bfv),{\sf q},\bfxi)+\nu\big(\Delta\bfB(\bfv)+\nabla\Div\bfB(\bfv)\big)\,,\\
\bfpzc{f}_{n\ell}(\bsf s,\bfv,\bfxi)&\!\!\!\!:=-\bfB(\bfv)\cdot\mathbb A\cdot\nabla\bfB(\bfv)-(\bfv+\bfB(\bfv)-\bfxi)\cdot\mathbb B\cdot\nabla\bfv-(\bfv-\bfxi)\cdot\mathbb B\cdot\nabla\bfB(\bfv)\ea
\eeq{6.11}
and
\be
\bfpzc{F}(\bsf{s},\bfv,{\sf q}):=-2\nu\int_{\varSigma}\mathbb D(\bfB(\bfv))\cdot\bfn+\bsf F(\bsf s,\bfv+\bfB(\bfv),{\sf q})\,.
\eeq{6.12}
In view of the properties of the operator $\bfB$, it follows that proving Theorem \ref{theo:6.1} is equivalent to prove the results there stated for the solution $(\bfv,{\sf q},\bfxi)$ to \Eqref{6.10}--\Eqref{6.12}.  
\subsection*{Step 2: Invading Domains Method}
To show existence for problem \Eqref{6.10}--\Eqref{6.12} we utilize the so-called invading domains technique. Precisely,
let $\{\varOmega_{R_k}\}$ be  a sequence of bounded domains such that
\be
\bar{\varOmega_{R_k}}\subset \varOmega_{R_{k+1}}\,,\ \ \mbox{for all $k\in\nat$}\,;\ \ 
\cup_{k\in\nat}\varOmega_{R_k}=\varOmega\,.
\eeq{invdo} 
We then consider the sequence of ``approximating problems" (for simplicity: $\varOmega_k\equiv\varOmega_{R_k}$, $B_k\equiv B_{R_k}$)
\be\ba{cc}\medskip\left.\ba{ll}\medskip
\partial_t{\bfv}+\partial_t\bfB(\bfv)+(\bfv+\bfB(\bfv)-\bfxi)\cdot\nabla\bfv+(\bfv-\bfxi)\cdot\nabla\bfB(\bfv)\\ 
\medskip
\hspace*{3.3cm}
=\nu\Delta\bfv-\nabla{\sf q}+\bfpzc{f}_\ell(\bsf{s},\bfv,{\sf q},\bfxi)+
\bfpzc{f}_{n\ell}(\bsf{s},\bfv,\bfxi)+\delta^2\bff_{\mbox{\tiny $\bfV$}}
\\
\Div\bfv=0\ea\right\}\,  \mbox{in $\varOmega_k\times\real$}\,,\\ \medskip
\bfv=\bfxi\,, \ \ \mbox{at $\varSigma\times \real$}\,;\ \ \bfv(x,t)=\0\,,\ \ \mbox{at $\partial B_k\times \real$}\,;\\ 
{\sf m}\dot{\bfxi}(t)=-\Int{\varSigma}{}\mathbb T(\bfv,{\sf q})\cdot\bfn{\rm d}\varSigma+ \bfpzc{F}(\bsf{s},\bfv,{\sf q})+\delta^2\bfF_{\mbox{\tiny $\bfV$}}\,,\ \ t\in\real\,.\ea
\eeq{6.3}
\subsection*{Step 3:  Unique Solvability of a Linearized Problem} We introduce a suitable {\em linearized} version of \Eqref{6.3}; see \Eqref{6.10_lin}--\Eqref{flin}. For the latter, we  can find $\delta_*>0$ such that, for any given $\delta\in (0,\delta_*)$, there exists a corresponding unique $T$-periodic solution $(\bfv,{\sf q},\bfxi)$ in the class $\mathscr W_T({\varOmega}_k)\times\calq(\varOmega_k)\times W_T^{2,2}(\real^3)$, provided the coefficients of the linearized equation satisfy appropriate conditions; see Proposition \ref{EUlin}. It is   of the utmost importance that, in this process, we are able to ``control" the constants entering in the estimates and show that they are independent of $k$.  Galerkin method typically satisfies this requirement. However, adapting the method to the case at hand is not effortless for the following reasons. In the first place, $\bfpzc{f}_\ell$ contains first and second order derivatives of $\bfv$ that prevent us from obtaining a straightforward uniform bound on the (kinetic) energy norm of the solution that, in the ``classical" approach, is crucial to establish the existence of time-periodic solutions \cite{GaKyH}. To overcome this issue we shall then employ  a number of suitable ``high-order" energy estimates in conjunction with \lemmref{Def}. Furthermore, both $\bfpzc{f}_\ell$ and $\bfpzc{F}$ include terms depending on the pressure field ${\sf q}$, which may not coincide with the pressure field ${\sf q}_*$ (say) recovered {\em a posteriori} by the Galerkin method. We thus employ an appropriate perturbation argument that shows ${\sf q}\equiv {\sf q}_*$.
\subsection*{{Step 4}: Resolution of Problem \Eqref{6.3} and Proof of Theorem \ref{theo:6.1}}  We couple the results of Step 3 with Schauder fixed-point theorem and prove the existence of a $T$-periodic solution ${\sf s}_k:=(\bfv_k,{\sf q_k},\bfxi_k)$ to the {\em nonlinear} problem \Eqref{6.3} in the class $\mathscr W_T(\varOmega_k)\times\calq(\varOmega_k)\times W_T^{2,2}(\real^3)$, along with  suitable corresponding estimates; see \propref{8.1}.  
As observed in Step 3, the crucial point is to secure that these  estimates involve constants that are  independent of $k$. In which case, by a routine procedure, we can show that as $k\to\infty$ the sequence $\{{\sf s}_k\}$ will tend, in suitable topologies, to a $T$-periodic solution to the original problem \Eqref{6.1}--\Eqref{6.2}, satisfying all the properties stated in the theorem. 
\medskip\par
The strategy just described will be implemented in the following sections. Precisely, in Section \ref{linear} we will accomplish Step 3, whereas Sections \ref{sec:Ok} and \ref{nonl} will be dedicated to the achievement of Step 4.
\setcounter{equation}{0} 
\section{$T$-periodic Solutions to a Linearized Version of \Eqref{6.3}}\label{linear}
As mentioned earlier on, our first objective is to show well-posedness of the $T$-periodic problem in the space $\mathscr W_T(\varOmega_k)$ for a linearized version of \Eqref{6.10} that we define next. Let $\tilde{\bfv}$ and $\tilde{\bfxi}$ be prescribed vector field and vector function, respectively, whose regularity is specified later on, and consider the following problem
\be\ba{cc}\medskip\left.\ba{ll}\medskip
\partial_t{\bfv}+\partial_t\bfB(\bfv)+(\tilde{\bfv}+\bfB(\tilde{\bfv})-\tilde{\bfxi})\cdot\nabla\bfv+(\tilde{\bfv}-\tilde{\bfxi})\cdot\nabla\bfB(\bfv)\\ \medskip
\hspace*{2.9cm}=\nu\Delta\bfv-\nabla{\sf q}+\bfpzc{f}_\ell(\bsf{s},\bfv,{\sf q},\bfxi)+
{\bfpzc{f}}(\bsf{s},\bfv,\tilde{\bfv},\tilde{\bfxi})+\delta^2{ \bf f}
\\
\Div\bfv=0\ea\right\}\,  \mbox{in $\varOmega_k\times\real$}\,,\\ \medskip
\bfv=\bfxi\,, \ \ \mbox{at $\varSigma\times (0,T)$}\,;\ \ \bfv(x,t)=\0\,,\ \ \mbox{at $\partial B_k\times \real$}\,;\\ 
{\sf m}\dot{\bfxi}(t)=-\Int{\varSigma}{}\mathbb T(\bfv,{\sf q})\cdot\bfn{\rm d}\varSigma+ \bfpzc{F}(\bsf{s},\bfv,{\sf q})+\delta^2{\bf F}\,,\ \ t\in\real\,,
\ea
\eeq{6.10_lin}
where
\be
{\bfpzc{f}}(\bsf s,\bfv,\tilde{\bfv},\tilde{\bfxi}):=-\bfB(\tilde{\bfv})\cdot\mathbb A\cdot\nabla\bfB(\bfv)-(\tilde{\bfv}+\bfB(\tilde{\bfv})-\tilde{\bfxi})\cdot\mathbb B\cdot\nabla\bfv-(\tilde{\bfv}-\tilde{\bfxi})\cdot\mathbb B\cdot\nabla\bfB(\bfv)\,,
\eeq{flin}
and ${\bf f}$, ${\bf F}$  are  given $T$-periodic functions to be specified later on.
\smallskip\par
For  $\epsilon_0,M>0$ we define
\be\ba{ll}\mathscr S_{\epsilon_0,M}:=\Big\{\tilde{\bfv}\in  W^{1,4}_{\rm loc}(\real;\cald^{1,2}(B_k)):\ 
\tilde{\bfv} \ \mbox{is $T$-periodic with $\tilde{\bfxi}=\tilde{\bfv}|_{\Omega_0}$};\\
\hspace*{5.5cm} \essup{t\in[0,T]}\|\mathbb D(\tilde{\bfv}(t))\|_2\le\epsilon_0\,;\ \ \|\tilde{\bfv}\|_{W^{1,4}(0,T;\cald^{1,2})}\le M\Big\} \,.
\ea
\eeq{Sepm}
Moreover, set
\be
\calf_1:=\|{\bf f}\|_2^2+\|{\bf f}\|^2_\frac65+|{\bf F}|^2\,,\ \ \ \calf_2:=\|\partial_t{\bf f}\|_2^2+|\dot{\bf F}|^2\,.
\eeq{cal}
and
\be
F:=\Big[\Int0T\big({\calf}_1+\calf_2+\calf_1^2\big){\rm d}t\Big]^\frac12\,,\ \ \mu:=M^2+1.
\eeq{fmu}
\par
Objective of this section is to show the following result.
\begin{proposition}{\sl Let $\tilde{\bfv}\in \mathscr S_{\epsilon_0,M}$, and suppose ${\bf f}$, ${\bf F}$ are given $T$-periodic functions with $F$ finite. Then,    
there are constants $C_i=C_i(\Omega,\nu,{\sf m},T)>0$, $i=0,1,2$,  such that, setting
\be
\delta_*:=\min\Big\{(C_1\mu)^{-1},(C_1 F)^{-\frac12},C_2\Big\} 
\eeq{assu}
the following properties hold. If $\varepsilon_0\le C_0$, then, for any $\delta\in(0,\delta_*)$ there exists one and only one  $T$-periodic solution 
$$(
\bfv,{\sf q},\bfxi)\in\mathscr W_T(\varOmega_k)\times\calq({\varOmega}_k)\times W^{2,2}(0,T)
$$ 
to problem \Eqref{6.10_lin}--\Eqref{flin} corresponding to the given ${\bf f},{\bf F}$. Moreover, the solution satisfies the following estimates}
\be\ba{ll}\medskip
\Max{t\in [0,T]}\big(\|\bfv(t)\|_6^2+|\bfxi(t)|^2\!+\!\|\mathbb D(\bfv(t))\|_2^2\big)+\Int0T\big(\|\mathbb D(\bfv)\|_2^2+\|\bfv_{t}\|_2^2+|\dot{\bfxi}|^2+\|D^2\bfv\|_2^2+\|\nabla{\sf q}\|_2^2\big){\rm d}t\\ \hspace*{5.5cm}\le 
c_0\,\delta^4\Int0T\calf_1(t)\,{\rm d}t,
\ea
\eeq{BaNa1}
{\sl and}
\be
\ba{ll}\medskip
\Max{t\in [0,T]}\big(\|\bfv_{t}(t)\|_2^2+|\dot{\bfxi}|^2\!+\!\|\mathbb D(\bfv_{t}(t))\|_2^2\big)+\Int0T\big(\|\mathbb D({\bfv}_{t})\|_2^2+\|\bfv_{tt}\|_2^2+|\ddot{\bfxi}|^2+\|D^2\bfv_{t}\|_2^2\big){\rm d}t
\\
\hspace*{5.5cm}
\le c_0\delta^4\Big[\Int0T\calf_1{\rm d}t\Int0T\|\mathbb D(\tilde{\bfv}_t)\|_2^4{\rm d}t+F^2\Big]\,,
\ea 
\eeq{BaNa2}
{\sl where  the constant $c_0$ and is independent of $R_k$.}   
\label{EUlin}
\end{proposition}
\par
We postpone the full proof of Proposition \ref{EUlin}  until Subsection \ref{propo}. 
As mentioned earlier on, we shall employ  Galerkin method applied to a suitable weak formulation of the problem. However, using this method {\em directly} on  \Eqref{6.10_lin}--\Eqref{flin} does not seem to be feasible, due to the presence of the pressure field terms $-\mathbb B\cdot\nabla{\sf q}$  in the function  ${\bfpzc{f}}_\ell$  defined in \Eqref{6.11}$_1$, and $-\int_\varSigma{\sf q}\bfn\,{\rm d}\varSigma$ in the function ${\bfpzc{F}}$ given in \Eqref{6.12}. In fact, as is well known, when one uses Galerkin method the existence of {\em some}  pressure field $p$ (say) can be obtained only a posteriori, once $\bfv$ and $\bfxi$ are determined and, in principle, it is not said that  $p\equiv {\sf q}$.  In order to overcome this difficulty, we proceed as follows. We prescribe ${\sf q}^\sharp$ as a $T$-periodic function in $\calq({\varOmega}_k)$ and show the existence of a corresponding solution $(\bfv,{\sf q},\bfxi)$ to \Eqref{6.10_lin}--\Eqref{flin} in the class $\mathscr W_T({\varOmega}_k)$ with ${\bfpzc{f}}_\ell$ and ${\bfpzc{F}}$ replaced by  
\be
\bfpzc{f}_\ell^\sharp:=\bfpzc{f}_\ell(\bsf{s},\bfv,{\sf q}^\sharp,\bfxi)\,,\ \ \bfpzc{F}^\sharp:=\bfpzc{F}(\bsf{s},\bfv,{\sf q}^\sharp)\,,
\eeq{fsha}
respectively.
Then, thanks to \lemmref{Def}, by a  perturbative argument we show ${\sf q}\equiv{\sf q}^\sharp$.

\subsection{Weak Formulation} \label{Section:WeakFormulation}
We begin to furnish a
weak formulation of \Eqref{6.10_lin}--\Eqref{flin} --with $\bfpzc{f}_\ell$ and $\bfpzc{F}$  replaced by $\bfpzc{f}_\ell^\sharp$ and $\bfpzc{F}^\sharp$, respectively,  given in \Eqref{fsha}-- that is amenable to  Galerkin method.
For this, we introduce the class   $\calc_T(B_k)$ of suitable ``test functions,"
constituted by the restriction to $[0,T]$ of functions $\bfphi \in C^{1}(B_k\times \mathbb{R})$, satisfying: 
\begin{itemize}
\item[{\rm (a)}] $\Div  \bfphi(x,t) = 0$ for $(x,t)\in B_k \times \mathbb{R}$\,; 
\item[{\rm (b)}] $\bfphi (x,t) = \hat\bfphi (t)$, some $\hat\bfphi \in C^{1}(\mathbb{R})$,   for $x$ in a neighborhood 
of $\Omega_0$ and $t \in \mathbb{R}$\,;
\item[{\rm (c)}] $\supp_x{\bfphi(x,t)}\subset B_k$ for all $t \in\real$\,;
\item[{\rm (d)}] $\bfphi(x,t+T)=\bfphi(x,t)$ for all $(x,t) \in B_k\times\mathbb{R}$\,.
\end{itemize}
\par
Testing \Eqref{6.10}$_1$ with $\bfpzc{f}_\ell\equiv \bfpzc{f}_\ell^\sharp$ by arbitrary $\bfphi\in \calc_T(B_k)$, integrating by parts over $\varOmega_k$ and employing \Eqref{6.10_lin}$_{2-4}$, we infer 
\be\ba{ll}\medskip
\langle\partial_t\bfv,\bfphi\rangle\!+\!(\partial_t\bfB(\bfv),\bfphi)+\big((\tilde{\bfv}+\bfB(\tilde{\bfv})-\tilde{ \bfxi})\cdot\nabla\bfv+(\tilde{\bfv}-\tilde{\bfxi})\cdot\nabla\bfB(\bfv)+ 2\nu(\mathbb D(\bfv),\mathbb D ({\bfphi}))\\
\hspace*{7.51cm}=(\bfpzc{f}_\ell^\sharp+\bfpzc{f}+\delta^2{\bf f},\bfphi)+\big(\bfpzc{F}^\sharp+\delta^2{\bf F}\big)\cdot  \hat\bfphi\,, 
\ea
\eeq{6.13}
where,  $(\cdot,\cdot)\equiv (\cdot,\cdot)_{\varOmega_k}$ is the $L^2(\varOmega)$-scalar product, and  
$$
\langle \bfphi,\bfpsi\rangle:=
{\sf m}\,\hat{\bfphi}\cdot\hat{\bfpsi}+(\bfphi,\bfpsi).
$$

\subsection{A Special Basis}
In order to apply  Galerkin method to our case, we need a suitable base whose properties are presented in the following lemma.
\Bl 
For any fixed $k\in\nat$, the problem  
\be
\begin{array}{l}
\left. 
\begin{array}{l}
\hspace{-0.2cm} \displaystyle  - \Div \mathbb T(\bfpsi,\phi)=\mu  \,\bfpsi \medskip \\ 
\hspace{-0.2cm} \displaystyle \Div   \bfpsi = 0 \medskip
\end{array}
\right\} \mbox{ in } \varOmega_k \,,\\  \medskip
\displaystyle \bfpsi  = \hat{\bfpsi} \ \   \mbox{ in } \Omega_0\,, \ \  
\displaystyle \bfpsi= 0\ \   \mbox{ at } \partial B_{R_k} \,, \\ 
\displaystyle \mu\,  \hat{\bfpsi}={\sf m}^{-1}\int_{\varSigma} \mathbb T(\bfpsi,\phi)\cdot \bfn\,{\rm d}\varSigma \,,
\end{array}
\eeq{3.4_0}
admits a denumerable number of positive eigenvalues $\{\mu_{i}\}$ clustering at infinity, 
and  corresponding eigenfunctions $\{\bfpsi_{i}\} \subset {\mathcal D}^{1,2}(B_{R_k})\cap W^{2,2}({\varOmega}_k)$    forming an orthonormal basis of ${\mathcal{H}}(B_{R_k})$ that is also orthogonal in $\cald^{1,2}({\varOmega}_k)$.  Precisely, 
\be
\langle\bfpsi_{j},\bfpsi_{i}\rangle:=(\bfpsi_{j},\bfpsi_{i})_{{\varOmega}_k}+{\sf m}\hat{\bfpsi}_{j}\cdot\hat{\bfpsi}_{i}=\delta_{ji}\,,
\eeq{ortg}
and
\be
\big(\mathbb D(\bfpsi_j),\mathbb D(\bfpsi_i)\big)_{{\varOmega}_k}=\mu_i\langle\bfpsi_{j},\bfpsi_{i}\rangle\,.
\eeq{poi}
Furthermore, the correspondent ``pressure" fields satisfy $\phi_{i}\in W^{1,2}({\varOmega}_k)$, $i\in\nat$.
Finally, let
$$
\bfPsi:=\sum_{i=1}^Nc_i\bfpsi_i\,,\ \ \Phi:=\sum_{i=1}^Nc_i\phi_i\,,\ \ \hat{\bfPsi}=\sum_{i=1}^n\hat{\bfpsi}_i\,,\ \ c_i\in\real\,, \ i=1,\ldots,N\,.
$$
Then,
there is a positive constant $c_0$ independent of $k\in\nat$ such that
\be
\|D^2\bfPsi\|_2+\|\nabla\Phi\|_2\le c_0\left(\|\Div\mathbb T(\bfPsi,\Phi)\|_2+\|\nabla\bfPsi\|_2\right)\,.
\eeq{hey}
\EL{Base}
{\em Proof.} The statements regarding the eigenvalue problem are proved in 
 \cite[Theorem 3.1]{GaSi1}. Concerning the validity of \Eqref{hey}, we notice that $(\bfPsi,\Phi)$ can be formally viewed as a solution to the following Stokes problem
\be
\begin{array}{l}
\left. 
\begin{array}{l}
\hspace{-0.2cm} \displaystyle   \Delta\bfPsi=\nabla\Phi+\bff  \medskip \\ 
\hspace{-0.2cm} \displaystyle \Div   \bfPsi = 0 \medskip
\end{array}
\right\} \mbox{ in } \varOmega_k \,,\\  
\displaystyle \bfPsi  = \hat{\bfPsi} \ \   \mbox{ at } \varSigma\,, \ \  
\displaystyle \bfPsi= 0\ \   \mbox{ at } \partial B_{R_k} \,, 
\end{array}
\eeq{3.4_00}
with $\bff:=\Div\mathbb T(\bfPsi,\Phi)$. 
Let
$$
\bfz(x):=x_3\hat{\Psi}_{2}\bfe_1+x_1\hat{\Psi}_{3}\bfe_2+x_2\hat{\Psi}_{1}\bfe_3
$$
and define
$$
\bfH(x):=\curl(\zeta(x)\bfz(x))\,,
$$
where $\zeta$ is a smooth function of bounded support in $\varOmega_1$, that is 1 in a neighborhood of $\varSigma$ and 0 away from $\varSigma$.
Since $\curl\bfz=\hat{\bfPsi}$, we get 
$$
\bfH(x)=\zeta(x)\hat{\bfPsi}-\bfz(x)\times\nabla\zeta(x)\,.
$$
It is at once established that $\bfH$ satisfies the following conditions
\be
\bfH=\hat{\bfPsi}\ \ \mbox{at $\varSigma$}\,;\ \ \ \|\bfH\|_{2,2}\le c_1\,|\hat{\bfPsi}|\le c_2\,\|\nabla\bfPsi\|_2\,,
\eeq{Psi}
with $c_2$ independent of $k\in\nat$, and where, in the last inequality, we have used \Eqref{1.10}, \Eqref{1.8} and \remref{2.1}. Setting $\bfpsi:=\bfPsi-\bfH$, from \Eqref{3.4_00} we get
$$
\begin{array}{l}\medskip
\left. 
\begin{array}{l}
\hspace{-0.2cm} \displaystyle   \Delta\bfpsi=\nabla\Phi+\bfg  \medskip \\ 
\hspace{-0.2cm} \displaystyle \Div   \bfpsi = 0 \medskip
\end{array}
\right\} \mbox{ in } \varOmega_k \,,\\  
\displaystyle \bfpsi=\0    \mbox{ at } \varSigma\cup \partial B_{R_k} \,, 
\end{array}
$$
with $\bfg:=\Div\mathbb T(\bfPsi,\Phi)-\Delta\bfH$. By \cite[Lemma 1]{Hey}, we then deduce that there is a constant $c$ independent of $k$ (and $N$) such that
$$
\|D^2\bfpsi\|_2\le c\,\left(\|\bfg\|_2+\|\nabla\bfpsi\|_2\right)\,.
$$
As a result, \Eqref{hey} follows from the latter and \Eqref{Psi}. \par\hfill$\square$
\par

\subsection{Approximated Solutions}
Let ${\sf q}^\sharp$ be a given function in $\calq({\varOmega}_k)$. 
We look for an approximated solution $(\bfv_N(x,t), \bfxi_N(t))$ to \Eqref{6.13} of the form
\be
\bfv_N(x,t)= \sum_{j=1}^Nc_{jN}(t)\bfpsi_{i}(x)\,,\ \  
\bfxi_N(t)= \sum_{j=1}^Nc_{iN}(t)\hat{\bfpsi}_{j}\,, 
\eeq{6.14}
where $\{\bfpsi_j\}$ is the basis introduced in \lemmref{Base},  and  the vector function $\bfc_N(t):=\{c_{1N}(t),\ldots c_{NN}(t)\}$ satisfies the following system of  linear equations, $i=1,\ldots,N$\,, 
\be\ba{ll}\medskip
\ode{}t\langle\bfv_N,\bfpsi_i\rangle+\ode{}{t}(\bfB(\bfv_N),\bfpsi_i)+\big((\tilde{\bfv}+\bfB(\tilde{\bfv})-\tilde{\bfxi})\cdot\nabla\bfv_N+(\tilde{\bfv}-\tilde{\bfxi})\cdot\nabla\bfB(\bfv_N)
+ 2\nu(\mathbb D(\bfv_N),\mathbb D ({\bfpsi_i}))
\\ 
\hspace*{.9cm}-\big(\bfpzc{f}_\ell(\bsf s,\bfv_N,{\sf q}^\sharp,\bfxi_N)+(\bfpzc{f}(\bsf s,\bfv_N,\tilde{\bfv},\tilde{\bfxi})+\delta^2{\bf f},\bfpsi_i\big)-\big(\bfpzc{F}(\bsf s,\bfv_N,{\sf q}^\sharp)+\delta^2{\bf F}\big)\cdot  \hat\bfpsi_i=0\,,
\ea
\eeq{6.15}
where
$(\cdot,\cdot)\equiv(\cdot,\cdot)_{{\varOmega}_k}$,  $\langle\cdot,\cdot\rangle\equiv\langle\cdot,\cdot\rangle_{B_k}$,   $\bfB$ is the operator introduced in \propref{Boo} and, in order to simplify notation, we set 
\be
\bfB(\tilde{\bfv})\equiv \bfB({-}\mathbb C^\top\cdot\tilde{\bfv})
\,,\ \  \bfB(\bfv_N)\equiv \bfB({-}\mathbb C^\top\cdot\bfv_N).
\eeq{6.16} 
These quantities are well defined. In fact, since
$$
\bfB(\bfv_N)=\Sum{j=1}N c_{jN}\bfB({-}\mathbb C^\top\cdot\bfpsi_j)\,,
$$
we at once obtain, from \lemmref{Def} and \lemmref{Base}, that $\bfg:={-}\mathbb C^\top\cdot\bfpsi_k\in W_T^{2,2}(W^{2,2}(\varOmega))$  and $\bfg\equiv\0$ for all $x$ in the exterior of $B_0$. Moreover, for each $k=1,\ldots,N$,
$$
\int_\varSigma\bfg\cdot\bfn
={-}\hat{\bfpsi}_{{k}}\cdot\int_\varSigma\bfn\cdot\mathbb C^\top=\hat{\bfpsi}_{{k}}\cdot\int_{\Omega_0}\Div\mathbb C^\top=0\,,
$$
where, in the last step, we have used \Eqref{PId}. Therefore, $\bfg\in\cals^{2,2}$ and \propref{Boo} applies. A similar proof shows that $\mathbb C^\top\cdot\tilde{\bfv}\in \cals^{1,2}$ so that $\bfB(\tilde{\bfv})$ is meaningful as well. 
Let us now prove that \Eqref{6.15} is a system of first order differential equations in normal form in the unknowns $\bfc_N$. To this end, we begin to observe that, in view of the linearity of $\bfB$,
$$
\ode{}t(\bfB(\bfv_N),\bfpsi_i)=\ode{}t\sum_{j=1}^nc_{jN}\left(\bfB({-}\mathbb C^\top\cdot\bfpsi_j),\bfpsi_i\right)=\sum_{j=1}^n\dot{c}_{jN}\left(\bfB({-}\mathbb C^\top\cdot\bfpsi_j),\bfpsi_i\right)+\left(\bfB({-}\partial_t\mathbb C^\top\cdot\bfv_N),\bfpsi_i\right)\,.
$$
As a consequence, we deduce that \Eqref{6.15} can be rewritten in the following form 
\be
\dot{c}_{jN}{\sf K}_{ji}= {\sf G}_i(\bfc_N,t)\,,\ \ i=1,\ldots,N\,,
\eeq{6.17}
where
$$
{\sf K}_{ji}:=\langle\bfpsi_j,\bfpsi_i\rangle+\left(\bfB({-}\mathbb C^\top\cdot\bfpsi_j),\bfpsi_i\right)\,,\ j,i=1,\ldots,N\,,
$$
and the function ${\sf G}$ is  linear in the the components $c_{jN}$, and $T$-periodic. Let us show that the matrix ${\sf K}$ is invertible for ``small" $\delta$. In fact, setting $\bfzeta=\zeta_m\bfpsi_m$, $\hat{\bfzeta}=\zeta_m\hat{\bfpsi}_m$,   and recalling that
$$
\langle\bfpsi_{j},\bfpsi_{i}\rangle:=(\bfpsi_{j},\bfpsi_{i})+{\sf m}\hat{\bfpsi}_{j}\cdot\hat{\bfpsi}_{i}=\delta_{ji}\,,
$$
we get
\be
\zeta_j{\sf K}_{ji}\zeta_i={\sf m}|\hat{\bfzeta}|^2+\|\bfzeta\|_2^2+(\bfB({-}\mathbb C^\top\cdot\bfzeta),\bfzeta)\,.
\eeq{6.18}
By \Eqref{Trunz}$_2$ with $k=0$, and \lemmref{Def} we infer
$$
\|\bfB({-}\mathbb C^\top\cdot\bfzeta)\|_2\le c\,\delta \left(\|\bfzeta\|_2+|\hat{\bfzeta}|\right)
$$
which combined with \Eqref{6.18} shows that there is $\delta_0>0$ such that ${\sf K}$ is invertible for any $\delta\in(0,\delta_0)$. This fact, in conjunction with the property that ${\sf G}$ is a continuous $T$-periodic function,   ensures that, corresponding to given initial conditions $\bfc_N(0)$, the linear system of differential equations \Eqref{6.17} has a unique solution in the time interval~$[0,\infty)$. 
\subsection{Uniform Estimates and Existence of Approximating $T$-periodic Solutions}
We begin to derive three basic ``energy equations," which are  formally obtained by testing \Eqref{6.10_lin}$_1$ with $\bfpzc{f}_\ell\equiv \bfpzc{f}_\ell^\sharp$ and $\bfpzc{F}\equiv \bfpzc{F}^\sharp$,  in the order, by $\bfv$, $\partial_t\bfv$ and $\Div\mathbb T(\bfv,{\sf q})$ and then integrating by parts.  
To make the above argument rigorous, we multiply both sides of \Eqref{6.15}$_1$ one time by by $c_{iN}$, the second time by $\dot{c}_{iN}$, the third time by $\mu_ic_{iN}$,  and sum over $i$. Thus, setting
\be
\bfH(\tilde{\bfv},\bfv):= -\big[\bfB_t(\bfv)+(\tilde{\bfv}-\tilde{\bfxi})\cdot\nabla\bfB(\bfv)+(\tilde{\bfv}+\bfB(\tilde{\bfv})-\tilde{\bfxi})\cdot\nabla\bfv\big]
\eeq{H}
we get, also with the help of \lemmref{Base}, the following three  relations (we omit the subscript $N$ and set $\partial_t(\cdot)=(\cdot)_t$ for simplicity)
\be\ba{ll}\medskip
\half\ode{}t\left(\|\bfv\|_2^2+|\bfxi|^2\right)+2\|\mathbb D(\bfv)\|_2^2=\!(\bfH(\tilde{\bfv},\bfv),\bfv)+\!(\bfpzc{f}_\ell^\sharp\!+\!\bfpzc{f}\!+\!\delta^2{\bf f},\bfv)\!+\!(\bfpzc{F}^\sharp+\delta^2{\bf F})\cdot\bfxi\,,\\ \medskip
\ode{}t\|\mathbb D(\bfv)\|_2^2+\|\bfv_t\|_2^2+|\dot{\bfxi}|^2=(\bfH(\tilde{\bfv},\bfv),\bfv_t\big)+(\bfpzc{f}_\ell^\sharp\!+\!\bfpzc{f}\!+\!\delta^2{\bf f},\bfv_t)+(\bfpzc{F}^\sharp+\delta^2{\bf F})\cdot\dot{\bfxi}\,\\ \smallskip
\ode{}t\|\mathbb D(\bfv)\|_2^2+\|\Div\mathbb T(\bfv,{\sf r})\|_2^2\!+\!|\bsf{T}|^2
=\big(\bfH(\tilde{\bfv},\bfv),\Div\mathbb T(\bfv,{\sf r})\big)
\\ \hspace*{5.6cm}
+(\bfpzc{f}_\ell^\sharp+\bfpzc{f}+\!\delta^2{\bf f},\Div\mathbb T(\bfv,{\sf r}))+\big(\bfpzc{F}^\sharp+\!\delta^2{\bf F}\big)\cdot{\bsf T}\,,
\ea
\eeq{6.19}
where
$$ 
{\sf r}:= \Sum{j=1}Nc_{jN}(t)\phi_{j}(x)\,,\ \
\bsf{T}:=\int_\varSigma\mathbb T(\bfv,{\sf r})\cdot\bfn{\rm d}\varSigma\,, 
$$
with $\{\phi_j\}$ ``pressures" given in \lemmref{Base}, and, again for simplicity since their actual value is irrelevant to the proof, we set ${\sf m}=\nu=1$.\smallskip
\par
We now estimate  the right-hand sides of \Eqref{6.19} suitably. To this end, 
in Appendix B the following inequalities are showed
\be\ba{ll}\medskip
|\big(\bfH(\tilde{\bfv},\bfv),\bfv\big)|\le c\,\delta\,\big[\|\mathbb D(\bfv)\|_2^2+(\|\mathbb D(\tilde{\bfv})\|_2+\|\mathbb D(\tilde{\bfv})\|_2^{\frac43})\|\mathbb D(\bfv)\|_2^{2}+\|\bfv_t\|_2^2+|\dot{\bfxi}|^2+\|\Div T(\bfv,{\sf r})\|_2^2\big]
\,,\\ 
\|\bfH(\tilde{\bfv},\bfv)\|_2^2\le c\big[\delta^2(\|\bfv_t\|_2^2+|\dot{\bfxi}|^2)\!+\!(1\!+\!\delta^2)(\|\mathbb D(\tilde{\bfv})\|_2^2+\|\mathbb D(\tilde{\bfv})\|_2^4)\|\mathbb D({\bfv})\|_2^2+\delta^2\|\mathbb D(\bfv)\|_2^2\big]
+\mbox{$\eta$}\|\Div\mathbb T(\bfv,{\sf r})\|_2^2
\ea
\eeq{6.20f}
and also 
\be\ba{rl}\medskip
|(\bfpzc f_\ell^\sharp,\bfa)|\equiv |(\bfpzc f_\ell^\sharp,\bfa)_{B_0}|&\!\!\!\le c \,\delta\,\big(\|\mathbb D(\bfv)\|_2^2+\|\Div\mathbb T(\bfv,{\sf r})\|_2^2+\delta\,\|\nabla{\sf q}^\sharp\|_2^2\big)+\mbox{$\eta$}\|\bfa\|_{2,B_0}^2
\,,
\\ \medskip
|(\bfpzc f,\bfa)|\equiv|(\bfpzc f,\bfa)_{B_0}| &\!\!\!\le c\,\delta\left[(\|\mathbb D(\tilde{\bfv})\|_2^2+\|\mathbb D(\tilde{\bfv})\|_2^4)\|\mathbb D(\bfv)\|_2^2+\|\Div\mathbb T(\bfv,{\sf r})\|_2^2\right]+\eta\|\bfa\|_{2,B_0}^2\,,
\\
|\bfpzc F^\sharp\cdot\bfb|&\!\!\!\le c\,\delta\left(\|\mathbb D(\bfv)\|_2^2+\|\Div\mathbb T(\bfv,{\sf r})\|_2^2+\delta\,\|\nabla{\sf q}^\sharp\|_2^2\right)+\mbox{$\eta$}|\bfb|^2
\,.
\ea
\eeq{6.21f}
where $\eta>0$ is arbitrary and $c$ is independent of $R_k$.  Finally, we observe that, by Lemma \ref{le:calf}, H\"older inequality  and \Eqref{1.9} we have
\be\ba{ll}\medskip
\delta^2|({\bf f},\bfv)|\le \delta^2\|{\bf f}\|_{\frac65}\|\bfv\|_6\le C\delta^4\calf_1+\frac{\eta}{\kappa_0}\,\|\bfv\|_6^2\le  C\delta^4\calf_1+\eta\,\|\mathbb D(\bfv)\|_2^2\,,\\ \medskip 
\delta^2|({\bf f},\bfv_t)|\le \|{\bf f}\|_2\|\bfv_t\|_2\le C\delta^4\calf_1+\eta\,\|\bfv_t\|_2^2\\
\delta^2\big|\big({\bf f},\Div\mathbb T(\bfv,{\sf r})\big)\big|\le\delta^2 \|{\bf f}\|_2\|\Div\mathbb T(\bfv,{\sf r})\|_2\le C\delta^4\calf_1+\eta\|\Div\mathbb T(\bfv,{\sf r})\|_2^2\,,
\ea
\eeq{727}
with $C$ independent of $R_k$.
\par
Choosing $\bfa\equiv\bfv$ in \Eqref{6.21f}$_1$, $\bfb\equiv\bfxi$ in \Eqref{6.21f}$_2$, and employing \Eqref{6.21}, \Eqref{6.22} in the appendix,  \Eqref{6.20f}$_1$ and \Eqref{727}$_1$, by taking $\delta$ and $\epsilon_0$ below a constant independent of $R_k$, from \Eqref{6.19}$_1$ it follows that 
\be\ba{ll}\medskip
\ode{}t\left(\|\bfv\|_2^2+|\bfxi|^2\right)+c_1\,\|\mathbb D(\bfv)\|_2^2\le 
c\,\delta\left(
\|\bfv_t\|_2^2+|\dot{\bfxi}|^2+\|\Div\mathbb T(\bfv,{\sf r})\|_2^2+\delta\|\nabla{\sf q}^\sharp\|_2^2\right)+C\,\delta^4\calf_1\,.
\ea
\eeq{6.39}
In the same way, taking $\bfa\equiv\bfv_t$ in \Eqref{6.21f}$_1$, $\bfb\equiv\dot{\bfxi}$ in \Eqref{6.21f}$_2$, and choosing $\delta,\eta$ and $\epsilon_0$  below a constant, with the help of \Eqref{6.20f}$_2$, \Eqref{727}$_2$, from \Eqref{6.19}$_2$ we show
\be
\ode{}t\|\mathbb D(\bfv)\|_2^2+c_1(\|\bfv_t\|_2^2+|\dot{\bfxi}|^2)\le c_2\,(\delta+\epsilon_0)\|\mathbb D(\bfv)\|_2^2+(\delta+\eta)\|\Div\mathbb T(\bfv,{\sf r})\|_2^2+c_3\,\delta^2\,\|\nabla{\sf q}^\sharp\|_2^2+C\delta^4\calf_1\,.
\eeq{naju3}
Likewise, taking $\bfa\equiv\Div\mathbb T(\bfv,{\sf r})$ in \Eqref{6.21f}$_1$, $\bfb\equiv{\sf T}$ in \Eqref{6.21f}$_2$  and using  \Eqref{6.20f}$_2$, \Eqref{727}$_3$, from \Eqref{6.19}$_3$ we infer
$$
\ode{}t\|\mathbb D(\bfv)\|_2^2+c_3(\|\Div\mathbb T(\bfv,{\sf r})\|_2^2+|\bsf{T}|^2)\le c\,(\delta+\epsilon_0)\|\mathbb D(\bfv)\|_2^2+\delta(\|\bfv_t\|_2^2+|\dot{\bfxi}|^2+\delta\,\|\nabla{\sf q}^\sharp\|_2^2)+C\delta^4\calf_1\,.
$$
In the above formulas, all constant involved are independent of the order of approximation $N$ and $k$. If we sum side-by-side the last three displayed inequalities and restrict the size of $\delta,\eta$ and $\epsilon_0$ below a quantity independent of $R_k$, we conclude, in particular,
\be
\ode{}t\left(\|\bfv\|_2^2+|\bfxi|^2+\|\mathbb D(\bfv)\|_2^2\right)+ c_0\,(\|\mathbb D(\bfv)\|_2^2+\|\bfv_t\|_2^2+|\dot{\bfxi}|^2+\|\Div\mathbb T(\bfv,{\sf r})\|_2^2)\le 
c_1\,\delta^2\,\|\nabla{\sf q}^\sharp\|_2^2+ C\delta^4\calf_1\,.
\eeq{6.41_0}
This inequality enables us to prove the following result.
\Bl  There exist constants $C_0=C_0(\Omega,{\sf m},\nu)>0$, $\delta_0=\delta_0(\Omega,{\sf m},\nu)>0$ such that if $\epsilon_0\le C_0$, then for all $\delta\in(0,\delta_0)$ the system of differential equations \Eqref{6.17} has one and only one $T$-periodic solution for each $N\in\nat$. Moreover, the 
Galerkin approximations \Eqref{6.19} satisfy the following uniform estimate
\be\ba{ll}\medskip
\Max{t\in [0,T]}\big(\|\bfv_N(t)\|_6^2+|\bfxi_N(t)|^2\!+\!\|\mathbb D(\bfv_N(t))\|_2^2\big)\\ \hspace*{1.1cm}+\Int0T\big(\|\mathbb D(\bfv_N)\|_2^2+\|(\bfv_{N})_t\|_2^2+|\dot{\bfxi}_N|^2+\|D^2\bfv_N\|_2^2+\|\nabla{\sf r}_N\|_2^2\big){\rm d}t\le 
C_0\,\delta^2\Big[\Int0T\big(\delta^2\calf_1+\|\nabla{\sf q}^\sharp\|_2^2\big){\rm d}t\Big]\,,
\ea
\eeq{intno}
where $C_0$ is independent of $R_k$ and $N$.
\EL{6.2}
{\em Proof.} For fixed $N\in\nat$, let
$$
E_N={\rm span}\,\{\bfpsi_1,\ldots,\bfpsi_N\}\,.
$$ 
We endow $E_N$ with the norm
$$
|\bfx|_\mu:=\left(\sum_{i=1}^N(1+\mu_i)|x_i|^2\right)^{\!\!\frac12}\,,\ \ \bfx\in E_N\,,
$$
where $\mu_i$, $i=1,\ldots N,$ are the first $N$ eigenvalues of the problem \Eqref{3.4_0}. Consider the map
$$
{\sf M}:\bfc_{N}(0)\in E_N\mapsto {\sf M}(\bfc_{N}(0))=\bfc_{N}(T)\in E_N\,,
$$
where $\bfc_N=\bfc_N(t)$, $t\ge 0$, is the unique solution to \Eqref{6.17} corresponding to initial data $\bfc_{N}(0)$. Clearly, the map ${\sf M}$ is continuous. We shall now show that, for a suitable choice of $R$, ${\sf M}$ transforms the ball of $E_N$ centered at the origin and of radius $R$, ${\sf B}_R$, into itself. To this end, we observe that from \Eqref{6.41_0} and \Eqref{1.9} it follows, in particular,
\be
\ode{} tE_N(t)+\kappa_kE_N(t)\le c_1\,\delta^2\|\nabla{\sf q}^\sharp\|_2^2+ C\,\delta^4\calf_1=:\calf\,,
\eeq{stbe}
where 
$$
E_N(t):=\|\bfv_N(t)\|_2^2+|\bfxi_N(t)|^2+\|\mathbb D(\bfv_N(t))\|_2^2\,.
$$ 
and $\kappa_k$ is a positive constant depending  on $R_k$.
Notice that from \Eqref{6.14}, \Eqref{ortg} and \Eqref{poi}, we infer
\be
E_N(t)=|\bfc_N(t)|_\mu^2\,.
\eeq{naju}
By integrating \Eqref{stbe} and taking into account \Eqref{naju}, it easily follows that
$$
|\bfc_N(T)|_\mu^2\le |\bfc_N(0)|_\mu^2{\rm e}^{-\kappa_kT}+\int_0^T\calf(t)\,{\rm d}t\,, 
$$
which ensures that ${\sf M}$ maps ${\sf B}_R$ into itself, provided we choose 
$$
R\ge \left(\frac{\int_0^T\calf(t)\,{\rm d}t}{1-{\rm e}^{-\kappa_kT}}\right)^{\!\!\frac12}\,.
$$
Thus, by Brouwer theorem,  ${\sf M}$ has a fixed point, which proves that \Eqref{6.17} has a time-periodic solution, that we continue to denote by $\bfc_N(t)$. Next, we integrate both sides of \Eqref{6.41_0} over $[0,T]$ and use $T$-periodicity to get
\be
\Int0T\Big(\|\mathbb D(\bfv_N)\|_2^2+\|(\bfv_{N})_t\|_2^2+|\dot{\bfxi}_N|^2+\|D^2\bfv_N\|_2^2\Big)\le 
\frac1{c_0}\int_0^T\calf(t)\,{\rm d}t\,.
\eeq{naju1}
Thus, by the mean-value theorem, there exists $t^*\in (0,T)$ such that
$$
\|\mathbb D(\bfv_N(t^*))\|_2^2=\frac1T\int_0^T\|\mathbb D(\bfv_N(t))\|_2^2{\rm d}t\,,
$$
which, in turn, by \Eqref{naju1} implies
\be
\|\mathbb D(\bfv_N(t^*))\|_2^2\le \frac1{c_0}\int_0^T\calf(t)\,{\rm d}t\,.
\eeq{naju4}
As a result, \Eqref{intno} follows after integrating both sides of \Eqref{naju3} over $[t^*,t^*+t]$, $t\in [0,T]$, and taking into account \Eqref{naju1}, \Eqref{naju4}, \Eqref{1.9} and $T$-periodicity. The uniqueness property is a consequence of \Eqref{naju1} and the linearity of the system \Eqref{6.17}.
\par\hfill$\square$
\subsection{Further Estimates}

We now derive some further estimates involving time-derivatives of the solution. They are formally deduced by first taking the time derivative of both sides of \Eqref{6.10_lin}$_1$ with $\bfpzc{f}_\ell\equiv \bfpzc{f}_\ell^\sharp$, and then testing the resulting equations, in the order, with $\partial_t\bfv$, $\partial^2_t\bfv$ and $\partial_t(\Div\mathbb T(\bfv,{\sf q}))$, $t>0$, and integrating by parts. In order to make such a procedure rigorous, 
we take the time derivative of both sides of \Eqref{6.15}$_1$, multiply the resulting equation one time by $\dot{c}_{iN}$, the second time by $\ddot{c}_{iN}$, the third time by $\mu_i\,\dot{c}_{iN}$, sum over $i$ and integrate by parts over $\varOmega_k$. Thus, setting
\be
\bfG(\tilde{\bfv},\bfv):=-\big[\bfB_{tt}(\bfv)+(\tilde{\bfv}_t+\bfB_t(\tilde\bfv)+\dot{\tilde{\bfxi}})\cdot\nabla\bfv
+(\tilde{\bfv}_t+\dot{\tilde{\bfxi}})\cdot\nabla\bfB(\bfv)+(\tilde{\bfv}-{\tilde{\bfxi}})\cdot\nabla\bfB_t(\bfv)\big]\,,
\eeq{6.34f}
we show the following ``energy equations" (again for simplicity, we set $\nu={\sf m}=1$ and suppress the subscript $N$)
\be 
\ba{ll}\medskip
\frac1{2}\ode{}t\left(\|\bfv_t\|_2^2+|\dot{\bfxi}|^2\right)+2\|\mathbb D(\bfv_t)\|_2^2= 
\half(\Div\bfB(\tilde{\bfv}),|\bfv_t|^2)
+\big(\bfG(\tilde{\bfv},\bfv),\bfv_t\big)
\\ \medskip
\hspace*{6.5cm}
+(\bfpzc{f}_{\ell t}^\sharp+\bfpzc{f}_t+\delta^2{\bf f}_t,\bfv_t)+(\bfpzc{F}^\sharp_t+\delta^2{\bf F}_t)\cdot\dot{\bfxi}
\\ \medskip
\ode{}t \|\mathbb D(\bfv_t)\|_2^2+\|\bfv_{tt}\|_2^2+\,|\ddot{\bfxi}|^2=-\big((\tilde{\bfv}+\bfB(\tilde{\bfv})-\tilde{\bfxi})\cdot\nabla\bfv_t,\bfv_{tt}\big)+\big(\bfG(\tilde{\bfv},\bfv),\bfv_{tt}\big)
\\ \medskip
\hspace*{6.5cm}
+(\bfpzc{f}_{\ell t}^\sharp+\bfpzc{f}_t+\delta^2{\bf f}_t,\bfv_{tt})+(\bfpzc{F}^\sharp_t+\delta^2{\bf F}_{t})\cdot\ddot{\bfxi}\\ \medskip
\ode{}t\big(\|\mathbb D(\bfv_t)\|_2^2\big)+\|\Div\bfT(\bfv_t,{\sf r}_t)\|_2^2+|\dot{\bsf{T}}|^2
=-\big((\tilde{\bfv}+\bfB(\tilde{\bfv})-\tilde{\bfxi})\cdot\nabla\bfv_t,\Div\mathbb T(\bfv_t,{\sf r}_t)\big)+\big(\bfG(\tilde{\bfv},\bfv),\Div\mathbb T(\bfv_t,{\sf r}_t)\big)\\
\hspace*{6.5cm}+(\bfpzc{f}_{\ell t}^\sharp+\bfpzc{f}_t+\delta^2{\bf f}_{t},\Div\mathbb T(\bfv_{t},{\sf r_t}))+(\bfpzc{F}^\sharp_t+\delta^2{\bf F}_{t})\cdot\dot{\bsf T}\,.
\ea
\eeq{6.36_f}
Set
\be
{\sf D}:=\Big[\Int0T\big(\delta^2\calf_1+\|\nabla{\sf q}^\sharp\|_2^2\big){\rm d}t\Big]^\frac12\,.
\eeq{sha}
In Appendix~C the following estimates are shown:
\be\ba{ll}\medskip
|(\bfG(\tilde{\bfv},\bfv),\bfa)|\le c_1\delta (\|\bfv_{tt}\|_2^2+|\ddot{\bfxi}|^2+\|\mathbb D(\bfv_t)\|_2^2+\|\Div\mathbb T(\bfv_t,{\sf r}_t)\|_2^2)+c_3\,\delta^2{\sf D}^2\|\mathbb D(\tilde{\bfv}_t)\|_2^4\\ \medskip
\hspace*{5.1cm}+c_4(1+{\sf D}^2)\|D^2\bfv\|_2^2+c_5\delta^2{\sf D}^2\big(1+{\sf D}^2\big)+\eta\|\bfa\|_2^2\\ \medskip
|(\Div\bfB(\tilde{\bfv}),|\bfv_t|^2)|\le \epsilon_0^\frac43\|\mathbb D(\bfv_t)\|_2^2+c_6\|\bfv_t\|_2^2
\\
\big|\big((\tilde{\bfv}+\bfB(\tilde{\bfv})-\tilde{\bfxi})\cdot\nabla\bfv_t,\bfa\big)\big|\le c\,\epsilon_0^2\left(\|\mathbb D(\bfv_t)\|_2^2+\|\Div\mathbb T(\bfv_t,{\sf r}_t)\|_2^2\right)+\eta\,\|\bfa\|^2_2\,,
\ea
\eeq{6.37}
and
\be\ba{ll}\medskip
|(\bfpzc{f}_t,\bfa)|\le c\delta\big(\|\mathbb D(\bfv_t)\|_2^2\!+\!\|\Div\mathbb T(\bfv_t,{\sf r}_t)\|_2^2\big)\!+\!c_1\delta^2{\sf D}^2\|\mathbb D(\tilde{\bfv}_t)\|_2^2\!+\!c_2(1\!+\!{\sf D}^2)\|D^2\bfv\|_2^2\!+\!c_3\delta^4({\sf D}^2\!+\!{\sf D}^4)\!+\!\eta\|\bfa\|_2^2\,,\\\medskip
|(\bfpzc{f}_{\ell t}^\sharp,\bfa)|\le c\, \delta\left(\|\mathbb D(\bfv_t)\|_2^2+\|\Div\mathbb T(\bfv_t,{\sf r}_t)\|_2^2\right)+c_1\!\|D^2\bfv\|_2^2+c_2\delta^2\big(\|\nabla{\sf q}^\sharp\|_2^2+\|\nabla{\sf q}^\sharp_t\|_2^2\big)+c_3\delta^4{\sf D}^2+\eta\,\|\bfa\|_2^2\\ 
|\bfpzc{F}_t\cdot\bfb|\le
c\,\delta\left(\|\mathbb D(\bfv_t)\|_2^2+\|\Div\mathbb T(\bfv_t,{\sf r}_t)\|_2^2\right)+c_1\|D^2\bfv\|_2^2+c_2\delta^2\|\nabla{\sf q}^\sharp_t\|_2^2+c_3\delta^4{\sf D}^2+\eta|\bfb|^2\,.
\ea
\eeq{6.38}
Set
\be
\cald:=
\|\bfv_t\|_2^2+|\dot{\bfxi}|^2+(1+{\sf D}^2)\|D^2\bfv\|_2^2+\delta^2{\sf D}^2(1+{\sf D}^2+\|\nabla\bfV_t\|_2^2+\|\nabla\bfV_t\|_3^2)+\delta^4\calf_2\\
+\delta^2(\|\nabla{\sf q}^\sharp\|_2^2+\|\nabla{\sf q}^\sharp_t\|_2^2)\,,
\eeq{6.39_00}
and use the inequality
\be
|({\bf f}_{t},\bfa)|+
|{\bf F}_{t}\cdot\bfb|\le c\,\calf_2+\eta\,(\|\bfa\|_2^2+|\bfb|^2)\,.
\eeq{742}
Choosing in \Eqref{6.37}$_1$--\Eqref{6.38}, and \Eqref{742} $\bfa=\bfv_t$ and $\bfb=\dot{\bfxi}$,  respectively, and using \Eqref{6.37}$_2$, from \Eqref{6.36_f}$_1$ we show, for $\epsilon_0, \delta$ sufficiently small, 
\be
\ode{}t\big(\|\bfv_t\|_2^2\!+\!|\dot{\bfxi}|^2\big)\!+\!c_1\|\mathbb D(\bfv_t)\|_2^2
\le c_2\,\!\delta\big(\|\bfv_{tt}\|_2^2+|\ddot{\bfxi}|^2\!+\!\|\Div\mathbb T(\bfv_t,{\sf r}_t)\|_2^2\big)+c_3\delta^2\,\|\mathbb D(\tilde{\bfv}_t)\|_2^4+c_4\cald\,.
\eeq{6.40}
In a similar fashion, from \Eqref{6.36_f}$_2$, and \Eqref{6.37}$_{1,3}$, \Eqref{6.38}, \Eqref{742}  with $\bfa=\bfv_{tt}$, $\bfb=\ddot{\bfxi}$, we obtain for $\delta,\eta$ sufficiently small,
\be
\ode{}t \|\mathbb D(\bfv_t)\|_2^2+c_1\big(\|\bfv_{tt}\|_2^2+\,|\ddot{\bfxi}|^2\big)\le c_2(\epsilon_0^2+\delta)\big(\|\mathbb D(\bfv_t)\|_2^2+\|\Div\mathbb T(\bfv_t,{\sf r}_t)\|_2^2\big)+c_3\delta^2\|\mathbb D(\tilde{\bfv}_t)\|_2^4+c_4\cald\,.
\eeq{6.41}
Finally,  from \Eqref{6.36_f}$_3$, \Eqref{6.37}$_{1,3}$, \Eqref{6.38} and \Eqref{742} with $\bfa=\|\Div\mathbb T(\bfv_{t},{\sf r}_t)\|_2$,  $\bfb=\dot{\bsf T}$, we deduce, again for $\delta,\eta$ sufficiently small,
\be
\ode{}t\|\mathbb D(\bfv_t)\|_2^2+c_1\big(\|\Div\bfT(\bfv_t,{\sf r}_t)\|_2^2+|\dot{\bsf{T}}|^2\big)\le c_2\,(\epsilon_0^2+\delta)\big(\|\bfv_{tt}\|_2^2+|\ddot{\bfxi}|_2^2+\|\mathbb D(\bfv_t)\|_2^2\big)+c_3\delta^2\|\mathbb D(\tilde{\bfv}_t)\|_2^4+c_4 \cald\,.
\eeq{6.42}
Thus, summing side-by-side \Eqref{6.40}--\Eqref{6.42}, and taking $\epsilon_0,\delta$ below certain constants we conclude, in particular,
\be
\ode{}t\big(\|\bfv_t\|_2^2+|\dot{\bfxi}|^2+\|\mathbb D(\bfv_t)\|_2^2\big)+
 c_1\,\big(\|\bfv_{tt}\|_2^2+|\ddot{\bfxi}|^2+\|\Div\mathbb T(\bfv_t,{\sf r}_t)\|_2^2\big)\le c_2\delta^2\,\|\mathbb D(\tilde{\bfv}_t)\|_2^4+c_3 \cald\,,
\eeq{6.43}
where all the constants involved are independent of $R_k$. Thanks to \Eqref{6.43}, we are now in a position to show the following result.
\Bl  There exist constants $C_0=C_0(\Omega,{\sf m},\nu)>0$, $\delta_0=\delta_0(\Omega,{\sf m},\nu)>0$ such that if $\epsilon_0\le C_0$, then for all $\delta\in(0,\delta_0)$ the 
Galerkin approximations \Eqref{6.19} satisfy the following uniform estimate
\be\ba{ll}\medskip
\Max{t\in [0,T]}\big(\|(\bfv_{N})_t(t)\|_2^2+|\dot{\bfxi}_N|^2\!+\!\|\mathbb D((\bfv_{N})_t(t))\|_2^2\big)
+\Int0T\big(\|\mathbb D((\bfv_{N})_t)\|_2^2+\|(\bfv_{N})_{tt}\|_2^2+|\ddot{\bfxi}_N|^2+\|D^2(\bfv_{N})_t\|_2^2\big){\rm d}t\\
\hspace*{6.5cm}\le 
c_0\,\delta^2\Big({\sf Q}^\sharp+{\sf D}^2\Int0T\|\mathbb D(\tilde{\bfv}_t)\|_2^4\,{\rm d}t\Big)\,,
\ea
\eeq{6.44}
where 
$$
{\sf Q}^\sharp:=\Int0T\big[\delta^2(\calf_1+\calf_2+\calf_1^2)+\|\nabla{\sf q}^\sharp\|_2^2+\|\nabla{\sf q}_t^\sharp\|_2^2\big]{\rm d}t+\big(\Int0T\big\|\nabla{\sf q}^\sharp\|_2^2\,{\rm d}t\big)^2\,,
$$
and 
$c_0$  is independent of $R_k$ and $N$.
\EL{6.3}
{\em Proof.} We begin to observe that from 
 \Eqref{6.39_00}, by employing \lemmref{6.2},  \propref{O_k} and recalling \Eqref{sha}, it follows that
$$\ba{ll}\medskip
\Int0T\mathcal D(t)\,{\rm d}t\le c_1\,\big[\delta^2({\sf D}^2 +{\sf D}^4)+\delta^2\Int0T\big(\delta^2\calf_2+\|\nabla{\sf q}^\sharp\|_2^2+\|\nabla{\sf q}_t^\sharp\|_2^2\big){\rm d}t\big]\le c_2\delta^2{\sf Q}^\sharp.
\,.
\ea
$$
Therefore, integrating both sides of \Eqref{6.43} from $0$ to $T$ and employing the $T$-periodicity of the approximations $\bfv_N(t)$, we show, in particular,
\be
\Int0T\big(\|\mathbb D((\bfv_{N})_t)\|_2^2+\|(\bfv_{N})_{tt}\|_2^2+|\ddot{\bfxi}_N|^2+\|D^2(\bfv_{N})_t\|_2^2\big)\le 
C_1\,\delta^2\Big({\sf Q}^\sharp+\Int0T\|\mathbb D(\tilde{\bfv}_t)\|_2^4\Big)\,.
\eeq{6.45}
By the mean-value theorem,  there is $t^*\in(0,T)$ such that
$$
\|(\bfv_{N})_t(t^*)\|_2^2+|\dot{\bfxi}_{N}(t^*)|^2+\|\mathbb D((\bfv_{N})_t(t^*))\|_2^2=\frac1T\int_0^T\left(\|(\bfv_{N})_t(t)\|_2^2+|\dot{\bfxi}_N(t)|^2+\|\mathbb D((\bfv_{N})_t)\|_2^2\right){\rm d}t\,,
$$
which, in turn, by \Eqref{6.45}, implies
\be
\|(\bfv_{N})_t(t^*)\|_2^2+|\dot{\bfxi}(t^*)|^2+\|\mathbb D((\bfv_{N})_t(t^*))\|_2^2\le \frac1T  
C_1\,\delta^2\Big({\sf Q}^\sharp+\Int0T\|\mathbb D(\tilde{\bfv}_t)\|_2^4\Big)\,.
\eeq{6.46}
As a result, integrating \Eqref{6.43} over $[t^*,t^*+t]$, $t\in [0,T]$, and using \Eqref{6.46} along with $T$-periodicity, in combination with \Eqref{6.45} proves \Eqref{6.44}.\par\hfill$\square$\par
\subsection{Proof of Proposition \ref{EUlin}}\label{propo}
We begin to show existence and uniqueness of $T$-periodic solutions to the problem \Eqref{6.10_lin} with $\bfpzc f_\ell$ and $\bfpzc F$  replaced by $\bfpzc f_{\ell}^\sharp$ and $\bfpzc F^\sharp$, respectively. Hereafter, we will  refer to this problem as  \Eqref{6.10_lin}$^\sharp$. 
From \lemmref{6.2} and \lemmref{6.3} we know that the Galerkin approximations are $T$-periodic and satisfy the uniform bound \Eqref{intno} and \Eqref{6.44}. Then, following a standard procedure \cite{GaSi}, we show the existence of a $T$-periodic pair $(\bfv,\bfxi)\in \mathscr W_T(\varOmega_k)\times W^{2,2}(0,T)$ such that (possibly, along a subsequence)
\be\ba{cc}\medskip
(\bfv_N, \bfxi_N)\to (\bfv,{\bfxi})\,,\ \ \mbox{weakly in $\mathscr W_T(\varOmega_k)\times W^{2,2}(0,T)$}\,,
\\
(\bfv_N, \bfxi_N)\to (\bfv,{\bfxi})\,,\ \ \mbox{strongly in $L^r(0,T; W^{1,2}(\varOmega_k))\times W^{1,2}(0,T)$\,,\ \ all $r\in[1,\infty)\,,$}
\ea
\eeq{6.50}
and which, in addition,  satisfies the estimates
\be\ba{ll}\medskip
\Max{t\in [0,T]}\big(\|\bfv(t)\|_6^2+|\bfxi(t)|^2\!+\!\|\mathbb D(\bfv(t))\|_2^2\big)+\Int0T\big(\|\mathbb D(\bfv)\|_2^2+\|\bfv_{t}\|_2^2+|\dot{\bfxi}|^2+\|D^2\bfv\|_2^2\big){\rm d}t\\ \hspace*{6.0cm}\le 
C_0\,\delta^2\!\Int0T\big(\delta^2\calf_1+\|\nabla{\sf q}^\sharp\|_2^2\big){\rm d}t\,,
\ea
\eeq{6.51}
and
\be
\ba{ll}\medskip
\Max{t\in [0,T]}\big(\|\bfv_{t}(t)\|_2^2+|\dot{\bfxi}|^2\!+\!\|\mathbb D(\bfv_{t}(t))\|_2^2\big)+\Int0T\big(\|\mathbb D({\bfv}_{t})\|_2^2+\|\bfv_{tt}\|_2^2+|\ddot{\bfxi}|^2+\|D^2\bfv_{t}\|_2^2\big){\rm d}t\\
\hspace*{1.2cm}\le 
c_0\delta^2\Big[{\sf D}^2\!\!\Int0T\|\mathbb D(\tilde{\bfv}_t)\|_2^4\,{\rm d}t\!+\!\!\Int0T\big[\delta^2(\calf_1+\calf_2+\calf_1^2)+\|\nabla{\sf q}^\sharp\|_2^2+\|\nabla{\sf q}_t^\sharp\|_2^2\big]{\rm d}t+\big(\Int0T\big\|\nabla{\sf q}^\sharp\|_2^2\,{\rm d}t\big)^2\Big]\,,
\ea 
\eeq{6.52}
where, we recall, ${\sf D}$ is defined in \Eqref{sha}.
Furthermore, setting
\be
\bfPhi(x,t):=\sum_{\ell=0}^m\gamma_\ell(t)\bfpsi_\ell(x)\,,\ \ \ \hat{\bfPhi}(x,t):=\sum_{\ell=0}^m\gamma_\ell(t)\hat{\bfpsi}_\ell(x)
\eeq{6.53}
with $\gamma_m$ arbitrary smooth $T$-periodic functions and $m\le N$, from \Eqref{6.15} we easily deduce that $(\bfv_N,\bfxi_N)$ satisfies
\be\ba{ll}\ms
\Int0T\Big(\langle\partial_t\bfv_N,\bfPhi\rangle+(\partial_t\bfB(\bfv_N),\bfPhi)+\big((\tilde{\bfv}+\bfB(\tilde{\bfv})-\tilde{\bfxi})\cdot\nabla\bfv_N+(\tilde{\bfv}-\tilde{\bfxi})\cdot\nabla\bfB(\bfv_N)\\ \medskip\hspace*{.9cm}+\delta(\bfV\!-\!\bfzeta)\cdot\nabla(\bfv_N\!+\!\bfB(\bfv_N))\!+\!\delta(\bfv_N+\bfB(\bfv_N)\!-\!\bfxi_N)\cdot\nabla\bfV\!+\!\delta^2(\bfV\!-\!\bfzeta)\cdot\nabla\bfV,\bfPhi\big)\!+\! 2\nu(\mathbb D(\bfv_N),\mathbb D ({\bfPhi}))\\ 
\hspace*{.9cm}-\big(\bfpzc{f}_\ell(\bsf s,\bfv_N,{\sf q}^\sharp,\bfxi_N)+\bfpzc{f}(\bsf s,\bfv_N,\tilde{\bfv},\tilde{\bfxi})+\delta^2{\bf f},\bfPhi\big)-(\bfpzc{F}(\bsf s,\bfv_N,{\sf q}^\sharp)+\delta^2{\bf F})\cdot  \hat{\bfPhi}\Big)\,{\rm d}t=0\,.
\ea
\eeq{6.54}
Thus, passing to the limit $N\to\infty$ in \Eqref{6.54} and using 
\Eqref{6.50} along with the assumptions on $\tilde{\bfv}$ and the properties of $\bfB$, one readily establishes that also $(\bfv,\bfxi)$ satisfies \Eqref{6.54} for all function $\bfPhi$ of the above form. Since every function $\bfphi\in \calc_T(B_k)$ can be approximated uniformly pointwise with its first spatial derivative by functions of the type \Eqref{6.53} \cite[Lemma 3.1]{GaSi}, by a standard argument we deduce that $(\bfv,\bfxi)$ obeys \Eqref{6.13}. We shall now show that, for such a pair $(\bfv,\bfxi)$, we can find a corresponding pressure field ${\sf q}\in \calq(\varOmega_k)$ such that the triple $(\bfv,{\sf q},\bfxi)$ solves the problem \Eqref{6.10_lin}--\Eqref{flin} along with the estimates \Eqref{BaNa1}--\Eqref{BaNa2}. To this end, choosing, in particular, in \Eqref{6.13} $\bfphi\in \calc_0(B_k)$ and integrating by parts we deduce
\be
\big(\partial_t\bfv-\bfH(\tilde{\bfv},\bfv)-\nu\Delta\bfv-\bfpzc{f}_\ell^\sharp-\bfpzc{f}-\delta^2{\bf f},\bfphi)=:\big(\bsf{h},\bfphi\big)=0\,, 
\eeq{6.55}
with, we recall, $\bfH(\tilde{\bfv},\bfv)$ given in \Eqref{H}. 
By employing the estimates \Eqref{H_2}--\Eqref{A24} in Appendix B along with Remark B.1 and the assumption $\|\mathbb D(\tilde{\bfv}(t))\|_2\le \epsilon_0$, we show that $\bsf h\in L^2(0,T;L^2(\varOmega))$ and, in addition, 
\be 
\|\bsf h\|_{L^2(L^2)}^2\le c\,\Int0T\big(\|\bfv_t\|_2^2+|\dot{\bfxi}|^2+\|\nabla\bfv\|_{1,2}^2+\delta^2\|\nabla{\sf q}^\sharp\|_2^2+ \delta^4\calf_1\big){\rm d}t\,.
\eeq{Ari}
Since $\bfphi$ is arbitrary in $\calc_0(B_k)$, from \Eqref{6.55} and classical results it follows that there exists ${\sf q}\in L^2(W^{1,2}({\varOmega}_k))$ such that 
\be
\bsf h=\nabla{\sf q}\,.
\eeq{6.58}
Therefore,  $(\bfv,{\sf q},\bfxi)$ obeys \Eqref{6.10_lin}$^\sharp_{1,2}$.  Moreover, from \Eqref{6.13} and \Eqref{6.10_lin}$^\sharp_{1,2}$ we show, by a standard argument,   that also \Eqref{6.10_lin}$^\sharp_4$ is satisfied.    
Combining \Eqref{6.58} with \Eqref{Ari} and
\Eqref{6.51}, we also establish the following estimate  
\be
\|\nabla{\sf q}\|_{L^2(L^2)}^2\le c \,\delta^2\!\Int0T\big(\delta^2\calf_1+\|\nabla{\sf q}^\sharp\|_2^2\big){\rm d}t\,,
\eeq{6.56}
where the constant $c$ is independent of $R_k$. 
Furthermore, from \Eqref{6.10_lin}$_1$, we get
$$
\nabla {\sf q}_t=\bfv_{tt}+(\tilde{\bfv}+\bfB(\tilde{\bfv})-\tilde{\bfxi})\cdot\nabla\bfv_t-\bfG(\tilde{\bfv},\bfv)-(\bfpzc{f}_{\ell t}^\sharp+\bfpzc{f}_t+\delta^2{\bf f}_{t})\,,
$$
with $\bfG(\tilde{\bfv},\bfv)$ given in \Eqref{6.34f}. Thus,   employing \Eqref{Gbe}, \Eqref{B.9}, \Eqref{ft},  \Eqref{fli} and Remark C.1 in Appendix~C, along with \propref{O_k}, we show 
$$ \ba{ll}\medskip
\|\nabla{\sf q}_t\|_{L^2(L^2)}^2
\le 
c\,\big\{\Int0T\!\!\big[\|\bfv_{tt}\|_2^2+|\ddot{\bfxi}|^2+\|\mathbb D(\bfv_t)\|_2^2+\|D^2\bfv_t\|_2^2+(1+{\sf D}^2)\|D^2\bfv\|_2^2\\
\hspace*{5.4cm}+\delta^2\big({\sf D}^2\|\mathbb D(\tilde{\bfv}_t)\|_2^4+\|\nabla{\sf q}^\sharp\|_2^2+\|\nabla{\sf q}^\sharp_t\|_2^2\big)+\calf_2\big]{\rm d} t+\delta^2{\sf D}^2\big(1+{\sf D}^2)\big\}\,.
\ea
$$
From \Eqref{6.52} and the latter, we obtain
\be
\|\nabla{\sf q}_t\|_{L^2(L^2)}^2\le c\,\delta^2\Big[{\sf D}^2\!\!\Int0T\|\mathbb D(\tilde{\bfv}_t)\|_2^4\,{\rm d}t\!+\!\!\Int0T\big[\delta^2(\calf_1+\calf_2+\calf_1^2)+\|\nabla{\sf q}^\sharp\|_2^2+\|\nabla{\sf q}_t^\sharp\|_2^2\big]{\rm d}t+\big(\Int0T\big\|\nabla{\sf q}^\sharp\|_2^2\,{\rm d}t\big)^2\Big]\,.
\eeq{6.57}
The solution $(\bfv,{\sf q},\bfxi)$ to \Eqref{6.10_lin}$^\sharp$  just constructed is also unique in its class of existence. In fact, since \Eqref{6.10_lin}$^\sharp$ is linear, it suffices to show the uniqueness of the null solution, namely, that \Eqref{6.10_lin}$^\sharp$ with ${\bf f}\equiv{\bf F}\equiv{\sf q}^\sharp\equiv\0$, has only the solution $\bfv\equiv\nabla {\sf q}\equiv\bfxi\equiv\0$. However, under these circumstances, this follows at once from  \Eqref{6.51}.
\par
With this result in hand, we may    employ (for example) the classical successive approximations method to prove the existence of a unique solution to the original problem \Eqref{6.10_lin} in the specified function class.
Consider the map   
$$
{\sf L}:{\sf q}^\sharp\in \calq({\varOmega}_k)\mapsto {\sf L}({\sf q}^\sharp)=(\bfv,{\sf q},\bfxi)\in \mathscr W_T({\varOmega}_k)\,,
$$
where $(\bfv,{\sf q},\bfxi)$ is the solution to \Eqref{6.10_lin}--\Eqref{flin} --with $\bfpzc{f}_\ell$ and $\bfpzc{F}$  replaced by $\bfpzc{f}_\ell^\sharp$ and $\bfpzc{F}^\sharp$, respectively,
The results we just obtained guarantee  that ${\sf L}$ is a well--defined linear map. Furthermore, it is continuous. In fact
setting
$$
\langle{\sf q}\rangle:=\|\nabla{\sf q}\|_{W^{1,2}(D^{1,2})}\,;\ \ \langle\langle(\bfv,{\sf q},\bfxi)\rangle\rangle:=\|\bfv\|_{\mathscr W_T({\varOmega}_k)}+\langle{\sf q}\rangle+\|\bfxi\|_{W^2(0,T)}\,, 
$$
and
$$
{\sf q}^\sharp_{1,2}:={\sf q}_{1}^\sharp-{\sf q}_{2}^\sharp\,,\ \  \mu:=M^2+1\,, 
$$
recalling that $\tilde{\bfv}\in\mathscr S_{\epsilon_0,M}$, 
by \Eqref{6.51}, \Eqref{6.52}, \Eqref{6.56}, and \Eqref{6.57} 
we deduce
\be
\langle\langle{\sf L}({\sf q}^\sharp_{1,2})\rangle\rangle\le \kappa\,\delta \big[\mu\langle{\sf q}_{1,2}^\sharp\rangle +\langle{\sf q}^\sharp_{1,2}\rangle^2\big]\,,
\eeq{7.61}
where $\kappa$ is a positive constant independent of $R_k$. 
Furthermore, if we define
$$
F:=\Big[\Int0T\big(\calf_1+\calf_2+\calf_1^2\big){\rm d}t\Big]^\frac12,
$$
and recall \Eqref{sha}, again by  \Eqref{6.51}, \Eqref{6.52}, \Eqref{6.56}, and \Eqref{6.57} we get
\be
\langle\langle {\sf L}({\sf q}^\sharp)\rangle\rangle\le \kappa \,\delta \left[\delta\,\mu\,F+\mu\langle{\sf q}^\sharp\rangle+\langle{\sf q}^\sharp\rangle^2\right].
\eeq{7.62}
 Consider the sequence $\{(\bfv_n,{\sf q}^\sharp_n,\bfxi_n)\}$ where
\be 
(\bfv_n,{\sf q}^\sharp_n,\bfxi_n)={\sf L}({\sf q}_{n-1}^\sharp)\,,\ n=1,2,\ldots;\ \ \ {\sf q}_0^\sharp=0
\,.
\eeq{CoN}
From \Eqref{7.62} we have 
$$
\langle\langle(\bfv_1,{\sf q}_1^\sharp,\bfxi_1)\rangle\rangle=\langle\langle {\sf L}({\sf q}_0^\sharp)\rangle\rangle\le \kappa\,\delta^2\mu\,F\,.
$$
We want to show that
\be
\langle\langle(\bfv_n,{\sf q}_n^\sharp,\bfxi_n)\rangle\rangle\le 2\kappa\,\delta^2\,\mu \,F\,,\ \ \mbox{for all $n\in\nat$}\,.
\eeq{7.63}
In fact, by induction, assuming
$$
\langle\langle {\sf L}({\sf q}_{n-1}^\sharp)\rangle\rangle\le 2\kappa\,\delta\,\mu \,F\,,
$$
from \Eqref{CoN} and \Eqref{7.62} we show
\be
\langle\langle(\bfv_n,{\sf q}_n^\sharp,\bfxi_n)\rangle\rangle\le \kappa\,\delta^2\mu\,F \left[1+2\kappa\,\delta\mu +4\kappa^2 \delta^3\mu\, F\right]\,.
\eeq{7.64}
Therefore, if  (for example)
\be
\delta\le \min\left\{(4\kappa\,\mu)^{-1},(4\kappa F)^{-\frac12}\right\}\,,
\eeq{c0}
then  \Eqref{7.63} holds. Combining \Eqref{7.63},  \Eqref{7.61} and \Eqref{CoN} we get
$$
\langle\langle(\bfv_{n+1}-\bfv_n,{\sf q}^\sharp_{n+1}-{\sf q}_n^\sharp,\bfxi_{n+1}-\bfxi_n\rangle\rangle\le \kappa\delta\,[(\mu+2\kappa\,\delta^2\mu\, F)\langle {\sf q}_{n}^\sharp-{\sf q}_{n-1}^\sharp\rangle]:=\alpha\,\langle {\sf q}_{n}^\sharp-{\sf q}_{n-1}^\sharp\rangle\,.
$$
Since, by \Eqref{c0}, it is $\alpha<3/8<1$, by a standard argument one shows that $\{(\bfv_n,{\sf q}^\sharp_n,\bfxi_n)\}$ is Cauchy in $\mathscr W_T(\varOmega_k)\times W^{1,2}(D^{1,2})\times W^{2,2}(0,T)$ converging there to some $\{(\bfv,{\sf q},\bfxi)\}$. Using the latter and letting $n\to\infty$ in  \Eqref{CoN}, we also deduce that 
$\{(\bfv,{\sf q},\bfxi)\}$ satisfies the original problem \Eqref{6.10_lin}. Likewise, writing \Eqref{6.51},  \Eqref{6.52}, \Eqref{6.56} and \Eqref{6.57} with  $(\bfv,{\sf q},\bfxi)$ and ${\sf q}^\sharp$ replaced by $(\bfv_n,{\sf q}_n,\bfxi_n)$ and ${\sf q}^\sharp_{n-1}$, respectively, and then passing to the limit $n\to\infty$, we establish that $(\bfv,{\sf q},\bfxi)$ satisfy the same inequalities with ${\sf q}$ in place of ${\sf q}^\sharp$.
Thus, using the latter and \Eqref{6.52} provides \Eqref{BaNa2}, which  completes the existence proof. With this in mind, we deduce that there exists a constant $C$ depending only on $\varOmega$, and the physical parameters such that if $\delta\le C$, the following properties hold. From \Eqref{6.56}, we obtain
\be
\|\nabla{\sf q}\|_{L^2(L^2)}^2\le c \,\delta^4\!\Int0T\calf_1\,{\rm d}t\,,
\eeq{7.67}
which, in turn, once used in \Eqref{sha}, furnishes
\be
{\sf D}^2\le c_1\delta^2\!\Int0T\calf_1\,{\rm d}t\,.
\eeq{7.68}
Similarly, from \Eqref{6.57}, with the help of \Eqref{7.67} and \Eqref{7.68}, we show 
\be
\|\nabla{\sf q}_t\|_{L^2(L^2)}^2\le c_2\delta^4\Big[\Int0T\calf_1{\rm d}t\Int0T\|\mathbb D(\tilde{\bfv}_t\|_2^4{\rm d}t+\Int0T(\calf_1+\calf_2+\calf_1^2){\rm d}t\Big]\,.
\eeq{7.69}
Employing \Eqref{7.67}--\Eqref{7.69} in \Eqref{6.51} and \Eqref{6.52} (with ${\sf q}^\sharp\equiv{\sf q}$) proves the validity of \Eqref{BaNa1} and \Eqref{BaNa2}. Concerning the uniqueness property, being \Eqref{6.10_lin} linear, it suffices to show the uniqueness of the null solution, namely, that \Eqref{6.10_lin} with ${\bf f}\equiv{\bf F}\equiv\0$, has only the solution $\bfv\equiv\nabla {\sf q}\equiv\bfxi\equiv\0$. However, under these circumstances, this follows at once from  \Eqref{BaNa1} or \Eqref{BaNa2}.
The proof of the proposition is thus completed.
\par\hfill$\square$
\setcounter{equation}{0}
\section{Solvability of the Approximating Nonlinear Problems \Eqref{6.3}}
\label{sec:Ok}
The main objective of this section is to prove existence of $T$-periodic solutions to the sequence of approximating problems  \Eqref{6.3} 
along with corresponding estimates, holding uniformly in $k$. To this end, we present some preliminary results. Let
$$
\calv_1:=\|\bff_{\mbox{\tiny $\bfV$}}\|_2^2+\|\bff_{\mbox{\tiny $\bfV$}}\|^2_\frac65+|\bfF_{\mbox{\tiny $\bfV$}}|^2\,,\ \ \ \calv_2:=\|\partial_t\bff_{\mbox{\tiny $\bfV$}}\|_2^2+|\dot{\bfF}_{\mbox{\tiny $\bfV$}}|^2\,.
$$
We have the following result,  whose proof is presented in Appendix A.
\begin{lemma}{\sl
$\calv_i\in L^1(0,T)$, $i=1,2$,   and there is a constant $C$ depending only on $\varOmega$, $T$ and the physical parameters such that}
\be
\|\calv_1\|_{L^1(0,T)}+\|\calv_2\|_{L^1(0,T)}\le C
\eeq{calf}
\label{le:calf}
\end{lemma}
We also have the following. 
\Bl
Let $\mathscr S_{\epsilon_0,M}$ be the set defined in \Eqref{Sepm},  pick $\epsilon_0$, $\delta_*$   as in {\rm Proposition \ref{EUlin}}, and consider the map
$$
\Phi:(\tilde{\bfv},\tilde{\bfxi})\in \mathscr S_{\epsilon_0,M}\subset W^{1,4}(0,T;\mathcal D^{1,2}(B_k))\mapsto (\bfv,\bfxi)\subset W^{1,4}(0,T;\mathcal D^{1,2}(B_k))
$$
where $(\bfv,{\sf q},\bfxi)$ is the unique solution given in {\rm Proposition \ref{EUlin}} corresponding to ${\bf f}\equiv\bff_{\mbox{\tiny$\bfV$}}$ and ${\bf F}\equiv\bfF_{\mbox{\tiny$\bfV$}}$. Then, the following properties hold.
\begin{itemize}   
\item[{\rm (a)}] $\mathscr S_{\epsilon_0,M}$ is  closed and convex;
\item[{\rm (b)}] $\Phi(\mathscr S_{\epsilon_0,M})$ is compact; 
\item[{\rm (c)}] There are  constants $C_i>0$, $i=1,2$, with $C_1$ depending only on $\varOmega$, $T$ and the physical parameters and $C_2$ depending also on $\epsilon_0$, such that if we choose $M=C_1F$, with $F$ given in \Eqref{cal}--\Eqref{fmu}, and take $\delta\in (0,\delta_1)$, with $\delta_1:=C_2F^{-\frac12}$, then $\Phi$ maps $\mathscr S_{\epsilon_0,M}$ into itself.
\end{itemize}
\EL{8.2}
{\em Proof.} Let 
$\{\tilde{\bfv}_n,\tilde{\bfxi}_n\}\subset \mathscr S_{\epsilon_0,M}$ be a sequence such that 
\be
(\tilde{\bfv}_n,\tilde{\bfxi}_n)\to(\tilde{\bfv},\tilde{\bfxi})\ \ \mbox{in $ W^{1,4}(0,T;\mathcal D^{1,2}(B_k))$}\,.
\eeq{8.2}
Then, clearly, $\|\tilde{\bfv}\|_{W^{1,4}(\cald^{1,2})}\le M$. Furthermore, we can select a subsequence converging weak star in $L^\infty(0,T;\cald^{1,2}(B_k))$ to some $\bfw\in L^\infty(0,T;\cald^{1,2}(B_k))$ satisfying the $\epsilon_0$ bound. However, because of \Eqref{8.2}, we must have $\bfw\equiv\tilde{\bfv}$, which shows the closedness property. Since convexity is obvious, the proof of (a) is completed. We next observe that, by \Eqref{1.9},  
for any sequence $\{\tilde{\bfv}_n,\tilde{\bfxi}_n\}\subset \mathscr S_{\epsilon_0,M}$, we have
$$
\|\tilde{\bfv}_n\|_{W^{1,4}(W^{1,2})}+\|\tilde{\bfxi}_n\|_{W^{1,4}(0,T)}\le C_0\,,
$$
for a suitable $C_0>0$. As a consequence, by Proposition \ref{EUlin},  \Eqref{BaNa1} and \Eqref{BaNa2}, we infer that the sequence of corresponding solutions $\{\bfv_n,{\sf q}_n,\bfxi_n\}$ satisfies    the following 
bound
$$
\|{\bfv}_n\|_{W^{1,\infty}(W^{1,2})}+\|{\bfv}_n\|_{W^{1,2}(W^{2,2})}+\|{\bfv}_n\|_{W^{2,2}(L^2)}+\|{\bfxi}_n\|_{W^{2,2}(0,T)}\le C_1\,,
$$
for another suitable $C_1>0$. Thus, $\{\bfxi_n\}$ is compact in $W^{1,4}(0,T)$. Moreover, setting $\bfw_n:=\partial_t\bfv_n$, from the above bound it follows  
$$
\{\bfw_n\}\ \mbox{bounded in}\ L^\infty(W^{1,2})\cap L^2(W^{2,2})\cap W^{1,2}(L^2)\,,  
$$
so that from \cite[Corollary 6]{Simo} $\{\bfw_n\}$ is compact in $ L^q(W^{1,2})$, for all $q\in[1,\infty)$ which furnishes, in particular, $\{\bfv_n\}$ compact in $W^{1,4}(\mathcal D^{1,2})$. The validity of property (b) is thus secured. In order to show (c), we begin to observe that from \Eqref{BaNa1} and \Eqref{BaNa2} it follows, in particular, 
\be\ba{rl}\medskip
\Max{t\in[0,T]}\|\mathbb D(\bfv(t))\|_2^2&\!\!\!\!\le c_0\delta^4F^2\\
\|\bfv\|^2_{W^{1,4}(\cald^{1,2}(B_k))}&\!\!\!\!\le c_0\delta^4F^2(M^4+1)\,.
\ea
\eeq{8.3}
Therefore, from \Eqref{8.3}$_1$ it follows $\max_{t\in[0,T]}\|\mathbb D(\bfv(t))\|_2\le \varepsilon_0$, provided $\delta\le c_0^{-\frac14}(\epsilon_0/F)^{\frac12}$.  Furthermore, using \Eqref{assu} in \Eqref{8.3}$_2$, furnishes
$$
\|\bfv\|^2_{W^{1,4}(\cald^{1,2}(B_k))}\le c_1\delta^2F^2\le C_1^2F^2\,,
$$
which completes the proof of the lemma.
\par\hfill$\square$\par
We are now in a position to prove the following proposition that establishes the existence of $T$-periodic solutions to the ``approximating" problems \Eqref{6.10}-\Eqref{6.12}, with corresponding bounds uniform with respect to $R_k$.
\Bp There is $\delta_0>0$ depending only on $\varOmega$, $T$ and the physical parameters such that, for every $\delta\in(0,\delta_0)$ there exists a corresponding $T$-periodic solution $(\bfv,{\sf q},\bfxi)\in \mathscr W_T(\varOmega_k)\times\calq(\varOmega_k)\times W_T^{2,2}(\real^3)$ to \Eqref{6.10}--\Eqref{6.12}.
Furthermore, this solution satisfies the estimate 
\be
\|\bfv\|_{\mathscr W_T(\varOmega_k)}+\|{\sf q}\|_{\calq(\varOmega_k)}+\|\bfxi\|_{W^{2,2}(0,T)}\le C_0\,\delta^2\,.
\eeq{belest}
where $C_0$ does not depend on $k$.
\EP{8.1}
{\em Proof.} Choose $\delta_0=\min\{\delta_*,\delta_1\}$ with $\delta_*$, $\delta_1$ as in  Proposition \ref{EUlin} and \lemmref{8.2}. Thus, in view of the latter, the stated properties will follow from Schauder fixed-point theorem, provided we show that the map $\Phi$ is continuous. To this end, let $\tilde{\bfv},\tilde{\bfv}_0\in \mathscr S_{\epsilon_0,M}$ with $\tilde{\bfxi}=\tilde{\bfv}|_{\Omega_0}$, $\tilde{\bfxi}_0=\tilde{\bfv}_0|_{\Omega_0}$, and denote by $(\bfv,{\sf q},\bfxi)$ and $(\bfv_0,{\sf q}_0,\bfxi_0)$ the corresponding solutions given in Proposition \ref{EUlin} with ${\bf f}\equiv\bff_{\mbox{\tiny$\bfV$}}$ and ${\bf F}\equiv\bfF_{\mbox{\tiny$\bfV$}}$. Thus, setting
$$
\bfw=\bfv-\bfv_0\,,\ \ {\sf r}:= {\sf q}-{\sf q}_0\,,\ \ \bfchi:=\bfxi-\bfxi_0\,, 
$$
from \Eqref{6.10_lin} we get
\be
\ba{cc}\medskip\left.\ba{ll}\medskip
\partial_t{\bfw}+\partial_t\bfB(\bfw)+(\tilde{\bfv}_0+\bfB(\tilde{\bfv}_0)-\tilde{\bfxi}_0)\cdot\nabla\bfw+(\tilde{\bfv}_0-\tilde{\bfxi}_0)\cdot\nabla\bfB(\bfw)\\ \medskip\hspace*{2.2cm}+\delta(\bfV-\bfzeta)\cdot\nabla(\bfw+\bfB(\bfw))+\delta(\bfw+\bfB(\bfw)-\bfchi)\cdot\nabla\bfV\\ \medskip
\hspace*{3cm}=\nu\Delta\bfw-\nabla{\sf r}+\bfpzc{f}_\ell(\bsf{s},\bfw,{\sf r},\bfchi)+{\bfpzc{f}}(\bsf s,\bfw,\tilde{\bfv}_0,\tilde{\bfxi}_0)+
{\bfcaln}
\\
\Div\bfw=0\ea\right\}\,  \mbox{in $\varOmega_k\times(0,T)$}\,,\\ \medskip
\bfw=\bfchi\,, \ \ \mbox{at $\varSigma\times (0,T)$}\,;\ \ \bfw(x,t)=\0\,,\ \ \mbox{at $\partial B_k\times (0,T)$}\,;\\ 
{\sf m}\dot{\bfchi}(t)=-\Int{\varSigma}{}\mathbb T(\bfw,{\sf r})\cdot\bfn{\rm d}\varSigma+ \bfpzc{F}(\bsf{s},\bfw,{\sf r})\,,\ \ t\in\real\,,
\ea
\eeq{8.4}
where
\be\ba{rl}\medskip
\bfcaln:=&\!\!\!\!-\big[(\tilde{\bfv}-\tilde{\bfv}_0)+\bfB(\tilde{\bfv}-\tilde{\bfv}_0)-(\tilde{\bfxi}-\tilde{\bfxi}_0)\big]\cdot\mathbb A\cdot\nabla\bfv -\big[(\tilde{\bfv}-\tilde{\bfv}_0)-(\tilde{\bfxi}-\tilde{\bfxi}_0)\big]\cdot\mathbb A\cdot\nabla\bfB(\bfv)\\
&\!\!\!\!-\bfB(\tilde{\bfv}-\tilde{\bfv}_0)\cdot\mathbb A\cdot\nabla\bfB(\bfv)\,.
\ea
\eeq{8.5}
We observe that, by \lemmref{Def}, we have
\be
\|\mathbb A\|_{W^{1,\infty}(L^\infty)}\le C\,.
\eeq{8.6} 
Employing \Eqref{Nato}, \Eqref{Nato1},  \Eqref{Nato2} in Appendix B, we deduce
\be\ba{ll}\medskip 
\|\big|(\tilde{\bfv}-\tilde{\bfv}_0)+\bfB(\tilde{\bfv}-\tilde{\bfv}_0)\big|\nabla\bfv\|_2^2+\|\big|\tilde{\bfv}-\tilde{\bfv}_0\big|\nabla\bfB(\bfv)\|_2^2+\|\big|\bfB(\tilde{\bfv}-\tilde{\bfv}_0)\big| \nabla\bfB(\bfv)\|_2^2\\
\hspace*{5.4cm} \le c\,\|\mathbb D(\tilde{\bfv}-\tilde{\bfv}_0)\|_2^2\left(\|\mathbb D(\bfv)\|_2^2+\|\mathbb D(\bfv)\|_2\|D^2\bfv\|_2\right)\,,
\ea
\eeq{8.7}
whereas, from \Eqref{6.20}, \Eqref{6.22} and \Eqref{B12},
\be
\|\big|\tilde{\bfxi}-\tilde{\bfxi}_0)\big|\nabla\bfv\|_2^2+\|\big|\tilde{\bfxi}-\tilde{\bfxi}_0\big|\nabla\bfB(\bfv)\|_2^2\le c\,\|\mathbb D(\tilde{\bfv}-\tilde{\bfv}_0)\|_2^2|\mathbb D(\bfv)\|_2^2\,.
\eeq{8.8}
Therefore, observing  that by, Proposition \ref{EUlin}, \lemmref{8.2}(c), \Eqref{BaNa1}, \Eqref{BaNa2} and classical embedding result,
we have
\be
\max_{t\in[0,T]}\left(\|\mathbb D(\bfv(t))\|_2+\|D^2\bfv\|_{2}\right)\le C
\eeq{8.9}
with $C$ depending only on $\delta$ and the data, and that ${\varOmega}_k$ is bounded, from \Eqref{8.5}--\Eqref{8.9} we infer
\be
c_k\|\bfcaln\|_{\frac65}+\|\bfcaln\|_{2}\le c\,\|\mathbb D(\tilde{\bfv}-\tilde{\bfv}_0)\|_{2}\,,
\eeq{8.10}
with $c_k$ depending on ${\varOmega}_k$.
Next, if we employ \Eqref{C6} and \Eqref{NaFr1} in Appendix C together with \Eqref{8.6}--\Eqref{8.8}, we obtain
\be\ba{ll}\medskip
\|\big\{\big[(\tilde{\bfv}-\tilde{\bfv}_0)+\bfB(\tilde{\bfv}-\tilde{\bfv}_0)-(\tilde{\bfxi}-\tilde{\bfxi}_0)\big]\cdot\mathbb A\cdot\nabla\bfv\big\}_t\|_2^2\\ \medskip
\hspace*{4cm}\le c\,\big[\big(\|\mathbb D(\tilde{\bfv}-\tilde{\bfv}_0)\|_2^2+\|\mathbb D(\tilde{\bfv}_t-\tilde{\bfv}_{0t})\|_2^2\big)\left(\|\mathbb D(\bfv)\|_2^2+\|\mathbb D(\bfv)\|_2\|D^2\bfv\|_2\right)\\
\hspace*{4.9cm}+\|\mathbb D(\tilde{\bfv}-\tilde{\bfv}_0)\|_2^2\left(\|\mathbb D(\bfv_t)\|_2^2+\|\mathbb D(\bfv_t)\|_2\|D^2\bfv_t\|_2\right)\big]\,.
\ea
\eeq{8.11}
Again from \Eqref{8.6}--\Eqref{8.8} and \Eqref{Sib},  \Eqref{sasi} we deduce
\be\ba{ll}\medskip
\|\big\{\big[(\tilde{\bfv}-\tilde{\bfv}_0)-(\tilde{\bfxi}-\tilde{\bfxi}_0)\big]\cdot\mathbb A\cdot\nabla\bfB(\bfv)\big\}_t\|_2^2\\ \medskip
\hspace*{4cm}\le c\,\big[\big(\|\mathbb D(\tilde{\bfv}-\tilde{\bfv}_0)\|_2^2+\|\mathbb D(\tilde{\bfv}_t-\tilde{\bfv}_{0t})\|_2^2\big)\left(\|\mathbb D(\bfv)\|_2^2+\|\mathbb D(\bfv)\|_2\|D^2\bfv\|_2\right)\\
\hspace*{4.9cm}+\|\mathbb D(\tilde{\bfv}-\tilde{\bfv}_0)\|_2^2\left(\|\mathbb D(\bfv_t)\|_2^2+\|D^2\bfv_t\|_2^2\right)\big]\,.
\ea
\eeq{8.12}
Finally,  \Eqref{8.6}--\Eqref{8.7},  \Eqref{Sibe1} and \Eqref{SiBe}, likewise imply 
\be\ba{ll}\medskip
\|\big\{\bfB(\tilde{\bfv}-\tilde{\bfv}_0)\cdot\mathbb A\cdot \nabla\bfB(\bfv)\big\}_t\|_2^2\\ \medskip
\hspace*{4cm}\le c\,\big[\big(\|\mathbb D(\tilde{\bfv}-\tilde{\bfv}_0)\|_2^2+\|\mathbb D(\tilde{\bfv}_t-\tilde{\bfv}_{0t})\|_2^2\big)\left(\|\mathbb D(\bfv)\|_2^2+\|\mathbb D(\bfv)\|_2\|D^2\bfv\|_2\right)\\
\hspace*{4.9cm}+\|\mathbb D(\tilde{\bfv}-\tilde{\bfv}_0)\|_2^2\left(\|\mathbb D(\bfv_t)\|_2^2+\|D^2\bfv_t\|_2^2\right)\big]\,.
\ea
\eeq{8.13}
By Proposition \ref{EUlin} and \Eqref{BaNa2} and \lemmref{8.2}(c),
we have
\be
\max_{t\in[0,T]}\|\mathbb D(\bfv_t(t))\|_2+\|D^2\bfv_t\|_{L^2(L^2)}\le C_1
\eeq{8.14}
where $C_1$ depends only on $\varOmega$, $T$ and the physical parameters.  As a result, collecting \Eqref{8.5}, \Eqref{8.11}--\Eqref{8.13} we gather
\be
\|\bfcaln_t\|_{L^2(L^2)}\le c\,\|\tilde{\bfv}-\tilde{\bfv}_0\|_{W^{1,4}(D^{1,2})}\,.
\eeq{8.15}
We now apply Proposition \ref{EUlin} to problem \Eqref{8.4} with $\bff:=\delta^{-2}\bfcaln$ and $\bfF\equiv\0$. By \Eqref{8.10} and \Eqref{8.15}, and the hypothesis on $\delta$, the assumptions of that proposition are satisfied. In particular, setting
$$
\caln_1:=\|\bfcaln\|_2^2+\|\bfcaln\|_{\frac65}^2\,,\ \ \caln_2:=\|\partial_t\bfcaln\|_2^2\,;\ \ N:=\Big[\int_0^T(\caln_1+\caln_2+\caln_1^2){\rm d}t\Big]^\frac12
$$
from \Eqref{BaNa1}, \Eqref{BaNa2} and \lemmref{8.2}(c), we get   
$$
\|\bfw\|_{W^{1,2}(D^{1,2})}^2\le c\,N^2\,,
$$
which, in turn, recalling the definition of $\bfw$ and using \Eqref{8.10} and \Eqref{8.15} furnishes
$$
\|\Phi(\tilde{\bfv})-\Phi(\tilde{\bfv}_0)\|_{W^{1,4}(D^{1,2})}^2\le c\,\left(\|\tilde{\bfv}-\tilde{\bfv}_0\|_{W^{1,4}(D^{1,2})}^2+\|\tilde{\bfv}-\tilde{\bfv}_0\|_{W^{1,4}(D^{1,2})}^4\right).
$$
This proves the continuity of the map $\Phi$ and concludes the proof of the proposition.
\par\hfill$\square$\par 
\setcounter{equation}{0}
\section{Proof of Theorem \ref{theo:6.1}}\label{nonl} With \propref{8.1} in hand, the proof is obtained by following a classical procedure \cite{GaKyH,GaSi2} that, for completeness, we sketch  here.
Let $(\bfv_k,{\sf q}_k,\bfxi_k)$ be the $T$-periodic solution to  \Eqref{6.3} determined in \propref{8.1}. Then, from \Eqref{belest} we deduce the following estimate valid for each $m<k$
\be
\|\bfv_k\|_{\mathscr W_T(\varOmega_m)}+\|{\sf q}_k\|_{\calq(\varOmega_m)}+\|\bfxi_k\|_{W^{2,2}(0,T)}\le C_0\,\delta^2\,.
\eeq{10.1}
As a result, for any fixed $m$, there is a subsequence of $(\bfv_k,{\sf q}_k,\bfxi_k)$ (that we continue to denote by the same symbols) 
converging weakly in the space $\mathscr W_T(\varOmega_m)\times\calq({\varOmega}_m)\times W_T^{2,2}(\real^3)$ to one of its elements. Furthermore, by known compactness theorems, $\{\bfv_k,\bfxi_k\}$ converge strongly in $W^{1,2}(\cali_T; W^{1,2}(\varOmega_m))\times W^{1,2}(\cali_T;\real^3)$, for instance, and $\{(\bfv_k,{\sf q}_k)\}$ converge strongly in $L^2(\cali_T; W^{1,2}(\varSigma))\cap L^2(\cali_T; L^{2}(\varSigma))$, with $\cali_T\subset\real$ arbitrary interval of length $T$. Recalling \Eqref{invdo} and using Cantor diagonalization procedure, we may thus select a subsequence, again denoted by $(\bfv_k,{\sf q}_k,\bfxi_k)$, satisfying the above convergence properties for every $m\in\nat$. Denote by $(\bfv,{\sf q},\bfxi)$ the weak limit of this subsequence. Then, from \Eqref{10.1} we obtain at once
\be
\|\bfv\|_{\mathscr W_T(\varOmega)}+\|{\sf q}\|_{\calq(\varOmega)}+\|\bfxi\|_{W^{2,2}(0,T)}\le C_0\,\delta^2\,.
\eeq{10.2}
It is also easy to show that $(\bfv,{\sf q},\bfxi)$ satisfies \Eqref{6.10}--\Eqref{6.12}. In fact, let $\psi=\varphi(x)\sigma(t)$, with $\varphi\in C_0(\varOmega)$, $\sigma\in  C_0(\real)$ arbitrary. From \Eqref{6.3} we thus get  
\be\ba{cc}\medskip\ba{ll}\medskip
\Int{\real}{}\Big\{\big(\partial_t{\bfv_k}+\partial_t\bfB(\bfv_k)+(\bfv_k+\bfB(\bfv_k)-\bfxi_k)\cdot\nabla\bfv_k+(\bfv_k-\bfxi_k)\cdot\nabla\bfB(\bfv_k)+\delta(\bfV-\bfzeta)\cdot\nabla(\bfv_k+\bfB(\bfv_k))\\ 
\medskip
\hspace*{1.3cm}
+\delta(\bfv_k+\bfB(\bfv_k)-\bfxi_k)\cdot\nabla\bfV
-\nu\Delta\bfv_k+\nabla{\sf q_k}-\bfpzc{f}_\ell(\bsf{s},\bfv_k,{\sf q}_k,\bfxi_k)-
\bfpzc{f}_{n\ell}(\bsf{s},\bfv_k,\bfxi_k)+\delta^2\bff_{\mbox{\tiny $\bfV$}},\psi\big)\Big\}{\rm d}t
\ea\,,\\  
{\sf m}\Int{\real}{}\dot{\bfxi_k}(t)\sigma(t){\rm d}t=-\Int{\real}{}\Big(\sigma(t)\Int{\varSigma}{}\mathbb T(\bfv_k,{\sf q}_k)\cdot\bfn{\rm d}\varSigma\Big){\rm d}t+ \Int{\real}{}\left(\bfpzc{F}(\bsf{s},\bfv_k,{\sf q}_k)+\delta^2\bfF_{\mbox{\tiny $\bfV$}}\right)\sigma(t){\rm d }t\,,\ea
\eeq{10.3}
where $(\cdot,\cdot):=(\cdot,\cdot)_\varOmega$.  
Using the above-mentioned convergence properties of the sequence $(\bfv_k,{\sf q}_k,\bfxi_k)$ along with the properties of the operator $\bfB$, we can pass to the limit $k\to\infty$ in \Eqref{10.3} and show, by a routine reasoning, that the limit functions $(\bfv,{\sf q},\bfxi)$ continue to satisfy \Eqref{10.3}. This, by the arbitrariness of $\varphi$ and $\sigma$, implies by a standard procedure that, in fact, $(\bfv,{\sf q},\bfxi)$ is a solution to  \Eqref{6.10}--\Eqref{6.12}. Finally, setting $\bsf{u}:=\bfv+\bfB(\bfv)$, in view of the properties of $\bfB$ and the argument presented in Section \ref{sub:7.1}, we conclude that $(\bsf{u},{\sf q},\bfxi)$ satisfies all the properties stated in Theorem  \ref{theo:6.1}, which is thus completely proved. 
\setcounter{equation}{0}
\section{Conditions for Self-Propulsion}\label{asy}
The objective of this section is  to provide the complete characterization of the average velocity, $\bar{\bfgamma}$, of the center of mass $G$ at the order of $\delta^2$, in an appropriate class of solutions. To this end, we give the following definition.
\Bd
A solution $(\bsf u,{\sf q}, \bfxi)$ to  \Eqref{6.1}--\Eqref{6.2} belongs to the class $\mathscr C$ if it meets the following properties:
\begin{itemize}
\item[{\rm (a)}] \ $(\bsf u,{\sf q}, \bfxi)\in \mathscr W_T(\varOmega)\times\calq({\varOmega})\times W_T^{2,2}(\real^3)$\,;
\item[{\rm (b)}]\  $(\bsf u,{\sf q}, \bfxi)$ satisfies the estimate \Eqref{goest}\,.  
\end{itemize}
\EED{11.1}
Theorem \ref{theo:6.1} ensures that $\mathscr C\neq\emptyset$.
Our  goal develops according to the following three logical steps. In the first one, we show that in the limit $\delta\to 0$, the average  $(\bar{\bsf{u}},\bar{\sf q},\bar{\bfxi})$ of any corresponding  solution in the class $\mathscr C$   tends, after suitable rescaling, to the unique solution,  $(\bsf v_0,{\sf r}_0,\bfsigma_0)$, of a Stokes-like (linear) steady-state problem; see \Eqref{11.2}. In the second step, we furnish necessary and sufficient conditions for $\bfsigma_0\neq \0$. In the third and final step, we combine these results with those of \propref{5.1} to provide an explicit characterization of the average velocity of $G$, $\bar{\bfgamma}$, up to the order of $\delta^2$; see \theoref{11.1}.     
\smallskip\par
Let $(\bfV_0,{\sf p}_0,\bfzeta_0)$ be the solution to the problem \Eqref{5.1} given in \lemmref{5.1}, and set
\be\ba{rl}\medskip
\bsf g_0(x,t)&\!\!\!\!:=(\nabla\bsf s)^\top:\nabla\, \mathbb T(\bfV_0,{\sf p}_0)+\nu\,\Div\big(\nabla\bsf{s}\cdot\nabla\bfV_0+(\nabla\bsf{s}\cdot\nabla\bfV_0)^\top\big)+(\bfV_0-\bfu_*-\bfzeta_0)\cdot\nabla\bfV_0\,,\\
\bsf G_0(t)&\!\!\!\!:=-\Int{\varSigma}{}\big[\mathbb T(\bfV_0,{\sf p}_0)\cdot\mathbb H^\top(\bsf{s})-\nu\,\big(\nabla\bsf{s}\cdot\nabla\bfV_0+(\nabla\bsf s\cdot\nabla\bfV_0)^\top\big)\big]\cdot\bfn\,{\rm d}\varSigma\,,
\ea
\eeq{11.1}
where $\mathbb H$ is defined in \Eqref{Hh}. Notice that both $\bsf{g}_0$ and $\bsf{G}_0$ are known  functions, depending only on $\bfu_*$, the physical parameters $\nu,{\sf m}$, the reference configuration $\Omega_0$ and the  displacement field $\bsf{s}$. Next, consider the boundary-value problem:
\be\ba{cc}\medskip\left.\ba{ll}\medskip
\Div\mathbb T({\bsf v}_0,{\sf r}_0)
=\bar{\bsf g_0}\\
\Div{{\bsf v}_0}=0\ea\right\}\ \ \mbox{in $\varOmega$}\,,\\ \medskip
{{\bsf v}_0}= {\bfsigma_0} \ \ \mbox{at $\varSigma$}\,;\ \ \Lim{|x|\to\infty}{{\bsf v}_0}(x)=\0\,;\\ 
\Int{\varSigma}{}\mathbb T({{\bsf v}_0},{{\sf r}_0})\cdot\bfn{\rm d}\varSigma=\bar{\bsf G_0}\,.
\ea
\eeq{11.2}
\Bd  A vector field ${\bsf v}_0$ is a weak solution to problem \Eqref{11.2} if {\rm (i)} ${\bsf v}_0\in\cald^{1,2}(\real^3)$, and {\rm (ii)} it satisfies the condition
\be
\big(\mathbb D({\bsf v}_0),\mathbb D(\bfphi)\big)=-(\bar{\bsf g_0},\bfphi)+\bar{\bsf G_0}\cdot\hat{\bfphi}\,,\ \ \mbox{for all $\bfphi\in\calk(\real^3)$}\,.
\eeq{11.5}
\EED{11.2}
\Br The equation in \Eqref{11.5} is formally obtained by dot-multiplying both sides of \Eqref{11.2}$_1$ by $\bfphi$, integrating by parts over $\varOmega$ and then taking into account \Eqref{11.2}$_5$. It is readily proved that, if ${\bsf v}_0$ is a weak solution that is also in $W^{2,2}(\varOmega_R)$ for all $R>R_*$, then there exists ${\sf r}_0\in W^{1,2}({\varOmega}_R)$ such that $({\bsf v}_0, {\sf r}_0,\bfsigma_0\equiv{\bsf v}_0|_{\Omega_0})$ is a solution to \Eqref{11.2}; see \cite[pp. 699-700]{Gah}.
\ER{11.1}
The following result holds.
\Bl Problem \Eqref{11.2} has one and only one weak solution. Furthermore,
$$
({\bsf v}_0,{\sf r}_0,\bfsigma_0)\in [L^{6}(\varOmega)\cap D^{1,2}(\varOmega)\cap D^{2,2}(\varOmega)]\times D^{1,2}(\varOmega)\times\real^3\,.
$$
\EL{11.1}
{\em Proof.}   Employing the properties of $\bsf s$ in \Eqref{Esse1} and those of $(\bfV_0,{\sf p}_0)$ in \lemmref{5.1}, we deduce
\be
\|\bar{\bfu_*\cdot\nabla\bfV_0-\nabla\bsf s:\nabla\, \mathbb T(\bfV_0,{\sf p}_0)}-\nu\,\Div\big(\bar{\nabla\bsf{s}\cdot\nabla\bfV_0+(\nabla\bsf{s}\cdot\nabla\bfV_0)^\top}\big)\|_{\frac65}\le c\,\left(\|\bfV_0\|_{\calw(\varOmega)}+\|{\sf p}_0\|_{\calq(\varOmega)}\right)\le C\,.
\eeq{11.6}
Likewise, by integration by parts, \Eqref{1.8} and again  \lemmref{5.1},
\be\ba{ll}\medskip
\big|\big(\bar{(\bfV_0-\bfzeta_0)\cdot\nabla\bfV_0},\bfphi\big)\big|=
\big|\big(\bar{(\bfV_0-\bfzeta_0)\otimes\bfV_0},\nabla\bfphi\big)\big|\\
\hspace*{3cm}\le c\,\left(\|\bfV_0\|^2_{L^4(L^4)}+\|\bfzeta_0\|_{W^{1,2}(0,T)}\|\bfV_0\|_{L^2(L^2)}\right)\|\nabla\bfphi\|_2\le C\,\|\mathbb D(\bfphi)\|_2\,,\ \  \bfphi\in\cald^{1,2}(\real^3)\,.
\ea
\eeq{11.7}
Using one more time \Eqref{Esse1},  \lemmref{5.1} and classical trace inequalities, we also show
\be
|\bar{\bsf G_0}|\le c\,\big(\|\bfV_0\|_{\calw({\varOmega})}+\|{\sf p}_0\|_{\calq({\varOmega})}\big)\le C\,.
\eeq{11.8}
As a result, from \Eqref{11.6}, \Eqref{11.7} and \Eqref{1.8}, \Eqref{1.9} we infer 
$$
|(\bar{\bsf g_0},\bfphi)|\le c\,(\|\bfphi\|_6+\|\nabla\bfphi\|_2)\le c\,\|\mathbb D(\bfphi)\|_2\,, 
$$
as well as, from \Eqref{11.8} and \Eqref{1.10}
$$
|\bar{\bsf G_0}\cdot\hat{\bfphi}|\le c\,|\hat{\bfphi}|\le c\,\|\mathbb D(\bfphi)\|_2\,.
$$
The last two displaced inequalities then entail that the right-hand side of \Eqref{11.5} defines a bounded linear functional on the Hilbert space $\cald^{1,2}(\real^3)$. Therefore, by Riesz theorem, there is one and only one ${\bsf v}_0\in\cald^{1,2}(\real^3)$ satisfying \Eqref{11.5}, which is equivalent to the existence of a unique weak solution. Moreover, from well-known regularity results \cite[Section V.4]{Gab}, we infer ${\bsf v}_0\in D^{2,2}({\varOmega})\cap \cald^{1,2}(\real^3)$. This furnishes, in particular, ${\bsf v}_0\in W^{2,2}(\varOmega_R)$, all $R>R_*$, which, by \remref{11.1}, completes the proof of the lemma. \par\hfill$\square$\par
We now introduce the vector $\bfpzc G_0\in\real^3$ with components
\be\ba{ll}\medskip
\mathpzc G_{0i}:=-\!\!\Int{\varOmega}{}\big\{[\bar{(\bfV_0-\bfu_*-\bfzeta_0)\cdot\nabla\bfV_0}+\bar{(\nabla\bsf s)^\top:\nabla\, \mathbb T(\bfV_0,{\sf p}_0)}]\cdot\bfh^{(i)}\\
\hspace*{2.1cm}-\nu\,\big(\bar{\nabla\bsf{s}\cdot\nabla\bfV_0+(\nabla\bsf{s}\cdot\nabla\bfV_0)^\top}\big):\mathbb D(\bfh^{(i)})\big\} 
-\bfe_i\cdot\Int{\varSigma}{}\bar{\mathbb T(\bfV_0,{\sf p}_0)\cdot\mathbb H^\top(\bsf{s})}\cdot\bfn{\rm d}\varSigma\,,\ \ i=1,2,3\,,
\ea
\eeq{G0}
where $\{\bfh^{(i)},p^{(i)}\}$ are the solutions to the Stokes problem \Eqref{4.32}. Notice that, taking into account the summability properties of the fields $\bfh^{(i)}$ and $\bfV_0$ given in  \Eqref{4.33} and \Eqref{5.3}, respectively, and the fact that both $\bfu_*$ and $\nabla{\bsf s}$ have bounded support, we obtain that the volume integral is well defined.
The following result holds.
\Bl Let $(\bsf v_0,\sf r_0,\bfsigma_0)$ be the solution to \Eqref{11.2} determined in \lemmref{11.1}. Then, $\bfsigma_0\neq 0$ if and only if $\bfpzc G_0\neq \0$. Precisely, we have 
\be
\bfsigma_0=\mathbb M^{-1}\cdot \bfpzc G_0\,,
\eeq{11.4}
with the matrix $\mathbb M$  defined in \Eqref{Matrix}\,.
\EL{11.2}
{\em Proof.} We dot-multiply \Eqref{11.2}$_1$ by $\bfh^{(i)}$ and integrate by parts over $\varOmega$, to get
\be
\big(\mathbb D({\bsf v}_0),\mathbb D(\bfh^{(i)})\big)=-(\bar{\bsf g}_0,\bfh^{(i)})+\bar{\bsf G}_0\cdot\bfe_i\,,\ \ i=1,2,3\,.
\eeq{11.9}
Analogously, if we dot-multiplying \Eqref{4.32}$_1$ by ${\bsf v}_0$ and integrate by parts over $\varOmega$, we deduce
\be
\big(\mathbb D({\bsf v}_0),\mathbb D(\bfh^{(i)})\big)=\mathbb M_{ji}\sigma_{0i}\,,\ \ i=1,2,3\,.
\eeq{11.10}
Thus, by operating with a further integration by parts in the integral on the right hand side of  \Eqref{11.9} and then combining the resulting equation with \Eqref{11.10} , we arrive at \Eqref{11.4}, thus completing the proof.\par\hfill$\square$\par
Let $(\bsf u,{\sf q},\bfxi)\in \mathscr C$, and define $(\bsf v,{\sf r},\bfsigma)$ as follows:
\be
\bsf u=\delta^2\bsf v\,,\ \ {\sf q}=\delta^2{\sf r}\,,\ \ \bfxi=\delta^2\bfsigma\,.
\eeq{11.11}
The following result holds.
\Bl Let $\delta_0$ be as in {\rm Theorem} {\rm \ref{theo:6.1}}. Then, for all $\delta\in (0,\delta_0)$, $\bar{\bfsigma}=\bar{\bfsigma}(\delta)$ satisfies
\be
\bar{\bfsigma}=\bfsigma_0+o\,(\delta)\,.
\eeq{11.13}
Thus, in particular
\be
\half|\bfsigma_0|\le |\bar{\bfsigma}|\le \mbox{$\frac32$}|\bfsigma_0|\,.
\eeq{11.14}
\EL{11.3}
{\em Proof.} For given $\delta\in (0,\delta_0)$ we denote by $(\bsf u(\delta),{\sf q}(\delta),\bfxi(\delta))$ the corresponding solution to \Eqref{6.1}, \Eqref{6.2} and by $(\bsf v(\delta),{\sf r}(\delta),\bfsigma(\delta))$ the associated rescaled solutions according to \Eqref{11.11}.
Let $\{\delta_k\}\subset (0,\delta_0)$ be an arbitrary, vanishing sequence.  We thus deduce that the sequence $\{(\bsf v_k,{\sf r}_k,\bfsigma_k)\}\equiv \{(\bsf v({\delta_k}),{\sf r}({\delta_k}),\bfsigma({\delta_k}))\}$ obeys the following problem
\be\ba{cc}\medskip\left.\ba{ll}\medskip
\Div\mathbb T(\bar{\bsf v}_k,\bar{\sf r}_k)=\delta_k\,\bfR_1{(\delta_k)}+\bfR_2(\delta_k)+\bfR(\delta_k)+{\bar{\bsf g}_0}\\
\Div\bar{\bsf v}_k=\bfR_3(\delta_k)\ea\right\}\ \ \mbox{in $\varOmega$}\\ \medskip
\bar{\bsf v}_k=\bar{\bfsigma}_k\ \ \mbox{at $\varSigma$}\,,\\
\Int{\varSigma}{}\mathbb T(\bar{\bsf v}_k,\bar{\sf r}_k)\cdot\bfn{\rm d}\varSigma=\bfS(\delta_k)+{\bar{\bsf G}_0}\,,
\ea
\eeq{11.15}
where
$$\ba{rl}\medskip
{\bfR_1(\delta_k)}&\!\!\!\!:=\delta_k\bar{(\bsf v_k-\bfsigma_k)\cdot\nabla\bsf v_k}
\\\medskip
\bfR_2(\delta_k)&\!\!\!\!:=-\bar{\bsf f_\ell(\bsf s,\bsf v_k,{\sf r}_k,\bfsigma_k)}-\delta^2_k\bar{\bff_{n\ell}(\bsf s,\bsf v_k,\bfsigma_k)}\\\medskip
\bfR_3(\delta_k)&\!\!\!\!:=-\Div(\bar{\mathbb C^\top(\delta_k)\cdot\bsf v_k})
\\
\bfR(\delta_k)&\!\!\!\!:=-\bar{\bff_{\mbox{\tiny$\bfV$}}}-\bar{\bsf g}_0\,;\ \ \bfS(\delta_k):={\bar{\bsf F(\bsf{s},\bsf v_k,{\sf r}_k)}}+\bar{\bfF_{\mbox{\tiny$\bfV$}}}-\bar{\bsf G}_0\,.
\ea
$$
From \Eqref{goest}, we also infer
\be
\|\bsf v_k\|_{\mathscr W_T(\varOmega)}+\|{\sf r}_k\|_{\calq(\varOmega)}+\|\bfsigma_k\|_{W^{2,2}(0,T)}\le C\,,
\eeq{11.16}
with $C$ independent of $k\in\nat$. This relation combined with \remref{norm} implies, in particular, the existence of $\bsf v_*\in D^{1,2}(\real^3)$ with $\bsf v_*|_{\Omega_0}=\bfsigma_*$ such that (along a subsequence)
\be
\nabla\bar{\bsf v}_k\to\nabla\bsf v_*\ \ \mbox{weakly in  $L^{2}(\real^3)$}\,;\ \ \bar{\bfsigma}_k\to\bfsigma_* \ \ \mbox{in $\real$.}\,.
\eeq{11.17}
We next dot-multiply both sides of \Eqref{11.15}$_1$ by $\bfphi\in\calk(\real^3)$, integrate by parts over $\varOmega$ and use \Eqref{11.15}$_4$ to show
\be
\big(\mathbb D(\bar{\bsf v}_k),\mathbb D(\bfphi)\big)=-(\bar{\bsf g},\bfphi)+\bar{\bsf G}\cdot\hat{\bfphi}-\big(\delta_k{\bfR_1(\delta_k)}+\bfR_2(\delta_k)+\bfR(\delta_k),\bfphi\big)-\bfS(\delta_k)\cdot\hat{\bfphi}\,,\ \ \mbox{for all $\bfphi\in\calk(\real^3)$}\,.
\eeq{11.18}
Employing \Eqref{11.16} along with H\"older inequality, it is not difficult to show that
\be
\|\bfR_1\|_{1,\varOmega_R}\le C_1\,,\ \ \mbox{all $R>R_*$}
\eeq{11.19}
where, here and in the rest of the proof, $C_i$, $i=1,2,\ldots$ denotes a constant possibly depending on $R$ but not on $k$. Furthermore, by using \lemmref{Def}, \Eqref{2.13} and \Eqref{11.16} we get
\be
\|\bfR_2(\delta)\|_{1,\varOmega_R}\le C\delta\,,\ \ \mbox{all $R>R_*$}.
\eeq{nova}
Moreover, setting
$$
\bfU:=\bfV-\bfV_0\,,\ \ p:={\sf p}-{\sf p}_0\,,\ \ \bfz=\bfzeta-\bfzeta_0\,,
$$
from \Eqref{6.2}$_2$ and \Eqref{11.1}$_1$ we have
\be
-\bff_{\mbox{\tiny$\bfV$}}-\bsf g_0=-\delta^{-1}\bff_\ell(\bsf s,\bfU,p)+\delta^{-1}\bff_{\mbox{\tiny$\bfV$}}^{(1)}+\bff_{\mbox{\tiny$\bfV$}}^{(2)}\,,
\eeq{11.20}
where
$$\ba{rl}\medskip
\bff_{\mbox{\tiny$\bfV$}}^{(1)}:=&\!\!\!\!-\delta\bfu_*\Cdot\mathbb B\Cdot\nabla\bfV_0+(\mathbb B+\delta\nabla\bsf s):\nabla\mathbb T(\bfV_0,{\sf p}_0)+\nu\mathbb A:\nabla\big [(\mathbb B+\delta\nabla\bsf s)\cdot\nabla\bfV_0+((\mathbb B+\delta\nabla\bsf s)\cdot\nabla\bfV_0)^\top\big]
\\
\bff_{\mbox{\tiny$\bfV$}}^{(2)}:=&\!\!\!\!-(\bfU-\bfz)\cdot\mathbb A\cdot\nabla\bfV-(\bfV-\bfzeta_0)\cdot\mathbb B\cdot\nabla\bfV_0-(\bfV_0-\bfzeta_0)\cdot\mathbb A\cdot\nabla\bfU\,.
\ea
$$
With the help of \lemmref{Def}, \lemmref{fF} and \lemmref{5.1}, it is not difficult to show that
$$
\delta^{-1}\left\|\bar{\bff_\ell(\bsf s,\bfU,p)}\right\|_{1,{\varOmega}_R}+\delta^{-1}\left\|\bar{\bff_{\mbox{\tiny$\bfV$}}^{(1)}}\right\|_{1,{\varOmega}_R}+\left\|\bar{\bff_{\mbox{\tiny$\bfV$}}^{(2)}}\right\|_{1,{\varOmega}_R}\le C_2\,\delta\,,\ \ \mbox{for all $R>R_*$}\,,
$$
which, by \Eqref{11.20}, thus implies
\be
\|\bfR(\delta)\|_{1,{\varOmega}_R}\le C_2\,\delta\,.
\eeq{11.21}
Likewise, from \Eqref{6.2}$_2$ and  \Eqref{11.1}$_2$ we deduce 
$$\ba{rl}\medskip
\bfF_{\mbox{\tiny $\bfV$}}-\bsf G_0=&\!\!\!\!\delta^{-1}\bfF(\bsf s,\bfU,p)-\delta^{-1}\Int{\varSigma}{}\mathbb T(\bfV_0,{\sf p}_0)\cdot(\mathbb C-\delta\mathbb H^\top(\bsf s))\cdot\bfn{\rm d}\varSigma\\\medskip&\!\!\!\!-\nu\delta^{-1}\Int \varSigma{}
J\,\big[(\mathbb B+\delta\nabla\bsf s)\cdot\nabla\bfV_0+(\nabla\bfV_0)^\top\cdot(\mathbb B+\delta\nabla\bsf s)^\top\big]\cdot\mathbb A\cdot\bfn\,{\rm d}\varSigma\\
&\!\!\!\!+\nu\Int \varSigma{}
\big[\nabla\bsf s\cdot\nabla\bfV_0+(\nabla\bfV_0)^\top\cdot(\nabla\bsf s)^\top\big]\cdot\mathbb C\cdot\bfn\,{\rm d}\varSigma\,.
\ea  
$$
Therefore, using again \lemmref{Def},  \lemmref{fF} and \lemmref{5.1}
\be
|\bfS(\delta)|\le C_3\delta\,.
\eeq{11.22}
Finally, for all $\varphi\in C_0^\infty(\varOmega)$, after integrating by parts, with the help of \lemmref{Def} and \Eqref{11.14}, we show 
\be
|(\bfR_3(\delta),\bfphi)|=|(\mathbb C^\top(\delta)\cdot\bsf v_k,\nabla\varphi)|\le c\,\delta\,\|\nabla\varphi\|_2
\eeq{11.22_1}
If we now pass to the limit $k\to\infty$ into \Eqref{11.18} and use \Eqref{11.17}, \Eqref{11.19}, \Eqref{nova}, \Eqref{11.21},   \Eqref{11.22} and \Eqref{11.22_1} we deduce 
\be
\big(\mathbb D({\bsf v}_*),\mathbb D(\bfphi)\big)=-(\bar{\bsf g}_0,\bfphi)+\bar{\bsf G}_0\cdot\hat{\bfphi}\,,\ \ \mbox{for all $\bfphi\in\calk(\real^3)$}\,;\ \ \Div\bsf v_*=0\,,
\eeq{11.23}
that is, $\bsf v_*$ is a weak solution to \Eqref{11.5}. However, from \lemmref{11.1} we know that this solution is unique, so that we must have $\bsf v_*\equiv \bsf v_0$, $\bsf v_*|_{\Omega_0}\equiv \bfsigma_0$, as well as  the convergence given in \Eqref{11.17} occurs not just along a sequence but as long as $\delta\to0$:
\be   
\bar{\bsf v}({\delta})\to \bsf v_0\ \ \mbox{as $\delta\to 0$}\,,\ \ \mbox{weakly in $\cald^{1,2}(\real^3)$}\,.
\eeq{11.24_0}
Next, if we dot-multiply both sides of \Eqref{4.32}$_1$ by $\bar{\bsf v}$ and integrate by parts over $\varOmega$, we find
$$
\big(\mathbb D(\bar{\bsf v}),\mathbb D(\bfh^{(i)})\big)=\mathbb M_{ji}\bar{\sigma}_{i}\,,\ \ i=1,2,3\,.
$$
Thus, subtracting side-by-side
\Eqref{11.10} from the latter, it follows that
$$
\mathbb M_{ji}(\bar{\sigma}_i({\delta})-\bfsigma_{0i})=\big(\mathbb D(\bar{\bsf v}_i(\delta)-\bsf v_0),\mathbb D(\bfh^{(i)})\big)\,.
$$
The statement in the lemma is thus a consequence of the last displayed equation and \Eqref{11.24_0}.\par\hfill$\square$\smallskip\par
Let
$$
\bfpzc G=\bfpzc G_0+\bfpzc G_1
$$
where $\bfpzc G_0$ and $\bfpzc G_1$ are given in \Eqref{G0} and \Eqref{5.6_1}, respectively. 
The following result, representing the main achievement of this paper, provides a complete characterization of $\bar{\bfgamma}$ up to the order of $\delta^2$ for solutions in the class $\mathscr C$.
\Bt Let $\delta\in(0,\delta_0)$ with $\delta_0$ as in \lemmref{11.1}, and suppose that $\bfu_*$ satisfies the assumption of \propref{O_k}. Let $(\bfu,p,\bfgamma)$ be a corresponding $T$-periodic solution to \Eqref{2.12}--\Eqref{2.13}  that, without loss of generality, can be written as
\be 
\bfu={\bsf u}+\delta\,\bfV\,, \ \ p={\sf q}+\delta\,{\sf p}\,,\ \ \bfgamma=\bfxi+\delta\,\bfzeta\,,
\eeq{11.24}
with $(\bfV,{\sf p},\bfzeta)$ given in \propref{O_k} and $(\bsf u,{\sf q},\bfxi)$ solution to the problem  \Eqref{6.10}--\Eqref{6.12}. Then, for any $(\bsf u,{\sf q},\bfxi)\in\mathscr C$ we have
\be
\bar{\bfgamma}=\delta^2\bfgamma_1+o\,(\delta^2)
\eeq{11.25}
with
$$
\bfgamma_1:=\mathbb M^{-1}\cdot\bfpzc G\,.
$$
Therefore, at the order of $\delta^2$, $\bar{\bfgamma}\neq\0$ if and only if $\bfpzc G\neq\0$.  
\ET{11.1}
{\em Proof.} Under the given assumptions, from \propref{5.1} we have
$$
\bar{\bfzeta}=\delta\,\mathbb M^{-1}\cdot\bfpzc G_1+O(\delta^2)\,,
$$
see \Eqref{5.5},
whereas, from \lemmref{11.2}, it follows that
$$
\bar{\bfxi}=\delta^2\mathbb M^{-1}\cdot\bfpzc G_0+o\,(\delta^2)\,,
$$
see \Eqref{11.4}. The proof of \Eqref{11.25} then follows from the last two displaced equations and \Eqref{11.24}$_3$.\par\hfill$\square$\par
\setcounter{equation}{0}
\section{An Application}\label{section:Application}
In this final section we will furnish an example where \theoref{11.1} is able to provide a {\em quantitative} evaluation of the thrust vector $\bfpzc G$ and, hence, of the associated propulsion velocity $\bar{\bfgamma}$ at the order of $\delta^2$. 
\par
Before proceeding in this direction, however, we offer some general comments regarding this evaluation. We recall (see \lemmref{1}) that the field $\bsf s$ is defined as $\bsf s(x,t):=\beta(x)\check{\bsf s}(x,t)$ where $\check{\bsf s}(x,t)$ is a suitable extension of the given displacement field $\hat{\bsf s}(x,t)$ and $\beta(x)$ is a smooth  ``cut-off" that is 1 in  $B_R\supset\bar{\Omega_0}$, and 0  in $\real^3\backslash\bar{B_{2R}}$.  
From 
\Eqref{G0} and \Eqref{5.6_1} it follows that 
$\bfpzc G$ depends linearly $\Div\bsf s$ and $\nabla\bsf s$. Thus, choosing $R:=c_0\delta^{-1}$, for appropriate choice of the constant $c_0$, by means of the properties \Eqref{beta} we easily show that  
$$
\bfpzc G(\bsf s)=\bfpzc G(\check{\bsf s})+o(1)\ \ \mbox{as $\delta\to 0$}\,.
$$
This means, in particular, that, in order to explicitly compute the thrust $\bfpzc G$, it is enough to perform this computation along the extension $\check{\bsf s}$.
\par
To simplify the notation, we will continue to denote this extension with $\bsf s$.
\par
In our example we take, as reference configuration,  the ball of radius $a$  centered at the origin. Denoting by $\omega$ the frequency of the oscillations, we scale velocities by $\omega\,a$, time by $\omega^{-1}$, length by $a$ and pressure by $\nu\,\omega$. In this way,  problem \Eqref{5.1} becomes, in dimensionless form,  
\be\ba{cc}\medskip\left.\ba{ll}\medskip
2h^2\,\partial_t\bfV_0=\Delta\bfV_0-\nabla p_0\\
\Div\bfV_0=0\ea\right\}\ \ \mbox{in $ \varOmega\times\real$}\\ \medskip
\bfV_0(x,t)=\bfu_{*}(x,t)+\bfzeta_0(t)\ \ (x,t)\in\varSigma\times\real\\
{\sf M}\dot{\bfzeta}_0+\Int\varSigma{}\mathbb T(\bfV_0,p_0)\cdot\bfn=\0\ \ \mbox{in $\real$}\,.
\ea
\eeq{V0}
where
$$
\varOmega=\{x\in\real^3:\ r:=|x|>1\}\,,\ \ \varSigma=\partial{\varOmega}\,,\ \ {\sf M}=2\frac{\sf m}{a^3}\,h^2\,, 
$$
and
$$
\mathbb T(\bfV_0,p_0)=\nabla\bfV_0+(\nabla\bfV_0)^\top-p_0\mathbb I:=2\mathbb D(\bfV_0)-p_0\mathbb I\,.
$$
Moreover,
\be
h=\sqrt{\frac{\omega}{2\nu}}a
\eeq{STOK}
denotes the Stokes number.
Concerning the solutions to the auxiliary problems \Eqref{4.32}, we scale $\bfh^{(i)}$ by 1 and $p^{(i)}$ by $\nu\,a^{-1}$. In doing so, we deduce \cite{Gasp},
$$
\mathbb M^{-1}=(6\pi)^{-1}\mathbb I\,.
$$ 
Taking into account all of the above and recalling \Eqref{5.6_1} and \Eqref{G0}, the propulsion  velocity $\bfgamma_1$ given in \theoref{11.1} 
takes the following non-dimensional form form
$$
\bfgamma_1=\frac{1}{6\pi}\left(\bfpzc G_0+\bfpzc G_1\right)
$$
with (in non-dimensional form)
\be
\bfpzc G_1= -\int_\varOmega \bar{(\nabla\bsf s)^\top:\nabla\bfV_0}\,p^{(i)}\bfe_i
\eeq{schi}
and $\bfpzc G_0=\mathpzc G_{0i}\,\bfe_i$  with components
\be\ba{rl}\medskip
\mathpzc G_{0i}:=&\!\!\!\!-\Int{\varOmega}{}\big\{[2h^2\,\bar{(\bfV_0-\bfu_*-\bfzeta_0)\cdot\nabla\bfV_0}+\bar{(\nabla\bsf s)^\top:\nabla\, \mathbb T(\bfV_0,{\sf p}_0)}]\cdot\bfh^{(i)}\\\medskip &\!\!\!\!\hspace*{1cm}-\big(\bar{\nabla\bsf{s}\cdot\nabla\bfV_0+(\nabla\bsf{s}\cdot\nabla\bfV_0)^\top}\big):\mathbb D(\bfh^{(i)})\big\}\\
&\!\!\!\!-\bfe_i\cdot\Int{\varSigma}{}\bar{\mathbb T(\bfV_0,{\sf p}_0)\cdot\mathbb H^\top(\bsf{s})}\cdot\bfn{\rm d}\varSigma\ \ i=1,2,3\,,
\ea
\eeq{2_2_2}
where the bar means average over the interval $[0,2\pi]$.
\par
Let $(\bfe_r,\bfe_\theta,\bfe_\varphi)$ be a local frame of spherical coordinates with zenith direction along $\bfe_3$, and set
\be\ba{ll}\medskip
\psi(x)=\Frac{\sin\theta\,\cos\varphi}{r^2}\,,\ \ \bfa(x):=\nabla\psi(x)\\
\medskip
\bfG_1=\Frac{{\rm e}^{-h(r-1)}}{r^2}\left(g_1(r)\cos[h(r-1)]+g_2(r)\sin[h(r-1)]\right)\sin\theta\bfe_\varphi:=G_1(r)\sin\theta\bfe_\varphi\\ \medskip
\bfG_2=\Frac{{\rm e}^{-h(r-1)}}{r^2}\left(g_1(r)\sin[h(r-1)]-g_2(r)\cos[h(r-1)]\right)\sin\theta\bfe_\varphi=G_2(r)\sin\theta\bfe_\varphi
\\ \medskip
g_1(r)=\Frac{1+h(r+1)+2h^2r}{1+2h+2h^2}\\
g_2(r)=\Frac{h(r-1)}{1+2h+2h^2}\,.
\ea
\eeq{00}
The displacement field is then chosen of the following form:
\be
\bsf s(x,t):=\bfa(x)\,\sin t -\!\!\Int{}{}\bfzeta_0(t)\,{\rm d}t+[\bfG_1(r,\theta,\varphi)\cos  t+\bfG_2(r,\theta,\varphi)\sin t]:=\bsf s_1(x,t)+\bsf s_2(x,t)\,.
\eeq{01}
From the physical viewpoint, $\bsf s$ is given, up to a function of time, by the sum of the  a time-periodic deformation associated to a dipole flow pattern ($\bsf s_1$) and a rigid oscillation  around  $\bfe_3$ ($\bsf s_2$). In other words, the body $\mathscr B$ deforms and oscillates appropriately.
Next, we introduce the following functions
\be\ba{ll}\medskip
\bfV_0(x,t):=\bfa(x)\cos t+[\bfG_2(r,\theta,\varphi)\cos t-\bfG_1(r,\theta,\varphi)\sin  t]:=\bfV_1(x,t)+\bfV_2(x,t)
\\
p_0(x,t):=-2h^2\,\psi(x)\,\sin  t\,,\ \ \bfzeta_0(t):=-\Frac1 {\sf M}\cos  t\int_{\varSigma}\psi\bfe_r\equiv -\Frac{4\pi} {3 \sf M}\cos  t\,\bfe_1\,.
\ea
\eeq{0_0}
Due to \cite[Lemma 9.1]{GaspU} we deduce that the fields
$$
\bfV_1(x,t)\,,\ \ p_0(x,t)\,,\ \ \bfzeta_0(t)\,,\ \ \bfu_{*1}:=\partial_t\bsf s_1
$$
satisfy the following Stokes problem
$$\ba{cc}\medskip\left.\ba{ll}\medskip
2h^2\,\partial_t\bfV_1=\Delta\bfV_1-\nabla p_0\\
\Div\bfV_1=0\ea\right\}\ \ \mbox{in $ {\varOmega}\times\real$}\\ \medskip
\bfV_0(x,t)=\bfu_{*1}(x,t)+\bfzeta_0(t)\ \ (x,t)\in\varSigma\times\real\\
{\sf M}\dot{\bfzeta}_0+\Int\varSigma{}\mathbb T(\bfV_1,p_0)\cdot\bfn=\0\ \ \mbox{in $\real$}\,.
\ea
$$
Furthermore,
from \cite[p. 220]{Berker} we know that
$$\ba{cc}\medskip\left.\ba{ll}\medskip
2h^2\,\partial_t\bsf s_2=\Delta\bsf s_2\\
\Div\bsf s_2=0\ea\right\}\ \ \mbox{in $ {\varOmega}\times\real$}\\ \medskip
\bsf s_2(x,t)=\bfG_1(1,\theta,\varphi)\cos t\ \ (x,t)\in\varSigma\times\real\\
\Int\varSigma{}\mathbb T(\bsf s_2,0)\cdot\bfn=\0\ \ \mbox{in $\real$}\,.
\ea
$$
As a result, observing that $\bfV_2=\partial_t\bsf s_2$, from all the above we conclude that $(\bfV_0,p_0,\bfzeta_0)$ is a solution to \Eqref{V0} with $\bfu_*(x,t):=\partial_t\bsf s(x,t)$\,.
\par
Concerning the auxiliary fields $\bfh^{(i)}$, $i=1,2,3$, it is well known that, in this case, they are given by \cite[Chapter 5]{happelbrenner}
\be
(\bfh^{(i)})_\ell=\mbox{$\frac34$}\,\Frac{x_ix_\ell}{r^3}\left(1-\Frac1{r^2}\right)+\mbox{$\frac14$}\left(\Frac3r+\Frac1{r^3}\right)\,\delta_{il}
\,,\ \ \ \ell=1,2,3\,.
\eeq{h}
\par
With the choice \Eqref{01}--\Eqref{h}, in the Appendix D it is shown that the the components of the thrust vector $\mathpzc G$ assume the following form
\be\ba{ll}\medskip
\mathpzc G_{i}= \Int\varOmega{} \left\{\nabla\bfa:(\nabla\bfG_1\cdot\nabla\bfh^{(i)})-\mathbb D(\bfG_1):(\nabla\bfa\cdot\nabla\bfh^{(i)})\right.\\ \hspace*{1.6cm}+\left.\half\big[\nabla\bfa\cdot(\nabla\bfG_1-(\nabla\bfG_1)^\top)-(\nabla\bfG_1-(\nabla\bfG_1)^\top)\cdot\nabla\bfa\big]:\mathbb D(\bfh^{(i)}\right\}-h^2\Int\varOmega{}\bfh^{(i)}\cdot\nabla\bfG_2\cdot\bfa + \mathcal R_i
\ea
\eeq{th}
where
\be
\ba{ll}\medskip
\mathcal R_i:=
\Int{\varOmega}{}\left\{\mathbb D(\bfG_2):(\nabla\bfG_1\cdot\nabla\bfh^{(i)})-\mathbb D(\bfG_1):(\nabla\bfG_2\cdot\nabla\bfh^{(i)})\right.\\
\hspace*{2cm}\left.+\half\big[\nabla\bfG_1\cdot\nabla\bfG_2+(\nabla\bfG_1\cdot\nabla\bfG_2)^\top-\big(\nabla\bfG_2\cdot\nabla\bfG_1+(\nabla\bfG_2\cdot\nabla\bfG_1)^\top\big)\big]:\mathbb D(\bfh^{(i)})\right\}\,.
\ea
\eeq{th1}

\begin{figure}[tbp] 
  \centering
  \includegraphics[width=4.93in,height=1.78in,keepaspectratio]{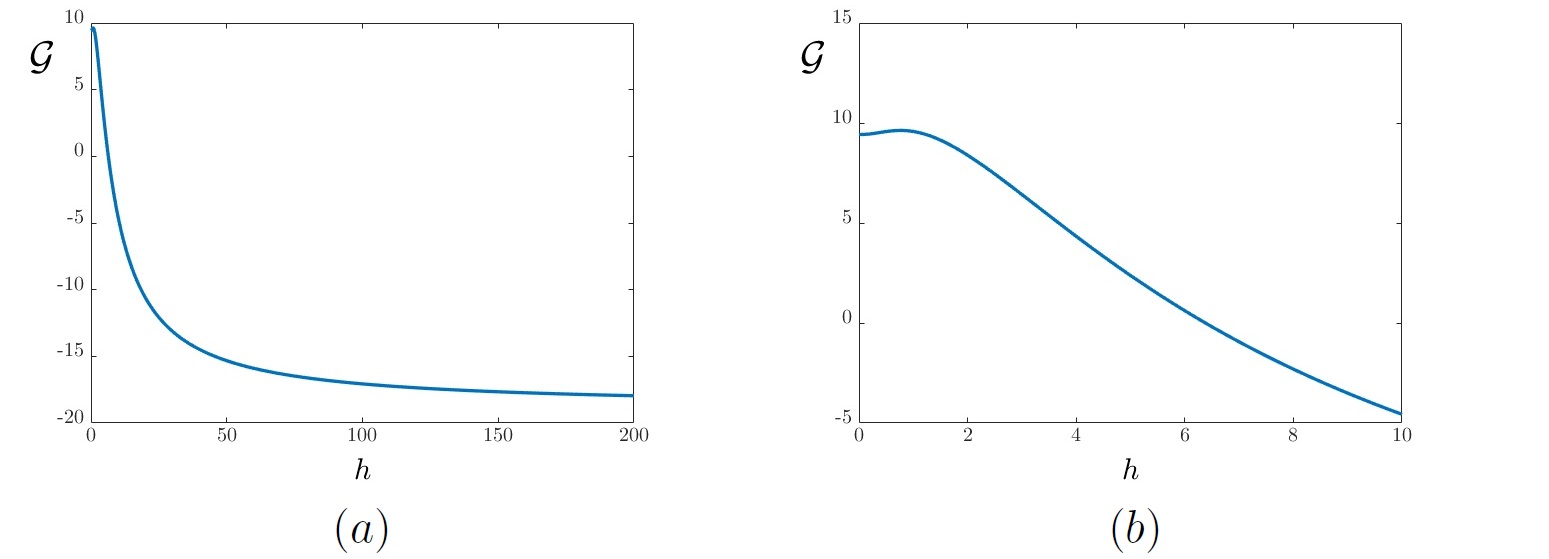}
  \caption{Dependence of the the thrust on the Stokes number $h$}
  \label{graph_h_G}
\end{figure}

From \Eqref{00}, \Eqref{h}--\Eqref{th1} it follows  that $\mathpzc G_{i}$ is a function only of the Stokes number $h$. Such a function can be explicitly evaluated by computing the integrals on the right-hand side of \Eqref{th}. This is done numerically.
Precisely, we fixed the range of variation of $h$ to be the interval $[0,200]$. Then, at each value of $h$ in this range, the components $\mathpzc G_{i}$, $i = 1, 2, 3,$ of the thrust vector were 
computed in MATLAB. The expressions for the integrands were then determined using the Symbolic 
Math Toolbox and then cast into MATLAB functions, which were then integrated numerically 
using the integral3 function with an absolute tolerance of $1\cdot{\rm e}^{-6}$. The computation shows that $\bfcalr= 0$ for all admissible values of $h$. This means --as expected-- that the propulsive thrust is the result only of the interaction of the periodic deformation of the ball with its oscillations around the $\bfe_3$ axis. Furthermore, it turns out that the components  $\calg_1$ and $\calg_3$ both vanish, and so  we conclude that
$$  
\mathpzc G= \calg(h)\,\bfe_2\,.
$$
The function $\calg(h)$ is plotted in {Figure \ref{graph_h_G}} and shows some interesting features.     
In fact, for small values of $h$, approximately in the interval $(0,6)$, the propulsion velocity $\bfgamma_1$ of the center of mass of the body  has the same orientation (and direction) as $\bfe_2$. However, $\bfgamma_1$ decreases in magnitude and when $h$ reaches a critical value $h_0\simeq 6$ it becomes 0; see {Figure \ref{graph_h_G}} (b). For $h>h_0$ the body starts moving in the opposite direction with $|\bfgamma_1|$ increasing with increasing $h$, until $|\bfgamma_1|$ reaches a constant value independent of $h$, for $h\gtrapprox150$. This, in particular, furnishes an ``optimal" frequency of oscillation  that maximize the speed,  given by $\omega_{opt}\simeq 2(150)^2\nu/a^2$. It is worth emphasizing that these are {\em entirely new features} that are not present in any of analogous researches performed on similar but  {\em linearized} models; see, e.g., \cite{FJ1,FJ2}. The reason is that the functions $\bsf g_0$ and $\bsf G_0$ in \Eqref{11.1}, which characterize the thrust, involve a certain number of terms that may be missing in a linearization procedure imposed {\em directly} on the starting equations.   

\setcounter{equation}{0}\renewcommand{\theequation}{A.\arabic{equation}}
\section*{Appendix A: Proof of Lemma \ref{le:calf}}\label{Ap:A} 
We begin to present some  properties of the fields $\bfV$ and $\mathbb A$ that will be useful also in a different context. By \propref{O_k} and classical embedding results we deduce, in particular,
\be\ba{c}\medskip
\|\bfzeta\|_{W^{1,\infty}(0,T)}+\|\bfV\|_{W^{1,\infty}(L^\infty(\varOmega))}\le C\,,\\ \medskip
\|\nabla\bfV\|_{L^\infty(L^q(\varOmega))}\le \|\nabla\bfV\|_{L^\infty(L^3(\varOmega_{\rho_*}))}+\|\nabla\bfV\|_{L^\infty(L^q({\varOmega}^{\rho_*}))}\le C\,,\ \ q\in(1,3]\,,\\
\|\partial_t\bfV\|_{L^\infty(W^{1,r})}\le C\,,\ \ r\in [2,6]\,,
\ea
\eeq{6.32}
as well as, by  \Eqref{AMPP_0}, \Eqref{west1_1} and again classical embedding,
\be
\|D^2\bfV\|_{L^\infty(L^2)}+\|\nabla p\|_{L^\infty(L^2)}\le C\,.
\eeq{6.36}
We also recall that from \Eqref{Si0} it follows, in particular, 
\be
\|\mathbb A\|_{W^{1,\infty}(L^\infty)}\le C\,.
\eeq{a3}
In the above inequalities, and in the rest of the proof, $C$ represents a constant depending at most on $\Omega_0$, $T$, the physical parameters and (the given) $\bsf s$.
Now, in view of \lemmref{fF} with $q=\frac65,2$ we deduce
$$
\|\partial_t^k\bff_\ell({\bsf s},\bfV,{\sf p})\|_q\le c\,\delta\left(\|\partial_t^k\nabla\bfV\|_{1,q,B_0}+\|\partial_t^k\nabla{\sf p}\|_{q,B_0}\right)\,,\ \ k=0,1\,,
$$
and therefore, thanks to the boundedness of $B_0$ and \propref{O_k},
\be
\|\bff_\ell({\bsf s},\bfV,{\sf p})\|_{W^{1,2}(L^2)}+\|\bff_\ell({\bsf s},\bfV,{\sf p})\|_{W^{1,2}( L^\frac65)}\le c\,\delta\left(\|\nabla\bfV\|_{W^{1,2}(L^2)}+\|\nabla{\sf p}\|_{W^{1,2}(L^2)}\right)\le C\,\delta\,.
\eeq{a4}
Furthermore, in view of \Eqref{6.32} and \Eqref{a3}, we show 
\be
\|(\bfV-\zeta)\cdot\mathbb A\cdot\nabla\bfV\|_q\le c\,\|\nabla\bfV\|_q\le C\,,\ \ q=\mbox{$\frac65,2$}.
\eeq{a5}
Likewise, we prove
\be
\|\partial_t[(\bfV-\zeta)\cdot\mathbb A\cdot\nabla\bfV]\|_{L^2(L^2)}\le c\,\|\nabla\bfV\|_{W^{1,2}(L^2)}\le C\,.
\eeq{a6}
Finally, since from \lemmref{fF} we have
$$
\|\bfF(\bsf s,\bfV,{\sf p})\|_{W^{1,2}(0,T)}\le c\,\delta \left(\|\nabla\bfV\|_{W^{1,2}(L^2)}+\|\nabla{\sf p}\|_{W^{1,2}(L^2)}\right)\,,
$$
by \propref{O_k} we get
\be
\|\bfF(\bsf s,\bfV,{\sf p})\|_{W^{1,2}(0,T)}\le C\,\delta\,.
\eeq{a7}
The property stated in \Eqref{calf} then follows from \Eqref{a4}--\Eqref{a7}
\setcounter{equation}{0}\renewcommand{\theequation}{B.\arabic{equation}}
\section*{Appendix B: Proof of \Eqref{6.20f} and \Eqref{6.21f}}\label{Ap:B}
We first recall or derive  some useful inequalities, which we shall state for $\bfv$, but with the observation that they are equally valid for $\tilde{\bfv}$. We begin to observe that, by \Eqref{1.8}, we have
\be
\|\nabla\bfv\|_2=\sqrt{2}\|\mathbb D(\bfv)\|_2
\eeq{6.20}
by \Eqref{1.9}
\be
\|\bfv\|_6\le \kappa_0\|\mathbb D(\bfv)\|_2
\eeq{6.21}
and, by \Eqref{1.10},
\be
|\bfxi|\le \kappa_1\|\mathbb D(\bfv)\|_2\,.
\eeq{6.22}
Moreover, from \Eqref{hey} and  \Eqref{6.20} we obtain 
\be
\|D^2\bfv\|_{2}+\|\nabla{\sf r}\|_{2}\le c\,\left(\|\Div\mathbb T(\bfv,{\sf r})\|_{2}+\|\mathbb D(\bfv)\|_2\right)\,,
\eeq{6.23}
where the constant $c$ depends solely on the regularity of $\varSigma$.
Thus, employing \Eqref{6.20}  in the  well-known inequality
\be
\|\nabla\bfw\|_{3,{\varOmega}_R}\le c\,\left(\|\nabla\bfw\|_{2,{\varOmega}_R}+\|\nabla\bfw\|_{2,{\varOmega}_R}^\frac12\|D^2\bfw\|_{2,{\varOmega}_R}^\frac12\right)\,,\ \ \bfw\in W^{2,2}(\varOmega_R)\,,
\eeq{ww}
with $c$ independent of $R$ \cite[p. 27]{Fri}, it follows 
\be
\|\nabla\bfv\|_3\le c\,\left(\|\mathbb D(\bfv)\|_{2}+\|\mathbb D(\bfv)\|_{2}^\frac12\|D^2\bfv\|_{2}^\frac12\right)
\eeq{vv}
or also, by \Eqref{6.23},
\be
\|\nabla\bfv\|_3
\le c\,\left(\|\mathbb D(\bfv)\|_2+\|\mathbb D(\bfv)\|_2^\frac12\|\Div\mathbb T(\bfv,{\sf r})\|_2^\frac12\right)
\eeq{A1}
From the latter in conjunction with H\"older inequality and \Eqref{6.20}, \Eqref{6.21}, we infer
\be
\||\tilde{\bfv}|\,\nabla\bfv\|_2\le c\,\|\tilde{\bfv}\|_6\|\nabla\bfv\|_3\le c\,\|\mathbb D(\tilde{\bfv})\|_2\left(\|\mathbb D(\bfv)\|_2+\|\mathbb D(\bfv)\|_2^\frac12\|D^2\bfv\|_2^\frac12\right)
\eeq{Nato}
as well as
\be
\||\tilde{\bfv}|\,\nabla\bfv\|_2\le c\,\|\mathbb D(\tilde{\bfv})\|_2\left(\|\mathbb D(\bfv)\|_2+\|\mathbb D(\bfv)\|_2^\frac12\|\Div\mathbb T(\bfv,{\sf r})\|_2^\frac12\right)\,.
\eeq{6.24}
We next derive estimates involving the operator $\bfB$. From \propref{Boo} and \Eqref{6.16} we deduce
\be
\|\bfB(\bfv)\|_{m,2}\le c\,\delta\,\|\bfv\|_{m,2,{\varOmega}_{\rho}}\,,\ \ m=1,2\,, 
\eeq{6.25}
which, with the help of \Eqref{6.25} and \Eqref{6.20}--\Eqref{6.21}, gives for $m=1$
\be
\|\bfB(\bfv)\|_{1,2}\le c\,\delta\,\|\mathbb D(\bfv)\|_2\,,
\eeq{B12}
whereas   for $m=2$ furnishes
\be
\|\bfB(\bfv)\|_{2,2}\le c\,\delta\,\left(\|\mathbb D(\bfv)\|_2+\|D^2\bfv\|_2\right)
\eeq{6.26_00}
or also, with the help of \Eqref{6.23},
\be\|\bfB(\bfv)\|_{2,2}
\le c\,\delta\,\left(\|\Div\mathbb T(\bfv,{\sf r})\|_2+\|\mathbb D(\bfv)\|_2\right)\,.
\eeq{6.26}
Likewise, from \lemmref{Def},  \propref{Boo} and \Eqref{6.21}, \Eqref{6.22} we show
\be\ba{ll}\medskip
\|\bfB(\bfv)\|_2\le c\,\delta\,\big(\|\bfv\|_{2,{\varOmega}_\rho}+|{\bfxi}|\big)\le c\,\delta\,\|\mathbb D(\bfv)\|_2\\
\|\bfB_t(\bfv)\|_2\le c\,\delta\,\big(\|\bfv\|_{2,{\varOmega}_\rho}+|{\bfxi}|+\|\bfv_t\|_2+|\dot{\bfxi}|\big)\le c\,\delta\,\big(\|\bfv_t\|_2+|\dot{\bfxi}|+\|\mathbb D(\bfv)\|_2\big) \,.
\ea
\eeq{6.27}
Moreover, 
using \Eqref{ww}, it follows that
$$ 
\|\nabla\bfB(\bfv)\|_3\le c\,\left(\|\nabla\bfB(\bfv)\|_2+
\|\nabla\bfB(\bfv)\|_2^\frac12\|\nabla\bfB(\bfv)\|_{2,2}^\frac12
\right)\,,
$$
which, in turn, 
by \Eqref{B12} and \Eqref{6.26}, 
furnishes
\be
\|\nabla\bfB(\bfv)\|_3\le c\,\delta\left(\|\mathbb D(\bfv)\|_2+\|\mathbb D(\bfv)\|_2^\frac12\|D^2\bfv\|_{2,2}^\frac12\right)
\eeq{Nuo}
as well as
\be
\|\nabla\bfB(\bfv)\|_3\le c\,\delta\left(\|\mathbb D(\bfv)\|_2+\|\mathbb D(\bfv)\|_2^\frac12\|\Div\mathbb T(\bfv,{\sf r})\|_{2,2}^\frac12\right)\,.
\eeq{Nuo1}
Thus, using \Eqref{Nuo1}, we get
$$
\||\bfB(\tilde{\bfv})|\,\nabla\bfB(\bfv)\|_2\le \|\bfB(\tilde{\bfv})\|_6\|\nabla\bfB(\bfv)\|_3\le c\,\delta\|\bfB(\tilde{\bfv})\|_{1,2}\left(\|\mathbb D(\bfv)\|_2+\|\mathbb D(\bfv)\|_2^\frac12\|\Div\mathbb T(\bfv,{\sf r})\|_{2,2}^\frac12\right)\,,
$$
which, after using \Eqref{B12}, produces
\be
\||\bfB(\tilde{\bfv})|\,\nabla\bfB(\bfv)\|_2\le c\,\delta^2 \|\mathbb D(\tilde{\bfv})\|_2\left(\|\mathbb D(\bfv)\|_2+\|\mathbb D(\bfv)\|_2^\frac12\|\Div\mathbb T(\bfv,{\sf r})\|_2^\frac12\right)\,,
\eeq{6.30}
as well as, by\Eqref{Nuo},
\be
\||\bfB(\tilde{\bfv})|\,\nabla\bfB(\bfv)\|_2\le c\,\delta^2 \|\mathbb D(\tilde{\bfv})\|_2\left(\|\mathbb D(\bfv)\|_2+\|\mathbb D(\bfv)\|_2^\frac12\|D^2\bfv\|_2^\frac12\right)\,,
\eeq{Nato1}
By a similar argument, one can prove also
\be
\||\bfB(\tilde{\bfv})|\nabla\bfv\|_2+\||\tilde{\bfv}|\nabla\bfB(\bfv)\|_2\le c\,\delta \|\mathbb D(\tilde{\bfv})\|_2\left(\|\mathbb D(\bfv)\|_2+\|\mathbb D(\bfv)\|_2^\frac12\|\Div\mathbb T(\bfv,{\sf r})\|_2^\frac12\right)\,,
\eeq{6.31}
as well as
\be
\||\bfB(\tilde{\bfv})|\nabla\bfv\|_2+\||\tilde{\bfv}|\nabla\bfB(\bfv)\|_2\le c\,\delta \|\mathbb D(\tilde{\bfv})\|_2\left(\|\mathbb D(\bfv)\|_2+\|\mathbb D(\bfv)\|_2^\frac12\|D^2\bfv\|_2^\frac12\right)\,.
\eeq{Nato2}
Thanks to the above inequalities, we shall next estimate the right-hand sides of \Eqref{6.19}.
Integrating by parts as necessary and recalling that $\bfB(\bfv_N)|_\varSigma=0$ we show
$$
(\bfH(\tilde{\bfv},\bfv),\bfv)=-(\bfB_t(\bfv),\bfv)+\half(\Div\bfB(\tilde{\bfv}),|\bfv|^2)+\big((\tilde{\bfv}-\tilde{\bfxi})\cdot\nabla\bfv,\bfB(\bfv)\big)\,.
$$
From \Eqref{6.19}$_1$ with the help of \Eqref{6.32}$_1$ and Schwarz inequality, we get
$$
|(\bfH(\tilde{\bfv},\bfv),\bfv)|\le c\,\big[\|\bfB_t(\bfv)\|_2\|\bfv\|_{2,{\varOmega}_\rho}+\|\nabla\bfB(\tilde{\bfv})\|_2\|\bfv\|_{4,{\varOmega}_\rho}^2+\|\bfB(\bfv)\|_2(\|\tilde{\bfv}\cdot\nabla\bfv\|_2
+|\tilde{\bfxi}|\,\|\nabla\bfv\|_2
)\big]\,.
$$
Thus, exploiting in the latter \Eqref{6.20}--\Eqref{6.22},  \Eqref{6.24}, \Eqref{6.27}, \Eqref{B12}  along with the Young inequality $ab\le r^{-1}a^r+s^{-1}b$, $a,b>0$, $ r^{-1}+ s^{-1}=1$ we obtain
\be\ba{rl}\medskip
|(\bfH(\tilde{\bfv},\bfv),\bfv)|&\!\!\!\le c\,\delta\left[\|\mathbb D(\bfv)\|_2^2+\|\mathbb D(\tilde{\bfv})\|_2(\|\mathbb D(\bfv)\|_2^2+\|\mathbb D(\bfv)\|_2^{\frac{3}2}\|\Div T(\bfv,{\sf r})\|_2^{\frac12})+\|\bfv_t\|_2^2+|\dot{\bfxi}|^2\right]
\\
&\!\!\!\le c\,\delta\,\left[\|\mathbb D(\bfv)\|_2^2+(\|\mathbb D(\tilde{\bfv})\|_2+\|\mathbb D(\tilde{\bfv})\|_2^{\frac43})\|\mathbb D(\bfv)\|_2^{2}+\|\bfv_t\|_2^2+|\dot{\bfxi}|^2+\|\Div T(\bfv,{\sf r})\|_2^2\right]
\,,\ea
\eeq{6.33}
which proves \Eqref{6.20f}$_1$. Next,  with the help of \Eqref{6.21} and \Eqref{6.22} (with $\bfv\equiv\tilde{\bfv}$) and Schwarz inequality, we show
$$
\|\bfH(\tilde{\bfv},\bfv)\|_2^2\!\le\! c\,\big[\|\bfB_t(\bfv)\|_2^2+\|\tilde{\bfv}\cdot\nabla\bfv\|_2^2\!+\!\|\bfB(\tilde{\bfv})\cdot\nabla\bfv\|_2^2+\|\tilde{\bfv}\cdot\nabla\bfB(\bfv)\|_2^2\!+\!\|\mathbb D(\tilde{\bfv})\|_2^2(\|\mathbb D(\bfv)\|_2^2\!+\!\|\nabla\bfB(\bfv)\|_2^2)\big]\,.
$$
If we use in the latter \Eqref{6.21}, \Eqref{6.24}, \Eqref{6.27}, \Eqref{B12}, \Eqref{6.31} and Young inequality, we arrive at the following estimate
\be\ba{ll}\medskip
\|\bfH(\tilde{\bfv},\bfv)\|_2^2\le c\big[\delta^2(\|\bfv_t\|_2^2+|\dot{\bfxi}|^2)+(1+\delta^2)\|\mathbb D(\tilde{\bfv})\|_2^2(\|\mathbb D(\bfv)\|_2^2+\|\mathbb D(\bfv)\|_2\|\Div\mathbb T(\bfv,{\sf r})\|_2)+\delta^2\|\mathbb D(\bfv)\|_2^2\big]\\ \medskip
\hspace*{.7cm}\le c\big[\delta^2(\|\bfv_t\|_2^2+|\dot{\bfxi}|^2)+(1+\delta^2)(\|\mathbb D(\tilde{\bfv})\|_2^2+\|\mathbb D(\tilde{\bfv})\|_2^4)\|\mathbb D({\bfv})\|_2^2+\delta^2\|\mathbb D(\bfv)\|_2^2\big]+\eta\|\Div\mathbb T(\bfv,{\sf r})\|_2^2\,,
\ea
\eeq{H_2}
with $\eta>0$  arbitrary and $c$ depends also on $\eta$, thus proving \Eqref{6.20f}$_2$.
\par
We now turn to the proof of \Eqref{6.21f}. From \lemmref{fF}, \Eqref{6.2} and \Eqref{6.11}, with the help of H\"older inequality we get 
$$\ba{ll}\medskip
\|\bfpzc f_\ell^\sharp\|_2^2\le c\,\delta^2\big[\|\nabla(\bfv+\bfB(\bfv)+\delta\,\bfV)\|_{1,2}^2+\|\nabla{\sf q}^\sharp\|_2^2+\|\Delta\bfB(\bfv)\|_2^2+\|\nabla\Div\bfB(\bfv)\|_2^2\big]\\ 
\hspace*{2cm}+\delta^2\big[(\|\bfV\|_\infty^2+|\bfzeta|^2)(\|\nabla(\bfv+\bfB(\bfv))\|_2^2)
+|\bfxi|^2\|\nabla\bfV\|_2^2+\|\bfv+\bfB(\bfv)\|_6^2\|\nabla\bfV\|_{\frac32}^2\big]\,.
\ea
$$
Using in the latter 
\Eqref{6.32}--\Eqref{6.23} and \Eqref{6.21},\Eqref{6.22},  \Eqref{B12} and \Eqref{6.26} we deduce ($\delta\le C$)
\be
\|\bfpzc f_\ell^\sharp\|_2^2\le c\,\delta \left(\|\mathbb D(\bfv)\|_2^2+\|\Div\mathbb T(\bfv,{\sf r})\|_2^2+\delta\,\|\nabla{\sf q}^\sharp\|_2^2\right)
\eeq{A23}
which, in conjunction with  Young inequality, proves \Eqref{6.21f}$_1$. Concerning \Eqref{6.21f}$_2$, with the help of \lemmref{fF}, we show
\vs
$$
\|\bfpzc f\|_2^2\le c\big[(1\!+\!\delta^2)\||\bfB(\tilde{\bfv})|\nabla\bfB(\bfv)\|_2^2
+\delta^2(\||\tilde{\bfv}|\,\nabla\bfv\|_2^2+\||\tilde{\bfv}|\,\nabla\bfB(\bfv)\|_2^2\!+\!\||\bfB(\tilde{\bfv})|\,\nabla\bfv\|_2^2+|\tilde{\bfxi}|^2\|\nabla(\bfv\!+\!\bfB(\bfv)\|_2^2)
$$
Therefore, employing \Eqref{6.22} (for $\tilde{\bfxi},\tilde{\bfv}$), \Eqref{6.23}, \Eqref{B12}--\Eqref{6.31}, from the previous inequality it easily follows that ($\delta\le C$)
\be
\|\bfpzc f\|_2^2 \le c\,\delta\left[(\|\mathbb D(\tilde{\bfv})\|_2^2+\|\mathbb D(\tilde{\bfv})\|_2^4)\|\mathbb D(\bfv)\|_2^2+\|\Div\mathbb T(\bfv,{\sf r})\|_2^2\right]\,,
\eeq{A24}
which leads to \Eqref{6.21f}$_2$.
Finally,  from \lemmref{fF}, \Eqref{6.2}, \Eqref{6.12},  \Eqref{6.20}, \Eqref{6.32} and \Eqref{6.36}, and classical trace theorems we readily deduce
$$
|\bfpzc F^\sharp|^2\le c\,\delta^2\left(\|\bfB(\bfv)\|_{2,2}^2+\|\mathbb D(\bfv)\|_2^2+\|D^2\bfv\|_2^2+\|\nabla{\sf q}^\sharp\|_2^2\right)\,,
$$
which, with the help of  \Eqref{6.23} and \Eqref{6.26}, entails $(\delta\le\delta_0$)
$$
|\bfpzc F^\sharp|^2\le c\,\delta\left(\|\mathbb D(\bfv)\|_2^2+\|\Div\mathbb T(\bfv,{\sf r})\|_2^2+\delta\,\|\nabla{\sf q}^\sharp\|_2^2\right) \,,
$$
and hence \Eqref{6.21f}$_3$ is recovered.\medskip\par\noindent
{\bf Remark B.1} If, in showing  inequalities \Eqref{H_2}, \Eqref{A23} and \Eqref{A24}, we use \Eqref{vv}, \Eqref{6.26_00} and \Eqref{Nuo} instead of
\Eqref{A1}, \Eqref{6.26} and \Eqref{Nuo1}, respectively, we get the same inequalities with $\|\Div\mathbb T(\bfv,{\sf r})\|_2$ replaced by $\|D^2\bfv\|_2$.

\setcounter{equation}{0}\renewcommand{\theequation}{C.\arabic{equation}}
\section*{Appendix C: Proof of \Eqref{6.37} and \Eqref{6.38}}\label{Ap:C}
We begin with some preliminary considerations. First of all, we observe that \Eqref{6.20}--\Eqref{6.22} and \Eqref{6.23} continue to hold if we replace $\bfv$ by $\bfv_t$ and ${\sf r}$ by ${\sf r}_t$. From \lemmref{6.2},   \Eqref{ww} and by the assumption in Proposition \ref{EUlin} we have the following two conditions 
\be\ba{ll}\medskip
\Max{t\in[0,T]}\|\nabla\bfv(t)\|_2\le c\,\delta\,{\sf D}\\
 \essup{t\in[0,T]}\|\nabla\tilde{\bfv}(t)\|_2\le\epsilon_0\,,
\ea
\eeq{B1}
with ${\sf D}$ given in \Eqref{sha}.
Combining \Eqref{B1}  with \Eqref{ww} and \Eqref{Nuo} implies, in particular,
\be\ba{rl}\medskip
\|\nabla{\bfv}\|_3^2&\!\!\!\!\le c_1\delta\,{\sf D}\|D^2\bfv\|_2+c_2\delta^2{\sf D}^2\\
\|\nabla\bfB({\bfv})\|_3^2&\!\!\!\!\le c\,\delta^2\left(\delta\,{\sf D}\|D^2\bfv\|_2+\delta^2{\sf D}^2)\right)
\,.
\ea
\eeq{B2}
From \Eqref{Trunz}, we deduce
$$
\|\bfB_{tt}(\bfv)\|_{2}\le c\,\delta\sum_{k=0}^2\left(\|\partial_t^k\bfv\|_{2,{\varOmega}_\rho}+\left|\frac{d^k\bfxi}{dt^k}\right|\right)\,,
$$
so that, with the help of \Eqref{6.20}--\Eqref{6.22} and \Eqref{B1}$_1$, we get
\be
\|\bfB_{tt}(\bfv)\|_{2}^2\le c\,\delta^2\left(\|\mathbb D(\bfv)\|_2^2+\|\mathbb D(\bfv_t)\|_2^2+\|\bfv_{tt}\|_2^2+|\ddot{\bfxi}|^2\right)\le c_1\,\delta^2\left(\|\mathbb D(\bfv_t)\|_2^2+\|\bfv_{tt}\|_2^2+|\ddot{\bfxi}|^2\right)+c_2\delta^4{\sf D}^2\,.
\eeq{B3}
Likewise, 
\be
\|\bfB_t(\bfv)\|_{1,2}^2\le c\,\delta^2\big(\|\mathbb D(\bfv)\|_2^2+\|\mathbb D(\bfv_t)\|_2^2\big)\le c_1\delta^2\|\mathbb D(\bfv_t)\|_2^2+c_2\delta^4{\sf D}^2\,\,.
\eeq{B4}
Thus, since by
\Eqref{6.20}--\Eqref{6.23}, and Sobolev embedding theorem, 
\be\ba{ll}\medskip
\|\bfB_t({\bfv})\|_6^2\le c\,\|\bfB_t({\bfv})\|_{1,2}^2\le c_1\delta^2\left(\|\mathbb D(\bfv_t)\|_2^2+\|\mathbb D(\bfv)\|_2^2\right)\\ 
\|\nabla\bfB_t(\bfv)\|_3^2\le c\,\|\bfB_t(\bfv)\|_{2,2}^2\le c\,\delta^2\left(\|\mathbb D(\bfv)\|_2^2+\|D^2\bfv\|_2^2+\|\mathbb D(\bfv_t)\|_2^2+\|D^2\bfv_t\|_2^2\right)
\,,
\ea\eeq{Sibe}
from \Eqref{B2} and \Eqref{B1}$_1$ we deduce
\be\ba{ll}\medskip
\|\bfB_t({\bfv})\|_6^2\le c\,\|\bfB_t({\bfv})\|_{1,2}^2\le c_1\delta^2\|\mathbb D(\bfv_t)\|_2^2+c_2\delta^4{\sf D}^2\\ \|\nabla\bfB_t(\bfv)\|_3^2\le c\,\delta^2\left(\|D^2\bfv\|_2^2+\|\mathbb D(\bfv_t)\|_2^2+\|\Div\mathbb T(\bfv_t,{\sf r}_t)\|_2^2\right)+c_2\delta^4{\sf D}^2
\,.
\ea\eeq{B5}
Finally, from \Eqref{6.20}, \Eqref{B12}, \Eqref{6.31},   \Eqref{B4} and \Eqref{B1}$_1$  we show at once
\be  
\delta^2\left(\|\,|\bfV_t-\dot{\bfzeta}|\,\nabla(\bfv+\bfB(\bfv))\|_2^2+\||\bfV-{\bfzeta}|\,\nabla(\bfv_t+\bfB_t(\bfv))\|_2^2\right)\le c_1\,\delta^2\|\mathbb D(\bfv_t)\|_2^2+c_2\delta^4{\sf D}^2\,.
\eeq{B11}
Analogously, by H\"older inequality, \Eqref{6.20}--\Eqref{6.22}, \Eqref{B1}$_1$, \Eqref{B5}$_1$ and \Eqref{6.32}
\be
\delta^2\left(\||{\bfv}_t+\bfB_t(\bfv)-\dot{{\bfxi}}|\,\nabla\bfV\|_2^2+\||{\bfv}+\bfB(\bfv)-{{\bfxi}}|\,\nabla\bfV_t\|_2^2\right)\le c_1\,\delta^2\|\mathbb D(\bfv_t)\|_2^2+c_2\delta^2{\sf D}^2\|\nabla\bfV_t\|_3^2\,.
\eeq{B.12}
We now pass to the proof of \Eqref{6.37} and \Eqref{6.38}.
In the estimates that follow, we will occasionally use the fact that $\delta\le C$, for some $C>0$. We begin to increase the term  involving $\bfG(\tilde{\bfv},\bfv)$.
From H\"older inequality, \Eqref{6.20}--\Eqref{6.22}, \Eqref{B1}, \Eqref{B2}, and \Eqref{B4} 
\be\ba{ll}\medskip
\||\tilde{\bfv}_t+\bfB_t(\tilde{\bfv})-\dot{\tilde{\bfxi}}|\nabla\bfv\|_2^2\le c\,\big[(\|\tilde{\bfv}_t\|_6^2+\|\bfB_t(\tilde{\bfv})\|_6^2)\|\nabla\bfv\|_3^2+\|\mathbb D(\tilde{\bfv}_t)\|_2^2 \|\mathbb D({\bfv})\|_2^2\big] \\
\hspace*{3.6cm}\le c\,\big[(\|\mathbb D(\tilde{\bfv})\|_2^2+\|\mathbb D(\tilde{\bfv}_t)\|_2^2)(\|\mathbb D({\bfv})\|_2^2+\|\mathbb D({\bfv})\|_2\|D^2\bfv\|_2)+\|\mathbb D(\tilde{\bfv}_t)\|_2^2 \|\mathbb D({\bfv})\|_2^2\big]\,,
\ea
\eeq{C6}
which, by \Eqref{B5}$_1$, furnishes
$$
\||\tilde{\bfv}_t+\bfB_t(\tilde{\bfv})-\dot{\tilde{\bfxi}}|\nabla\bfv\|_2^2\le c\,\big[(\delta^2{\sf D}^2+\|\mathbb D(\tilde{\bfv}_t)\|_2^2)(\delta^2{\sf D}^2+\delta\,{\sf D}\|D^2\bfv\|_2) +\delta^2{\sf D}^2\|\mathbb D(\tilde{\bfv}_t)\|_2^2\big]
\,.
$$
The latter, with the help of Young inequality, implies
\be
\||\tilde{\bfv}_t+\bfB_t(\tilde{\bfv})-\dot{\tilde{\bfxi}}|\nabla\bfv\|_2^2\le c_1\,\delta^2{\sf D}^2\|\mathbb D(\tilde{\bfv}_t)\|_2^4+ c_2(1+{\sf D}^2)\|D^2\bfv\|_2^2+c_3\,\delta^2({\sf D}^2+{\sf D}^4)
\,.
\eeq{B6} 
Next, by H\"older inequality, \Eqref{6.21} and  \Eqref{B12},   we get
$$
\||\tilde{\bfv}_t-\dot{\tilde{\bfxi}}|\nabla\bfB(\bfv)\|_2^2\le \|\tilde{\bfv}_t\|_6^2\|\nabla\bfB(\bfv)\|_3^2+c_1\|\mathbb D(\bfv)\|_2^2\|\mathbb D(\tilde{\bfv}_t)\|_2^2\,,
$$
which furnishes, on the one hand, by  \Eqref{6.21} and \Eqref{Nuo},
\be
\||\tilde{\bfv}_t-\dot{\tilde{\bfxi}}|\nabla\bfB(\bfv)\|_2^2\le  c\,\|\mathbb D(\tilde{\bfv}_t)\|_2^2\left(\|\mathbb D(\bfv)\|_2^2+\|\mathbb D(\bfv)\|_2\|D^2\bfv\|_{2,2}\right)
\eeq{Sib} 
and, on the other hand, by \Eqref{6.21}, \Eqref{B1}$_1$ and \Eqref{B2}$_2$,
$$
\||\tilde{\bfv}_t-\dot{\tilde{\bfxi}}|\nabla\bfB(\bfv)\|_2^2\le \|\mathbb D(\tilde{\bfv}_t)\|_2^2\|\nabla\bfB(\bfv)\|_3^2+c_1\,\delta^2{\sf D}^2\|\mathbb D(\tilde{\bfv}_t)\|_2^2
\le c_2\,\delta^3{\sf D}\|\mathbb D(\tilde{\bfv}_t)\|_2^2\|D^2\bfv\|_2+c_3\delta^2{\sf D}^2\|\mathbb D(\tilde{\bfv}_t)\|_2^2\,,
$$
which, upon using Young inequality,
entails
\be
\||\tilde{\bfv}_t-\dot{\tilde{\bfxi}}|\nabla\bfB(\bfv)\|_2^2\le c_4\,\delta^2{\sf D}^2\|\mathbb D(\tilde{\bfv}_t)\|_2^4+c_5\|D^2\bfv\|_2^2+c_6\delta^2{\sf D}^2\,.
\eeq{B7}
Moreover, by H\"older inequality, \Eqref{6.21},  \Eqref{B4} and \Eqref{Sibe}$_2$ we show
\be\ba{rl}\medskip
\||\tilde{\bfv}-{\tilde{\bfxi}}|\,\nabla\bfB_t(\bfv)\|_2^2&\!\!\!\le  \|\tilde{\bfv}\|_6^2\|\nabla\bfB_t(\bfv)\|_3^2+c\,\|\mathbb D(\tilde{\bfv})\|_2^2\|\bfB_t(\bfv)\|_{1,2}^2
\\
&\!\!\!\le c\,\delta^2\|\mathbb D(\tilde{\bfv})\|_2^2\big(\|\mathbb D(\bfv)\|_2^2+\|D^2\bfv\|_2^2+\|\mathbb D(\bfv_t)\|_2^2+\|D^2\bfv_t\|_2^2\big)\,,\ea
\eeq{sasi}
which, in view of \Eqref{B1}$_2$ and \Eqref{6.23}
furnishes, in particular,
\be
\||\tilde{\bfv}-{\tilde{\bfxi}}|\,\nabla\bfB_t(\bfv)\|_2^2\le c\,\delta^2\big(\|\mathbb D(\bfv_t)\|_2^2+\|\Div\mathbb T(\bfv_t,{\sf r}_t)\|_2^2\big)+ c_1\,\|D^2\bfv\|_2^2+c_2\delta^2{\sf D}^2\,.
\eeq{B8}
Collecting \Eqref{B3}, and \Eqref{B6}--\Eqref{B8}
we conclude ($\delta\le C$)
\be\ba{ll}\medskip
\|\bfG(\tilde{\bfv},\bfv)\|_2^2\le 
c_1\delta (\|\bfv_{tt}\|_2^2+|\ddot{\bfxi}|^2+\|\mathbb D(\bfv_t)\|_2^2+\|\Div\mathbb T(\bfv_t,{\sf r}_t)\|_2^2)+c_3\,\delta^2{\sf D}^2\|\mathbb D(\tilde{\bfv}_t)\|_2^4\\
\hspace*{7.5cm}+c_4(1+{\sf D}^2)\|D^2\bfv\|_2^2+c_5\delta^2{\sf D}^2\big(1+{\sf D}^2\big)\,.
\ea
\eeq{Gbe}
By the well-known inequality \cite[Theorem 2.1]{CrMa}
$$
\|\bfv_t\|_4\le c\,\|\bfv_t\|_2^\frac14\|\nabla\bfv_t\|_2^\frac34
$$
with $c$ independent of $R_k$, we find
$$
|(\Div\bfB(\tilde{\bfv}),|\bfv_t|^2)|\le c\,\|\mathbb D(\tilde{\bfv})\|_2\|\bfv_t\|_4^2\le c\,\|\mathbb D(\tilde{\bfv})\|_2\|\bfv_t\|_2^\frac12\|\mathbb D(\bfv_t)\|_2^\frac32\,,
$$
so that, by \Eqref{B2}$_2$ and Young inequality we obtain
\be
|(\Div\bfB(\tilde{\bfv}),|\bfv_t|^2)|\le c\,\epsilon_0^\frac43\|\mathbb D(\bfv_t)\|_2^2+ c_1\|\bfv_t\|_2^2\,.
\eeq{B10}
Again by H\"older inequality, \Eqref{6.20}--\Eqref{6.22},  Sobolev theorem and  \Eqref{B12}, we show
\be\ba{rl}\medskip
\||\tilde{\bfv}+\bfB(\tilde{\bfv})-\tilde{\bfxi}|\,\nabla\bfv_t\|_2^2&\!\!\!\le \big(\|\tilde{\bfv}\|_6^2+\|\bfB(\tilde{\bfv})\|_6^2\big)\|\nabla\bfv_t\|_3^2+c\,\|\mathbb D(\tilde{\bfv})\|_2^2\|\mathbb D(\bfv_t)\|_2^2\\
&\!\!\!\le c\,\|\mathbb D(\tilde{\bfv})\|_2^2\left(\|\nabla\bfv_t\|_3^2+\|\mathbb D(\bfv_t)\|_2^2\right)\,,
\ea
\eeq{NaFr}
which, in turn, by \Eqref{vv}, implies
\be
\||\tilde{\bfv}+\bfB(\tilde{\bfv})-\tilde{\bfxi}|\,\nabla\bfv_t\|_2^2\le c\,\|\mathbb D(\tilde{\bfv})\|_2^2\left(\|\mathbb D(\bfv_t)\|_2^2+\|\mathbb D(\bfv_t)\|_2\|D^2\bfv_t\|_2\right)\,.
\eeq{NaFr1}
Also, using in \Eqref{NaFr} inequalities \Eqref{B1}$_2$ and \Eqref{A1} we conclude
\be
\||\tilde{\bfv}+\bfB(\tilde{\bfv})-\tilde{\bfxi}|\,\nabla\bfv_t\|_2^2\le c\,\epsilon_0^2\left(\|\mathbb D(\bfv_t)\|_2^2+\|\Div\mathbb T(\bfv_t,{\sf r}_t)\|_2^2\right)\,.
\eeq{B.9}
The validity of \Eqref{6.37} is then a consequence of   \Eqref{Gbe}, \Eqref{B10} and \Eqref{B.9}.  
We now pass to the proof of \Eqref{6.38}. From \Eqref{flin} it follows that
\be\ba{rl}\medskip
{\bfpzc{f}}_t=&\!\!\!-\bfB_t(\tilde{\bfv})\cdot\nabla\bfB(\bfv)-\bfB(\tilde{\bfv})\cdot\nabla\bfB_t(\bfv)+(\tilde{\bfv}_t+\bfB_t(\tilde{\bfv})-\dot{\tilde{\bfxi}})\cdot\mathbb B\cdot\nabla(\bfv+\bfB(\bfv))
\\
&\!\!\!+(\tilde{\bfv}+\bfB(\tilde{\bfv})-\tilde{\bfxi})\cdot\mathbb B_t\cdot\nabla(\bfv+\bfB(\bfv))+(\tilde{\bfv}+\bfB(\tilde{\bfv})-\tilde{\bfxi})\cdot\mathbb B\cdot\nabla(\bfv_t+\bfB_t(\bfv))\,.
\ea
\eeq{B14}
By H\"older inequality, \Eqref{Sibe}$_1$ and \Eqref{Nuo} we get
\be\ba{rl}\medskip
\||\bfB_t(\tilde{\bfv})|\,\nabla\bfB(\bfv)\|_2^2\le \|\bfB_t(\tilde{\bfv})\|_6^2\|\nabla\bfB(\bfv)\|_{3}^2&\!\!\!\!\le c\,\delta^4\big(\|\mathbb D(\tilde{\bfv})\|_2^2+\|\mathbb D(\tilde{\bfv}_t)\|_2^2\big)\big(\|\mathbb D({\bfv})\|_2^2+\|\mathbb D({\bfv})\|_2\|D^2\bfv\|_2\big)\,,
\ea
\eeq{Sibe1}
as well as, by  \Eqref{B2}$_2$, \Eqref{B5}$_1$ and Young's inequality ($\delta\le c_0$)
\be\ba{rl}\medskip
\||\bfB_t(\tilde{\bfv})|\,\nabla\bfB(\bfv)\|_2^2&\!\!\!\!\le c\,\delta^4\big(\delta^4{\sf D}^2+\|\mathbb D(\tilde{\bfv}_t)\|_2^2\big)\big(\delta^2{\sf D}^2+\delta\,{\sf D}\|D^2\bfv\|_2\big)\\
&\!\!\!\!\le c_1\,\delta^2{\sf D}^2\|\mathbb D(\tilde{\bfv}_t)\|_2^4+ c_2(1+{\sf D}^2)\|D^2\bfv\|_2^2+c_3\,\delta^2({\sf D}^2+{\sf D}^4)
\,.
\ea
\eeq{B15}
Again from H\"older and Sobolev inequalities, \Eqref{B12} and \Eqref{Sibe}$_2$ we show
\be
\||\bfB(\tilde{\bfv})|\,\nabla\bfB_t(\bfv)\|_2^2\le \|\bfB(\tilde{\bfv})\|_6^2\|\nabla\bfB_t(\bfv)\|_3^2\le c_1\,\delta\,\|\mathbb D(\tilde{\bfv})\|_2^2\big(\|\mathbb D(\bfv)\|_2^2+\|D^2\bfv\|_2^2+\|\mathbb D(\bfv_t)\|_2^2+\|D^2\bfv_t\|_2^2\big)\,,
\eeq{SiBe}
as well as, from \Eqref{B1}$_2$, \Eqref{B12} and \Eqref{B5}, 
\be
\||\bfB(\tilde{\bfv})|\,\nabla\bfB_t(\bfv)\|_2^2\le \|\bfB(\tilde{\bfv})\|_6^2\|\nabla\bfB_t(\bfv)\|_3^2\le c_1\,\delta\big(\|\mathbb D(\bfv_t)\|_2^2+\|\Div\mathbb T(\bfv_t,{\sf r}_t)\|_2^2\big)+c_2\|D^2\bfv\|_2^2+c_3\delta^4{\sf D}^2\,.
\eeq{B16}
Moreover, from \Eqref{B6}, \Eqref{B7} \Eqref{B15} and \lemmref{Def} we at once obtain
\be
\|(\tilde{\bfv}_t+\bfB_t(\tilde{\bfv})-\dot{\tilde{\bfxi}})\cdot\mathbb B\cdot\nabla(\bfv+\bfB(\bfv))\|_2^2
\le  c_1\,\delta^2{\sf D}^2\|\mathbb D(\tilde{\bfv}_t)\|_2^4+ c_2(1+{\sf D}^2)\|D^2\bfv\|_2^2+c_3\,\delta^2({\sf D}^2+{\sf D}^4)
\,,
\eeq{B17}
whereas from \Eqref{Nato}, \Eqref{Nato1}, \Eqref{Nato2}, \Eqref{B1}  and again \lemmref{Def} we infer
\be
\|(\tilde{\bfv}+\bfB(\tilde{\bfv})-\tilde{\bfxi})\cdot\mathbb B_t\cdot\nabla(\bfv+\bfB(\bfv))\|_2^2\le c_1\|D^2\bfv\|_2^2+c_2\delta^4{\sf D}^2\,.
\eeq{B18}
Finally, from \Eqref{B8}, \Eqref{B.9} and \Eqref{B16} we easily show
\be
\|(\tilde{\bfv}+\bfB(\tilde{\bfv})-\tilde{\bfxi})\cdot\mathbb B\cdot\nabla(\bfv_t+\bfB_t(\bfv))\|_2^2\le c\,\delta\big(\|\mathbb D(\bfv_t)\|_2^2+\|\Div\mathbb T(\bfv_t,{\sf r}_t)\|_2^2)+c_1\|D^2\bfv\|_2^2+c_2\delta^4{\sf D}^2\,.
\eeq{B19} 
From \Eqref{B14}--\Eqref{B19} we deduce
\be
\|\bfpzc{f}_{t}\|_2^2\le c\delta\big(\|\mathbb D(\bfv_t)\|_2^2\!+\!\|\Div\mathbb T(\bfv_t,{\sf r}_t)\|_2^2\big)\!+\!c_1\delta^2{\sf D}^2\|\mathbb D(\tilde{\bfv}_t)\|_2^4\!+\!c_2(1\!+\!{\sf D}^2)\|D^2\bfv\|_2^2\!+\!c_3\delta^4({\sf D}^2\!+\!{\sf D}^4)\,,
\eeq{ft}
which implies  \Eqref{6.38}$_1$. 
In order to prove \Eqref{6.38}$_2$, we begin to observe that from   \Eqref{6.2}$_1$ and \Eqref{6.11}$_1$, we have
\be\ba{rl}\medskip
\bfpzc{f}_{\ell t}^\sharp=&\!\!\!\!\big[\bff_{\ell t}(\bsf s,\bfv+\bfB(\bfv),{\sf q}^\sharp)+\nu\big(\Delta\bfB_t(\bfv)+\nabla\Div\bfB_t(\bfv)\big)\big]\! +\!\delta\,\big[(\bfV-\bfzeta)\cdot\mathbb A\cdot\nabla(\bfv+\bfB(\bfv))\big]_t\!\\
&\!\!\!\!+\delta\big[(\bfv+\bfB(\bfv)-\bfxi)\cdot\mathbb A\cdot\nabla\bfV\big]_t:= \bff_1+\bff_2+\bff_3\,.
\ea
\eeq{B20}
From \lemmref{fF} it follows that
$$
\|\bff_1\|_2^2\le c\,\delta^2\big(\|\nabla\bfv\|_{1,2}^2+\|\nabla\bfv_t\|_{1,2}^2+\|\bfB(\bfv)\|_{2,2}^2+\|\bfB_t(\bfv)\|_{2,2}^2+\|\nabla{\sf q}^\sharp\|_2^2+\|\nabla{\sf q}^\sharp_t\|_2^2\big)\,.
$$
Thus, using \Eqref{6.23}, \Eqref{6.26}, \Eqref{6.32}$_2$, \Eqref{B1}$_1$ and \Eqref{B5}$_2$ we infer
\be
\|\bff_{1}\|_2^2\le c\, \delta\left(\|\mathbb D(\bfv_t)\|_2^2+\|\Div\mathbb T(\bfv_t,{\sf r}_t)\|_2^2\right)+c_1\|D^2\bfv\|_2^2+c_2\delta^2\big(\|\nabla{\sf q}^\sharp\|_2^2+\|\nabla{\sf q}^\sharp_t\|_2^2\big)+c_3\delta^4{\sf D}^2\,.
\eeq{B21}
Furthermore, from \Eqref{6.20}, \Eqref{6.32}, \Eqref{a3}, and \Eqref{B11},
\be  
\|\bff_2\|_2^2\le c\,\delta^2\big(\|\mathbb D(\bfv_t)\|_2^2+\delta^2{\sf D}^2\big)\,.
\eeq{B22}
Employing H\"older inequality along with \Eqref{6.22} and \lemmref{Def}, we show
$$\ba{ll}\medskip
\|\bff_3\|_2^2\le c\,\delta^4\left(\|\bfv\|_6^2+\|\bfv_t\|_6^2+\|\bfB(\bfv)\|_6^2+\|\bfB_t(\bfv)\|_6^2\right)\|\nabla\bfV\|_3^2 +c\,\delta^4(\|\mathbb D(\bfv)\|_2^2+\|\mathbb D(\bfv_t)\|_2^2)\|\nabla\bfV\|_2^2\\
\hspace*{6.3cm}+c\,\delta^4\left(\|\bfv\|_6^2+\|\bfB(\bfv)\|_6^2\right)\|\nabla\bfV_t\|_3^2+c\,\delta^4\|\mathbb D(\bfv)\|_2^2\|\nabla\bfV_t\|_2^2\,.
\ea
$$
Therefore, with the help of \Eqref{6.32}, \Eqref{B12}, \Eqref{B1}$_1$ and \Eqref{B5}$_1$, we conclude  
\be
\|\bff_3\|_2^2\le c_1\,\delta\,\|\mathbb D(\bfv_t)\|_2^2+c_2\delta^4{\sf D}^2\,.
\eeq{B23}
From \Eqref{B20}--\Eqref{B23} it follows that
\be\ba{ll}\medskip
\|\bfpzc{f}_{\ell t}^\sharp\|_2^2\le c\, \delta\left(\|\mathbb D(\bfv_t)\|_2^2+\|\Div\mathbb T(\bfv_t,{\sf r}_t)\|_2^2\right)+c_1\!\|D^2\bfv\|_2^2+c_2\delta^2\big(\|\nabla{\sf q}^\sharp\|_2^2+\|\nabla{\sf q}^\sharp_t\|_2^2\big)+c_3\delta^4{\sf D}^2
\ea
\eeq{fli}
which entails \Eqref{6.38}$_2$.
Finally, from \lemmref{fF}, \Eqref{6.2}$_3$ and \Eqref{6.12} we get
$$
|\bfpzc{F}^\sharp_t|^2\le c\,\delta^2\left[\|\nabla\bfv\|_{1,2}^2+\|\nabla\bfv_t\|_{1,2}^2+\|\bfB(\bfv)\|_{2,2}^2+\|\bfB_t(\bfv)\|_{2,2}^2+\|\nabla{\sf q}^\sharp\|_2^2\right]\,.
$$
Thus, proceeding as in the proof of \Eqref{B21}, we derive
\be
|\bfpzc{F}^\sharp_t|^2\le c\,\delta\left(\|\mathbb D(\bfv_t)\|_2^2+\|\Div\mathbb T(\bfv_t,{\sf r}_t)\|_2^2\right)+c_1\|D^2\bfv\|_2^2+c_2\delta^2\|\nabla{\sf q}^\sharp_t\|_2^2+c_3\delta^4{\sf D}^2\,,
\eeq{C.26}
which  proves \Eqref{6.38}$_3$.
\medskip\par\noindent
{\bf Remark C.1} Arguing as in Remark B.1, in the inequality  \Eqref{Gbe}, \Eqref{B.9}, \Eqref{ft}, \Eqref{fli}  and \Eqref{C.26}  we can replace $\|\Div\mathbb T(\bfv_t,{\sf r}_t)\|_2$ with $\|D^2\bfv_t\|_2$
\renewcommand{\theequation}{D.\arabic{equation}}\setcounter{equation}{0}
\section*{Appendix D: Proof of \Eqref{th}--\Eqref{th1}}\label{Ap:D}
We write $\bfpzc G$ as the sum of five contributions: 
\be\ba{rl}\medskip
\mathpzc G_{i}:=&\!\!\!\!\mathpzc G_{1i}-\Int{\varOmega}{}\big\{[2h^2\,\bar{(\bfV_0-\bfu_*-\bfzeta_0)\cdot\nabla\bfV_0}+\bar{(\nabla\bsf s)^\top:\nabla\, \mathbb T(\bfV_0,{\sf p}_0)}]\cdot\bfh^{(i)}\\\medskip &\!\!\!\!\hspace*{1.4cm}-\big(\bar{\nabla\bsf{s}\cdot\nabla\bfV_0+(\nabla\bsf{s}\cdot\nabla\bfV_0)^\top}\big):\mathbb D(\bfh^{(i)})\big\}\\
&\!\!\!\!-\bfe_i\cdot\Int{\varSigma}{}\bar{\mathbb T(\bfV_0,{\sf p}_0)\cdot\mathbb H^\top(\bsf{s})}\cdot\bfn{\rm d}\varSigma:=\mathpzc G_{1i}+\bar{\cali_i^1}+\bar{\cali_i^2}+\bar{\cali_i^3}+\bar{\cali_i^4}\,,\ \ i=1,2,3\,,
\ea
\eeq{2_2_2_2}
with $\bfpzc G_1$ given in \Eqref{schi}.
We begin to notice that, by its own definition, it is
\be
\bfV_0(x,t)=\partial_t{\bsf s}(x,t)+\bfzeta_0(t)\equiv \bfu_*(x,t)+\bfzeta_0(t)\,.
\eeq{scippo}
In view of \Eqref{scippo} we get
$$
(\nabla\bsf s)^\top:\nabla\bfV_0=(\nabla\bsf s)^\top:(\partial_t\nabla\bsf s)=\half \partial_t\big[(\nabla\bsf s)^\top:\nabla\bsf s\big]\,.
$$
Therefore, by the $2\pi$-periodicity property of $\bsf s(x,t)$, we conclude
\be
\bfpzc G_1=\0\,.
\eeq{G_1}
We next evaluate  
$\bfpzc G_0$.
Again by \Eqref{scippo}, it follows at once that
\be
\cali_i^1=0\,.
\eeq{I1}
Furthermore, integrating by parts and taking into account that $\Div\bsf s=0$, we get
$$\ba{rl}\medskip
\cali_i^2=&\!\!\!-\Int\varOmega{}\partial_j{\sf s}_\ell\partial_\ell\mathbb T_{jk}h^{(i)}_k=-\Int\varOmega{}\partial_\ell\big(\partial_j{\sf s}_\ell\mathbb T_{jk}h^{(i)}_k\big) +\Int\varOmega{}\partial_j{\sf s}_\ell\mathbb T_{jk}\partial_\ell h_k^{(i)}\\
=&\!\!\!\bfe_i\cdot\Int{\varSigma}{}{\mathbb T(\bfV_0,{\sf p}_0)\cdot\mathbb H^\top(\bsf{s})}\cdot\bfn{\rm d}\varSigma+\Int\varOmega{}\mathbb T(\bfV_0,p_0):(\nabla\bsf s\cdot\nabla\bfh^{(i)})\,,
\ea
$$
which, in turn, furnishes
\be
\bar{\cali_i^2}+\bar{\cali_i^4}=\Int\varOmega{}\bar{\mathbb T(\bfV_0,p_0):(\nabla\bsf s\cdot\nabla\bfh^{(i)})}\,.
\eeq{I}
We next observe that
$$\ba{rl}\medskip
-\Int\varOmega{}{p_0\,\mathbb I:(\nabla\bsf s\cdot\nabla\bfh^{(i)})}&\!\!\!=-\Int\varOmega{}p_0\,\partial_j{\sf s}_\ell\partial_\ell h^{(i)}_j=-\Int\varOmega{}\partial_\ell\big(p_0\partial_j{\sf s}_\ell h_j^{(i)}\big)+\Int\varOmega{}p_0\partial_j{\sf s}_\ell h_j^{(i)}\\
&\!\!\!=-\Int\varSigma{}p_0\,\bfe_i\cdot\nabla\bsf s\cdot\bfn+\Int\varOmega{}\bfh^{(i)}\cdot\nabla\bsf s\cdot\nabla p_0\,.
\ea
$$
From \Eqref{01}-\Eqref{0_0} we deduce
$$
\bar{p_0\,\bfe_i\cdot\nabla\bsf s\cdot\bfn}=-2h^2\bar{\sin^2  t}\left[\psi\,\partial_i\nabla\psi\cdot\bfn-\bfe_i\cdot\nabla\bfG_2\cdot\bfn\right]\,.
$$
Hence, since $(\nabla\psi)^2|_\varSigma=1$, we show
\be
\Int\varSigma{}\psi\,\partial_i\nabla\psi\cdot\bfn=\Int\varOmega{}\partial_\ell(\psi\partial_i\partial_\ell\psi)=\half\Int\varSigma{}(\nabla\psi)^2\bfe_i\cdot\bfn=0\,.
\eeq{6}
Moreover, from \Eqref{00} it follows that in the frame $(\bfe_r,\bfe_\theta,\bfe_\varphi)$ we have
$$\nabla\bfG_2(1,\theta,\varphi)=\left(\ba{ccc}\medskip 0 & 0 & G_2^\prime(1)\sin\theta\\ \medskip 
0 & 0 & G_2^\prime(1)\cos\theta
\\
-G_2(1)\sin\theta & -G_2(1)\cos\theta&0\ea\right)\,,
$$
so that, recalling that $\bfn=-\bfe_r$, and $G_2(1)=0$
$$
\nabla\bfG_2\cdot\bfn|_\varSigma=G_2(1)\sin\theta\,\bfe_\varphi=\0\,.
$$
Next, we have
$$
\bar{\nabla\bsf s\cdot\nabla p_0}=-2h^2\bar{\sin^2  t}\left(\nabla\bfa\cdot\nabla\psi+\nabla\bfG_2\cdot\nabla\psi\right)=-  h^2\left(\nabla\nabla\psi\cdot\nabla\psi+\nabla\bfG_2\cdot\nabla\psi\right)
$$
and, hence,
$$
\Int\varOmega{}\bfh^{(i)}\cdot\bar{\nabla\bsf s\cdot\nabla p_0}=- h^2\big(\Int\varOmega{}h^{(i)}_k\partial_k\partial_\ell\psi\partial_\ell \psi+\Int\varOmega{}\bfh^{(i)}\cdot\nabla\bfG_2\cdot\nabla\psi\big):=-h^2\left(I_1+I_2\right)\,.
$$
By integrating by parts and using $\Div\bfh^{(i)}=0$ and \Eqref{6}, we infer
$$
I_1=\half\int_\varSigma (\nabla\psi)^2\bfe_i\cdot\bfn=0\,,
$$
which furnishes
\be
-\int_\varOmega\bar{p_0\,\mathbb I:(\nabla\bsf s\cdot\nabla\bfh^{(i)})}=- h^2\Int\varOmega{}\bfh^{(i)}\cdot\nabla\bfG_2\cdot\bfa\,.
\eeq{7}
Successively, from \Eqref{0} and \Eqref{01} we obtain
$$\ba{rl}\medskip
\bar{\mathbb D(\bfV_0):(\nabla\bsf s\cdot\nabla\bfh^{(i)})}=&\!\!\!\left\{\bar{\cos^2  t}\left[\mathbb D(\bfa):(\nabla\bfG_1\cdot\nabla\bfh^{(i)})+\mathbb D(\bfG_2):(\nabla\bfG_1\cdot\nabla\bfh^{(i)})\right]\right.\\
&\!\!\!\hspace*{1cm}-\bar{\sin^2  t}\left.\left[\mathbb D(\bfG_1):(\nabla\bfa\cdot\nabla\bfh^{(i)})+\mathbb D(\bfG_1):(\nabla\bfG_2\cdot\nabla\bfh^{(i)})\right]\right\}
\ea
$$
and so, observing that  
$\mathbb D(\bfa)=\nabla\bfa$, and $\bar{\cos^2  t}=\bar{\sin^2  t}=\half$,  
combining \Eqref{7} with the latter and \Eqref{I}, it follows that
\be\ba{rl}\medskip
\bar{\cali_i^2}+\bar{\cali_i^4}=&\!\!\! \Int\varOmega{} \left[\nabla\bfa:(\nabla\bfG_1\cdot\nabla\bfh^{(i)})-\mathbb D(\bfG_1):(\nabla\bfa\cdot\nabla\bfh^{(i)})\right.
\\
&\!\!\!\left.\hspace*{1.2cm}+\mathbb D(\bfG_2):(\nabla\bfG_1\cdot\nabla\bfh^{(i)})-\mathbb D(\bfG_1):(\nabla\bfG_2\cdot\nabla\bfh^{(i)})\right]
-h^2\Int\varOmega{}\bfh^{(i)}\cdot\nabla\bfG_2\cdot\bfa\,.
\ea
\eeq{I24}
Finally, since
$$\ba{rl}\medskip
\bar{\nabla\bsf{s}\cdot\nabla\bfV_0+(\nabla\bsf{s}\cdot\nabla\bfV_0)^\top}=&\!\!\!\bar{\cos^2 t}\,\nabla\bfa\cdot(\nabla\bfG_1-(\nabla\bfG_1)^\top)-\bar{\sin^2t}\,(\nabla\bfG_1-(\nabla\bfG_1)^\top)\cdot\nabla\bfa\\ \medskip
&\!\!\!+\bar{\cos^2t}\,(\nabla\bfG_1\cdot\nabla\bfG_2+(\nabla\bfG_1\cdot\nabla\bfG_2)^\top)\\&\!\!\!-\bar{\sin^2t}\,(\nabla\bfG_2\cdot\nabla\bfG_1+(\nabla\bfG_2\cdot\nabla\bfG_1)^\top)\,.
\ea
$$
we deduce
\be\ba{ll}\medskip
\cali_i^3=\half\Int{\varOmega}{}\big[\nabla\bfa\cdot(\nabla\bfG_1-(\nabla\bfG_1)^\top)-(\nabla\bfG_1-(\nabla\bfG_1)^\top)\cdot\nabla\bfa\\
\hspace*{2cm}+\nabla\bfG_1\cdot\nabla\bfG_2+(\nabla\bfG_1\cdot\nabla\bfG_2)^\top-\big(\nabla\bfG_2\cdot\nabla\bfG_1+(\nabla\bfG_2\cdot\nabla\bfG_1)^\top\big)\big]:\mathbb D(\bfh^{(i)})\,.
\ea
\eeq{I2}
Consequently, collecting  \Eqref{2_2_2_2}, \Eqref{G_1}, \Eqref{I1}, \Eqref{I24} and \Eqref{I2} we show the validity of \Eqref{th}--\Eqref{th1}.
\newpage
\section*{Notation}
\renewcommand{\arraystretch}{1.5}
\begin{longtable}	
	{| p{.15\textwidth} | p{.55\textwidth} | p{.20\textwidth} |} 
	\hline
	\textbf{Label} & \textbf{Description} &
	\textbf{definition/1st appearance} 
	\\ \hline\hline 
	$\Omega(t)$, $\Omega_0$
	&
	$\Omega(t)\subset \real^3$ is a bounded, sufficiently smooth domain representing the configuration of the body at time $t$ ,  $\Omega_0=\Omega(0)$
	&
	Section \ref{Section: ProblemFormulation}
	\\ \hline
	$\mathscr E(t)$, $S(t)$
	&
	$\mathscr E(t)=\real^3\setminus \Omega(t)$ is a fluid domain,  $S(t):=\partial\Omega(t)\, (\equiv\partial\mathcal E(t))$
	&
	Section \ref{Section: ProblemFormulation}
	\\ \hline
	$\bfv$,${\sf p}$ $\bfgamma$
	&
	$\bfv$ and $\rho\,{\sf p}$ are velocity and pressure field of the liquid, $\bfgamma$ is the velocity of the center of mass
	&
	Section \ref{Section: ProblemFormulation}, equation \Eqref{Eu}
	\\ \hline
	$B_R$, $R_*$
	& $B_R$ is the open ball  with the origin in $\Omega_0$, and radius $R>R_*:={\rm diam}\,(\Omega_0)$
	& Section \ref{Section:Notation_and_FunctionSpaces}
	\\ \hline
	$\varOmega$, $\varOmega_R$, $\varOmega^R$
	& $\varOmega=\real^3\backslash\bar{\Omega_0}$, $\varOmega_R:=\varOmega\cap B_R$, $\varOmega^R:=\real^3\backslash\bar{B_R}$
	&  Section \ref{Section:Notation_and_FunctionSpaces}
	\\ \hline
	$\varOmega_{R_1,R_2}$
	& $\varOmega_{R_1,R_2}:=B_{R_2}\backslash\bar{B_{R_1}}$
	& Section \ref{Section:Notation_and_FunctionSpaces}
	\\ \hline
	$L^q$, $W^{m,2}$
	&
	$L^q=L^q(A)$, $W^{m,2}=W^{m,2}(A)$ are Lebesgue and  Sobolev spaces with norm  $\|\cdot\|_{q,A}$, and $\|\cdot\|_{m,2,A}$. 
	& Section \ref{Section:Notation_and_FunctionSpaces}
	\\ \hline
	$D^{m,q}$
	&
	$D^{m,q}=D^{m,q}(A)$ is the homogeneous Sobolev space with semi-norm $\sum_{|l|=m}\|D^lu\|_{q,A}$.
	& Section \ref{Section:Notation_and_FunctionSpaces}
	\\ \hline
	$(\,\ ,\,\ )_A$
	&
	$L^2(A)$-scalar product 
	& Section \ref{Section:Notation_and_FunctionSpaces}
	\\ \hline
	$\calk$
	&
	$\calk=\mathcal K({\sf A}):=\big\{\bfphi\in C_0^\infty({\sf A}) : 
	\exists\,\hat{\bfphi}\in\real^3 \mbox{ s.t. }\bfphi(x)=\hat{\bfphi} \mbox{ in a neighborhood of }\Omega_0\big\}$
	& Section \ref{Section:Notation_and_FunctionSpaces}, equation \Eqref{kcc}
	\\ \hline
	$\calc$, $\calc_0$
	&
	$\calc=\mathcal C({\sf A}):=\{\bfphi\in\calk({\sf A}):\ \Div\bfphi=0\ \mbox{in } {\sf A}\},$ 
	$\calc_0=\mathcal C_0({\sf A}):=\{\bfphi\in\calc({\sf A}): \hat{\bfphi}=\0\}$
	& Section \ref{Section:Notation_and_FunctionSpaces}, equation \Eqref{kcc}
	\\ \hline
	$\langle \bfphi,\bfpsi\rangle_{\sf A}$
	&
	$\langle \bfphi,\bfpsi\rangle_{\sf A}:=
	{\sf m}\,\hat{\bfphi}\cdot\hat{\bfpsi}+(\bfphi,\bfpsi)_{{\sf A}\cap\varOmega}\,,\ \ \bfphi,\bfpsi\in\calk$
	& equation \Eqref{0.0}
	\\ \hline
	$\call^2(\real^3)$
	&
	$\call^2(\real^3)=\{\bfu\in L^2(\real^3): \ \bfu=\hat{\bfu}\ \mbox{in}\ \Omega_0, \ \mbox{for some } \hat{\bfu}\in \real^3\}$
	& equation \Eqref{spazi}
	\\ \hline
	$\calh(\real^3)$
	&
	$\calh(\real^3)=\{\bfu\in \call^2(\real^3): \ \Div\bfu=0\,\}$
	& equation \Eqref{spazi}
	\\ \hline
	$\calg(\real^3)$
	&
	$\calg(\real^3):=\big\{\bfh\in \call^2(\real^3): \, \exists\, p\in D^{1,2}(\varOmega) \mbox{ s.t. } \bfh=\nabla p \mbox{ in } \varOmega,
	\mbox{and}\ \bfh=-\varpi\int_{\partial\varOmega}p\,\bfn   \ \mbox{in}\ \Omega_0\,\big\}\,.$
	& equation \Eqref{spazi}
	\\ \hline
	$\cald^{1,2}$
	&
	${\cald^{1,2}}=\big\{\bfu\in L^6(\real^3)\cap D^{1,2}(\real^3)\,;\ \Div\bfu=0\,;\,
	\bfu=\hat{\bfu} \ \mbox{in } \Omega_0\,,\ \mbox{for some } \hat{\bfu}\in\real^3 \big\}$
	& Section \ref{Section:Notation_and_FunctionSpaces}, equation \Eqref{1.7777}
	\\ \hline 
	$\call^2(B_R)$
	&
	$
	\call^2(B_R):= \{\bfphi\in L^2(B_R): \, \bfphi|_{\Omega_0}=\hat{\bfphi}\,\ \mbox{for some } \hat{\bfphi}\in\real^3\}
	$
	& Section \ref{Section:Notation_and_FunctionSpaces}
	\\ \hline 
	$\calh(B_R)$
	&
	$\calh(B_R):=\{\bfphi\in \call^2(B_R):\, \Div\bfphi=0\,,\ \bfphi\cdot\bfn|_{\partial B_R}=0\}$
	& Section \ref{Section:Notation_and_FunctionSpaces}
	\\ \hline 
	$\cald^{1,2}(B_R)$
	&
	$\cald^{1,2}(B_R):= \{\bfphi\in W^{1,2}(B_R): \, \Div\bfphi=0\,,\ \bfphi|_{\Omega_0}=\hat{\bfphi}\,\ \mbox{for some } \hat{\bfphi}\in\real^3\,,\ \bfphi|_{\partial B_R}=\0 \}$
	& Section \ref{Section:Notation_and_FunctionSpaces}
	\\ \hline 
	$\cald^{1,2}_0(B_R)$
	&
	$\cald^{1,2}_0(B_R):= \{\bfphi\in \cald^{1,2}(B_R): \, \hat{\bfphi}=\0 \}$
	& Section \ref{Section:Notation_and_FunctionSpaces}
	\\ \hline
	$W_T^{m,r}(B)$
	& 
	$
	W_T^{m,r}(B):=\{u\in W^{m,r}(I,B),\ \mbox{for all } I\subset \real:\ u\, \,\mbox{is } T\mbox{-periodic}\} 
	$
	& Section \ref{Section:Notation_and_FunctionSpaces}
	\\ \hline
	$\hat{\bsf{s}}$, $\hat{\bfchi}$, $\delta$
	& 
	$\hat{\bfchi}:(\bfx,t)\in \Omega_0\times\real\mapsto \bfx+\delta\,\hat{\bsf{s}}(x,t)\equiv\bfy\in \Omega(t)$
	& Section \ref{Section:ProblemInReferenceConfiguration}, \Eqref{Esse}-\Eqref{map}
	\\ \hline
	$\bfchi$, $\bsf{s}$
	& 
	$\bfchi:(\bfx,t)\in\real^3\times\real\mapsto \bfx+\delta\,\bsf{s}(x,t)\equiv \bfz\in \real^3$
	(extension of $\hat{\bfchi}$)
	& \lemmref{1}, equation \Eqref{7_}
	\\ \hline
	$\mathbb{A}$, $\mathbb B$, $\mathbb C$, $J$
	& 
	$\mathbb{A}:=\nabla\bfchi^{-1}\,,\ \mathbb B:=\mathbb A-\mathbb I\,,\ \ \mathbb C:=J\mathbb A-\mathbb I\,,\ \ J:={\rm det}\,\nabla\bfchi$
	& equation \Eqref{2.9}
	\\ \hline
	$B_0$
	&  $B_0=\supp \mathbb B = \supp \mathbb C$
	& \lemmref{Def}
	\\ \hline
	$\bfu$, $p$, $\bfu_*$
	&
	$\bfu$, $p$ are the fluid velocity and fluid pressure in the reference configuration,
	$ \bfu_*(x,t)=  \partial_t{\bsf{s}}(x,t),\, (x,t)\in\Omega_0\times\real$  
	& equation \Eqref{2.8}
	\\ \hline
	$\varSigma$
	& $\varSigma=\partial\Omega_0$ 
	&
	\\ \hline
	$\bff_{\ell}({\bsf{s},\bfu,p})$
	& 
	$\bff_{\ell}({\bsf{s},\bfu,p}):=\delta\bfu_*\!\cdot\!\mathbb{A}\cdot\nabla\bfu
	+\mathbb B^\top:\nabla\,\mathbb T(\bfu,p)+\nu\,\mathbb A^\top:\nabla\big(\mathbb B\cdot\nabla\bfu+(\mathbb B\cdot\nabla\bfu)^\top\big)$
	& equation \Eqref{2.13}
	\\ \hline
	$\bff_{n\ell}({\bsf{s},\bfu,\bfgamma})$
	& 
	$\bff_{n\ell}({\bsf{s},\bfu,\bfgamma}):=-(\bfu-\bfgamma)\cdot\mathbb{B}\!\cdot\!\nabla\bfu$
	& equation \Eqref{2.13}
	\\ \hline
	$\bfF(\bsf{s},\bfu,p)$
	& 
	$\bfF(\bsf{s},\bfu,p):=-\Int \varSigma{}
	\mathbb T(\bfu,{p})\cdot\mathbb C\cdot\bfn\,{\rm d}\varSigma-\nu\Int \varSigma{}J\,
	\big(\mathbb B\cdot\nabla\bfu+(\mathbb B\cdot\nabla\bfu)^\top\big)\cdot\mathbb A\cdot\bfn\,{\rm d}\varSigma$ 
	& equation \Eqref{2.13}
	\\ \hline
	$\cals^{k,m}$
	& 
	$
	\cals^{k,m}:=\big\{\bfg\in W^{k,2}_{T}(W^{m,2}(\varOmega)):\, \int_\varSigma \bfg(t)\cdot\bfn=0\,,\ \mbox{all } t\in\real; \ \bfg(x,t)=\0 \,\ \mbox{for all } |x|\ge \rho_0\ \mbox{and}\ t\in\real\big\}
	$
	& Section \ref{Section:DivOperator}
	\\ \hline
	$\cals_0^{k,m}$
	& 
	$
	\cals_0^{k,m}:=\big\{\bfw\in W^{k,2}_{T}(W^{m,2}(\varOmega)):\, \bfw|_{\varSigma}=\0\,; \ \bfw(x,t)=\0 \,\ \mbox{for all } |x|\ge \rho\big\}\,. 
	$
	& Section \ref{Section:DivOperator}
	\\ \hline
	$\bfB$
	& 
	$\bfB:\bfg\in\cals^{k,m}\mapsto \cals_0^{k,m}$ such that
	$
	\Div\bfB(\bfg)+\Div(\mathbb C^\top\cdot\bfB(\bfg))=\Div\bfg\ \ \mbox{in } \varOmega_\rho
	$
	& \propref{Boo}
	\\ \hline
	$\bfV$, ${\sf p}$, $\bfzeta$
	& 
	$T$-periodic solutions to the linear problem \Eqref{3.1_0}
	& equation \Eqref{3.1_0} 
	\\ \hline
	$\bsf V$, ${\sf q}$, $\bfchi$
	& 
	$\bfV=\bar{\bfV}+\bsf V\,,\ \ {\sf p}=\bar{\sf p}+{\sf q}\,,\ \ \bfzeta=\bar{\bfzeta}+\bfchi$
	& equation \Eqref{Split} 
	\\ \hline
	$\cald^s(\varOmega)$
	& 
	$\cald^s(\varOmega):=L^{\frac{3s}{3-2s}}(\varOmega)\cap D^{1,\frac{3s}{3-s}}(\varOmega)\cap D^{2,s}(\varOmega)\cap D^{2,2}(\varOmega)$, $s\in (1,\frac32)$
	& equation \Eqref{calw}
	\\ \hline
	$\calp^s(\varOmega)$
	& 
	$\calp^s(\varOmega):= L^{\frac{3s}{3-s}}(\varOmega)\cap D^{1,s}(\varOmega)\cap D^{1,2}(\varOmega)$, $s\in (1,\frac32)$
	& equation \Eqref{calw}
	\\ \hline
	$\calw({\varOmega})$
	& 
	$\calw({\varOmega}):=W^{3,2}_T(0,T;L^2({\varOmega}))\cap W^{2,2}_T(0,T;W^{2,2}({\varOmega}))$
	& equation \Eqref{calw}
	\\ \hline
	$\calq(\varOmega)$
	& 
	$\calq(\varOmega):=W^{2,2}_T(0,T;D^{1,2}({\varOmega})\cap L^6({\varOmega}))$
	&  equation \Eqref{calw}
	\\ \hline
	$\calv$
	& 
	$\calv:=W_T^{3,2}(W^{\frac32,2}(\varSigma))$
	&  equation \Eqref{calw}
	\\ \hline
	$\bfh^{i)},$ ${p}^{(i)}$
	& 
	solution to the Stokes problem \Eqref{4.32}
	&  equation \Eqref{4.32}
	\\ \hline
	$\mathbb M$
	& 
	$\mathbb M_{ji}=\Int{\varSigma}{}\mathbb T_{jk}(\bfh^{(i)},{ p}^{(i)})\,n_k{\rm d}\varSigma\,,\ \ i,j=1,2,3,$
	&  equation \Eqref{Matrix}
	\\ \hline
	$\bfV_0,$ ${\sf p}_0,$ $\bfzeta_0$
	& 
	$T$-periodic solution to the linearized problem \Eqref{5.1}
	&  equation \Eqref{5.1}
	\\ \hline
	$\bfpzc G_1$
	& 
	$\bfpzc G_1:=\sum_{i=1}^3\big[\int_\varOmega p^{(i)}\bar{\mathbb H(\bsf s):\nabla\bfV_0}\big]\bfe_i$
	& equation \Eqref{5.6_1}
	\\ \hline
	$\mathbb H(\bsf s)$
	& 
	$\mathbb H(\bsf s):= \Div\bsf s\,\mathbb I-(\nabla\bsf s)^\top$
	& equation \Eqref{Hh}
	\\ \hline
	$\bsf u,$ ${\sf q}$, $\bfxi$
	& 
	solution to problem \Eqref{6.1}, $\bfu=\bsf u+\delta\,\bfV\,,\,\ p= {\sf q}+\delta\,{\sf p}\,,\,\ \bfgamma=\bfxi+\delta\,\bfzeta,\,$
	& Section \ref{Section:NonlinearProblem}
	\\ \hline
	$\bsf f_{\ell}(\bsf{s},\bsf u,{\sf q},\bfxi)$
	& 
	$\bsf f_{\ell}(\bsf{s},\bsf u,{\sf q},\bfxi):=\bff_{\ell}(\bsf{s},\bsf u,0)-\mathbb B\cdot\nabla{\sf q}-\delta(\bfV-\bfzeta)\cdot\mathbb A\cdot\nabla\bsf u-\delta\,(\bsf u-\bfxi)\cdot\mathbb A\cdot\nabla\bfV$
	& equation \Eqref{6.2}
	\\ \hline
	$\bff_{\tiny \bfV}$
	& 
	$\bff_{\tiny \bfV}:=\delta^{-1}\,\bff_\ell({\bsf s},\bfV,{ \sf p})-(\bfV-\bfzeta)\cdot\mathbb A\cdot\nabla\bfV$
	& equation \Eqref{6.2}
	\\ \hline
	$\bsf F(\bsf{s},\bsf u,{\sf q})$
	& 
	$\bsf F(\bsf{s},\bsf u,{\sf q}):=\bfF(\bsf{s},\bsf u,0)-\Int\varSigma{}{\sf q}\,\mathbb  C\cdot\bfn\,{\rm d}\varSigma$
	& equation \Eqref{6.2}
	\\ \hline
	$\bfF_{\tiny \bfV}$
	& 
	$\bfF_{\tiny \bfV}:=\delta^{-1}\bfF(\bsf{s},\bfV,{\sf p})$
	& equation \Eqref{6.2}
	\\ \hline
	$\mathscr W_T({\sf A})$
	& 
	$\mathscr W_T({\sf A}):=\big\{ \bfu\in W^{1,2}_T(L^6({\sf A})\cap D^{1,2}({\sf A})\cap D^{2,2}({\sf A}));\,\ \partial_t\bfu\in W^{1,2}_T(L^2({\sf A}))\big\}$
	& Section \ref{Section:NonlinearProblem}
	\\ \hline
	$\bfv$
	& solenoidal velocity field, 
	$\bsf u=\bfv +\bfB(\bfv)$
	& equation \Eqref{6.8}
	\\ \hline
	$\bfpzc{f}_\ell(\bsf s,\bfv,{\sf q},\bfxi)$
	& 
	$\bfpzc{f}_\ell(\bsf s,\bfv,{\sf q},\bfxi):=\bsf f_\ell(\bsf s,\bfv+\bfB(\bfv),{\sf q},\bfxi)+\nu\big(\Delta\bfB(\bfv)+\nabla\Div\bfB(\bfv)\big)$
	& equation \Eqref{6.11}
	\\ \hline
	$\bfpzc{f}_{n\ell}(\bsf s,\bfv,\bfxi)$
	& 
	$\bfpzc{f}_{n\ell}(\bsf s,\bfv,\bfxi):=-\bfB(\bfv)\cdot\mathbb A\cdot\nabla\bfB(\bfv)-(\bfv+\bfB(\bfv)-\bfxi)\cdot\mathbb B\cdot\nabla\bfv-(\bfv-\bfxi)\cdot\mathbb B\cdot\nabla\bfB(\bfv)$
	& equation \Eqref{6.11}
	\\ \hline
	$\bfpzc{F}(\bsf{s},\bfv,{\sf q})$
	& 
	$\bfpzc{F}(\bsf{s},\bfv,{\sf q}):=-2\nu\int_{\varSigma}\mathbb D(\bfB(\bfv))\cdot\bfn+\bsf F(\bsf s,\bfv+\bfB(\bfv),{\sf q})$
	& equation \Eqref{6.12}
	\\ \hline
	$\varOmega_{k}$, $B_k$
	& $\varOmega_{k}\equiv \varOmega_{R_k}$, $B_k\equiv B_{R_k}$ such that
	$\bar{\varOmega_{R_k}}\subset \varOmega_{R_{k+1}}\,,\ \ \mbox{for all $k\in\nat$}\,;\ \ 
	\cup_{k\in\nat}\varOmega_{R_k}=\varOmega$
	& equation \Eqref{invdo}
	\\ \hline
	$\tilde{\bfv}$, $\tilde{\bfxi}$
	& 
	$\tilde{\bfv}$ and $\tilde{\bfxi}$ are prescribed vector field and vector function
	& Problem \Eqref{6.10_lin}
	\\ \hline
	${\bf f}$, ${\bf F}$
	& 
	${\bf f}$, ${\bf F}$  are  given $T$-periodic functions
	& Problem \Eqref{6.10_lin}
	\\ \hline
	${\bfpzc{f}}(\bsf s,\bfv,\tilde{\bfv},\tilde{\bfxi})$
	& 
	${\bfpzc{f}}(\bsf s,\bfv,\tilde{\bfv},\tilde{\bfxi}):=-\bfB(\tilde{\bfv})\cdot\mathbb A\cdot\nabla\bfB(\bfv)-(\tilde{\bfv}+\bfB(\tilde{\bfv})-\tilde{\bfxi})\cdot\mathbb B\cdot\nabla\bfv-(\tilde{\bfv}-\tilde{\bfxi})\cdot\mathbb B\cdot\nabla\bfB(\bfv)$
	& equation \Eqref{flin}
	\\ \hline
	$S_{\epsilon_0,M}$
	& 
	$ S_{\epsilon_0,M}:=\Big\{\tilde{\bfv}\in  W^{1,4}_{\rm loc}(\real;\cald^{1,2}(B_k)):\ 
	\tilde{\bfv} \ \mbox{is } T\mbox{-periodic with } \tilde{\bfxi}=\tilde{\bfv}|_{\Omega_0}; \essup{t\in[0,T]}\|\mathbb D(\tilde{\bfv}(t))\|_2\le\epsilon_0\,;\ \ \|\tilde{\bfv}\|_{W^{1,4}(0,T;\cald^{1,2})}\le M\Big\}$
	& equation \Eqref{Sepm}
	\\ \hline
	$\calf_1$, $\calf_2$
	& 
	$\calf_1:=\|{\bf f}\|_2^2+\|{\bf f}\|^2_\frac65+|{\bf F}|^2\,,\ \ \ \calf_2:=\|\partial_t{\bf f}\|_2^2+|\dot{\bf F}|^2$
	& equation \Eqref{cal}
	\\ \hline
	$F$, $\mu$
	& 
	$F:=\Big[\Int0T\big({\calf}_1+\calf_2+\calf_1^2\big){\rm d}t\Big]^\frac12\,,\ \ \mu:=M^2+1$
	& equation \Eqref{fmu}
	\\ \hline
	${\sf q}^\sharp$, $\bfpzc{f}_\ell^\sharp$, $\bfpzc{F}^\sharp$
	& 
	${\sf q}^\sharp$ is a prescribed $T$-periodic function in $\calq(\varOmega_k)$,
	$\bfpzc{f}_\ell^\sharp:=\bfpzc{f}_\ell(\bsf{s},\bfv,{\sf q}^\sharp,\bfxi)\,,\ \ \bfpzc{F}^\sharp:=\bfpzc{F}(\bsf{s},\bfv,{\sf q}^\sharp)$
	& equation \Eqref{fsha}
	\\ \hline
	$\calc_T(B_k)$
	& 
	the class of suitable "test functions", constituted by the restriction to $[0,T]$ of functions $\bfphi \in C^{1}(B_k\times \mathbb{R})$
	& Section \ref{Section:WeakFormulation}
	\\ \hline
	$\bfpsi_{i}$
	& 
	$\{\bfpsi_{i}\} \subset {\mathcal D}^{1,2}(B_{R_k})\cap W^{2,2}(\varOmega_k)$ is an orthonormal basis of ${\mathcal{H}}(B_{R_k})$ that is also orthogonal in $\cald^{1,2}(\varOmega_k)$. 
	& equations \Eqref{ortg}, \Eqref{poi}
	\\ \hline
	$\bfv_N, \bfxi_N$
	& 
	approximated solution to \Eqref{6.13}, i.e. solution to \Eqref{6.15}
	& equation \Eqref{6.14}
	\\ \hline
	$\bfH(\tilde{\bfv},\bfv)$
	& 
	$\bfH(\tilde{\bfv},\bfv):= -\big[\bfB_t(\bfv)+(\tilde{\bfv}-\tilde{\bfxi})\cdot\nabla\bfB(\bfv)+(\tilde{\bfv}+\bfB(\tilde{\bfv})-\tilde{\bfxi})\cdot\nabla\bfv\big]$
	& equation \Eqref{H}
	\\ \hline
	$\bfG(\tilde{\bfv},\bfv)$
	& 
	$\bfG(\tilde{\bfv},\bfv):=-\big[\bfB_{tt}(\bfv)+(\tilde{\bfv}_t+\bfB_t(\tilde\bfv)+\dot{\tilde{\bfxi}})\cdot\nabla\bfv
	+(\tilde{\bfv}_t+\dot{\tilde{\bfxi}})\cdot\nabla\bfB(\bfv)+(\tilde{\bfv}-{\tilde{\bfxi}})\cdot\nabla\bfB_t(\bfv)\big]$
	& equation \Eqref{6.34f}
	\\ \hline
	${\sf D}$
	& 
	${\sf D}:=\Big[\Int0T\big(\delta^2\calf_1+\|\nabla{\sf q}^\sharp\|_2^2\big){\rm d}t\Big]^\frac12$
	& equation \Eqref{sha}
	\\ \hline
	$\cald$
	&
	$\cald:=
	\|\bfv_t\|_2^2+|\dot{\bfxi}|^2+(1+{\sf D}^2)\|D^2\bfv\|_2^2+\delta^2{\sf D}^2(1+{\sf D}^2+\|\nabla\bfV_t\|_2^2+\|\nabla\bfV_t\|_3^2)+\delta^4\calf_2
	+\delta^2(\|\nabla{\sf q}^\sharp\|_2^2+\|\nabla{\sf q}^\sharp_t\|_2^2)$
	& equation \Eqref{6.39_00}
	\\ \hline
	$\calv_1$, $\calv_2$
	&
	$
	\calv_1:=\|\bff_{\tiny \bfV}\|_2^2+\|\bff_{\tiny \bfV}\|^2_\frac65+|\bfF_{\tiny \bfV}|^2\,,\ \ \ \calv_2:=\|\partial_t\bff_{\tiny \bfV}\|_2^2+|\dot{\bfF}_{\mbox{\tiny $\bfV$}}|^2\,.
	$
	& Section \ref{sec:Ok}
	\\ \hline
	$\mathscr C$
	&
	& \defref{11.1}
	\\ \hline
	$\bsf g_0$
	&
	$\bsf g_0(x,t):=(\nabla\bsf s)^\top:\nabla\, \mathbb T(\bfV_0,{\sf p}_0)+\nu\,\Div\big(\nabla\bsf{s}\cdot\nabla\bfV_0+(\nabla\bsf{s}\cdot\nabla\bfV_0)^\top\big)+(\bfV_0-\bfu_*-\bfzeta_0)\cdot\nabla\bfV_0$
	& equation \Eqref{11.1}
	\\ \hline
	$\bsf G_0$
	&
	$\bsf G_0(t):=-\Int{\varSigma}{}\big[\mathbb T(\bfV_0,{\sf p}_0)\cdot\mathbb H^\top(\bsf{s})-\nu\,\big(\nabla\bsf{s}\cdot\nabla\bfV_0+(\nabla\bsf s\cdot\nabla\bfV_0)^\top\big)\big]\cdot\bfn\,{\rm d}\varSigma$
	& equation \Eqref{11.1}
	\\ \hline
	${\bsf v}_0,{\sf r}_0,\bfsigma_0$
	& solution to \Eqref{11.2}
	& Problem \Eqref{11.2}
	\\ \hline
	$\bfpzc G_0$
	&
	$
	\mathpzc G_{0i}:=-\Int{\varOmega}{}\big\{[\bar{(\bfV_0-\bfu_*-\bfzeta_0)\cdot\nabla\bfV_0}+\bar{(\nabla\bsf s)^\top:\nabla\, \mathbb T(\bfV_0,{\sf p}_0)}]\cdot\bfh^{(i)}
	-\nu\,\big(\bar{\nabla\bsf{s}\cdot\nabla\bfV_0+(\nabla\bsf{s}\cdot\nabla\bfV_0)^\top}\big):\mathbb D(\bfh^{(i)})\big\} 
	-\bfe_i\cdot\Int{\varSigma}{}\bar{\mathbb T(\bfV_0,{\sf p}_0)\cdot\mathbb H^\top(\bsf{s})}\cdot\bfn{\rm d}\varSigma\,,\ \ i=1,2,3$
	& equation \Eqref{G0}
	\\ \hline
	$\bsf v,{\sf r},\bfsigma$
	&
	$\bsf u=\delta^2\bsf v\,,\ \ {\sf q}=\delta^2{\sf r}\,,\ \ \bfxi=\delta^2\bfsigma$
	& equation \Eqref{11.11}
	\\ \hline
	$\bfpzc G$
	&
	$\bfpzc G = \bfpzc G_0 + \bfpzc G_1$
	& Section \ref{asy}
	\\ \hline
\end{longtable}

\section*{Acknowledgements} The authors wish to thank Messrs Carter Gassler and Marc Karakouzian for their careful computation of the thrust vector. G.P.~Galdi was partially supported by the National Science Foundation Grant
DMS-2307811, while B.~Muha and A.~Rado\v{s}evi\' c were supported by the Croatian
Science Foundation under project number IP-2022-10-2962. Finally, G.P.~Galdi wishes to thank Mr. J.A.~Wein for helpful conversations.

\ed

Finally,

first part of the lemma. In order to show the second part,
\begin{thebibliography}{99}
\bibitem{Berker}Berker, R., Int\'egration des \'equations du mouvement d'un fluide visqueux incompressible.  Handbuch der Physik  Bd. VIII/2. Springer, Berlin (1963) 1--384 
\bibitem{Bore} Borelli, G., {\em On the Movement of Animals}, Springer, New York, 1989; English Translation of {\em De Motu Animalium}, Rome: A. Bernabo  (1680)
\bibitem{Chil}Childress, S.,  {\em Mechanics of Swimming and Flying}, vol. 2 of Cambridge Studies in
Mathematical Biology. Cambridge University Press, Cambridge (1981)
\bibitem{Court}Court, S.,  Existence of 3D strong solutions for a system modeling a deformable solid inside a viscous incompressible fluid, {\em J. Dynam. Differential Equations} {\bf 29} (2017)  737--782
\bibitem{Court1} Court, S.,  Stabilization of a fluid-solid system, by the deformation of the self-propelled solid. Part I: The linearized system. {\em Evol. Equ. Control Theory} {\bf 3} (2014) 59--82
\bibitem{Court2}Court, S.,  Stabilization of a fluid-solid system, by the deformation of the self-propelled solid. Part II: The nonlinear system. {\em Evol. Equ. Control Theory} {\bf 3} (2014) 83--118
\bibitem{CrMa} Crispo F., and Maremonti P., An interpolation inequality in exterior domains, {\em Rend. Sem.
Mat. Univ. Padova} {\bf 112} (2004), 11--39
\bibitem{DHR} Davis, R.H.,  Deep diving and submarine operations (6th ed.). Tolworth, Surbiton, Surrey: Siebe Gorman \& Company Ltd. (1955) p. 693.
\bibitem{LdV} da Vinci, L.  Sul volo degli uccelli, {\em Bibiotheca Regia Taurinensi} no. 23560 (1501)
\bibitem{FJ1} Felderhof, B.U. and Jones, R.B.,  Inertial effects in small-amplitude swimming
of a finite body, {\em Physica} A, {\bf 202} (1994) 94--118.
\bibitem{FJ2} Felderhof, B.U. and Jones, R.B.,  Small-amplitude swimming of a sphere,
{\em Physica} A, {\bf 202} (1994)  119--144
\bibitem{FJ3} Felderhof, B.U. and Jones, R.B.,   Swimming of a uniform deformable sphere in a viscous incompressible fluid with inertia, {\em Eur. J. Mech. B} Fluids {\bf 85} (2021) 58--67
{\bibitem{Fri} A. Friedman, {\em Partial differential equations}, Holt, Rinehart and Winston, Inc., New York-Montreal, Que.-London (1969)}
\bibitem{Gasp}Galdi, G.P., On the steady self-propelled motion of a body in a viscous incompressible fluid. {\em Arch. Ration. Mech. Anal.
} {\bf 148}, 53--88 (1999)
\bibitem{Gah} Galdi, G.P., On the motion of a rigid body in a viscous liquid: A mathematical analysis with 
applications, {\it Handbook of Mathematical Fluid Mechanics}, Elsevier Science, 653--791 (2002)
\bibitem{Gab}Galdi, G.P.,
     \newblock {\em An Introduction to the Mathematical Theory of the Navier--Stokes Equations: Steady-State Problems},
     \newblock 2$^{nd}$ edition, Springer-Verlag, New York  (2011)
\bibitem{GaZ}Galdi, G.P., Viscous flow past a body translating by
time--periodic motion with zero average, {\em Arch. Rational Mech. Anal.} {\bf 237} (2020) 1237--1269
\bibitem{GaspU}Galdi, G.P., On the Self-propulsion of a Rigid Body in a Viscous Liquid by Time-Periodic
Boundary Data, {\em J. Math. Fluid Mech.} (2020) 22-61
\bibitem{GaKyH} Galdi, G.P., and Kyed, M., Time-periodic solutions to the Navier--Stokes equations. In: Giga, Y., Novotn\'y, A. (eds.) Handbook
of Mathematical Analysis in Mechanics of Viscous Fluids, pp. 509-578. Springer, Cham (2018)
\bibitem{6}Galdi, G.P., Kyed, M.,  Time--periodic flow of a viscous liquid past a body. {Partial differential equations in fluid mechanics}, 20--49, {\em London Math. Soc. Lecture Note Ser.}, {\bf 452}, Cambridge Univ. Press, Cambridge (2018)
\bibitem{7} Galdi, G.P., Kyed, M.,  Time-periodic solutions to the Navier-Stokes equations in the three-dimensional whole-space with a non-zero drift term: asymptotic profile at spatial infinity. {\em Mathematical analysis in fluid mechanics--selected recent results}, 121--144, Contemp. Math., {\bf 710}, Amer. Math. Soc., Providence, RI (2018)
\bibitem{GaSi1}Galdi, G.P. and Silvestre, A.L., Strong solutions to the problem of motion of a rigid body in
a Navier-Stokes liquid under the action of prescribed forces and torques, Nonlinear Problems in
Mathematical Physics and Related Topics, I, Int. Math. Ser. (N. Y.), vol. 1, Kluwer/Plenum,
New York, 2002, 121--144
\bibitem{GaSi2} Galdi G.P., and Silvestre, A.L., Existence of time-periodic solutions to the Navier-Stokes
equations around a moving body, {\em Pacific J. Math.}, {\bf 223} (2006), 251--267
\bibitem{GaSi}Galdi, G.P., and Silvestre, A.L., On the motion of a rigid body in a Navier-Stokes liquid under the action of a time--periodic force. {\it Indiana Univ. Math. J.} {\bf 58} (2009) 2805--2842
\bibitem{Gray2} Gray, J.,  Study in animal locomotion IV -- the propulsive powers of the dolphin. {\em J. Exp. Biol.} {\bf 10} (1932), 192--199
\bibitem{Gray1} Gray, J.,  Study in animal locomotion I -- the movement of a fish with special reference
to the eel. {\em J. Exp. Biol.} {\bf 13} (1936), 88--104 
\bibitem{Gray}Gray, J.,  {\it Animal Locomotion},  
Weidenfeld and Nicolson, London. 
(1968)
\bibitem{happelbrenner} Happel, V. and Brenner, H., {\em Low Reynolds Number Hydrodynamics}, Prentice Hall (1965)
\bibitem{Hey}Heywood, J.G., The Navier-Stokes equations: on the existence, regularity and decay
of solutions,  {\em Indiana Univ. Math. J.} {\bf 29}, (1980) 639--681 
\bibitem{Hishida1}Hishida, T., Silvestre, A.L. and Takahashi, T.,  Optimal boundary control for steady motions of a self-propelled body in a Navier-Stokes liquid. ESAIM {\em Control Optim. Calc. Var.} {\bf 26} Paper No. 92, (2008) 42 pp
\bibitem{Hishida} Hishida, T., Silvestre, A. L., and Takahashi, T.,  A boundary control problem for the steady self-propelled motion of a rigid body in a Navier-Stokes fluid. {\em  Ann. Inst. H. Poincar\'{e} Anal. Non Lin\'{e}aire} {\bf 34} (2017)  1507--1541
\bibitem{HiMa}Hishida, T. and Maremonti, P., 
Navier-Stokes flow past a rigid body: attainability of steady solutions as limits of unsteady weak solutions, starting and landing cases, {\em J. Math. Fluid Mech.} {\bf 20} (2018) 771--800
\bibitem{KaAk}Kantorovich, L.V., and   Akilov, G.P.,   {\em Functional Analysis in Normed Spaces}, Pergamon  (1964)
\bibitem{Ke}Keller, S.R. and Wu, T.Y.,  A porous prolate-spheroidal model for ciliated
micro-organisms, {\em J. Fluid Mech.}, {\bf 80} (1977) 259--278
\bibitem{Kapa} Khapalov, A., Cannarsa, P., Priuli, F.S, and Floridia, G., Well-posedness of 2-D and 3-D swimming models in incompressible fluids governed by Navier-Stokes equations. {\em J. Math. Anal. Appl.} {\bf 429} (2015)  1059--1085
\bibitem{Ligh}Lighthill, J., {\em Mathematical Biofluiddynamics} Society for Industrial and Applied
Mathematics, Philadelphia, (1975)
\bibitem{MS}M\'{a}cha, V.,  Ne\v{c}asov\'{a}, \v{S}.,   Self-propelled motion in a viscous compressible fluid.{\em  Proc. Roy. Soc. Edinburgh} Sect. A ( 2016) {\bf 146}  415--433
\bibitem{Macha}M\'{a}cha, V.,  Ne\v{c}asov\'{a}, \v{S}.,  Self-propelled motion in a viscous compressible fluid--unbounded domains. {\em Math. Models Methods Appl. Sci}. {\bf 26} (2016) 627--643
\bibitem{Medk}Medkov\'a, D., {\it The Laplace equation. Boundary value problems on bounded and unbounded Lipschitz domains}. Springer, Cham, (2018)
\bibitem{Necasova} Ne\v{c}asov\'{a}, \v{S}., Takahashi, T., and Tucsnak, M.,  Weak solutions for the motion of a self-propelled deformable structure in a viscous incompressible fluid. {\em Acta Appl. Math.} {\bf 116} (2011) 329--352
\bibitem{Nau} Ne\v{c}asov\'a, \v{S}. Ramaswamy, M., Roy, A. and  Schl\"{o}merkemper, A., Self-propelled motion of a rigid body inside a
 density dependent incompressible fluid, {\em Mathematical Modelling of Natural Phenomena}, {\bf 16} (2021) p. 9
\bibitem{Pu}Pukhnachev, V.V., Stokes approximation in a problem of the flow around
a self-propelled body, {\em Boundary Value Problems in Mathematical Physics}, Naukova
Dumka, Kiev, ( 1990) 65--73 (in Russian).
\bibitem{Raymond} Raymond, J.-P., Vanninathan, M.,  A fluid-structure model coupling the Navier-Stokes equations and the Lam\'{e} system. {\em J. Math. Pures Appl.} {\bf 102} (2014) 546--596
\bibitem{Tuc}San Martin, J., Scheid, J-F, Takahashi, T., and M. Tucsnak,  An initial and boundary value problem modeling of fish-like swimming,  {\em Arch. Ration. Mech. Anal.}  {\bf 188} (2008) 429--455
\bibitem{SanM}San Martin, J.A., Starovoitov,V., and Tucsnak, M., 2002, Global Weak Solutions for the Two Dimensional
Motion of Several Rigid Bodies in an Incompressible Viscous Fluid {\em Arch. Rational Mech. Anal.},{ \bf 161} 
93--112
\bibitem{Se}Sennitskii, V. L.,  Self-propulsion of a body in a fluid, {\em J. Appl. Mech. Tech.
Phys.}, {\bf 31} (1990) 266--272
\bibitem{GiW}Serchi, F-G. and Weymouth, G.D.,  Underwater soft robotics, the benefit
 of body-shape variations in aquatic
 propulsion, {\em Soft Robotics: Trends, applications and challenges}, Biosystems \& Biorobotics
 {\bf 17} (2016) 36--47
\bibitem{SW1}Shapere, A. and Wilczek, F.,  
Geometry of self-propulsion at low Reynolds number.
{\em J. Fluid Mech.} {\bf 198}  (1989) 557--585
\bibitem{SW2} Shapere, A. and Wilczek, F.,  Efficiencies of self-propulsion at low Reynolds number, {\em J. Fluid Mech.} {\bf 198} (1989) 587--599
\bibitem{ALS}Silvestre, A.L., On the self-propelled motion of a rigid body in a viscous liquid and on the attainability
of steady symmetric self-propelled motions, {\em J. Math. Fluid Mech}. {\bf 4} (2002) 285--326 
\bibitem{AS2}Silvestre, A.L.,  On the self-propelled motion of a rigid body in a viscous liquid and on the attainability of steady symmetric self-propelled motions. {\em J. Math. Fluid Mech.} {\bf 4} (2002) 285--326
\bibitem{Staro}Starovoitov, V.N.,  Solvability of the problem of the self-propelled motion of several rigid bodies in a viscous incompressible fluid. {\em  Comput. Math. Appl.} {\bf 53} (2007)   413--435
\bibitem{Simo} Simon, J., Compact sets in the space $L^p(0,T;B)$, {\em Ann. Mat. Pura Appl.} {\bf 146} (1987), 65--96
\bibitem{EP}Purcell, E.M.,  Life at low Reynolds number, {\em Ame. J. Phys.} {\bf 45} (1977), 3--11
\bibitem{Ta}Taylor, G.I.,  Analysis of the Swimming of Microscopic
Organisms, {\it Proc. Royal Soc. London A,} {\bf 209} (1951) 447-461
\bibitem{Tay}Taylor, G.I., {\em Low-Reynolds-Number Flow},  Videotape,  33min, Encyclopaedia Britannica Educational Corporation
\bibitem{WBB} Wu, T.Y-T., Brokav, C.J., and C.~Brennen (Eds), {\em Swimming and Flying in Nature}, Vol.1 and Vol.2, Plenum Press, N.Y. and London (1975)
\bibitem{NaU} Zoppello, M. and Cardin, F., Swim-like motion of bodies immersed in an ideal fluid. ESAIM {\em Control Optim. Calc. Var.} {\bf 25} (2019), Paper No. 16, 38 pp


\end{thebibliography}
